\documentclass[11pt]{article}
\usepackage{amsmath}

\usepackage{color}

\usepackage{amsfonts}
\usepackage{amssymb,amsmath,amsthm, amsfonts}
\allowdisplaybreaks[4]
\usepackage{url}

\def\sqr#1#2{{\vcenter{\vbox{\hrule height.#2pt
				\hbox{\vrule width.#2pt height#1pt \kern#1pt \vrule width.#2pt}
				\hrule height.#2pt}}}}
\usepackage[
bookmarks=true,         %generates bookmarks for all entries in the "Table of Contents"
bookmarksnumbered=true, %includes chapter-/header-/section-/subsection-/... -
colorlinks=true, pdfstartview=FitV, linkcolor=blue, citecolor=blue,
urlcolor=blue]{hyperref}

\def\sqr#1#2{{\vcenter{\vbox{\hrule height.#2pt
				\hbox{\vrule width.#2pt height#1pt \kern#1pt \vrule width.#2pt}
				\hrule height.#2pt}}}}
\def\3n{\negthinspace \negthinspace \negthinspace }
\def\2n{\negthinspace \negthinspace }
\def\1n{\negthinspace }

\def\={\buildrel \triangle \over =}

\def\exp{\mathop{\rm exp}}
\def\sup{\mathop{\rm sup}}
\def\inf{\mathop{\rm inf}}

\def\({\Big (}
\def\){\Big )}
\def\[{\Big[}
\def\]{\Big]}
\def\be{\begin{equation}}
	
	\def\ee{\end{equation}}

\def\square#1{\vbox{\hrule\hbox{\vrule height#1%
			\kern#1\vrule}\hrule}}
\def\rectangle#1#2{\vbox{\hrule\hbox{\vrule height#1%
			\kern#2\vrule}\hrule}}

\font\tenbb=msbm10 \font\sevenbb=msbm7 \font\fivebb=msbm5

\newfam\bbfam
\scriptscriptfont\bbfam=\fivebb \textfont\bbfam=\tenbb
\scriptfont\bbfam=\sevenbb

\newtheorem{lemma}{Lemma}[section]
\newtheorem{remark}{Remark}[section]

\newtheorem{theorem}{Theorem}[section]
\newtheorem{corollary}{Corollary}[section]

\newtheorem{definition}{Definition}[section]
\newtheorem{proposition}{Proposition}[section]

\makeatletter

\@addtoreset{equation}{section}
\makeatother

\begin{document}

\title{A stochastic maximum principle for partially observed general mean-field control problems with only weak solution\footnotemark[1]}
\author{Juan Li$^{1,**}$,\,\,Hao Liang$^{1}$,\,\, Chao Mi$^{1,**}$ \\
 {\small School of Mathematics and Statistics, Shandong University, Weihai, Weihai 264209, P.~R.~China.$^1$}\\
 {\small{\it E-mails: juanli@sdu.edu.cn,\,\ haoliang@mail.sdu.edu.cn,\,\ michao94@mail.sdu.edu.cn.}}
 \date{November 21, 2021.}
	}
\renewcommand{\thefootnote}{\fnsymbol{footnote}}
\footnotetext[1]{The work is supported by the NSF of P.R. China (NOs. 12031009, 11871037), National Key R and D Program of China (NO. 2018YFA0703900), and NSFC-RS (No. 11661130148; NA150344).\\
\mbox{ }\ \ \ \ $^{**}$Corresponding authors. }
\maketitle

\textbf{Abstract}. In this paper we focus on a general type of mean-field stochastic control problem with partial observation, in which the coefficients depend in a non-linear way not only on the state process $X_t$ and its control $u_t$ but also on the conditional law $E[X_t|\mathcal{F}_t^Y]$ of the state process conditioned with respect to the past of observation process $Y$. We first deduce the well-posedness of the controlled system by showing weak existence and uniqueness in law. Neither supposing convexity of the control state space nor differentiability of the coefficients with respect to the control variable,
we study Peng's stochastic maximum principle for our control problem. The novelty and the difficulty of our work stem from the fact that, given an admissible control $u$, the solution of the associated control problem is only a weak one. This has as consequence that also the probability measure in the solution $P^{u}=L^{u}_TQ$ depends on $u$ and has a density $L^{u}_T$ with respect to a reference measure $Q$. So characterizing an optimal control leads to the differentiation of non-linear functions $f(P^{u}\circ\{E^{P^{u}}[X_t|\mathcal{F}_t^Y]\}^{-1})$ with respect to $(L^{u}_T,X_t)$. This has as consequence for the study of Peng's maximum principle that we get a new type of first and second order variational equations and adjoint backward stochastic differential equations, all with new mean-field terms and with coefficients which are not Lipschitz. For their estimates and for those for the Taylor expansion new techniques have had to be introduced and rather technical results have had to be established. The necessary optimality condition we get extends Peng's one with new, non-trivial terms.


\section{Introduction}
Let $(\Omega,\mathcal{F},\mathbb{F}=\{\mathcal{F}_t\}_{t\geq0},P)$ be a filtered probability space endowed with an $(\mathbb{F},P)$ Brownian motion $B=(B^1,B^2)$.
In this paper, we focus on the following partially observed stochastic control system of mean-field type:
\begin{equation}\label{Intro}
\bigg\{
\begin{tabular}{l}
$dX_t=\sigma(t,Y_{\cdot\wedge t},X_t,\mu_t^{X|Y},u_t)dB_t^1,\ X_0=x_0\in\mathbb{R};$ \\
$dY_t=h(t,Y_{\cdot\wedge t},X_t,\mu_t^{X|Y},u_t)dt+dB_t^2,\ Y_0=0,\ t\in[0,T].$
\end{tabular}
\end{equation}
Here $\mu_t^{X|Y}=P_{E[X_t\,|\,\mathcal{F}_t^Y]}=P\circ\{E[X_t\,|\,\mathcal{F}_t^Y]\}^{-1}$ is the law of $E[X_t\,|\,\mathcal{F}_t^Y]$ under $P$, and $\mathbb{F}^Y=\{\mathcal{F}_t^Y\}_{t\geq 0}$ stands for the filtration generated by $Y$, and $Y_{\cdot\wedge t}$ is the process $Y$ stopped at $t$.

Special cases of the above controlled system have been studied by different authors, mainly in applications concerning partial information problems: The general interpretation here is that the controlled state process $X$ can only be observed through the process $Y$, which expresses, in particular, through the conditional expectation  $E[X_t\,|\,\mathcal{F}_t^Y]$. Relating the rather classical control problem with partial observation with a mean-field problem, in which the coefficients can depend on the law of $E[X_t\,|\,\mathcal{F}_t^Y]$, Buckdahn, Li and Ma  \cite{BLM17} studied Pontryagin's stochastic maximum principle (SMP) for a special case of the above control problem, where, firstly, $h$ does not depend on the control nor on $\mu^{X|Y}$ and, secondly, $\sigma$ depends on $\mu^{X|Y}$ only linearly. The by far non-trivial extensions we study here are multiple: The coefficients can depend in a general manner, i.e., non-linearly, on the law  $\mu^{X|Y}$, and as no convexity assumption on the control state space nor any differentiability assumption with respect to (w.r.t.) the control variable is made, we study here Peng's SMP. The dynamic programming principle as another method to solve stochastic optimization problems cannot be used here, as the stochastic control system involves non-linearly the (conditional) law of the state process which makes that the control is highly time inconsistent.

Since Kushner's pioneering works \cite{K65,K72} on the SMP, ground-breaking contributions were made namely by Bensoussan \cite{B81}, Bismut \cite{B78}, Haussmann \cite{H86}. The work of Peng \cite{P90} on the SMP was the first to overcome the restrictive assumption of convexity of control space and the regularity of the coefficients w.r.t. the control variable.

On the other hand, as explained above, SDE (\ref{Intro}) unites the problem of control under partial observation with that of optimal control in the mean-field case. With their seminal paper \cite{LL07} Lasry and Lions  were  the first to study mean-field games. Their work has been the starting point for a lot of different works related with by different authors; by citing Buckdahn, Djehiche, Li and Peng \cite{BDLP09} and \cite{BLP09}, Buckdahn, Djehiche, Li \cite{BDL11}, Carmona, Delarue, Lachapelle \cite{CD13,CDL13}, we restrict here to works related with our work.
Furthermore, also in the framework of mean-field dynamics, Carmona and Zhu \cite{CZ16} studied a kind of mean-field SDE, called conditional mean-field SDE, in which the coefficients depend not only on the path of the state but also on its conditional expectation w.r.t. a given filtration. Finally, Buckdahn, Li, Ma \cite{BLM17}--the work the most related with the present one- has been already mentioned above.

Let us come back to our problem. The first point we have to study is the well-posedness of (\ref{Intro}). The fact that (\ref{Intro})--with the sought solution process $(X,Y)$--involves (the law of) the conditional expectation $E[X_t\,|\,\mathcal{F}_t^Y]$ makes that we can not expect the strong existence of a solution, but only the weak one. It will be constructed with the help of a Girsanov transformation argument, which allows to reduce the study of weak existence for (\ref{Intro}) to that of a strong solution of the controlled SDE
\begin{equation}\label{trans}
\Bigg\{
\begin{tabular}{l}
$dX_t^u=\sigma(t,X_t^u,\mu_t^u,u_t)dB_t^1,\ X_0^u=x;$ \\
$dL_t^u=L_t^u h(t,X_t^u,\mu_t^u,u_t)dY_t,\ L_0^u=1,\ t\in[0,T],$  \\
\end{tabular}
\end{equation}
where $(B^1,Y)$ is a Brownian motion under a reference probability measure $Q$. As we show, the weak solution which turns out to be unique in law, is then obtained by putting $P^u=L_T^u Q$, and $\mu_t^u=\mu_t^{X^u|Y}=P^u_{E^{P^{u}}[X_t^u\,|\,\mathcal{F}_t^Y]}$.

Our main objective consists in minimizing the cost functional
$$ J(u)=E^Q\Big[\Phi(X_T^u,\mu_T^u)+\int_0^T f(t,X_t^u,\mu_t^u,u_t)dt\Big], $$
for $f$ and $\Phi$ two given functions satisfying appropriate conditions. In addition to the already above listed new, main difficulties which we have to deal in our work--non-linear control state space, non-differentiability of the coefficients w.r.t. the control variable, non-linear dependence of the coefficients of the law of the conditional expectation, we meet as main problem that the underlying probability measure $P^u=L_T^u Q$ in the weak solution depends itself through the density $L_T^{u}$  on the control $u$. This typical property of a weak solution has as consequence that we have to consider the derivative of functionals of the form
$F_Q(L):=f\big((LQ)_{\xi}\big),\ L\in\mathcal{L}^Q$, where $\mathcal{L}^Q$ denotes the collection of densities under a given probability $Q$; $\xi$ is an arbitrary given random variable; $f:\mathcal{P}(\mathbb{R}^d)\rightarrow\mathbb{R}$ is a function differentiable in P.L. Lion's sense (\cite{Lions}, or see Cardaliaguet \cite{C12}). Such derivatives have been recently studied in \cite{Dwrtd}. So considering, for simplicity, the coefficient of the running cost $f(t,x,\mu,u)=f(\mu)$ only as depending (non-linearly) on the probability law, we will see that we have to take the derivative w.r.t. $(L^{u}_t,X^{u}_t)$ of
$$\displaystyle f(\mu_t^{u})=f\big(P^{u}_{E^{u}[X_t^{u}|\mathcal{F}_t^Y]}\big)
=f\bigg((L^{u}_tQ)\circ\big\{E^Q\big[L^{u}_tX_t^{u}\frac{1}{E^Q[L^{u}_t|\mathcal{F}^Y_t]}\big|\mathcal{F}_t^Y\big]\big\}^{-1}\bigg),$$
and later even second order derivatives. This leads not only to a new type of first order and second order coupled variational equations but also to new terms in the adjoint coupled backward stochastic differential equations of first and second order. All these equations are of mean-field type. It is also worth to point out that these equations have coefficients which do not have a Lipschitz constant independent of $\omega\in\Omega$. Additional difficulties come from the derivative of the coefficients depending on $\mu_t^{u}$ as already explained above: The terms of the derivative also appear in the variational equations and in the adjoint equations. For the estimate of the error of the Taylor expansion of $(X^\varepsilon,L^\varepsilon)$ given as the sum of
$(X,L)$ and the solution couples of the first and the second order adjoint equations, in particular the technical result given by  Proposition \ref{tech} turn out to be crucial. Here $(X^\varepsilon,L^\varepsilon)$ denotes the solution couple  of  (\ref{trans})  governed by a spike variation $u^\varepsilon$ of an optimal control $u$, and $(X,L)$ denotes the solution of  (\ref{trans}) for the control $u$. Just as crucial and highly technical is the also new, associated Taylor expansion for the conditional expectation $E^{u^\varepsilon}[X_t^{\varepsilon}|\mathcal{F}_t^Y].$

Finally, comparing with earlier works on Peng's SMP, the necessary optimality condition we get contains additional terms, which have their origin in the mean-field terms of the control problem, and, in particular, in the derivative w.r.t. the density, which comes from the fact that we have to do with weak solutions of our control problem. To highlight these new terms we discuss briefly in the appendix the case, where we have in the control problem the law of $\varphi(X_t,Y_{\cdot\wedge t})$, i.e., $P\circ\{\varphi(X_t,Y_{\cdot\wedge t}) \}^{-1}$, instead of the law $P\circ\{E^P[X_t|\mathcal{F}_t^Y]\}^{-1}$,  that is, we consider the general full observation, but still more general than classical cases, to have a comparison with Peng's SMP we get for $U_t=E^P[X_t\,|\,\mathcal{F}_t^Y]$. Although easier, also this is a novel problem, and also here the variational equations and adjoint backward equations of first and second order have new terms, but the necessary optimality condition we get here resembles at first glance to the classical one, obtained by Peng.

This paper is organized as follows. In Section 2, we recall some notations and the definition of derivative with respect to a probability measure. Section 3 proves the well-posedness (weak existence and uniqueness in law) of the state-observation dynamics we focus on and allows to deal with the system given in Section 4. We introduce our stochastic optimization problems in Subsection 4.1; Subsections 4.2 and 4.3 study the variational equations and adjoint equations, respectively. We make the computations and conclude the stochastic maximum principle in Subsection 4.4. Last but not least, the Appendix is devoted to the discussion of the case, where in the control problem the law $P\circ\{E^P[X_t|\mathcal{F}_t^Y]\}^{-1}$ is replaced by that of $\varphi(X_t,Y_{\cdot\wedge t})$.

\section{Preliminaries}
Let $(E,d)$ be a complete metric space, for which $\mathcal{P}(E)$ denotes the space of all probability measures. For $p\geq 1$, we consider
 the space of probability measures on $(E,d,\mathcal{B}(E))$ with finite $p$-th moment, denoted by $\mathcal{P}_p(E)$, and we endow it with
  the $p$-Wasserstein metric
\begin{equation*}
W_{p}(\mu,\nu):=\inf\Big\{\left(\int_{E\times E}\big(d(z,z')\big)^p\rho(dzdz')\right)^{\frac{1}{p}}\Big|\ \rho\in\mathcal{P}_p(E\times E)
\mbox{ with }
                      \rho(\cdot\times E)=\mu,\ \rho(E\times\cdot)=\nu\Big\},
\end{equation*}
where $\mu,\nu\in\mathcal{P}_p(E)$. Note that also $\big(\mathcal{P}_p(E),W_{p}(\cdot,\cdot)\big)$ is a complete metric space.
%In this paper, we consider $p=2$. It's easy to see that $(\mathcal{P}_2(C_T),W_2)$ is a separatable complete metric space.

Let now $E=\mathbb{R}^k,\, k\geq1$. For a complete probability space $(\Omega,\mathcal{F},P)$ and a filtration $\mathbb{F}$ satisfying the usual assumptions, we first introduce the following spaces: For any sub-$\sigma$-field $\mathcal{G}$ of $\mathcal{F}$ and any subfiltration $\mathbb{G}$ of $\mathbb{F}$, $p\geq1$,
\begin{itemize}
  \item $L^p(\mathcal{G},P;\mathbb{R}^k)$ is the set of $\mathbb{R}^k$-valued, $\mathcal{G}$-measurable random variables $\xi$ over $\Omega$ with $E^P[|\xi|^p]<\infty$. Here $E^P[\,\cdot\,]$ denotes the expectation with respect to $P$.
  \item $S^p_\mathbb{G}([0,T],P;\mathbb{R}^k)$ denotes the set of $\mathbb{R}^k$-valued, $\mathbb{G}$-adapted continuous stochastic processes $X$ on $[0,T]\times\Omega$, with
      $\displaystyle E^P\big[\sup_{t\in[0,T]}|X_t|^p\big]<\infty. $
  \item $L^p_\mathbb{G}([0,T],P;\mathbb{R}^k)$ denotes the set of $\mathbb{R}^k$-valued, $\mathbb{G}$-progressively measurable stochastic processes $X$ on $[0,T]\times\Omega$, with
      $\displaystyle E^P\Big[\Big(\int_0^T|X_t|^2dt\Big)^{\frac{p}{2}}\Big]<\infty. $
\end{itemize}

For $E=\mathbb{R}$ we omit $\mathbb{R}$ in the notations of the above spaces.
We also assume that $(\Omega,\mathcal{F},P)$ is ``rich enough'' in the sense that $\mathcal{P}_2(\mathbb{R}^k)=\{P_\vartheta, \vartheta\in L^2(\mathcal{F},P;\mathbb{R}^k)\}$, $k\geq1$. The space   $L^{2}(\mathcal{F},P)$ can be regarded as a Hilbert space with inner product $(\xi,\eta)_{L^{2}}=E^P[\xi\cdot\eta],\ \xi,\eta\in L^{2}(\mathcal{F},P)$, and norm $|\xi|_{L^{2}}=(\xi,\xi)_{L^{2}}^{\frac{1}{2}}$. For a given function $f:\mathcal{P}_{2}(\mathbb{R})\rightarrow\mathbb{R}$, we define its ``lifted'' function $\widetilde{f}(\vartheta):=f(P_{\vartheta}),\ \vartheta\in L^{2}(\mathcal{F},P)$. The function $f:\mathcal{P}_{2}(\mathbb{R})\rightarrow\mathbb{R}$ is said to be differentiable at $\mu\in\mathcal{P}_{2}(\mathbb{R})$, if there exists $\vartheta_{0}\in L^{2}(\mathcal{F},P)$ with $P_{\vartheta_{0}}=\mu$, such that the function $\widetilde{f}:L^{2}(\mathcal{F},P)\rightarrow\mathbb{R}$ is Fr\'{e}chet differentiable at $\vartheta_{0}$, i.e., there exists a continuous linear mapping $D\widetilde{f}(\vartheta_{0}):L^{2}(\mathcal{F},P)\rightarrow \mathbb{R}$ (i.e., $D\widetilde{f}(\vartheta_{0})\in L(L^{2}(\mathcal{F},P);\mathbb{R})$) such that
\begin{equation}
\widetilde{f}(\vartheta_{0}+\eta)-\widetilde{f}(\vartheta_{0})=D\widetilde{f}(\vartheta_{0})(\eta)+o(|\eta|_{L^{2}}),
\end{equation}
with $|\eta|_{L^{2}}\rightarrow0$, for $\eta\in L^{2}(\mathcal{F},P)$.

Noting that $D\widetilde{f}(\vartheta_{0})\in L(L^{2}(\mathcal{F},P);\mathbb{R})$,
we obtain from Riesz' Representation Theorem the existence of a ($P$-a.s.) unique random variable $\theta_{0}\in L^{2}(\mathcal{F},P)$ such that $D\widetilde{f}(\vartheta_{0})(\eta)=(\theta_{0},\eta)_{L^{2}}=E^P[\theta_{0}\eta]$,
for all $\eta\in L^{2}(\mathcal{F},P)$. P.L. Lions (\cite{Lions}, or see the notes of Cardaliaguet \cite{C12}) proved the existence of a Borel function $h_{0}:\mathbb{R}\rightarrow\mathbb{R}$ such that $\theta_{0}=h_{0}(\vartheta_{0})$, $P$-a.s. Moreover, the function $h_0$ was shown to depend on $\vartheta_0$ only through its law $P_{\vartheta_0}$, but not on the particular choice of $\vartheta_0$. Then we can write
\begin{equation}\label{equ2.4}
f(P_{\vartheta})-f(P_{\vartheta_{0}})=E^P[h_{0}(\vartheta_{0})\cdot(\vartheta-\vartheta_{0})]+o(|\vartheta-\vartheta_{0}|_{L^{2}}),
\end{equation}
with $|\vartheta-\vartheta_{0}|_{L^{2}}\rightarrow0$, for $\vartheta\in L^{2}(\mathcal{F},P)$.

We define $\partial_{\mu}f(P_{\vartheta_{0}},y):=h_{0}(y),\,y\in\mathbb{R}$, and call it the derivative of $f:\mathcal{P}_{2}(\mathbb{R})\rightarrow\mathbb{R}$ at $P_{\vartheta_{0}}$. We will consider the case that $f:\mathcal{P}_{2}(\mathbb{R})\rightarrow\mathbb{R}$ is differentiable over the whole space, and we have the derivative $\partial_{\mu}f(P_{\vartheta},y)$ defined $P_{\vartheta}(dy)$-a.e., for all $\vartheta\in L^{2}(\mathcal{F},P)$,
through the relation $D\widetilde{f}(\vartheta)(\eta)=E[\partial_{\mu}f(P_{\vartheta},\vartheta)\eta]$, for all $\eta\in L^{2}(\mathcal{F},P)$.

\section{Well-posedness of the state-observation dynamics}

Through this paper, we consider the canonical space $(\Omega,\mathcal{F}):=(C_T^2,\mathcal{B}(C_T^2))$, where $C_T^2=C([0,T];\mathbb{R}^2)$.  Let $\mathbb{F}$ be the natural filtration generated by the coordinate process on $\Omega$. Let $P\in\mathcal{P}(\Omega)$ be a probability over $(\Omega,\mathcal{F})$. The filtered probability space $(\Omega,\mathcal{F},\mathbb{F},P)$ is supposed to satisfy the usual condition.

Now we give the dynamics of the state and the observation process which we focus on in this paper:
\begin{equation}\label{A4}
\bigg\{
\begin{tabular}{l}
$dX_t=\sigma(t,Y_{\cdot\wedge t},X_t,\mu_t^{X|Y})dB_t^1,\ X_0=x_0\in\mathbb{R};$ \\
$dY_t=h(t,Y_{\cdot\wedge t},X_t,\mu_t^{X|Y})dt+dB_t^2,\ Y_0=0,\ t\in[0,T],$
\end{tabular}
\end{equation}
where $(B^1,B^2)$ is an $(\mathbb{F},P)$-Brownian motion. In the above system, $X$ is the state process, while $Y$ is the observation process, defined on $(\Omega,\mathcal{F},P)$. Let $U_t^{X|Y}:=E^{P}[X_t\,|\,\mathcal{F}_t^Y],\ t\in[0,T]$, denote the ``filtered'' state process and $\mu_t^{X|Y}$ its law under $P$, i.e., $\mu_t^{X|Y}:=P_{U_t^{X|Y}}$.
We will consider the well-posedness of \eqref{A4} under the following Assumption (H1).\\% and let $U_t^{X|Y}=E^{P}[X_t\,|\,\mathcal{F}_t^Y],\ t\in[0,T]$ denote the ``filtered'' state process and $\mu_t^{X|Y}$ its law under $P$. \\
\textbf{Assumption (H1)}:\\
(i) The functions $\sigma, h:[0,T]\times C_T\times \mathbb{R}\times\mathcal{P}_2(\mathbb{R})\rightarrow\mathbb{R}$ are Borel measurable and bounded;\\
(ii) For all $(t,y)\in[0,T]\times C_T$, $x,x'\in\mathbb{R}$, $\gamma,\gamma'\in\mathcal{P}_2(\mathbb{R}):$
$$|\phi(t,y_{\cdot\wedge t},x,\gamma)-\phi(t,y_{\cdot\wedge t},x',\gamma')|\leq C\big(|x-x'|+W_1(\gamma,\gamma')\big),$$
for $\phi=\sigma,h$.
\begin{remark}
In \eqref{A4} we have assumed that the drift coefficient $b=0$. Indeed, the extension of our discussion to the case with a drift does not add additional difficulties.
\end{remark}
We first consider the well-posedness of the following type of SDE. Let $(B^1,Y)$ be the coordinate process on $\Omega=C_T^2$, $(B_t^1(\omega),Y_t(\omega))=(\omega_1(t),\omega_2(t))$, $\omega=(\omega_1,\omega_2)\in\Omega$, $t\in[0,T]$ and let $Q$ be the Wiener measure over $(\Omega,\mathcal{F})=(C_T^2,\mathcal{B}(C_T^2))$. By $\mathbb{F}=\mathbb{F}^{B^1,Y}$ we denote the filtration generated by $(B^1,Y)$ and augmented by all $Q$-null sets; also $\mathcal{F}$ is considered to be completed w.r.t. $Q$. In particular, $(B^1,Y)$ is an $(\mathbb{F},Q)$-Brownian motion. For $t\in[0,T]$,
\begin{equation}\label{A1}
\bigg\{
\begin{tabular}{l}
$dX_t=\sigma(t,Y_{\cdot\wedge t},X_t,\mu_t^{X|Y})dB_t^1,\ X_0=x_0\in\mathbb{R};$ \\
$dL_t=h(t,Y_{\cdot\wedge t},X_t,\mu_t^{X|Y})L_tdY_t,\ L_0=1.$
\end{tabular}
\end{equation}
Notice that also $P=L_T Q$ is a probability.
\begin{proposition}\label{existofstr}
Under Assumption (H1) equation \eqref{A1} possesses a unique strong solution.
\end{proposition}
\begin{proof}
Given any $V\in S_{\mathbb{F}}^2([0,T],Q)$, and $K\in\mathcal{K}_{\mathbb{F}}^2([0,T],Q):=\big\{K\in S_{\mathbb{F}}^2([0,T],Q)\,\big|\, K_T\geq0,\,E^Q[K_T]\\=1,\, K_t=E^Q[K_T\,|\, \mathcal{F}_t],\, t\in[0,T]\big\}$, and putting $P:=K_T Q$, and $\mu_t:=P_{E^{P}[V_t\,|\,\mathcal{F}_t^Y]}$, $t\in[0,T]$, we consider the following SDE:
\begin{equation}\label{A2}
\bigg\{
\begin{tabular}{l}
$d\overline{X}_t=\sigma(t,Y_{\cdot\wedge t},\overline{X}_t,\mu_t)dB_t^1,\ \overline{X}_0=x_0\in\mathbb{R};$ \\
$d\overline{L}_t=h(t,Y_{\cdot\wedge t},\overline{X}_t,\mu_t)\overline{L}_tdY_t,\ \overline{L}_0=1.$
\end{tabular}
\end{equation}
We deduce from the classical theory of SDEs that \eqref{A2} admits a unique solution $(\overline{X},\overline{L})\in S_{\mathbb{F}}^2([0,T],Q)\\
\times \mathcal{K}_{\mathbb{F}}^2([0,T],Q)$. This allows to define by putting $\Phi(V,K):=(\overline{X},\overline{L})$ a mapping from $S_{\mathbb{F}}^2([0,T],Q)\times \mathcal{K}_{\mathbb{F}}^2([0,T],Q)$ to itself.
As $\sigma,h$ are bounded, we have that for all $p\geq 1$, there exists a constant $C_p(\sigma,h)\in\mathbb{R}_+$ (only depending on $p$, and the bounds of $\sigma$ and $h$), such that for all $(V,K)\in S_{\mathbb{F}}^2([0,T],Q)\times \mathcal{K}_{\mathbb{F}}^2([0,T],Q)$, $\displaystyle E^Q\Big[\sup_{t\in[0,T]}\Big(|\overline{X}_t|^p+|\overline{L}_t|^p+|\frac{1}{\overline{L}_t}|^p\Big)\Big]\leq C_p(\sigma,h). $

We define $\mathcal{H}:=\big\{(V,K)\in S_{\mathbb{F}}^2([0,T],Q)\times \mathcal{K}_{\mathbb{F}}^2([0,T],Q):\, E^Q[\sup_{t\in[0,T]}(|V_t|^p+|K_t|^p+|\frac{1}{K_t}|^p)]\leq C_p(\sigma,h),\, p>1\big\}$. Let $(V^i,K^i)\in\mathcal{H},\, i=1,2$, and $(\overline{X}^i,\overline{L}^i):=\Phi(V^i,K^i),\, i=1,2$. We set $(\widehat{\overline{X}},\widehat{\overline{L}}):=(\overline{X}^1-\overline{X}^2,\overline{L}^1-\overline{L}^2)$ and $(\widehat{V},\widehat{K}):=(V^1-V^2,K^1-K^2)$. For $i=1,2$, let $\mu_t^i:=P^i_{E^{P^i}[V^i_t\,|\,\mathcal{F}_t^Y]}$, $t\in[0,T]$, for $P^i=K_T^i Q$. Now we first estimate $W_1(\mu_t^1,\mu_t^2)$.

Let $\phi\in \mbox{Lip}_1(\mathbb{R})$ (the space of real Lipschitz continuous functions defined over $\mathbb{R}$ with Lipschitz constant $1$) with $\phi(0)=0$. To use the Kantorovich-Rubinstein duality, we note that
\begin{equation*}
\begin{split}
&\Big|\int_\mathbb{R}\phi d\mu_t^1-\int_\mathbb{R}\phi d\mu_t^2\Big|=\Big| E^{P^1}\big[\phi\big(E^{P^1}[V_t^1\,|\,\mathcal{F}_t^Y]\big)\big]-E^{P^2}\big[\phi\big(E^{P^2}[V_t^2\,|\,\mathcal{F}_t^Y]\big)\big]\Big|\\
&=\Big| E^{Q}\big[K_t^1\phi\big(E^{P^1}[V_t^1\,|\,\mathcal{F}_t^Y]\big)-K_t^2\phi\big(E^{P^2}[V_t^2\,|\,\mathcal{F}_t^Y]\big)\big]\Big|\\
&\leq E^Q\big[\big|E^{P^1}[V_t^1\,|\,\mathcal{F}_t^Y]\big|\cdot|K_t^1-K_t^2|\big]+E^Q\big[E^{Q}[K_t^2\,|\,\mathcal{F}_t^Y]\cdot\big|E^{P^1}[V_t^1\,|\,\mathcal{F}_t^Y]-E^{P^2}[V_t^2\,|
\,\mathcal{F}_t^Y]\big|\big]\\
&=:I_1+I_2.
\end{split}
\end{equation*}
We observe that, since $(V^i,K^i)\in \mathcal{H}$, and $\frac{1}{E^{Q}[K_t^i\,|\,\mathcal{F}_t^Y]}\leq E^{Q}[\frac{1}{K_t^i}\,|\,\mathcal{F}_t^Y],\, i=1,2$,
\begin{equation*}
\begin{split}
I_1&=E^Q\Big[\frac{E^{Q}[K_t^1 V_t^1\,|\,\mathcal{F}_t^Y]}{E^{Q}[K_t^1\,|\,\mathcal{F}_t^Y]}|K_t^1-K_t^2|\Big]\\
&\leq \Big(E^Q\Big[\big(E^{Q}[K_t^1 V_t^1\,|\,\mathcal{F}_t^Y]E^{Q}[\frac{1}{K_t^1}\,|\,\mathcal{F}_t^Y]\big)^2\Big]\Big)^{\frac{1}{2}}\big(E^{Q}[|K_t^1-K_t^2|^2]\big)^{\frac{1}{2}}\\
&\leq C(\sigma,h)\big(E^{Q}[|K_t^1-K_t^2|^2]\big)^{\frac{1}{2}},
\end{split}
\end{equation*}
for $C(\sigma,h)=C_8(\sigma,h)$,
and
\begin{equation*}
\begin{split}
I_2&=E^Q\Big[E^{Q}[K_t^2\,|\,\mathcal{F}_t^Y]\cdot\Big|\frac{E^{Q}[K_t^1 V_t^1\,|\,\mathcal{F}_t^Y]}{E^{Q}[K_t^1\,|\,\mathcal{F}_t^Y]}-\frac{E^{Q}[K_t^2 V_t^2\,|\,\mathcal{F}_t^Y]}{E^{Q}[K_t^2\,|\,\mathcal{F}_t^Y]}\Big|\Big]\\
&\leq E^Q\big[\big|E^Q[K_t^1 V_t^1-K_t^2 V_t^2\,|\,\mathcal{F}_t^Y]\big|\big]+E^Q\Big[\big|E^{Q}[K_t^1 V_t^1\,|\,\mathcal{F}_t^Y]\big|\cdot\Big|1-\frac{E^Q[K_t^2\,|\,\mathcal{F}_t^Y]}{E^Q[K_t^1\,|\,\mathcal{F}_t^Y]}\Big|\Big]\\
&\leq \big(E^Q[|K_t^1|^2]\big)^{\frac{1}{2}}\big(E^Q[|V_t^1-V_t^2|^2]\big)^{\frac{1}{2}}+\big(E^Q[|V_t^2|^2]\big)^{\frac{1}{2}}\big(E^Q[|K_t^1-K_t^2|^2]\big)^{\frac{1}{2}}\\
&\qquad +\Big(E^Q\Big[\Big| E^{Q}[K_t^1 V_t^1\,|\,\mathcal{F}_t^Y]\,E^{Q}[\frac{1}{K_t^1}\,|\,\mathcal{F}_t^Y]\Big|^2\Big]\Big)^{\frac{1}{2}}\big(E^Q[|K_t^1-K_t^2|^2]\big)^{\frac{1}{2}}\\
&\leq C(\sigma,h)\Big(\big(E^{Q}[|V_t^1-V_t^2|^2]\big)^{\frac{1}{2}}+\big(E^{Q}[|K_t^1-K_t^2|^2]\big)^{\frac{1}{2}}\Big).
\end{split}
\end{equation*}
Consequently, for all $(V^1,K^1),\,(V^2,K^2)\in\mathcal{H}$, thanks to the Kantorovich-Rubinstein duality,
\begin{equation}\label{difmu}
\begin{split}
W_1(\mu_t^1,\mu_t^2)&=\sup\Big\{\Big|\int_\mathbb{R}\phi d\mu_t^1-\int_\mathbb{R}\phi d\mu_t^2\Big|\,:\,\phi\in\mbox{Lip}_1(\mathbb{R}),\,\phi(0)=0\Big\}\\
&\leq C(\sigma,h)\Big(\big(E^{Q}[|V_t^1-V_t^2|^2]\big)^{\frac{1}{2}}+\big(E^{Q}[|K_t^1-K_t^2|^2]\big)^{\frac{1}{2}}\Big).
\end{split}
\end{equation}
Now let us estimate $E^Q[|\overline{L}_t^1-\overline{L}_t^2|^2]$. From the assumptions on $h$, we deduce
\begin{equation*}
\begin{split}
&E^Q[\sup_{s\leq t}|\overline{L}_s^1-\overline{L}_s^2|^2]\leq C\int_0^t E^Q\big[\big|h(s,Y_{\cdot\wedge s},\overline{X}^1_s,\mu^1_s)\overline{L}^1_s-h(s,Y_{\cdot\wedge s},\overline{X}^2_s,\mu^2_s)\overline{L}^2_s\big|^2\big]ds\\
&\leq C\int_0^t E^Q[|\overline{L}_s^1-\overline{L}_s^2|^2]ds+C\int_0^t E^Q\big[(\overline{L}^2_s)^2\big|h(s,Y_{\cdot\wedge s},\overline{X}^1_s,\mu^1_s)-h(s,Y_{\cdot\wedge s},\overline{X}^2_s,\mu^2_s)\big|^2\big]ds\\
&\leq C\int_0^t E^Q[|\overline{L}_s^1-\overline{L}_s^2|^2]ds+C\int_0^t \Big(E^Q\big[\big|h(s,Y_{\cdot\wedge s},\overline{X}^1_s,\mu^1_s)-h(s,Y_{\cdot\wedge s},\overline{X}^2_s,\mu^2_s)\big|^4\big]\Big)^{\frac{1}{2}}ds\\
&\leq C\int_0^t\Big( E^Q[|\overline{L}_s^1-\overline{L}_s^2|^2]+ \big(E^Q[|\overline{X}_s^1-\overline{X}_s^2|^4]\big)^{\frac{1}{2}}+W_1^2(\mu_s^1,\mu_s^2)\Big)ds,\ t\in[0,T].
\end{split}
\end{equation*}
Hence, due to Gronwall's inequality, it holds that
\begin{equation}\label{L}
E^Q[\sup_{s\leq t}|\overline{L}_s^1-\overline{L}_s^2|^2]\leq C\int_0^t\Big( \big(E^Q[|\overline{X}_s^1-\overline{X}_s^2|^4]\big)^{\frac{1}{2}}+W_1^2(\mu_s^1,\mu_s^2)\Big)ds,\ t\in[0,T].
\end{equation}
Next we come to the estimate of $E^Q[\sup\limits_{s\leq t}|\overline{X}_s^1-\overline{X}_s^2|^4]$. Thanks to the Lipschitz condition (H1)-(ii), we have
\begin{equation*}
\begin{split}
E^Q[\sup_{s\leq t}|\overline{X}_s^1-\overline{X}_s^2|^4]&\leq C E^Q\Big[\Big(\int_0^t \big|\sigma(s,Y_{\cdot\wedge s},\overline{X}^1_s,\mu^1_s)-\sigma(s,Y_{\cdot\wedge s},\overline{X}^2_s,\mu^2_s)\big|^2ds\Big)^2\Big]\\
&\leq C\int_0^t\Big(E^Q[|\overline{X}_s^1-\overline{X}_s^2|^4]+W_1^4(\mu_s^1,\mu_s^2)\Big)ds,\ t\in[0,T],
\end{split}
\end{equation*}
and again from Gronwall's inequality,
\begin{equation}\label{X}
E^Q[\sup_{s\leq t}|\overline{X}_s^1-\overline{X}_s^2|^4]\leq C\int_0^t W_1^4(\mu_s^1,\mu_s^2)ds,\ t\in[0,T].
\end{equation}
Consequently, the above estimates yield that
\begin{equation}\label{A3}
\begin{split}
E^Q[\sup_{s\leq t}|\overline{X}_s^1-\overline{X}_s^2|^4]+\big(E^Q[\sup_{s\leq t}|\overline{L}_s^1-\overline{L}_s^2|^2]\big)^2\leq C \int_0^t W_1^4(\mu_s^1,\mu_s^2)ds\\
\leq C \int_0^t\Big( E^Q[|V_s^1-V_s^2|^4]+\big(E^Q[|K_s^1-K_s^2|^2]\big)^2\Big)ds,\ t\in[0,T],
\end{split}
\end{equation}
for all $(V^i,K^i)\in\mathcal{H},\,i=1,2$. As $\Phi(V^i,K^i)=(\overline{X}^i,\overline{L}^i)$, \eqref{A3} combined with standard arguments, in particular the contraction mapping arguments, proves the existence and uniqueness of a strong solution of \eqref{A1}, given the $(\mathbb{F},Q)$-Brownian motion $(B^1,Y)$.
\end{proof}

The existence of a strong solution $(X,L)$ of SDE \eqref{A1} implies, in particular, that of a weak solution of \eqref{A4}. We recall the definition of a weak solution here.
\begin{definition}\label{WeSol}
A six-tuple $(\Omega,\mathcal{F},\mathbb{F},P,(B^1,B^2),(X,Y))$ is called a weak solution of \eqref{A4}, if\\
\mbox{ }\ \ {\rm i)}\ $(\Omega,\mathcal{F},\mathbb{F},P)$ is a filtered probability space satisfying the usual assumptions;\\
\mbox{ }\ \ {\rm ii)}\ $(B^1,B^2)$ is an $(\mathbb{F},P)$-Brownian motion;\\
\mbox{ }\ \ {\rm iii)}\ All terms in \eqref{A4} are well-defined, $(X,Y)$ is an $\mathbb{F}$-adapted process and equation \eqref{A4} holds\\
\mbox{ }\ \ \ \ \ \ \ true, for all $t\in[0,T]$,\ $P$-a.s.
\end{definition}
From the Girsanov theorem, we know that, given a strong solution $(X,L)$ of \eqref{A1} with driving $(\mathbb{F},Q)$-Brownian motion $(B^1,Y)$, $(\Omega,\mathcal{F},\mathbb{F},P,(B^1,B^2),(X,Y))$ is a weak solution of \eqref{A4}, where $P=L_T Q$ and
$$ B_t^2=Y_t-\int_0^t h(s,Y_{\cdot\wedge s},X_s,\mu_s^{X|Y})ds,\ t\in[0,T]. $$
As a conclusion, under Assumption (H1), the dynamics \eqref{A4} admits at least one solution in the sense of Definition \ref{WeSol}.

In addition to the existence of the weak solution of \eqref{A4} we also have uniqueness in law. Before showing this, let us make the following remark.
\begin{remark}\label{ContofU}
Note that  $\displaystyle U_t^{X|Y}:=E^{P}[X_t\,|\,\mathcal{F}_t^Y]=\frac{E^Q[L_t X_t\,|\,\mathcal{F}_t^Y]}{E^Q[L_t\,|\,\mathcal{F}_t^Y]},\ Q\mbox{-a.s.}, \ t\in[0,T]. $\
Furthermore, as $L_t$ and $X_t$ are both $\mathcal{F}_t^{B^1,Y}$-measurable and thus independent of $\sigma\{Y_s-Y_t, s\in[t,T]\}$, we also have
$$ U_t^{X|Y}=\frac{E^Q[L_t X_t\,|\,\mathcal{F}_t^Y]}{E^Q[L_t\,|\,\mathcal{F}_t^Y]}=\frac{E^Q[L_t X_t\,|\,\mathcal{F}_T^Y]}{E^Q[L_t\,|\,\mathcal{F}_T^Y]},\ Q\mbox{-a.s.}, \ t\in[0,T]. $$
From \eqref{A1}, %and for simplicity we assume now that $b\equiv0$,
it follows that
\begin{equation*}
    \begin{split}
E^Q[L_t\,|\,\mathcal{F}_t^Y]&=1+\int_0^t E^Q[L_s h(s,Y_{\cdot\wedge s},X_s,\mu_s^{X|Y})\,|\,\mathcal{F}_s^Y]dY_s, \ t\in [0,T],
    \end{split}
\end{equation*}
and applying It\^{o}'s formula to \eqref{A1} before taking conditional expectation gives that
\begin{equation*}
E^Q[X_t L_t\,|\,\mathcal{F}_t^Y]=x_0+\int_0^t E^Q[X_s L_s h(s,Y_{\cdot\wedge s},X_s,\mu_s^{X|Y})\,|\,\mathcal{F}_s^Y]dY_s,\ t\in [0,T].
\end{equation*}
Thus, applying It\^{o}'s formula to $U_t^{X|Y}=\frac{E^Q[L_t X_t\,|\,\mathcal{F}_t^Y]}{E^Q[L_t\,|\,\mathcal{F}_t^Y]}$ we deduce the so-called Fujisaki-Kallianpur-Kunita (FKK) equation: For $t\in[0,T],\ Q$-a.s.,
\begin{equation}\label{FKK}
\begin{split}
d U_t^{X|Y}=&d E^{P}[X_t\,|\,\mathcal{F}_t^Y]\\
=&\big\{E^{P}[X_t h(t,Y_{\cdot\wedge t},X_t,\mu_t^{X|Y})\,|\,\mathcal{F}_t^Y]-E^{P}[X_t\,|\,\mathcal{F}_t^Y]E^{P}[h(t,Y_{\cdot\wedge t},X_t,\mu_t^{X|Y})\,|\,\mathcal{F}_t^Y]\big\}dY_t\\
&+\Big\{E^{P}[X_t\,|\,\mathcal{F}_t^Y]\big(E^{P}[h(t,Y_{\cdot\wedge t},X_t,\mu_t^{X|Y})\,|\,\mathcal{F}_t^Y]\big)^2\\
&\quad \ \ -E^{P}[X_t h(t,Y_{\cdot\wedge t},X_t,\mu_t^{X|Y})\,|\,\mathcal{F}_t^Y]E^{P}[h(t,Y_{\cdot\wedge t},X_t,\mu_t^{X|Y})\,|\,\mathcal{F}_t^Y]\Big\}dt.
\end{split}
\end{equation}
Equation \eqref{FKK} shows especially that $U^{X|Y}$ has a continuous version with which we identify $U^{X|Y}$.
\end{remark}

\begin{theorem}\label{weak sol equa}
Under the Assumption (H1), let $(\Omega^i,\mathcal{F}^i,\mathbb{F}^i,P^i,(B^{1,i},B^{2,i}),(X^i,Y^i)),\ i=1,2,$ be two weak solutions of \eqref{A4}. Then it holds that
\begin{equation}
P^1_{((B^{1,1},B^{2,1}),(X^1,Y^1))}=P^2_{((B^{1,2},B^{2,2}),(X^2,Y^2))}.
\end{equation}
\end{theorem}

\begin{proof}
For $i=1,2$ and $t\in[0,T]$, we denote $U_t^i=E^{P^i}[X_t^i\,|\,\mathcal{F}_t^{Y^i}]$ and $\mu_t^i=P^i_{U_t^i}$. We define
$$ \widetilde{L}^i_t=\exp\Big\{-\int_0^t h(s,Y^i_{\cdot\wedge s},X^i_s,\mu_s^i)dB_s^{2,i}-\frac{1}{2}\int_0^t |h(s,Y^i_{\cdot\wedge s},X^i_s,\mu_s^i)|^2 ds\Big\},\ i=1,2. $$
Then $(B^{1,i},Y^i)$ is an $(\mathbb{F}^i,Q^i)$ Brownian motion where $dQ^i=\widetilde{L}^i_T dP^i$. Setting $L_t^i=\frac{1}{\widetilde{L}^i_t}$, it holds that $dP^i=L^i_T dQ^i,\ i=1,2$. Then we consider the following SDE on the probability space $(\Omega^1,\mathcal{F}^1,Q^1)$:
\begin{equation}\label{Q1}
\left\{
\begin{aligned}
&dX^1_t=\sigma(t,Y^1_{\cdot\wedge t},X^1_t,\mu_t^1)dB_t^{1,1},\ t\in[0,T],\ X^1_0=x_0;  \\
&dL^1_t= L^1_t h(t,Y^1_{\cdot\wedge t},X^1_t,\mu_t^1)dY^1_t,\ t\in[0,T],\ L^1_0=1.  \\
\end{aligned}\right.
\end{equation}
Putting $\widetilde{\phi}(t,y_{\cdot\wedge t},x):=\phi(t,y_{\cdot\wedge t},x,\mu_t^1)$, for $\phi=\sigma,\,h$, the following intermediate SDE on $(\Omega^2,\mathcal{F}^2,Q^2)$ is a classical one:
\begin{equation}\label{InQ2}
\left\{
\begin{aligned}
&d\widetilde{X}^2_t=\sigma(t,Y^2_{\cdot\wedge t},\widetilde{X}^2_t,\mu_t^1)dB_t^{1,2},\ t\in[0,T],\ \widetilde{X}^2_0=x_0;  \\
&d\widetilde{L}^2_t=\widetilde{L}^2_t h(t,Y^2_{\cdot\wedge t},\widetilde{X}^2_t,\mu_t^1)dY^2_t,\ t\in[0,T],\ \widetilde{L}^2_0=1.   \\
\end{aligned}\right.
\end{equation}
Then from the classical SDE theory (cf. \cite{IW81}) it follows that there is a unique measurable functional $\Phi:C_T^2:=C_T\times C_T\rightarrow C_T^2$ such that $(X^1,L^1)=\Phi(B^{1,1},Y^1),\ Q^1$-a.s.
and $(\widetilde{X}^2,\widetilde{L}^2)=\Phi(B^{1,2},Y^2),\ Q^2$-a.s. Since $Q^1_{(B^{1,1},Y^1)}=Q^2_{(B^{1,2},Y^2)}$ are the Wiener measure on $(C_T^2,\mathcal{B}(C_T^2))$, we deduce that
\begin{equation}\label{id}
Q^1_{(B^{1,1},Y^1,X^1,L^1)}=Q^2_{(B^{1,2},Y^2,\widetilde{X}^2,\widetilde{L}^2)}.
\end{equation}
Now our task is to show that $(\widetilde{X}^2,\widetilde{L}^2)$ coincides with $(X^2,L^2)$ under $Q^2$. Using the solution $(\widetilde{X}^2,\widetilde{L}^2)$ of \eqref{InQ2}, we define
 $\widetilde{P}^2=\widetilde{L}^2 Q^2$, $\widetilde{U}_t^2=E^{\widetilde{P}^2}[\widetilde{X}_t^2\,|\,\mathcal{F}_t^{Y^2}]$ and $\widetilde{\mu}_t^2=\widetilde{P}^2_{\widetilde{U}_t^2}=\widetilde{P}^2_{E^{\widetilde{P}^2}[\widetilde{X}_t^2\,|\,\mathcal{F}_t^{Y^2}]}$, $t\in[0,T]$. Due to the preceding remark, the process $\widetilde{U}^2$ can be considered to be continuous. On the other hand, recalling that
\begin{equation}
U_t^1=E^{P^1}[X_t^1|\mathcal{F}_t^{Y^1}]=\frac{E^{Q^1}[L_t^1 X_t^1\,|\,\mathcal{F}_t^{Y^1}]}{E^{Q^1}[L_t^1\,|\,\mathcal{F}_t^{Y^1}]},\ t\in[0,T],
\end{equation}
we know there exists a measurable functional $\Psi:[0,T]\times C_T\rightarrow\mathbb{R}$ such that $U_t^1=\Psi(t,Y^1_{\cdot\wedge t}),\ Q^1$-a.s., $t\in[0,T]$. For all $t\in[0,T]$ and any bounded Borel measurable function $f:C_T\rightarrow\mathbb{R}$, from \eqref{id} it follows that
\begin{equation*}
\begin{split}
&E^{\widetilde{P}^2}[\widetilde{U}_t^2 f(Y^2_{\cdot\wedge t})]=E^{\widetilde{P}^2}[\widetilde{X}_t^2 f(Y^2_{\cdot\wedge t})]
=E^{Q^2}[\widetilde{L}_t^2\widetilde{X}_t^2 f(Y^2_{\cdot\wedge t})]=E^{Q^1}[L_t^1 X_t^1 f(Y^1_{\cdot\wedge t})]\\
=&E^{P^1}[X_t^1 f(Y^1_{\cdot\wedge t})]=E^{P^1}[U_t^1 f(Y^1_{\cdot\wedge t})]=E^{Q^1}[L_t^1 U_t^1 f(Y^1_{\cdot\wedge t})]
=E^{Q^1}[L_t^1 \Psi(t,Y^1_{\cdot\wedge t})f(Y^1_{\cdot\wedge t})]\\
=&E^{Q^2}[\widetilde{L}_t^2 \Psi(t,Y^2_{\cdot\wedge t})f(Y^2_{\cdot\wedge t})]=E^{\widetilde{P}^2}[\Psi(t,Y^2_{\cdot\wedge t})f(Y^2_{\cdot\wedge t})].
\end{split}
\end{equation*}
This implies $\widetilde{U}_t^2=\Psi(t,Y^2_{\cdot\wedge t}),\ \widetilde{P}^2$-a.s., $t\in[0,T]$. Then, for any bounded Borel measurable function $g:\mathbb{R}\rightarrow\mathbb{R}$, we have
\begin{equation*}
\begin{split}
\int_{\mathbb{R}}g(x)\widetilde{\mu}_t^2(dx)&=E^{\widetilde{P}^2}[g(\widetilde{U}_t^2)]=E^{Q^2}[\widetilde{L}_t^2 g(\Psi(t,Y^2_{\cdot\wedge t}))]=E^{Q^1}[L_t^1 g(\Psi(t,Y^1_{\cdot\wedge t}))]\\
&=E^{P^1}[g(U_t^1)]=\int_{\mathbb{R}}g(x)\mu_t^1(dx),
\end{split}
\end{equation*}
which means $\widetilde{\mu}_t^2=\mu_t^1,\ t\in[0,T]$. Hence, we can rewrite \eqref{InQ2} as follow
\begin{equation}\label{InQ2'}
\left\{
\begin{split}
&d\widetilde{X}^2_t=\sigma(t,Y^2_{\cdot\wedge t},\widetilde{X}^2_t,\widetilde{\mu}_t^2)dB_t^{1,2},\ t\in[0,T],\ \widetilde{X}^2_0=x_0;  \\
&d\widetilde{L}^2_t=\widetilde{L}^2_t h(t,Y^2_{\cdot\wedge t},\widetilde{X}^2_t,\widetilde{\mu}_t^2)dY^2_t,\ t\in[0,T],\ \widetilde{L}^2_0=1.  \\
\end{split}\right.
\end{equation}
Now we compare \eqref{InQ2'} with
\begin{equation}\label{Q2}
\left\{
\begin{split}
&dX^2_t=\sigma(t,Y^2_{\cdot\wedge t},X^2_t,\mu_t^2)dB_t^{1,2},\ t\in[0,T],\ X^2_0=x_0;  \\
&dL^2_t=L^2_t h(t,Y^2_{\cdot\wedge t},X^2_t,\mu_t^2)dY^2_t,\ t\in[0,T],\ L^2_0=1.   \\
\end{split}\right.
\end{equation}
Applying the same method as that for the proof of \eqref{difmu}, we deduce that
\begin{equation}
\begin{split}
W_1(\mu_t^2,\widetilde{\mu}_t^2)&=\sup\Big\{\Big|\int_{\mathbb{R}}\varphi(x)\mu^2_t(dx)-\int_{\mathbb{R}}\varphi(x)\widetilde{\mu}^2_t(dx)\Big|:\ \varphi\in \mbox{Lip}_1(\mathbb{R}),\ \varphi(0)=0\Big\}\\
&\leq C \Big(\big(E^{Q^2}[|L_t^2-\widetilde{L}_t^2|^2]\big)^{\frac{1}{2}}+\big(E^{Q^2}[|X_t^2-\widetilde{X}_t^2|^2]\big)^{\frac{1}{2}}\Big).
\end{split}
\end{equation}
Combining this estimate with the corresponding one of the form \eqref{A3} (But, of course, now with $((\overline{X}^1,\overline{X}^2),(\overline{L}^1,\overline{L}^2))=((V^1,V^2),(K^1,K^2))=((X^2,\widetilde{X}^2),(L^2,\widetilde{L}^2))$), we have
\begin{equation}
\begin{split}
W_1^4(\mu_t^2,\widetilde{\mu}_t^2)\leq C\int_0^t W_1^4(\mu_s^2,\widetilde{\mu}_s^2) ds,\ t\in[0,T].
\end{split}
\end{equation}
\iffalse
On the other hand, it holds that $t\mapsto W_1(\mu_t^2,\widetilde{\mu}_t^2)$ is continuous. Indeed, for $s,t\in[0,T]$,
\begin{equation*}
\begin{split}
W_1(\mu_t^2,\widetilde{\mu}_t^2)&\leq W_1(\mu_t^2,\mu_s^2)+ W_1(\mu_s^2,\widetilde{\mu}_t^2) \\
&\leq W_1(\mu_t^2,\mu_s^2)+ W_1(\mu_s^2,\widetilde{\mu}_s^2)+ W_1(\widetilde{\mu}_s^2,\widetilde{\mu}_t^2).
\end{split}
\end{equation*}
There is symmetry and then we have
\begin{equation*}
\begin{split}
|W_1(\mu_t^2,\widetilde{\mu}_t^2)-W_1(\mu_s^2,\widetilde{\mu}_s^2)|&\leq W_1(\mu_t^2,\mu_s^2)+W_1(\widetilde{\mu}_s^2,\widetilde{\mu}_t^2)\\
&=W_1(P^2_{U_t^2},P^2_{U_s^2})+W_1(\widetilde{P}^2_{\widetilde{U}_s^2},\widetilde{P}^2_{\widetilde{U}_s^2})\\
&\leq E^{P^2}[|U_t^2-U_s^2|]+E^{\widetilde{P}^2}[|\widetilde{U}_t^2-\widetilde{U}_s^2|].
\end{split}
\end{equation*}
\fi
On the other hand, as the processes $U_t^2=E^{P^2}[X_t^2\,|\,\mathcal{F}_t^{Y^2}]$ and $\widetilde{U}_t^2=E^{\widetilde{P}^2}[\widetilde{X}_t^2\,|\,\mathcal{F}_t^{Y^2}]$, $t\in[0,T]$, are continuous ($L^1(P^2)$- and $L^1(\widetilde{P}^2)$-continuous, respectively), also $\mu_t^2=P^2_{U_t^2}$, $\widetilde{\mu}_t^2=\widetilde{P}^2_{\widetilde{U}_t^2}$, $t\in[0,T]$, are continuous in $\mathcal{P}_1(\mathbb{R})$. Hence, also $t\mapsto W_1(\mu_t^2,\widetilde{\mu}_t^2)$ is continuous, and from Gronwall's Lemma we deduce that, for $t\in[0,T]$,
\begin{equation}
W_1(\mu_t^2,\widetilde{\mu}_t^2)=0.
\end{equation}
Thus, $\mu_t^2=\widetilde{\mu}_t^2,\ t\in[0,T]$. By the pathwise uniqueness of SDE \eqref{Q2}, we conclude that $(\widetilde{X}^2,\widetilde{L}^2)$ coincides with $(X^2,L^2)$ under $Q^2$. Together with \eqref{id} we complete the proof.

\end{proof}

\section{Stochastic Optimization Problem}
\subsection{Problem Formulation}
The well-posedness obtained in the preceding section for \eqref{A4} and \eqref{A1}, respectively, allows to consider an associate control problem.
From now on we let $Q$ be the reference probability measure on $(\Omega,\mathcal{F})$, under which the coordinate process $(B^1,Y)$ is a Brownian motion. Recall that $\mathbb{F}:=\mathbb{F}^{B^1,Y}$ is the filtration generated by $(B^1,Y)$. The $\sigma$-field $\mathcal{F}$ as well as the filtration $\mathbb{F}$ are considered as complete under $Q$. We study the following controlled stochastic system
\begin{equation}\label{system}
\Bigg\{
\begin{tabular}{l}
$dX_t^u=\sigma(t,X_t^u,\mu_t^u,u_t)dB_t^1,\ X_0^u=x;$ \\
$dL_t^u=L_t^u h(t,X_t^u,\mu_t^u,u_t)dY_t,\ L_0^u=1,\ t\in[0,T],$  \\
\end{tabular}
\end{equation}
where $P^u=L_T^u Q$, $\mu_t^u=\mu_t^{X^u|Y}=P^u_{E^u[X_t^u\,|\,\mathcal{F}_t^Y]}$ and $E^u[\cdot]:=E^{P^u}[\cdot]$ is the expectation under $P^u$. Here $u$ is an admissible control process. To be more precise,
for an arbitrarily fixed nonempty subset $U\subset\mathbb{R}^k$ (the control state space) the control $u$ runs the set of admissible controls\  $\mathcal{U}_{ad}=L_{\mathbb{F}^Y}^0([0, T], Q;U)$, where $$L_{\mathbb{F}^Y}^0([0, T], Q;U):=\Big\{v\,\Big|\,\mbox{the stochastic process}\ v=(v_t)_{t\in[0,T]}\mbox{ is }U\mbox{-valued, and }\mathbb{F}^Y\mbox{-adapted}\Big\}.$$

Let us define the cost functional:
\begin{equation}\label{cost} J(u):=E^Q\Big[\Phi(X_T^u,\mu_T^u)+\int_0^T f(t,X_t^u,\mu_t^u,u_t)dt\Big],\ u\in\mathcal{U}_{ad}. \end{equation}
The goal of the control problem is to minimize the cost functional.

 We say a control $u^*\in\mathcal{U}_{ad}$ is optimal if $J(u^*)=\inf\limits_{u\in\mathcal{U}_{ad}}J(u)$.\\

Let us make the following standard assumptions.\\
\smallskip
\textbf{Assumption (H2)}: For $\phi=\sigma,\,h,\,f,\,\Phi$, we suppose\\
(i) The function $\phi:[0,T]\times\mathbb{R}\times\mathcal{P}_2(\mathbb{R})\times U\rightarrow\mathbb{R}$ is Borel measurable, and to simplify the computations, we also suppose the boundedness;\\
(ii) For all $t\in[0,T]$, $\mu\in\mathcal{P}_2(\mathbb{R})$ and $v\in U$, the function $\phi(t,\cdot,\mu,v)$ is in $C^2_b(\mathbb{R})$;\\
(iii) For all $t\in[0,T]$, $x\in\mathbb{R}$ and $v\in U$, the function $\phi(t,x,\cdot,v)$ is differentiable on $\mathcal{P}_2(\mathbb{R})$; $\partial_\mu\phi(t,x,\mu,v;y)$ is bounded and also differentiable w.r.t. $\mu\in\mathcal{P}_2(\mathbb{R})$ and $x, y\in\mathbb{R}$, and the derivatives, denoted by $\partial_\mu(\partial_\mu\phi)$, $\partial_x(\partial_\mu\phi)$ and $\partial_z(\partial_\mu\phi)$, respectively, are bounded. Moreover, we have the following continuity conditions: For $t\in[0,T]$, $v\in U$, $\mu,\mu'\in\mathcal{P}_2(\mathbb{R})$ and $x,x',y,y',z,z'\in\mathbb{R}$,
\begin{equation*}
\begin{split}
&\ \ |\phi(t,x,\mu,v)-\phi(t,x,\mu',v)|\leq CW_1(\mu,\mu'), \\
&\ \ |\partial_\mu\phi(t,x,\mu,v;y)-\partial_\mu\phi(t,x',\mu',v;y')|\leq C\big(W_1(\mu,\mu')+|x-x'|+|y-y'|\big),\\
&\ \ |\partial_x(\partial_\mu\phi)(t,x,\mu,v;y)-\partial_x(\partial_\mu\phi)(t,x',\mu',v;y')|\leq C\big(W_1(\mu,\mu')+|x-x'|+|y-y'|\big),\\
&\ \ |\partial_z(\partial_\mu\phi)(t,x,\mu,v;y)-\partial_z(\partial_\mu\phi)(t,x',\mu',v;y')|\leq C\big(W_1(\mu,\mu')+|x-x'|+|y-y'|\big),\\
&\ \ |\partial_\mu(\partial_\mu\phi)(t,x,\mu,v;y,z)-\partial_\mu(\partial_\mu\phi)(t,x',\mu',v;y',z')|\leq C\big(W_1(\mu,\mu')+|x-x'|+|y-y'|+|z-z'|\big).
\end{split}
\end{equation*}

Under the Assumption (H2) we know from Proposition \ref{weak sol equa} that, for all $u\in\mathcal{U}_{ad}$, SDE \eqref{system} admits a unique solution
 $(X^u,L^u)\in S^2_{\mathbb{F}}([0,T],Q)\times S^2_{\mathbb{F}}([0,T],Q)$. Moreover, $X^u,\,L^u,\,U^u$ are in all $S^p_{\mathbb{F}}([0,T],Q)$, for $p\geq1$.

\begin{remark}
For all $p\geq1$, we have $\mu_t^u\in\mathcal{P}_p(\mathbb{R}),\ t\in[0,T]$. Indeed,
$\!\displaystyle \int_\mathbb{R}\!|x|^p\mu_t^u(dx)\!=\! E^u[|U_t^u|^p]<\infty.$
\end{remark}

\begin{remark}
In \cite{BLM17}, the setting of the coefficients is $\displaystyle \phi(t,x,\gamma,u):=\int \phi(t,x,z,u)\gamma(dz)$, for $\phi=\sigma, f$;
$\displaystyle h(t,x,\gamma,u):=h(t,x);\  \Phi(x,\gamma):=\int\Phi(x,z)\gamma(dz), \ (t,x,\gamma,u)\in[0,T]\times\mathbb{R}\times\mathcal{P}_2(\mathbb{R})\times\mathcal{U}_{ad}.$

 That is, the coefficients $\sigma, f$ and $\Phi$ all depend on $\gamma$ only linearly, and $h$ does not depend on the control nor on $\gamma$. Moreover, the SMP studied there is the Pontryagin one.
\end{remark}

We suppose the existence of an optimal control $u^*\in\mathcal{U}_{ad}$, and for shortness we write $u:=u^*$. Our objective is to derive a necessary optimality condition for $u$,
using Peng's SMP.

In the following calculations throughout this paper, we let $(\widetilde{\Omega},\widetilde{\mathcal{F}},\widetilde{Q})$ be a copy of $(\Omega,\mathcal{F},Q)$. Furthermore, for each
 $\xi\in L^0(\Omega,\mathcal{F},Q)$, we denote by $\widetilde{\xi}\in L^0(\widetilde{\Omega},\widetilde{\mathcal{F}},\widetilde{Q})$ an independent copy of $\xi$, i.e.,
 $\widetilde{\xi}$ under $\widetilde{Q}$ has the same law as $\xi$ under $Q$ and they are independent. In the same spirit we can consider another copy $(\widehat{\Omega},\widehat{\mathcal{F}},\widehat{Q})$,
  and the copied random variable $\widehat{\xi}$ on $(\widehat{\Omega},\widehat{\mathcal{F}},\widehat{Q})$, which is independent of both $\widetilde{\xi}$ and $\xi$, and of the same law.

\subsection{Variational Equation}
As the control state set $U$ is not supposed to be convex, we shall consider Peng's stochastic maximum principle here. Assume that $u\in\mathcal{U}_{ad}$ is an optimal control, and consider an arbitrary but fixed $v\in\mathcal{U}_{ad}$. For $\varepsilon>0$, we let $E_{\varepsilon}\in\mathcal{B}([0,T])$ with $|E_{\varepsilon}|=\varepsilon$, and we put
$$ u^{\varepsilon}:=u \mathbf{1}_{E_\varepsilon^c}(t)+v\mathbf{1}_{E_\varepsilon}(t),\ t\in[0,T]. $$
The process $u^{\varepsilon}\in\mathcal{U}_{ad}$ is a so-called spike variation of the optimal control $u$.

For simplicity we also introduce the following notations: For $\phi=\sigma,h,f$ and $\Phi$ we set
\begin{equation}
\begin{array}{ll}
\phi(t):=\phi(t,X_t^u,\mu_t^u,u_t), & \delta\phi(t):=\phi(t,X_t^u,\mu_t^u,v_t)-\phi(t,X_t^u,\mu_t^u,u_t),\\
\phi_x(t):=\partial_x \phi(t,X_t^u,\mu_t^u,u_t), & \phi_{xx}(t):=\partial^2_{xx} \phi(t,X_t^u,\mu_t^u,u_t),\\
\phi_\mu(t,y):=\partial_\mu \phi(t,X_t^u,\mu_t^u,u_t;y), & \widetilde{\phi}_\mu(t):=\phi_\mu(t,\widetilde{U}_t^u)=\partial_\mu \phi(t,X_t^u,\mu_t^u,u_t;\widetilde{U}_t^u),\\
\widetilde{\phi}^*_\mu(t):=\partial_\mu \phi(t,\widetilde{X}_t^u,\mu_t^u,\widetilde{u}_t;U_t^u), & \phi^*_\mu(t,y):=\partial_\mu \phi(t,\widetilde{X}_t^u,\mu_t^u,\widetilde{u}_t;y)\\
\phi_{z\mu}(t,y):=\partial_z(\partial_\mu\phi)(t,X_t^u,\mu_t^u,u_t;y), &\phi_{z\mu}^*(t,y):=\partial_z(\partial_\mu\phi)(t,\widetilde{X}_t^u,\mu_t^u,\widetilde{u}_t;y),\\
\widetilde{\phi}_{z\mu}(t):=\partial_z(\partial_\mu\phi)(t,X_t^u,\mu_t^u,u_t;\widetilde{U}_t^u), &
\widetilde{\phi}^*_{z\mu}(t):=\partial_z(\partial_\mu\phi)(t,\widetilde{X}_t^u,\mu_t^u,\widetilde{u}_t;U_t^u).
\end{array}
\end{equation}

Moreover, we denote $(X,L):=(X^u,L^u)$, $P:=L_T Q(=P^u)$, $\mu:=\mu^u$, $U_t=U_t^u:=E^P[X_t\,|\,\mathcal{F}_t^Y]$, $t\in[0,T]$. In the same spirit we set $(X^\varepsilon,L^\varepsilon):=(X^{u^\varepsilon},L^{u^\varepsilon})$, $P^\varepsilon:=P^{u^\varepsilon}$, $\mu^\varepsilon:=\mu^{u^\varepsilon}$ and $U_t^\varepsilon:=E^{P^\varepsilon}[X_t^\varepsilon\,|\,\mathcal{F}_t^Y]$, $t\in[0,T]$.

For $\varepsilon>0$, the state-observation dynamics is as follows:
\begin{equation}\label{SMP1}
\left\{
\begin{split}
&dX_t^\varepsilon=\sigma(t,X_t^\varepsilon,\mu_t^\varepsilon,u^\varepsilon_t)dB_t^1,\ X_0^\varepsilon=x; \\
&dL_t^\varepsilon=L_t^\varepsilon h(t,X_t^\varepsilon,\mu_t^\varepsilon,u^\varepsilon_t)dY_t,\ L_0^\varepsilon=1,\ t\in[0,T];  \\
&\mu_t^\varepsilon=P^\varepsilon_{U_t^\varepsilon}, \mbox{ with } P^\varepsilon=L_T^\varepsilon Q, \ U_t^\varepsilon=E^{P^\varepsilon}[X_t^\varepsilon\,|\,\mathcal{F}_t^Y]=\frac{E^{Q}[L_t^\varepsilon X_t^\varepsilon\,|\,\mathcal{F}_t^Y]}{E^{Q}[L_t^\varepsilon\,|\,\mathcal{F}_t^Y]}.
\end{split}\right.
\end{equation}
For $\varepsilon=0$, we put $(X^0,L^0,U^0,\mu^0,u^0,P^0):=(X,L,U,\mu,u,P)$.

Formally, we should derive \eqref{SMP1} with respect to $\varepsilon$ at $\varepsilon=0$, but as $\phi=\sigma,h$, is not differentiable in the control variable, we take $\delta\phi(t)=\phi(t,X_t,\mu_t,v_t)-\phi(t,X_t,\mu_t,u_t)$ instead of $\partial_\varepsilon\big[\phi(t,X_t,\mu_t,u^\varepsilon_t)\big]_{|\varepsilon=0}$.
For this we recall that, if $f:\mathcal{P}_2(\mathbb{R})\rightarrow\mathbb{R}$ is continuously differentiable and $\varepsilon\rightarrow(X^\varepsilon,L^\varepsilon,U^\varepsilon)$ were differentiable in $\varepsilon=0$, then from Theorem 3.2 in \cite{Dwrtd} we would have
\begin{equation*}
\begin{split}
&\partial_\varepsilon f(\mu_t^\varepsilon)_{|\varepsilon=0}=\partial_\varepsilon \big[f\big((L_T^\varepsilon Q)_{U_t^\varepsilon}\big)\big]_{|\varepsilon=0}=\partial_\varepsilon \big[f\big((L_T^\varepsilon Q)_{U_t}\big)\big]_{|\varepsilon=0}+\partial_\varepsilon \big[f\big((L_T Q)_{U_t^\varepsilon}\big)\big]_{|\varepsilon=0}\\
=&E^Q\Big[\int_0^{U_t}\partial_\mu f\big((L_T Q)_{U_t},y\big)dy\cdot{\partial_\varepsilon L_T^\varepsilon}_{|\varepsilon=0}\Big]+E^{L_T Q}\big[\partial_\mu f\big((L_T Q)_{U_t},U_t\big)\cdot{\partial_\varepsilon U_t^\varepsilon}_{|\varepsilon=0}\big].
\end{split}
\end{equation*}
As $\displaystyle \int_0^{U_t}\partial_\mu f\big((L_T Q)_{U_t},y\big)dy$ is $\mathcal{F}_t$-measurable and $L^\varepsilon$ is an $(\mathbb{F},Q)$-martingale, this would give
\begin{equation*}
\begin{split}
\partial_\varepsilon f(\mu_t^\varepsilon)_{|\varepsilon=0}=&E^Q\Big[\int_0^{U_t}\partial_\mu f(\mu_t,y)dy\cdot{\partial_\varepsilon L_t^\varepsilon}_{|\varepsilon=0}\Big]+E^{Q}\big[\partial_\mu f(\mu_t,U_t)L_t\cdot{\partial_\varepsilon U_t^\varepsilon}_{|\varepsilon=0}\big]\\
\bigg(=&E^{Q}\Big[\partial_\varepsilon \Big(L_t^\varepsilon\int_0^{U^\varepsilon_t}\partial_\mu f(\mu_t,y)dy\Big)_{|\varepsilon=0}\Big]\bigg),
\end{split}
\end{equation*}
with
\begin{equation*}
\begin{split}
&{\partial_\varepsilon U^\varepsilon_t}_{|\varepsilon=0}=\frac{E^Q[X_t {\partial_\varepsilon L_t^\varepsilon}_{|\varepsilon=0}+L_t {\partial_\varepsilon X_t^\varepsilon}_{|\varepsilon=0}\,|\,\mathcal{F}_t^Y]}
{E^Q[L_t\,|\,\mathcal{F}_t^Y]}-\frac{E^{Q}[L_t X_t\,|\,\mathcal{F}_t^Y]}{(E^{Q}[L_t\,|\,\mathcal{F}_t^Y])^2} E^Q[{\partial_\varepsilon L_t^\varepsilon}_{|\varepsilon=0}\,|\,\mathcal{F}_t^Y]\\
&=E^P\big[X_t \partial_\varepsilon[\ln L_t^\varepsilon]_{|\varepsilon=0}\,\big|\,\mathcal{F}_t^Y\big]+E^P\big[{\partial_\varepsilon X_t^\varepsilon}_{|\varepsilon=0}\,\big|\,\mathcal{F}_t^Y\big]-E^P[X_t\,|\,
\mathcal{F}_t^Y]E^P\big[\partial_\varepsilon[\ln L_t^\varepsilon]_{|\varepsilon=0}\,\big|\,\mathcal{F}_t^Y\big].
\end{split}
\end{equation*}
But, however, the derivatives ${\partial_\varepsilon X^\varepsilon}_{|\varepsilon=0}$ and ${\partial_\varepsilon L^\varepsilon}_{|\varepsilon=0}$ don't exist. They will be replaced by the solution of the first order variational equation $Y^{1,\varepsilon}=(Y^{1,\varepsilon})_{t\in[0,T]}$ and $K^{1,\varepsilon}=(K^{1,\varepsilon})_{t\in[0,T]}$, respectively. Together with the classical dependence of the coefficients $\phi=\sigma,h$ on $X^\varepsilon$ this suggests the following first order variational equations whose choice will have to be confirmed by the fact that $X_t^\varepsilon-(X_t+Y_t^{1,\varepsilon})=O(\varepsilon)$ and $L_t^\varepsilon-(L_t+K_t^{1,\varepsilon})=O(\varepsilon)$, in $L^2([0,T],Q)$, uniformly in $t\in[0,T]$, as $\varepsilon\searrow0$. For $\varepsilon>0$,
\begin{equation}\label{SMP2}
\left\{
\begin{split}
&dY_t^{1,\varepsilon}=\Big\{\sigma_x(t)Y_t^{1,\varepsilon}+\widetilde{E}^Q\Big[\int_0^{\widetilde{U}_t}\sigma_\mu(t,y)dy\cdot\widetilde{K}_t^{1,\varepsilon}\Big]+\widetilde{E}^Q\big
[\widetilde{\sigma}_\mu(t)\widetilde{L}_t\widetilde{V}_t^{1,\varepsilon}\big]+\delta\sigma(t)\mathbf{1}_{E_\varepsilon}(t)\Big\}dB_t^1,\\ &Y_0^{1,\varepsilon}=0; \\ &dK_t^{1,\varepsilon}=\Big\{h(t)K_t^{1,\varepsilon}+\Big(h_x(t)Y_t^{1,\varepsilon}+\widetilde{E}^Q\Big[\int_0^{\widetilde{U}_t}h_\mu(t,y)dy\cdot\widetilde{K}_t^{1,\varepsilon}\Big]\\
&\qquad\qquad\qquad\qquad\qquad\qquad+\widetilde{E}^Q\big[\widetilde{h}_\mu(t)\widetilde{L}_t\widetilde{V}_t^{1,\varepsilon}\big]+\delta h(t)\mathbf{1}_{E_\varepsilon}(t)\Big)L_t\Big\}dY_t,\\
&K_0^{1,\varepsilon}=0;  \\
&V_t^{1,\varepsilon}=\frac{E^Q[L_t Y_t^{1,\varepsilon}+X_t K_t^{1,\varepsilon}\,|\,\mathcal{F}_t^Y]}{E^Q[L_t\,|\,\mathcal{F}_t^Y]}-\frac{E^{Q}[L_t X_t\,|\,\mathcal{F}_t^Y] E^{Q}[K_t^{1,\varepsilon}\,|\,\mathcal{F}_t^Y]}{(E^{Q}[L_t\,|\,\mathcal{F}_t^Y])^2},\ t\in[0,T].
\end{split}\right.
\end{equation}
\begin{remark}
In \cite{BLM17}, Pontryagin's stochastic maximum principle is studied for the case $h=h(t,x)$ and $\sigma(t,x,\mu,u)=\int\sigma(t,x,y,u)\mu(dy)$, and it holds that $h_\mu(t,y)=\partial_\mu h(t,X_t,y)=0$, $\widetilde{h}_\mu(t)=0$ and $\sigma_\mu(t,y)=\partial_y\sigma(t,X_t,y,u_t)$, $\widetilde{\sigma}_\mu(t)=\sigma_\mu(t,\widetilde{U}_t)$, and
\begin{equation*}
\begin{split}
\widetilde{E}^Q\Big[\int_0^{\widetilde{U}_t}\sigma_\mu(t,y)dy\cdot\widetilde{K}_t^{1,\varepsilon}\Big]&=\widetilde{E}^Q\big[\big(\sigma(t,X_t,\widetilde{U}_t,u_t)-\sigma(t,X_t,0,u_t)
\big)\widetilde{K}_t^{1,\varepsilon}\big]\\
&=\widetilde{E}^Q\big[\sigma(t,X_t,\widetilde{U}_t,u_t)\widetilde{K}_t^{1,\varepsilon}\big].
\end{split}
\end{equation*}
Indeed, observe that $E^Q[K_t^{1,\varepsilon}]=0,\ t\in[0,T]$.
\end{remark}
\begin{proposition}\label{Wellof1stVar}
Under the Assumption (H2), \eqref{SMP2} has a unique solution $(Y^{1,\varepsilon},K^{1,\varepsilon})\in S_{\mathbb{F}}^2([0,T],Q)\times S_{\mathbb{F}}^2([0,T],Q)$. Moreover, $Y^{1,\varepsilon},\ K^{1,\varepsilon}$, $V^{1,\varepsilon}\in S_{\mathbb{F}}^p([0,T],Q)$ for all $p\geq 1$.
\end{proposition}
\begin{proof}
For $p\geq 2$, let $\zeta\in S_{\mathbb{F}}^{2p}([0,T],Q)$, $\eta\in S_{\mathbb{F}}^{p}([0,T],Q)$ and define
\begin{equation}\label{theta}\theta_t:=\theta_t(\zeta_t,\eta_t)=\frac{E^Q[L_t \zeta_t+X_t \eta_t\,|\,\mathcal{F}_t^Y]}{E^Q[L_t\,|\,\mathcal{F}_t^Y]}-\frac{E^{Q}[L_t X_t\,|\,\mathcal{F}_t^Y] E^{Q}[\eta_t\,|\,\mathcal{F}_t^Y]}{(E^{Q}[L_t\,|\,\mathcal{F}_t^Y])^2}. \end{equation}
Standard arguments shows that $\theta$ is continuous, and using that $X,\,L$ and $L^{-1}\in S_{\mathbb{F}}^{p}([0,T],Q)$, for all $p\geq 2$, we see that $\theta\in S_{\mathbb{F}}^{p}([0,T],Q)$. By putting, for $\varphi=\sigma,\,h$,
\begin{equation*}
\begin{split}
&F_1^\varphi(t,\zeta_t,\eta_t):=\varphi_x(t)\zeta_t+\widetilde{E}^Q\Big[\int_0^{\widetilde{U}_t}\varphi_\mu(t,y)dy\cdot\widetilde{\eta}_t\Big]+\widetilde{E}^Q\big[\widetilde{\varphi}
_\mu(t)\widetilde{L}_t \widetilde{\theta_t(\zeta_t,\eta_t)}\big]+\delta\varphi(t)\mathbf{1}_{E_\varepsilon}(t);\\
&F_2^h(t,\zeta_t,\eta_t):=h(t)\eta_t+L_t F_1^h(t,\zeta_t,\eta_t),
\end{split}
\end{equation*}
we define the functional $\big(F_1^\varphi(t,\cdot,\cdot),\, F_2^h(t,\cdot,\cdot)\big):L^{2p}(\mathcal{F}_t,Q)\times L^{p}
(\mathcal{F}_t,Q)\longrightarrow L^{2p}(\mathcal{F}_t,Q)\times L^{p}(\mathcal{F}_t,Q)$. Let us show that this functional is Lipschitz. As
$F_1^\sigma(t,0,0)=\delta\sigma(t)\mathbf{1}_{E_\varepsilon}(t)$, and $F_2^h(t,0,0)=L_t\delta h(t)\mathbf{1}_{E_\varepsilon}(t)$, a direct
consequence of this Lipschitz property is that the functional maps $L^{2p}(\mathcal{F}_t,Q)\times L^{p}(\mathcal{F}_t,Q)$ into itself. For
 the proof of Lipschitz property, let us proceed in several steps:\\
(1) Using the fact that $X,L,L^{-1}$ are in $S_{\mathbb{F}}^{p'}([0,T],Q)$, for all $p'\geq 1$, standard estimates using the definition
of $\theta_t$ (\ref{theta}) yield
$$ \big(E^Q[|\theta_t(\zeta_t,\eta_t)|^p]\big)^{\frac{1}{p}}\leq C_{p,q}\big(E^Q[|\zeta_t|^q+|\eta_t|^q]\big)^{\frac{1}{q}}, \mbox{ for all } 1\leq p<q, $$
where $C_{p,q}$ only depends on $p,q$.\\
(2) From (1), for $(\zeta,\eta),\,(\zeta',\eta')\in S_{\mathbb{F}}^{2p}([0,T],Q)\times S_{\mathbb{F}}^{p}([0,T],Q)$, recalling that $\widetilde{\sigma}_\mu(\cdot)$ is bounded, we have
\begin{equation}\label{Estoftheta}
\begin{split}
&\Big(E^Q\Big[\Big|\widetilde{E}^Q\big[\widetilde{\sigma}_\mu(t)\widetilde{L}_t\big(\widetilde{\theta_t(\zeta_t,\eta_t)}-\widetilde{\theta_t(\zeta'_t,\eta'_t)}\big)\big]\Big|^{2p}\Big]
\Big)^{\frac{1}{2p}}\leq C E^Q\big[L_t\big|\theta_t(\zeta_t,\eta_t)-\theta_t(\zeta'_t,\eta'_t)\big|\big]\\
&\leq C \Big(E^Q\big[\big|\theta_t(\zeta_t,\eta_t)-\theta_t(\zeta'_t,\eta'_t)\big|^{\frac{3}{2}}\big]\Big)^{\frac{2}{3}}\leq C\Big(E^Q\big[\big|\zeta_t-\zeta'_t\big|^2+\big|\eta_t-\eta'_t\big|^2\big]\Big)^{\frac{1}{2}}.
\end{split}
\end{equation}
Moreover, from the fact that $U\in S_{\mathbb{F}}^{p}([0,T],Q)$, for all $p\geq1$, we have
\begin{equation*}
\begin{split}
&\Big(E^Q\Big[\Big|\widetilde{E}^Q\big[\int_0^{\widetilde{U}_t}\sigma_\mu(t,y)dy(\widetilde{\eta}_t-\widetilde{\eta}'_t)\big]\Big|^{2p}\Big]\Big)^{\frac{1}{2p}}\\
&\leq C_p E^Q\big[|U_t|\cdot|\eta_t-\eta'_t|\big]\leq C_p\big(E^Q[|\eta_t-\eta'_t|^2]\big)^{\frac{1}{2}}.
\end{split}
\end{equation*}
Consequently, as $\sigma_x(\cdot)$ and $\delta\sigma(\cdot)$ are bounded, for all $p\geq2$, it holds that
$$ E^Q\big[\big|F_1^\sigma(t,\zeta_t,\eta_t)-F_1^\sigma(t,\zeta'_t,\eta'_t)\big|^{2p}\big]\leq C_p \Big(E^Q[|\zeta_t-\zeta'_t|^{2p}]+\big(E^Q[|\eta_t-\eta'_t|^p]\big)^2\Big), $$
and the same estimates holds true for $F_1^h$.\\
(3) From (2), as $h(t),\, t\in[0,T]$, is bounded, it follows that
\begin{equation*}
\begin{split}
&\Big(E^Q\big[\big|F_2^h(t,\zeta_t,\eta_t)-F_2^h(t,\zeta'_t,\eta'_t)\big|^{p}\big]\Big)^{\frac{1}{p}}\\
&\leq C_p \big(E^Q[|\eta_t-\eta'_t|^p]\big)^{\frac{1}{p}}+C \Big(E^Q\big[\big(L_t|F_1^h(t,\zeta_t,\eta_t)-F_1^h(t,\zeta'_t,\eta'_t)|\big)^{p}\big]\Big)^{\frac{1}{p}}\\
&\leq C_p \big(E^Q[|\eta_t-\eta'_t|^p]\big)^{\frac{1}{p}}+C \Big(E^Q\big[\big|F_1^h(t,\zeta_t,\eta_t)-F_1^h(t,\zeta'_t,\eta'_t)\big|^{2p}\big]\Big)^{\frac{1}{2p}}\\
&\leq C_p \Big(E^Q[|\zeta_t-\zeta'_t|^{2p}]+\big(E^Q[|\eta_t-\eta'_t|^{p}]\big)^2\Big)^{\frac{1}{2p}}.
\end{split}
\end{equation*}
Hence, from (2) and (3), for all $p\geq2$, for all $(\zeta,\eta),\,(\zeta',\eta')\in S_{\mathbb{F}}^{2p}([0,T],Q)\times S_{\mathbb{F}}^{p}([0,T],Q)$,

\begin{equation*}
\begin{split}
&E^Q\big[\big|F_1^\sigma(t,\zeta_t,\eta_t)-F_1^\sigma(t,\zeta'_t,\eta'_t)\big|^{2p}\big]\leq C_p \Big(E^Q[|\zeta_t-\zeta'_t|^{2p}]+\big(E^Q[|\eta_t-\eta'_t|^p]\big)^2\Big);\\
&E^Q\big[\big|F_2^h(t,\zeta_t,\eta_t)-F_2^h(t,\zeta'_t,\eta'_t)\big|^{p}\big]\leq C_p \Big(\big(E^Q[|\zeta_t-\zeta'_t|^{2p}]\big)^{\frac{1}{2}}+E^Q[|\eta_t-\eta'_t|^p]\Big).
\end{split}
\end{equation*}
Now by using standard arguments, this type of Lipschitz condition allows to show with Picard's iteration method that equation \eqref{SMP2}
\begin{equation*}\left\{
\begin{split}
dY_t^{1,\varepsilon}=F_1^\sigma(t,Y_t^{1,\varepsilon},K_t^{1,\varepsilon})dB_t^1,\ t\in[0,T],\ Y_0^{1,\varepsilon}=0;\\
dK_t^{1,\varepsilon}=F_2^h(t,Y_t^{1,\varepsilon},K_t^{1,\varepsilon})dY_t,\ t\in[0,T],\ K_0^{1,\varepsilon}=0,\\
\end{split}\right.
\end{equation*}
has a unique solution $(Y^{1,\varepsilon},K^{1,\varepsilon})$, and $(Y^{1,\varepsilon},K^{1,\varepsilon})\in S_{\mathbb{F}}^{2p}([0,T],Q)\times S_{\mathbb{F}}^{p}([0,T],Q)$ for all $p\geq2$.
\end{proof}
\begin{remark}
The computations in the above proof also show the Lipschitz property of $\theta_t:L^{q}(\mathcal{F}_t,Q)\times L^{q}(\mathcal{F}_t,Q)\rightarrow L^{p}(\mathcal{F}_t,Q)$, uniformly in $t\in[0,T]$, for
all $1\leq p<q$. For this it suffices to use (1) of the proof and the fact that $\theta_t$ is linear.
\end{remark}
Let us now give the following useful estimates.
\begin{proposition}\label{EstofXL}
For all $k\geq1$, there exists $C_k\in\mathbb{R}_+$, such that,
\begin{equation*}
\begin{array}{lll}
&\mbox{ }\ \ \ \ \ {\rm (i)}\  E^Q\big[\sup_{t\in[0,T]}\big(|X_t^\varepsilon|^{2k}+|L_t^\varepsilon|^{2k}\big)\big]\leq C_k;\\
&\mbox{ }\ \ \ \ \  {\rm (ii)}\ E^Q\big[\sup_{t\in[0,T]}\big(|X_t^\varepsilon-X_t|^{2k}+|L_t^\varepsilon-L_t|^{2k}\big)\big]\leq C_k \varepsilon^{k},\ \varepsilon>0;\\
&\mbox{ }\ \ \ \ \  {\rm (iii)}\ E^Q\big[\sup_{t\in[0,T]}\big(|Y_t^{1,\varepsilon}|^{2k}+|K_t^{1,\varepsilon}|^{2k}\big)\big]\leq C_k \varepsilon^{k},\ \varepsilon>0;\\
&\mbox{ }\ \ \ \ \  {\rm (iv)}\ E^Q\big[\sup_{t\in[0,T]}\big(|X_t^\varepsilon-(X_t+Y_t^{1,\varepsilon})|^{2k}+|L_t^\varepsilon-(L_t+K_t^{1,\varepsilon})|^{2k}\big)\big]\leq C_k \varepsilon^{2k},\ \varepsilon>0. \hspace{5cm}
\end{array}
\end{equation*}
\end{proposition}
\begin{proof}
We restrict here to the proof for (iv), that of (i), (ii) and (iii) is easier and more standard, so omitted here. To simplify the proof for (iv) we concentrate on the essential difficulties by assuming that $\sigma(t,x,\gamma,u)=\sigma(\gamma,u)$, $h(t,x,\gamma,u)=h(x,u)$. This allows to rewrite the equations as follows: For $t\in[0,T]$,
\begin{equation}\label{jj}
\begin{split}
\left\{
\begin{split}
&dX_t^\varepsilon=\sigma(\mu_t^\varepsilon,u^\varepsilon_t)dB_t^1,\ X_0^\varepsilon=x; \\
&dX_t=\sigma(\mu_t,u_t)dB_t^1,\ X_0=x; \\
&dY_t^{1,\varepsilon}=\Big\{\widetilde{E}^Q\Big[\int_0^{\widetilde{U}_t}\partial_\mu\sigma(\mu_t,u_t;y)dy\cdot\widetilde{K}_t^{1,\varepsilon}\Big]+\widetilde{E}^Q\big[\partial_\mu
\sigma(\mu_t,u_t;\widetilde{U}_t)\widetilde{L}_t\widetilde{V}_t^{1,\varepsilon}\big]+\delta\sigma(t)\mathbf{1}_{E_\varepsilon}(t)\Big\}dB_t^1,\\
&Y_0^{1,\varepsilon}=0.
\end{split}\right.\\
\left\{
\begin{split}
&dL_t^\varepsilon=h(X_t^\varepsilon,u^\varepsilon_t)L_t^\varepsilon dY_t,\ L_0^\varepsilon=1, \\
&dL_t=h(X_t,u_t)L_t dY_t,\ L_0=1, \\
&dK_t^{1,\varepsilon}=\Big\{h(X_t,u_t)K_t^{1,\varepsilon}+L_t h_x(X_t,u_t)Y_t^{1,\varepsilon}+L_t \delta h(t)\mathbf{1}_{E_\varepsilon}(t)\Big\}dY_t,\ K_0^{1,\varepsilon}=0. \hspace{2.7cm}
\end{split}\right.
\end{split}
\end{equation}

Recall the definition of $V_t^{1,\varepsilon}$ in \eqref{SMP2}. We put $(X^{\varepsilon,\lambda},L^{\varepsilon,\lambda},U^{\varepsilon,\lambda}):=(1-\lambda)(X,L,U)+\lambda(X^{\varepsilon},L^{\varepsilon},U^{\varepsilon}),\ \lambda\in[0,1]$, and we remark that, due to Theorem 3.2 in \cite{Dwrtd}, for $\mu_t^{\varepsilon,\lambda}:=(L_t^{\varepsilon,\lambda}Q)_{U_t^{\varepsilon,\lambda}}$,
\begin{equation*}
\begin{split}
\partial_\lambda \sigma(\mu_t^{\varepsilon,\lambda},u_t^\varepsilon)=&\widetilde{E}^Q\Big[\int_0^{\widetilde{U}_t^{\varepsilon,\lambda}}\partial_\mu \sigma(\mu_t^{\varepsilon,\lambda},u_t^\varepsilon;y)dy\cdot \partial_\lambda \widetilde{L}_t^{\varepsilon,\lambda}\Big]+\widetilde{E}^{Q}\big[\partial_\mu \sigma(\mu_t^{\varepsilon,\lambda},u_t^\varepsilon;\widetilde{U}_t^{\varepsilon,\lambda})\widetilde{L}_t^{\varepsilon,\lambda}\cdot\partial_\lambda\widetilde{U}_t^{\varepsilon,\lambda}\big]\\
=&\widetilde{E}^Q\Big[\int_0^{\widetilde{U}_t^{\varepsilon,\lambda}}\partial_\mu \sigma(\mu_t^{\varepsilon,\lambda},u_t^\varepsilon;y)dy(\widetilde{L}_t^\varepsilon-\widetilde{L}_t)\Big]+\widetilde{E}^{Q}\big[\partial_\mu \sigma(\mu_t^{\varepsilon,\lambda},u_t^\varepsilon;\widetilde{U}_t^{\varepsilon,\lambda})\widetilde{L}_t^{\varepsilon,\lambda}(\widetilde{U}_t^\varepsilon-\widetilde{U}_t)\big].\\
\end{split}
\end{equation*}
Thus, as $\displaystyle
\sigma(\mu_t^{\varepsilon},u_t^\varepsilon)-\sigma(\mu_t,u_t^\varepsilon)=\int_0^1 \partial_\lambda[\sigma(\mu_t^{\varepsilon,\lambda},u_t^\varepsilon)]d\lambda$,
\begin{equation}\label{jtri}
\begin{split}
&\sigma(\mu_t^{\varepsilon},u_t^\varepsilon)-\sigma(\mu_t,u_t^\varepsilon)-\Big\{\widetilde{E}^Q\Big[\int_0^{\widetilde{U}_t}\partial_\mu\sigma(\mu_t,u_t^\varepsilon;y)dy\cdot
\widetilde{K}_t^{1,\varepsilon}\Big]+\widetilde{E}^Q\big[\partial_\mu\sigma(\mu_t,u_t^\varepsilon;\widetilde{U}_t)\widetilde{L}_t\widetilde{V}_t^{1,\varepsilon}\big]\Big\}\\
=&\int_0^1\widetilde{E}^Q\Big[\int_0^{\widetilde{U}_t^{\varepsilon,\lambda}}\partial_\mu\sigma(\mu_t^{\varepsilon,\lambda},u_t^\varepsilon;y)dy(\widetilde{L}_t^\varepsilon-
\widetilde{L}_t)-\int_0^{\widetilde{U}_t}\partial_\mu\sigma(\mu_t,u_t^\varepsilon;y)dy\cdot\widetilde{K}_t^{1,\varepsilon}\Big]d\lambda\\
&+\int_0^1\widetilde{E}^Q\big[\partial_\mu \sigma(\mu_t^{\varepsilon,\lambda},u_t^\varepsilon;\widetilde{U}_t^{\varepsilon,\lambda})\widetilde{L}_t^{\varepsilon,\lambda}(\widetilde{U}_t^\varepsilon-\widetilde{U}_t)-
\partial_\mu\sigma(\mu_t,u_t^\varepsilon;\widetilde{U}_t)\widetilde{L}_t\widetilde{V}_t^{1,\varepsilon}\big]d\lambda\\
\end{split}
\end{equation}
\begin{equation*}
\begin{split}
=: &\int^1_0\big( I_1^{\varepsilon,\lambda}(t)+I_2^{\varepsilon,\lambda}(t)\big)d\lambda,\quad\quad\hskip9.5cm\
\end{split}
\end{equation*}
where $ \displaystyle I_1^{\varepsilon,\lambda}(t)=\widetilde{E}^Q\Big[\int_0^{\widetilde{U}_t}\partial_\mu\sigma(\mu_t,u_t^\varepsilon;y)dy\big(\widetilde{L}_t^\varepsilon-(\widetilde{L}_t+
\widetilde{K}_t^{1,\varepsilon})\big)\Big]+R_1^{\varepsilon,\lambda}(t),$
with
$$ R_1^{\varepsilon,\lambda}(t)=\widetilde{E}^Q\Big[\Big(\int_0^{\widetilde{U}_t^{\varepsilon,\lambda}}\partial_\mu\sigma(\mu_t^{\varepsilon,\lambda},u_t^\varepsilon;y)dy-
\int_0^{\widetilde{U}_t}\partial_\mu\sigma(\mu_t,u_t^\varepsilon;y)dy\Big)(\widetilde{L}_t^\varepsilon-\widetilde{L}_t)\Big].$$
Obviously,
\begin{equation*}
\begin{split}
&\partial_\mu\sigma(\mu_t^{\varepsilon,\lambda},u_t^\varepsilon;y)-\partial_\mu\sigma(\mu_t,u_t^\varepsilon;y)=\int_0^\lambda\partial_\rho\big[(\partial_\mu\sigma)(\mu_t^{\varepsilon,\rho},
u_t^\varepsilon;y)\big]d\rho\\
=&\int_0^\lambda\Big\{\widetilde{E}^Q\Big[\int_0^{\widetilde{U}_t^{\varepsilon,\rho}}\partial_\mu(\partial_\mu\sigma)(\mu_t^{\varepsilon,\rho},u_t^\varepsilon;y,y')dy'(
\widetilde{L}_t^\varepsilon-\widetilde{L}_t)\Big]\\
&\ \ \ \quad \ \ +\widetilde{E}^{Q}\big[\partial_\mu(\partial_\mu \sigma)(\mu_t^{\varepsilon,\rho},u_t^\varepsilon;y,\widetilde{U}_t^{\varepsilon,\rho})\widetilde{L}_t^{\varepsilon,\rho}(\widetilde{U}_t^\varepsilon-\widetilde{U}_t)\big]\Big\}d\rho,
\end{split}
\end{equation*}
and, as $\partial_\mu(\partial_\mu\sigma)$ is bounded,
\begin{equation}\label{j}
\big|\partial_\mu\sigma(\mu_t^{\varepsilon,\lambda},u_t^\varepsilon;y)-\partial_\mu\sigma(\mu_t,u_t^\varepsilon;y)\big|\leq C \int_0^\lambda \Big(\widetilde{E}^Q
\big[|\widetilde{U}_t^{\varepsilon,\rho}|\cdot|\widetilde{L}_t^\varepsilon-\widetilde{L}_t|+|\widetilde{L}_t^{\varepsilon,\rho}|\cdot|\widetilde{U}_t^\varepsilon
-\widetilde{U}_t|\big]\Big)d\rho.
\end{equation}
Recall that $U_t^{\varepsilon,\rho}=U_t+\rho(U_t^\varepsilon-U_t)$, and
$$ U_t^\varepsilon-U_t=\frac{E^{Q}[L_t^\varepsilon X_t^\varepsilon\,|\,\mathcal{F}_t^Y]}{E^{Q}[L_t^\varepsilon\,|\,\mathcal{F}_t^Y]}-\frac{E^{Q}[L_t X_t\,|\,\mathcal{F}_t^Y]}{E^{Q}[L_t\,|\,\mathcal{F}_t^Y]}. $$
From standard estimates using that for all $p'\geq1$, $X^\varepsilon,\, L^\varepsilon,\,(L^\varepsilon)^{-1},\,X,\,L,\,L^{-1}$ are bounded in $S^{p'}_{\mathbb{F}}([0,T],Q)$, uniformly with respect to $\varepsilon>0$, we deduce that, for all $p\geq1$:
\begin{equation*}
\begin{split}
&E^Q[|U_t^{\varepsilon,\rho}|^p]\leq C_p,\ t\in[0,T],\ \varepsilon>0,\ \rho\in[0,1];\\
&E^Q[|U_t^{\varepsilon}-U_t|^p]\leq C_p\big(E^Q[|X_t^{\varepsilon}-X_t|^{2p}+|L_t^{\varepsilon}-L_t|^{2p}]\big)^{\frac{1}{2}},\ t\in[0,T],\ p\geq 1.
\end{split}
\end{equation*}
Consequently, as due to (i) $E^Q[|L_t^{\varepsilon,\rho}|^{p}]\leq C_p,\ t\in[0,T],\ \varepsilon>0,\ \rho\in[0,1],\ p\geq1$, we get from \eqref{j} combined with (ii),
\begin{equation}\label{SMPa}
\big|\partial_\mu\sigma(\mu_t^{\varepsilon,\lambda},u_t^\varepsilon;y)-\partial_\mu\sigma(\mu_t,u_t^\varepsilon;y)\big|\leq C\big(E^Q[|X_t^{\varepsilon}-X_t|^{4}+|L_t^{\varepsilon}-L_t|^{4}]\big)^{\frac{1}{4}}\leq C\sqrt{\varepsilon},
\end{equation}
for $t\in[0,T],\ \lambda\in[0,1],\ \varepsilon>0$. Hence,
$$ |R_1^{\varepsilon,\lambda}(t)|\leq C\widetilde{E}^Q\big[|\widetilde{L}_t^\varepsilon-\widetilde{L}_t|\big(|\widetilde{U}_t^\varepsilon-\widetilde{U}_t|+|\widetilde{U}_t|\sqrt{\varepsilon}\big)\big]\leq C\varepsilon,\ t\in[0,T],\ \varepsilon>0,\ \lambda\in[0,1]. $$
For $I_2^{\varepsilon,\lambda}(t)$ we have that
\begin{equation*}
\begin{split}
I_2^{\varepsilon,\lambda}(t)&=\widetilde{E}^Q\big[\partial_\mu \sigma(\mu_t^{\varepsilon,\lambda},u_t^\varepsilon;\widetilde{U}_t^{\varepsilon,\lambda})\widetilde{L}_t^{\varepsilon,\lambda}(\widetilde{U}_t^\varepsilon-\widetilde{U}_t)-\partial_\mu
\sigma(\mu_t,u_t^\varepsilon;\widetilde{U}_t)\widetilde{L}_t\widetilde{V}_t^{1,\varepsilon}\big]\\
&=\widetilde{E}^Q\big[\partial_\mu\sigma(\mu_t,u_t^\varepsilon;\widetilde{U}_t)\widetilde{L}_t\big(\widetilde{U}_t^\varepsilon-(\widetilde{U}_t+\widetilde{V}_t^{1,\varepsilon})\big)\big]
+R_2^{\varepsilon,\lambda}(t),
\end{split}
\end{equation*}
with $R_2^{\varepsilon,\lambda}(t):=\widetilde{E}^Q\big[\big(\partial_\mu \sigma(\mu_t^{\varepsilon,\lambda},u_t^\varepsilon;\widetilde{U}_t^{\varepsilon,\lambda})\widetilde{L}_t^{\varepsilon,\lambda}-\partial_\mu\sigma(\mu_t,u_t^\varepsilon;\widetilde{U}_t)
\widetilde{L}_t\big)(\widetilde{U}_t^\varepsilon-\widetilde{U}_t)\big]$. Thus, using (ii) and \eqref{SMPa} as well as that $\partial_z(\partial_\mu\sigma)$ is bounded,
\begin{equation*}
\begin{split}
|R_2^{\varepsilon,\lambda}(t)|\leq &C\widetilde{E}^Q\big[\sqrt{\varepsilon}\widetilde{L}_t^{\varepsilon,\lambda}|\widetilde{U}_t^\varepsilon-\widetilde{U}_t|+
|\widetilde{L}_t^{\varepsilon,\lambda}-\widetilde{L}_t||
\widetilde{U}_t^\varepsilon-\widetilde{U}_t|+\widetilde{L}_t^{\varepsilon,\lambda}|\widetilde{U}_t^\varepsilon-\widetilde{U}_t||\widetilde{U}_t^{\varepsilon,\lambda}-\widetilde{U}_t|\big]\\
\leq &C\varepsilon,\ \varepsilon>0,\ t\in[0,T],\ \lambda\in[0,1].
\end{split}
\end{equation*}
Consequently, from \eqref{jj} and \eqref{jtri}, for $R^{\varepsilon,\lambda}(s)=R_1^{\varepsilon,\lambda}(s)+R_2^{\varepsilon,\lambda}(s)$,
\begin{equation}\label{SMPb}
\begin{split}
X_t^\varepsilon-(X_t&+Y_t^{1,\varepsilon})=\int_0^t\Big\{\widetilde{E}^Q\Big[\int_0^{\widetilde{U}_s}\partial_\mu\sigma(\mu_s,u_s^\varepsilon;y)dy\big(\widetilde{L}_s^\varepsilon-
(\widetilde{L}_s+\widetilde{K}_s^{1,\varepsilon})\big)\Big]\\
&\qquad\quad+\widetilde{E}^Q\big[\partial_\mu\sigma(\mu_s,u_s^\varepsilon;\widetilde{U}_s)\widetilde{L}_s\big(\widetilde{U}_s^\varepsilon-(\widetilde{U}_s+\widetilde{V}_s^{1,\varepsilon})
\big)\big]\Big\}dB_s^1+\int_0^t\int_0^1 R^{\varepsilon,\lambda}(s)d\lambda dB_s^1\\
&+\int_0^t\Big(\widetilde{E}^Q\Big[\int_0^{\widetilde{U}_s}\delta(\partial_\mu\sigma)(s,y)dy\cdot\widetilde{K}_s^{1,\varepsilon}\Big]+\widetilde{E}^Q\big[\delta(\partial_\mu\sigma)(s,
\widetilde{U}_s)\widetilde{L}_s\widetilde{V}_s^{1,\varepsilon}\big]\Big)\mathbf{1}_{E_\varepsilon}(s)dB_s^1,\\
\end{split}
\end{equation}
where $\delta(\partial_\mu\sigma)(s,y):=\partial_\mu\sigma(\mu_s,u_s^\varepsilon;y)-\partial_\mu\sigma(\mu_s,u_s;y)$.\\
From (1) of the proof of Proposition \ref{Wellof1stVar} combined with (iii) we have for $V_t^{1,\varepsilon}=\theta_t(Y_t^{1,\varepsilon},K_t^{1,\varepsilon})$,
$$ \big(E^Q[|V_t^{1,\varepsilon}|^2]\big)^{\frac{1}{2}}\leq C\big(E^Q[|Y_t^{1,\varepsilon}|^4+|K_t^{1,\varepsilon}|^4]\big)^{\frac{1}{4}}\leq C\sqrt{\varepsilon}. $$
Hence, $\widetilde{E}^Q\big[\big|\delta(\partial_\mu\sigma)(s,\widetilde{U}_s)\widetilde{L}_s\widetilde{V}_s^{1,\varepsilon}\big|\big]\leq C\sqrt{\varepsilon}$.\ On the other hand, again from (iii) we have
$$ \widetilde{E}^Q\Big[\Big|\int_0^{\widetilde{U}_s}\delta(\partial_\mu\sigma)(s,y)dy\cdot\widetilde{K}_s^{1,\varepsilon}\Big|\Big]\leq C\sqrt{\varepsilon}. $$
Consequently, combining with \eqref{SMPb} we see that
\begin{equation}\label{SMPc}
\begin{split}
X_t^\varepsilon-(X_t+&Y_t^{1,\varepsilon})=\int_0^t\Big\{\widetilde{E}^Q\Big[\int_0^{\widetilde{U}_s}\partial_\mu\sigma(\mu_s,u_s^\varepsilon;y)dy\big(\widetilde{L}_s^\varepsilon-
(\widetilde{L}_s+\widetilde{K}_s^{1,\varepsilon})\big)\Big]\\
&+\widetilde{E}^Q\big[\partial_\mu\sigma(\mu_s,u_s^\varepsilon;\widetilde{U}_s)\widetilde{L}_s\big(\widetilde{U}_s^\varepsilon-(\widetilde{U}_s+\widetilde{V}_s^{1,\varepsilon})
\big)\big]+\int_0^1 R^{\varepsilon,\lambda}(s)d\lambda\Big\}dB_s^1,\ t\in[0,T],\\
\end{split}
\end{equation}
where $|R^{\varepsilon,\lambda}(s)|\leq C\big(\varepsilon+\sqrt{\varepsilon}\mathbf{1}_{E_\varepsilon}(s)\big)$, $s\in[0,T]$, $\lambda\in[0,1]$.

Let us now estimate $L_t^\varepsilon-(L_t+K_t^{1,\varepsilon})$. As $h,\,\partial_x h,\,\partial^2_{xx} h$ are bounded,
$$ h(X_t^\varepsilon,u_t^\varepsilon)L_t^\varepsilon-h(X_t,u_t^\varepsilon)L_t=h(X_t,u_t^\varepsilon)(L_t^\varepsilon-L_t)+L_t\partial_x h(X_t,u_t^\varepsilon)(X_t^\varepsilon-X_t)+R_t^\varepsilon, $$
where, thanks to (ii) $\big(E^Q[|R_t^\varepsilon|^p]\big)^{\frac{1}{p}}\leq C_p\varepsilon,\ \varepsilon>0,\ t\in[0,T],\ p\geq1$. Then, using \eqref{jj},
\begin{equation}\label{SMPd}
\begin{split}
L_t^\varepsilon-(L_t+K_t^{1,\varepsilon})=&\int_0^t\Big\{h(X_s,u_s^\varepsilon)\big(L_s^\varepsilon-(L_s+K_s^{1,\varepsilon})\big)+L_s\partial_x h(X_s,u_s^\varepsilon)\big(X_s^\varepsilon-(X_s+Y_s^{1,\varepsilon})\big)\\
&\qquad+R_s^\varepsilon+\big(\delta h(s)K_s^{1,\varepsilon}+L_s\delta(\partial_x h)(s)Y_s^{1,\varepsilon}\big)\mathbf{1}_{E_\varepsilon}(s)\Big\}dY_s,\ t\in[0,T].\\
\end{split}
\end{equation}
Obviously, for $\widetilde{R}_s^\varepsilon=R_s^\varepsilon+\big(\delta h(s)K_s^{1,\varepsilon}+L_s\delta(\partial_x h)(s)Y_s^{1,\varepsilon}\big)\mathbf{1}_{E_\varepsilon}(s)$,
$$ \big(E^Q[|\widetilde{R}_s^\varepsilon|^p]\big)^{\frac{1}{p}}\leq C_p\big(\varepsilon+\sqrt{\varepsilon}\mathbf{1}_{E_\varepsilon}(s)\big),\ s\in[0,T],\ \varepsilon>0,\ p\geq2. $$
Now we need to compute $U_t^\varepsilon-(U_t+V_t^{1,\varepsilon})$. A straightly-forward computations yield
\begin{equation*}
\begin{split}
&U_t^\varepsilon-U_t=\frac{E^{Q}[L_t^\varepsilon X_t^\varepsilon\,|\,\mathcal{F}_t^Y]}{E^{Q}[L_t^\varepsilon\,|\,\mathcal{F}_t^Y]}-\frac{E^{Q}[L_t X_t\,|\,\mathcal{F}_t^Y]}{E^{Q}[L_t\,|\,\mathcal{F}_t^Y]}=\int_0^1\partial_\lambda\Big(\frac{E^{Q}[L_t^{\varepsilon,\lambda} X_t^{\varepsilon,\lambda}\,|\,\mathcal{F}_t^Y]}{E^{Q}[L_t^{\varepsilon,\lambda}\,|\,\mathcal{F}_t^Y]}\Big)d\lambda\\
&=\int_0^1\Big(\frac{E^{Q}[L_t^{\varepsilon,\lambda} (X_t^\varepsilon-X_t)+X_t^{\varepsilon,\lambda}(L_t^\varepsilon-L_t)\,|\,\mathcal{F}_t^Y]}{E^{Q}[L_t^{\varepsilon,\lambda}\,|\,\mathcal{F}_t^Y]}-\frac{E^{Q}[L_t^{\varepsilon,\lambda} X_t^{\varepsilon,\lambda}\,|\,\mathcal{F}_t^Y]E^{Q}[L_t^\varepsilon-L_t\,|\,\mathcal{F}_t^Y]}{(E^{Q}[L_t^{\varepsilon,\lambda}\,|\,\mathcal{F}_t^Y])^2}\Big)d\lambda\\
&=\frac{E^{Q}[L_t(X_t^\varepsilon-X_t)+X_t(L_t^\varepsilon-L_t)\,|\,\mathcal{F}_t^Y]}{E^{Q}[L_t\,|\,\mathcal{F}_t^Y]}-\frac{E^{Q}[L_t X_t\,|\,\mathcal{F}_t^Y]}{(E^{Q}[L_t\,|\,\mathcal{F}_t^Y])^2}E^{Q}[L_t^\varepsilon-L_t\,|\,\mathcal{F}_t^Y]+\vartheta_t^\varepsilon,
\end{split}
\end{equation*}
with $\big(E^Q[|\vartheta_t^\varepsilon|^p]\big)^{\frac{1}{p}}\leq C_p\varepsilon,\ \varepsilon>0,\ t\in[0,T],\ p\geq1$. Consequently, recalling the definition of $\theta_t$ (\ref{theta}),
\begin{equation}\label{4.14-1}
\begin{split}
U_t^\varepsilon-(U_t+V_t^{1,\varepsilon})=&\frac{E^{Q}[L_t\big(X_t^\varepsilon-(X_t+Y_t^{1,\varepsilon})\big)+X_t\big(L_t^\varepsilon-(L_t+K_t^{1,\varepsilon})\big)\,|\,
\mathcal{F}_t^Y]}{E^{Q}[L_t\,|\,\mathcal{F}_t^Y]}\\
&-\frac{E^{Q}[L_t X_t\,|\,\mathcal{F}_t^Y]}{(E^{Q}[L_t\,|\,\mathcal{F}_t^Y])^2}E^{Q}[L_t^\varepsilon-(L_t+K_t^{1,\varepsilon})\,|\,\mathcal{F}_t^Y]+\vartheta_t^\varepsilon\\
=&\theta_t\big(X_t^\varepsilon-(X_t+Y_t^{1,\varepsilon}),L_t^\varepsilon-(L_t+K_t^{1,\varepsilon})\big)+\vartheta_t^\varepsilon,
\end{split}
\end{equation}
and from (1) in the proof of Proposition \ref{Wellof1stVar}, for $1\leq p<q<\infty$,
\begin{equation*}
\begin{split}
\Big(E^Q\big[\big|U_t^\varepsilon-(U_t+V_t^{1,\varepsilon})\big|^p\big]\Big)^{\frac{1}{p}}&\leq\Big(E^Q\big[\big|\theta_t\big(X_t^\varepsilon-(X_t+Y_t^{1,\varepsilon}),L_t^\varepsilon
-(L_t+K_t^{1,\varepsilon})\big)\big|^p\big]\Big)^{\frac{1}{p}}+\big(E^Q[|\vartheta_t^\varepsilon|^p]\big)^{\frac{1}{p}}\\
&\leq C_{p,q}\Big(E^Q\big[\big|X_t^\varepsilon-(X_t+Y_t^{1,\varepsilon})\big|^q+\big|L_t^\varepsilon-(L_t+K_t^{1,\varepsilon})\big|^q\big]\Big)^{\frac{1}{q}}+C_p\varepsilon.
\end{split}
\end{equation*}
Hence, for all $p\geq2$,
$$ \widetilde{E}^Q\big[\widetilde{L}_t\big|\widetilde{U}_t^\varepsilon-(\widetilde{U}_t+\widetilde{V}_t^{1,\varepsilon})\big|\big]\leq C_p\Big(E^Q\big[\big|X_t^\varepsilon-(X_t+Y_t^{1,\varepsilon})\big|^p+\big|L_t^\varepsilon-(L_t+K_t^{1,\varepsilon})\big|^p\big]\Big)^{\frac{1}{p}}+C_p\varepsilon,\ t\in[0,T], $$
and thanks to \eqref{SMPc} and \eqref{SMPd}, using Gronwall's inequality we get
$$ E^Q\big[\sup_{t\in[0,T]}\big(|X_t^\varepsilon-(X_t+Y_t^{1,\varepsilon})|^{p}+|L_t^\varepsilon-(L_t+K_t^{1,\varepsilon})|^{p}\big)\big]\leq C_p \varepsilon^{p},\ \varepsilon>0,\ p\geq2. $$
Furthermore, the above estimate allows to deduce from (\ref{4.14-1}) (with $\mathcal{F}_t^Y$ which can be replaced by $\mathcal{F}_T^Y$) by using standard methods that $\displaystyle E^Q\big[\sup_{t\in[0,T]}|U_t^\varepsilon-(U_t+V_t^{1,\varepsilon})|^{p}\big]\leq C_p \varepsilon^{p},\ \varepsilon>0,\ p\geq2. $
\end{proof}

In the proof of the above proposition we also have shown the following important estimates.
\begin{corollary}\label{EstofUV}
For all $k\geq1$, there exists a constant $C_k>0$ such that,
\begin{equation*}
\begin{array}{lll}
&\mbox{ }\qquad\qquad\quad{\rm (i)}\  E^Q\big[\sup_{t\in[0,T]}|U_t^\varepsilon|^{2k}\big]\leq C_k;\\
&\mbox{ }\qquad\qquad\quad{\rm (ii)}\ E^Q\big[\sup_{t\in[0,T]}|U_t^\varepsilon-U_t|^{2k}\big]\leq C_k \varepsilon^{k},\ \varepsilon>0;\\
&\mbox{ }\qquad\qquad\quad{\rm (iii)}\ E^Q\big[\sup_{t\in[0,T]}|V_t^{1,\varepsilon}|^{2k}\big]\leq C_k \varepsilon^{k},\ \varepsilon>0;\\
&\mbox{ }\qquad\qquad\quad{\rm (iv)}\ E^Q\big[\sup_{t\in[0,T]}|U_t^\varepsilon-(U_t+V_t^{1,\varepsilon})|^{2k}\big]\leq C_k \varepsilon^{2k},\ \varepsilon>0. \hspace{8cm}
\end{array}
\end{equation*}
\end{corollary}
Now we present a very subtle and useful estimate, whose proof applies an idea which was first introduced in \cite{BCL}.
\begin{proposition}\label{tech}
For all $\theta=(\theta^1,\theta^2)\in L^2_{\mathbb{F}}([0,T],Q;\mathbb{R}^2)$ %(=L^2_{\mathbb{F}}([0,T]\times\Omega,dsdQ,\mathbb{R}^2))$
with $\displaystyle E^Q\Big[\int_0^T\big(|\theta^1_t|^2+|L_t\theta^2_t|^2\big)dt\Big]<+\infty$, and $(\theta^1_t,L_t\theta^2_t)\in L^2(\mathcal{F}_t,Q;\mathbb{R}^2)$ for all $t\in[0,T]$, there exists $\rho:[0,T]\times\mathbb{R}_+\rightarrow\mathbb{R}_+$ such that
$$ \big|E^Q[\theta^1_t Y_t^{1,\varepsilon}+\theta^2_t K_t^{1,\varepsilon}]\big|\leq \rho_t(\varepsilon)\sqrt{\varepsilon},\ \varepsilon\in(0,1],\ t\in[0,T], $$
with $\rho_t(\varepsilon)\rightarrow0\, (\varepsilon\searrow0)$, $t\in[0,T]$, and\ $\displaystyle \rho_t(\varepsilon)\leq CE^Q\big[|\theta^1_t|^2+|L_t\theta^2_t|^2\big],\ \varepsilon\in(0,1],\ t\in[0,T]. $
\end{proposition}
\begin{proof}
For simplicity of the arguments let us keep the same simplified setting as in the proof of Proposition \ref{EstofXL}. The same arguments--for the price of much longer computations--can be used for the general case. Recall from \eqref{jj} that
\begin{equation*}
\begin{split}
&dK_t^{1,\varepsilon}=\big(h(t)K_t^{1,\varepsilon}+L_t h_x(t)Y_t^{1,\varepsilon}+L_t \delta h(t)\mathbf{1}_{E_\varepsilon}(t)\big)dY_t,\ K_0^{1,\varepsilon}=0,\\
\end{split}
\end{equation*}
and notice that $$ dL^{-1}_t=-h(t)L^{-1}_t dY_t+|h(t)|^2 L^{-1}_tdt,\ L^{-1}_0=1. $$
Then we define $N_t^{1,\varepsilon}:=L^{-1}_t K_t^{1,\varepsilon}$, and from It\^{o}'s formula we have
$$ dN_t^{1,\varepsilon}=\big(h_x(t)Y_t^{1,\varepsilon}+\delta h(t)\mathbf{1}_{E_\varepsilon}(t)\big)dY_t-h(t)\big(h_x(t)Y_t^{1,\varepsilon}+\delta h(t)\mathbf{1}_{E_\varepsilon}(t)\big)dt. $$
Recall also from \eqref{jj} that
$$ dY_t^{1,\varepsilon}=\Big\{\widetilde{E}^Q\Big[\int_0^{\widetilde{U}_t}\partial_\mu\sigma(\mu_t,u_t;y)dy\cdot\widetilde{K}_t^{1,\varepsilon}\Big]+\widetilde{E}^Q\big[\partial_\mu\sigma
(\mu_t,u_t;\widetilde{U}_t)\widetilde{L}_t\widetilde{V}_t^{1,\varepsilon}\big]+\delta\sigma(t)\mathbf{1}_{E_\varepsilon}(t)\Big\}dB_t^1.\\
$$
Moreover, for $(\zeta, \eta)\in L^2(\Omega, {\mathcal{F}}, Q; \mathbb{R}^2)$ let us put
\begin{equation}\label{sqr}
\begin{split}
\theta_t^*(\zeta, \eta):=&\frac{E^Q[L_t(\zeta+X_t \eta)\,|\,\mathcal{F}_t^Y]}{E^Q[L_t\,|\,\mathcal{F}_t^Y]}-\frac{E^{Q}[L_t X_t\,|\,\mathcal{F}_t^Y] }
{E^{Q}[L_t\,|\,\mathcal{F}_t^Y]}\frac{E^{Q}[L_t \eta\,|\,\mathcal{F}_t^Y]}{E^{Q}[L_t\,|\,\mathcal{F}_t^Y]}\\
=&E^P[\zeta+X_t \eta\,|\,\mathcal{F}_t^Y]-E^{P}[X_t\,|\,\mathcal{F}_t^Y]E^{P}[\eta\,|\,\mathcal{F}_t^Y].
\end{split}
\end{equation}
Then, obviously we have $\theta^*_t (Y_t^{1,\varepsilon},N_t^{1,\varepsilon})= \theta_t(Y_t^{1,\varepsilon},K_t^{1,\varepsilon})= V_t^{1,\varepsilon}.$\\
Furthermore, setting now
\begin{equation*}
    \begin{split}
  &\alpha_t^1(u):=L_t\int_0^{U_t}\partial_\mu\sigma(\mu_s,u;y)dy,\ \alpha_t^2(u):=L_t\cdot\partial_\mu\sigma(\mu_t,u;U_t),\ u\in U,\ t\in[0,T]; \\
  &\widetilde{\alpha}_t^1(u):=\widetilde{L}_t\int_0^{\widetilde{U}_t}\partial_\mu\sigma(\mu_s,u;y)dy,\ \widetilde{\alpha}_t^2(u):=\widetilde{L}_t\cdot\partial_\mu\sigma(\mu_t,u;\widetilde{U}_t),\ u\in U,\ t\in[0,T].
\end{split}
\end{equation*}
Inspired by \cite{BCL} we consider the Hilbert space $H:=L^2_{\mathbb{F}}([0,T],Q;\mathbb{R}^2)$, and for $Z=(Z^1,Z^2)\in H$, we define the mappings
\begin{equation*}
\begin{split}
&D_1[Z](t):=\int_0^t\widetilde{E}^Q\big[\widetilde{\alpha}_s^1(u_s)\widetilde{Z}_s^2+\widetilde{\alpha}_s^2(u_s)\widetilde{\theta_s^*(Z_s^1,Z_s^2)}\big]dB_s^1,\\
&D_2[Z](t):=\int_0^t h_x(s)Z_s^1 dY_s-\int_0^t h(s)h_x(s)Z_s^1 ds,\quad D[Z]:=(D_1[Z],D_2[Z]).
\end{split}
\end{equation*}
Point (i) of the proof of Proposition \ref{Wellof1stVar} allows to get that $E^Q\big[(1+L_s)\big|\theta_s^*(Z_s^1,Z_s^2)\big|\big]\leq C\big(E^Q[|Z_s|^2]\big)^{\frac{1}{2}}$, $Z\in H$. Hence, using standard arguments we see that $D\in L(H)$ is a bounded and linear operator over $H$, and
$$ E^Q\big[\big|D[Z](t)\big|^2\big]\leq C\int_0^t E^Q[|Z_s|^2]ds,\ t\in[0,T],\ Z\in H. $$
Consequently, by iteration,
$$ E^Q\big[\big|D^N[Z](t)\big|^2\big]\leq \frac{(Ct)^{N-1}}{(N-1)!}E^Q\Big[\int_0^T|Z_s|^2ds\Big]=\frac{(Ct)^{N-1}}{(N-1)!}|Z|_H^2, $$
with $\displaystyle |Z|_H^2=E^Q\Big[\int_0^T|Z_s|^2ds\Big]$,
and, thus, $\displaystyle \big|D^N[Z]\big|_H^2\leq \frac{(CT)^N}{N!}|Z|_H^2,\ Z\in H,\ N\geq1, $\ i.e., for the operator norm in $L(H)$ we have $||D^N||^2\leq\frac{(CT)^N}{N!},\ N\geq 1$. Consequently, $\displaystyle A:=\sum_{n\geq0}D^n\in L(H)$ converges in $L(H)$, and the adjoint operator $A^*$ of $A$ satisfies
$$ \big(A^*[Z'],Z\big)_H=\big(Z',AZ\big)_H,\ Z,Z'\in H. $$
Let us point out that, for all $\theta\in H$, $t\in[0,T]$,
\begin{equation*}
\begin{split}
& A^*[\theta\mathbf{1}_{[0,t]}]=A^*[\theta\mathbf{1}_{[0,t]}]\mathbf{1}_{[0,t]},\ \ \mbox{i.e.,} \ \ A^*[\theta\mathbf{1}_{[0,t]}](s)=0,\ s\in[t,T],\ \mbox{dsdQ-a.e.}
\end{split}
\end{equation*}
Indeed, let $Z\in H$. Then, taking into account the definition of $D[Z]$, and the fact that $\theta_s^*(0,0)=0$, $s\in[0,T]$, we see that $D[Z](t\wedge s)=D[Z\mathbf{1}_{[0,t]}](s)$, and, thus, also
$$ D[Z]\mathbf{1}_{[0,t]}=D[Z](t\wedge\cdot)\mathbf{1}_{[0,t]}=D[Z\mathbf{1}_{[0,t]}]\mathbf{1}_{[0,t]}. $$
Hence, for all \ $Z\in H$,\ \  $\displaystyle (\theta\mathbf{1}_{[0,t]},D[Z])_H=(\theta\mathbf{1}_{[0,t]},D[Z\mathbf{1}_{[0,t]}])_H, $\ \ and
\begin{equation*}
\begin{split}
&(\theta\mathbf{1}_{[0,t]},D^2[Z])_H=(\theta\mathbf{1}_{[0,t]},D(D[Z]\mathbf{1}_{[0,t]}))_H=(\theta\mathbf{1}_{[0,t]},D(D[Z\mathbf{1}_{[0,t]}]\mathbf{1}_{[0,t]}))_H\\
=&(\theta\mathbf{1}_{[0,t]},D(D[Z\mathbf{1}_{[0,t]}]))_H=(\theta\mathbf{1}_{[0,t]},D^2[Z\mathbf{1}_{[0,t]}])_H.
\end{split}
\end{equation*}
Iterating this argument we have $\displaystyle (\theta\mathbf{1}_{[0,t]},D^n[Z])_H=(\theta\mathbf{1}_{[0,t]},D^n[Z\mathbf{1}_{[0,t]}])_H,\ n\geq0, $\
and, consequently,
$$ (\theta\mathbf{1}_{[0,t]},A[Z])_H=(\theta\mathbf{1}_{[0,t]},A[Z\mathbf{1}_{[0,t]}]\mathbf{1}_{[0,t]})_H, $$
and \  $\displaystyle (A^*(\theta\mathbf{1}_{[0,t]}),Z)_H=(A^*(\theta\mathbf{1}_{[0,t]})\mathbf{1}_{[0,t]},Z)_H,\ \mbox{for all } Z\in H. $\ \
This proves
$$ A^*(\theta\mathbf{1}_{[0,t]})=A^*(\theta\mathbf{1}_{[0,t]})\mathbf{1}_{[0,t]}. $$
Let us define
$$ Z^{1,\varepsilon}:=\Big(\begin{array}{c}Y^{1,\varepsilon}\\ N^{1,\varepsilon}\end{array}\Big),\ \ \ \ \ G_t^\varepsilon:=\int_0^t\Big(\begin{array}{c}\delta\sigma(s)\mathbf{1}_{E_\varepsilon}(s)dB_s^1\\ \delta h(s)\mathbf{1}_{E_\varepsilon}(s)dY_s-h(s)\delta h(s)\mathbf{1}_{E_\varepsilon}(s)ds\end{array}\Big). $$
Then, from the definition of $D$,
$$ Z^{1,\varepsilon}=D[Z^{1,\varepsilon}]+G^\varepsilon=G^\varepsilon+D\big[G^\varepsilon+D[Z^{1,\varepsilon}]\big]=...=\sum_{j=0}^N D^j[G^\varepsilon]+D^{N+1}[Z^{1,\varepsilon}],\ N\geq1. $$
Consequently, as $\big|D^{N+1}[Z^{1,\varepsilon}]\big|_H\leq ||D^{N+1}||\cdot|Z^{1,\varepsilon}|_H\rightarrow0\ (N\rightarrow\infty)$, and $\displaystyle \sum_{j=0}^N D^j\xrightarrow[]{L(H)}A$ but also $A=I+D\circ A$, we have
$$ Z^{1,\varepsilon}=A[G^\varepsilon]=G^\varepsilon+D\big[A[G^\varepsilon]\big]. $$
Given any $\theta=(\theta^1,\theta^2)\in H$, we want to compute $E[\theta_t Z_t^{1,\varepsilon}]=E\big[\theta_t\big(G_t^\varepsilon+D\big[A[G^\varepsilon]\big]_t\big)\big]$. For this notice that,\\
1) In analogy to \eqref{sqr} we see
$$ \theta_s^*(Z_s^1,Z_s^2)=E^P[Z_s^1+X_s Z_s^2\,|\,\mathcal{F}_s^Y]-E^{P}[X_s\,|\,\mathcal{F}_s^Y]E^{P}[Z_s^2\,|\,\mathcal{F}_s^Y]. $$
Thus, for all $Z=(Z^1,Z^2)\in H:$
\begin{equation*}
\begin{split}
&E^Q\big[\alpha_s^2(\widetilde{u}_s)\theta_s^*(Z_s^1,Z_s^2)\big]=E^P\big[\partial_\mu\sigma(\mu_s,\widetilde{u}_s;U_s)\big(E^P[Z_s^1+X_s Z_s^2\,|\,\mathcal{F}_s^Y]-E^P[X_s\,|\,\mathcal{F}_s^Y]E^P[Z_s^2\,|\,\mathcal{F}_s^Y]\big)\big]\\
=&E^P\big[E^P[\partial_\mu\sigma(\mu_s,\widetilde{u}_s;U_s)\,|\,\mathcal{F}_s^Y]\big(Z_s^1+(X_s-E^P[X_s\,|\,\mathcal{F}_s^Y])Z_s^2\big)\big],
\end{split}
\end{equation*}
i.e., $\displaystyle D_1[Z](t)=\int_0^t\widetilde{E}^Q\big[\widetilde{\alpha}_s^1(u_s)\widetilde{Z}_s^2+\widetilde{\alpha}_s^2(u_s)\widetilde{\theta_s^*(Z_s^1,Z_s^2)}\big]dB_s^1
=\int_0^t\widetilde{E}^Q\big[\widetilde{\beta}_s^1(u_s)\widetilde{Z}_s^1+\widetilde{\beta}_s^2(u_s)\widetilde{Z}_s^2\big]dB_s^1,
$\
with $\displaystyle  \beta_s^1(v):=E^P[\partial_\mu\sigma(\mu_s,v;U_s)\,|\,\mathcal{F}_s^Y]L_s$, $\displaystyle \beta_s^2(v):=\alpha_s^1(v)+E^P[\partial_\mu\sigma(\mu_s,v;U_s)\,|\,\mathcal{F}_s^Y]L_s(X_s-E^P[X_s\,|\,\mathcal{F}_s^Y]).$ Hence,
$$ D_1[Z](t)=\int_0^t\widetilde{E}^Q\big[\widetilde{\beta}_s^1(u_s)\widetilde{Z}_s^1+\widetilde{\beta}_s^2(u_s)\widetilde{Z}_s^2\big]dB_s^1. $$
2) Now, for our given $\theta=(\theta^1,\theta^2)$, as $\theta_t^i\in L^2(\mathcal{F}_t,Q)$, $t\in[0,T]$, $i=1,2$, there are jointly measurable random fields $\theta^{i,j}=\{\theta_{t,s}^{i,j},\,0\leq s\leq t\leq T\}$ with $\theta_{t,\cdot}^{i,j}\in L^2_{\mathbb{F}}([0,t],Q),\ t\in[0,T]$, such that, $dtdQ$-a.e.,
$$ \theta_t^i=E^Q[\theta_t^i]+\int_0^t \theta_{t,s}^{i,1}dB_s^1+\int_0^t \theta_{t,s}^{i,2}dY_s. $$
Then,
\begin{equation}\label{4.15-1}
\begin{split}
&E^Q\big[\theta_t^1 D_1[Z](t)\big]=E^Q\Big[\int_0^t\theta_{t,s}^{1,1}\widetilde{E}^Q\big[\widetilde{\beta}_s^1(u_s)\widetilde{Z}_s^1+\widetilde{\beta}_s^2(u_s)\widetilde{Z}_s^2\big]ds\Big]\\
=&E^Q\Big[\int_0^T\Big(\widetilde{E}^Q\big[\widetilde{\theta}_{t,s}^{1,1}\beta_s^1(\widetilde{u}_s)\big]Z_s^1+\widetilde{E}^Q\big[\widetilde{\theta}_{t,s}^{1,1}\beta_s^2(\widetilde{u}_s)
\big]Z_s^2\Big)\mathbf{1}_{[0,t]}(s)ds\Big]=(\gamma_t^\theta,Z)_H,
\end{split}
\end{equation}
for $\gamma_{t,\cdot}^\theta=(\gamma_{t,\cdot}^{\theta,1},\gamma_{t,\cdot}^{\theta,2})$ with $\gamma_{t,s}^{\theta,1}:=\widetilde{E}^Q\big[\widetilde{\theta}_{t,s}^{1,1}\beta_s^1(\widetilde{u}_s)\big]\mathbf{1}_{[0,t]}(s)$, $\gamma_{t,s}^{\theta,2}:=\widetilde{E}^Q\big[\widetilde{\theta}_{t,s}^{1,1}\beta_s^2(\widetilde{u}_s)\big]\mathbf{1}_{[0,t]}(s)$. Note that
\begin{equation*}
\begin{split}
&|\gamma_{t,s}^{\theta,1}|\leq CE^Q[|\theta_{t,s}^{1,1}|]L_s\mathbf{1}_{[0,t]}(s),\\
&|\gamma_{t,s}^{\theta,2}|\leq CE^Q[|\theta_{t,s}^{1,1}|]L_s(|X_s|+E^{P}[|X_s|\,|\,\mathcal{F}_s^Y]+|U_s|)\mathbf{1}_{[0,t]}(s),
\end{split}
\end{equation*}
and we obtain
$$ E^Q\Big[\int_0^T|\gamma_{t,s}^{\theta,i}|^2 ds\Big]\leq CE^Q\Big[\int_0^t|\theta_{t,s}^{1,1}|^2 ds\Big]\leq CE^Q[|\theta_t^1|^2],\ i=1,2,\ t\in[0,T]. $$
3) We define now $\widehat{\gamma}_{t,s}^\theta=(\widehat{\gamma}_{t,s}^{\theta,1},\widehat{\gamma}_{t,s}^{\theta,2})$, $0\leq s\leq t\leq T$, by putting $\widehat{\gamma}_{t,\cdot}^\theta
=A^*[\gamma_{t,\cdot}^\theta]$. Note that $\widehat{\gamma}_{t,\cdot}^\theta=A^*[\gamma_{t,\cdot}^\theta\mathbf{1}_{[0,t]}]=A^*[\gamma_{t,\cdot}^\theta\mathbf{1}_{[0,t]}]\mathbf{1}_{[0,t]}=
\widehat{\gamma}_{t,\cdot}^\theta
\mathbf{1}_{[0,t]}$. Obviously, for $t\in[0,T]$,
$$ E^Q\Big[\int_0^t|\widehat{\gamma}_{t,s}^{\theta}|^2 ds\Big]=E^Q\Big[\int_0^T|\widehat{\gamma}_{t,s}^{\theta}|^2 ds\Big]\leq CE^Q\Big[\int_0^T|\gamma_{t,s}^{\theta}|^2 ds\Big]=CE^Q\Big[\int_0^t|\gamma_{t,s}^{\theta}|^2 ds\Big]\leq CE^Q[|\theta_t^1|^2]. $$
Let now $\gamma^{\theta,i,j}=\{\gamma_{t,s,r}^{\theta,i,j},\, 0\leq r\leq s\leq t\leq T\}$ be the jointly measurable field with $\gamma_{t,s,\cdot}^{\theta,i,j}\in L^2_{\mathbb{F}}([0,s],Q)$,
$0\leq s\leq t\leq T$, $j=1,2$, such that $\displaystyle \widehat{\gamma}_{t,s}^{\theta,i}=E^Q[\widehat{\gamma}_{t,s}^{\theta,i}]+\int_0^s \widehat{\gamma}_{t,s,r}^{\theta,i,1}dB_s^1+\int_0^s
\widehat{\gamma}_{t,s,r}^{\theta,i,2}dY_s,$ $dsdtdQ$-a.e., $i=1,2$. Remark that $\widehat{\gamma}_{t,\cdot}^\theta$ does not depend on $\varepsilon>0$. Then, from \eqref{4.15-1},
\begin{equation*}
\begin{split}
&E^Q\big[\theta_t^1 D_1\big[A[G^\varepsilon]\big](t)\big]=\big(\gamma_t^\theta,A[G^\varepsilon]\big)_H=\big(A^*[\gamma_t^\theta],G^\varepsilon\big)_H=\big(\widehat{\gamma}_t^\theta,G^\varepsilon\big)_H=
E^Q\Big[\int_0^t\widehat{\gamma}_{t,s}^{\theta}G_s^\varepsilon ds\Big]\\
=&\int_0^t E^Q\Big[\widehat{\gamma}_{t,s}^{\theta,1}\int_0^s\delta\sigma(r)\mathbf{1}_{E_\varepsilon}(r)dB_r^1 +\widehat{\gamma}_{t,s}^{\theta,2}\Big(\int_0^s\delta h(r)\mathbf{1}_{E_\varepsilon}(r)dY_r-\int_0^s h(r)\delta h(r)\mathbf{1}_{E_\varepsilon}(r)dr\Big)\Big]ds\\
=&\int_0^t E^Q\Big[\int_0^s\widehat{\gamma}_{t,s,r}^{\theta,1,1}\delta \sigma(r)\mathbf{1}_{E_\varepsilon}(r)dr+\int_0^s\widehat{\gamma}_{t,s,r}^{\theta,2,2}\delta h(r)\mathbf{1}_{E_\varepsilon}(r)dr-\widehat{\gamma}_{t,s}^{\theta,2}\int_0^s h(r)\delta h(r)\mathbf{1}_{E_\varepsilon}(r)dr\Big]ds\\
=&E^Q\Big[\int_0^t\int_0^s\big(\widehat{\gamma}_{t,s,r}^{\theta,1,1}\delta \sigma(r)+\widehat{\gamma}_{t,s,r}^{\theta,2,2}\delta h(r)\big)\mathbf{1}_{E_\varepsilon}(r)drds\Big]-E^Q\Big[\int_0^t\widehat{\gamma}_{t,s}^{\theta,2}\Big(\int_0^s h(r)\delta h(r)\mathbf{1}_{E_\varepsilon}(r)dr\Big)ds\Big].
\end{split}
\end{equation*}
Hence,
\begin{equation*}
\begin{split}
&\big|E^Q\big[\theta_t^1 D_1\big[A[G^\varepsilon]\big](t)\big]\big|\\
\leq&\ C\Big(E^Q\Big[\int_0^t\int_0^s\big(|\widehat{\gamma}_{t,s,r}^{\theta,1,1}|^2+|\widehat{\gamma}_{t,s,r}^{\theta,2,2}|^2\big)\mathbf{1}_{E_\varepsilon}(r)drds\Big]\Big)^{\frac{1}{2}}
\sqrt{\varepsilon}+C\Big(E^Q\Big[\int_0^t|\widehat{\gamma}_{t,s}^{\theta,2}|^2 ds\Big]\Big)^{\frac{1}{2}}\varepsilon\\
=&\ \rho_t^1(\varepsilon)\sqrt{\varepsilon},\ \sqrt{\varepsilon}>0,
\end{split}
\end{equation*}
where
\begin{equation*}
\begin{split}
\rho_t^1(\varepsilon):&=C\Big(E^Q\Big[\int_0^t\int_0^s\big(|\widehat{\gamma}_{t,s,r}^{\theta,1,1}|^2+|\widehat{\gamma}_{t,s,r}^{\theta,2,2}|^2\big)\mathbf{1}_{E_\varepsilon}(r)drds\Big]
\Big)^{\frac{1}{2}}+C\sqrt{\varepsilon}\Big(E^Q\Big[\int_0^t|\widehat{\gamma}_{t,s}^{\theta,2}|^2 ds\Big]\Big)^{\frac{1}{2}}\\
&\leq C\Big(E^Q\Big[\int_0^t\int_0^s|\widehat{\gamma}_{t,s,r}^{\theta,1,1}|^2drds\Big]+E^Q\Big[\int_0^t\int_0^s|\widehat{\gamma}_{t,s,r}^{\theta,2,2}|^2drds\Big]\Big)^{\frac{1}{2}}+C\Big
(E^Q\Big[\int_0^t|\widehat{\gamma}_{t,s}^{\theta,2}|^2 ds\Big]\Big)^{\frac{1}{2}}\\
&\leq C\Big(E^Q\Big[\int_0^t|\widehat{\gamma}_{t,s}^{\theta}|^2ds\Big]\Big)^{\frac{1}{2}}\leq C\big(E^Q[|\theta_t^1|^2]\big)^{\frac{1}{2}},\ t\in[0,T],\ \varepsilon\in(0,1],
\end{split}
\end{equation*}
and, by the dominated convergence theorem, it holds that $$ \rho_t^1(\varepsilon)\rightarrow0\ (\varepsilon\searrow0),\ t\in[0,T]. $$
4) We now make the computation for $E^Q\big[\theta_t^2 D_2\big[A[G^\varepsilon]\big](t)\big]$. Remark that, with the notation \\ $A[G^\varepsilon]=\big(A[G^\varepsilon]^1,A[G^\varepsilon]^2\big)$,
\begin{equation*}
\begin{split}
&E^Q\big[\theta_t^2 D_2\big[A[G^\varepsilon]\big](t)\big]=E^Q\Big[\theta_t^2\Big(\int_0^t h_x(s)A[G^\varepsilon]^1(s)dY_s-\int_0^t h(s)h_x(s)A[G^\varepsilon]^1(s)ds\Big)\Big]\\
=&E^Q\Big[\int_0^t\theta_{t,s}^{2,2}h_x(s)A[G^\varepsilon]^1(s)ds-\theta_t^2\int_0^t h(s)h_x(s)A[G^\varepsilon]^1(s)ds\Big]\\
=&E^Q\Big[\int_0^t h_x(s)\big(\theta_{t,s}^{2,2}-E^Q[\theta_t^2\,|\,\mathcal{F}_s]h(s)\big)A[G^\varepsilon]^1(s)ds\Big]=\big(\overline{\gamma}_t^\theta,A[G^\varepsilon]\big)_H,
\end{split}
\end{equation*}
for $\overline{\gamma}_t^\theta=(\overline{\gamma}_t^{\theta,1},\overline{\gamma}_t^{\theta,2})$, $\overline{\gamma}_{t,s}^{\theta,1}:=h_x(s)\big(\theta_{t,s}^{2,2}-E^Q[\theta_t^2\,|\,\mathcal{F}_s]h(s)\big)
\mathbf{1}_{[0,t]}(s)$, $\overline{\gamma}_{t,s}^{\theta,2}:=0$, $s\in[0,T]$. Hence, for $\overline{\overline{\gamma}}_t^\theta:=A^*[\overline{\gamma}_t^\theta]$, as $\overline{\overline{\gamma}}_t^
\theta=A^*[\overline{\gamma}_t^\theta\mathbf{1}_{[0,t]}]=A^*[\overline{\gamma}_t^\theta\mathbf{1}_{[0,t]}]\mathbf{1}_{[0,t]}=\overline{\overline{\gamma}}_t^\theta\mathbf{1}_{[0,t]}$,
$$ E^Q\big[\theta_t^2 D_2\big[A[G^\varepsilon]\big](t)\big]=\big(\overline{\overline{\gamma}}_t^\theta,G^\varepsilon\big)_H=E^Q\Big[\int_0^t \overline{\overline{\gamma}}_{t,s}^\theta G_s^\varepsilon ds\Big], $$
and in analogy to 3),
$$ \big|E^Q\big[\theta_t^2 D_2\big[A[G^\varepsilon]\big](t)\big]\big|\leq \rho_t^2(\varepsilon)\sqrt{\varepsilon},\ \varepsilon\in(0,1], $$
with $\rho_t^2(\varepsilon)\rightarrow0\ (\varepsilon\searrow0)$, $\rho_t^2(\varepsilon)\leq C \big(E^Q[|\theta_t^2|^2]\big)^{\frac{1}{2}},\ \varepsilon\in(0,1],\ t\in[0,T]$.\\
5) Recall that
\begin{equation*}
\begin{split}
&E^Q[\theta_t Z_t^{1,\varepsilon}]=E^Q\big[\theta_t\big(G_t^\varepsilon+D\big[A[G^\varepsilon]\big](t)\big)\big]\\
=&E^Q[\theta_t G_t^\varepsilon]+E^Q\big[\theta_t^1 D_1\big[A[G^\varepsilon]\big](t)\big]+E^Q\big[\theta_t^2 D_2\big[A[G^\varepsilon]\big](t)\big],
\end{split}
\end{equation*}
and
\begin{equation*}
\begin{split}
&E^Q[\theta_t G_t^\varepsilon]=E^Q\Big[\theta_t^1\int_0^t\delta\sigma(s)\mathbf{1}_{E_\varepsilon}(s)dB_s^1\Big]+E^Q\Big[\theta_t^2\Big(\int_0^t\delta h(s)\mathbf{1}_{E_\varepsilon}(s)dY_s\!-\!\int_0^t h(s)\delta h(s)\mathbf{1}_{E_\varepsilon}(s)ds\Big)\Big]\\
=&E^Q\Big[\int_0^t\theta_{t,s}^{1,1}\delta\sigma(s)\mathbf{1}_{E_\varepsilon}(s)ds\Big]+E^Q\Big[\int_0^t\theta_{t,s}^{2,2}\delta h(s)\mathbf{1}_{E_\varepsilon}(s)ds\Big]-E^Q\Big[\int_0^t\theta_t^2 h(s)\delta h(s)\mathbf{1}_{E_\varepsilon}(s)ds\Big]\\
=&E^Q\Big[\int_0^t\big(\theta_{t,s}^{1,1}\delta\sigma(s)+\theta_{t,s}^{2,2}\delta h(s)-E^Q[\theta_t^2\,|\,\mathcal{F}_s]h(s)\delta h(s)\big)\mathbf{1}_{E_\varepsilon}(s)ds\Big].
\end{split}
\end{equation*}
Hence, $\big|E^Q[\theta_t G_t^\varepsilon]\big|\leq \rho_t^0(\varepsilon)\sqrt{\varepsilon}$, where
\begin{equation*}
\begin{split}
\rho_t^0(\varepsilon)&=C\Big(E^Q\Big[\int_0^t\big(|\theta_{t,s}^{1,1}|^2+|\theta_{t,s}^{2,2}|^2+|\theta_{t}^{2}|^2\big)\mathbf{1}_{E_\varepsilon}(s)ds\Big]\Big)^{\frac{1}{2}}\\
&\leq C\Big(E^Q\big[|\theta_t^1|^2+|\theta_t^2|^2\big]\Big)^{\frac{1}{2}}=C\Big(E^Q\big[|\theta_t|^2\big]\Big)^{\frac{1}{2}},\ \varepsilon>0,
\end{split}
\end{equation*}
and from the dominated convergence theorem, it follows that $\rho_t^0(\varepsilon)\rightarrow0\ (\varepsilon\searrow0)$.

Summarizing, we see, for all $\theta\in H$, there exists a function $\rho:[0,T]\times\mathbb{R}_+\rightarrow\mathbb{R}_+$ with
$$ \big|E^Q[\theta_t Z_t^{1,\varepsilon}]\big|\big(=\big|E^Q[\theta_t^1 Y_t^{1,\varepsilon}]+E^Q[\theta_t^2 N_t^{1,\varepsilon}]\big|\big)\leq C\rho_t(\varepsilon)\sqrt{\varepsilon},\ \ \varepsilon\in(0,1], $$
such that $\displaystyle {\rm i)}\  \rho_t(\varepsilon)\rightarrow 0\, (\varepsilon\searrow0);\ \mbox{and}\
{\rm ii)}\  \rho_t(\varepsilon)\leq C\big(E^Q[|\theta_t|^2]\big)^{\frac{1}{2}},\ t\in[0,T],\ \varepsilon\in(0,1].$
\end{proof}

Now we want to deduce the second order variational equation, and for this we also apply the notations introduced in the proof of Proposition \ref{EstofXL}. To better emphasize our approach without overcharging the computations, we still restrict to the case $\sigma=\sigma(\gamma,v),\ h=h(x,v)$, $(\gamma,x,v)\in\mathcal{P}_2(\mathbb{R})\times\mathbb{R}\times U$. Then, for all $t\in[0,T]$, $\varepsilon>0$,
\begin{equation*}
\begin{split}
&X_t^\varepsilon-X_t=\int_0^t\big(\sigma(\mu_s^\varepsilon,u_s^\varepsilon)-\sigma(\mu_s,u_s^\varepsilon)\big)dB_s^1+\int_0^t\delta\sigma(s)\mathbf{1}_{E_\varepsilon}(s)dB_s^1,\\ &Y_t^{1,\varepsilon}=\int_0^t\Big\{\widetilde{E}^Q\Big[\int_0^{\widetilde{U}_s}\partial_\mu\sigma(\mu_s,u_s;y)dy\cdot\widetilde{K}_s^{1,\varepsilon}\Big]+\widetilde{E}^Q\big[\partial_\mu
\sigma(\mu_s,u_s;\widetilde{U}_s)\widetilde{L}_s\widetilde{V}_s^{1,\varepsilon}\big]+\delta\sigma(s)\mathbf{1}_{E_\varepsilon}(s)\Big\}dB_s^1,\\
&K_t^{1,\varepsilon}=\int_0^t\big(h(s)K_s^{1,\varepsilon}+L_s h_x(s)Y_s^{1,\varepsilon}+L_s \delta h(s)\mathbf{1}_{E_\varepsilon}(s)\big)dY_s,\ t\in[0,T].\\
\end{split}
\end{equation*}
Recalling that $(X^{\varepsilon,\lambda},L^{\varepsilon,\lambda},U^{\varepsilon,\lambda}):=(1-\lambda)(X,L,U)+\lambda(X^{\varepsilon},L^{\varepsilon},U^{\varepsilon}),\ \lambda\in[0,1]$, we get that, for $\mu_t^{\varepsilon,\lambda}:=(L_t^{\varepsilon,\lambda}Q)_{U_t^{\varepsilon,\lambda}}$,
\begin{equation*}
\begin{split}
&\sigma(\mu_t^\varepsilon,u_t^\varepsilon)-\sigma(\mu_t,u_t^\varepsilon)=\int_0^1\partial_\lambda \big[\sigma(\mu_t^{\varepsilon,\lambda},u_t^\varepsilon)\big]dt\\
=&\int_0^1\Big\{\widetilde{E}^Q\Big[\int_0^{\widetilde{U}_t^{\varepsilon,\lambda}}\partial_\mu \sigma(\mu_t^{\varepsilon,\lambda},u_t^\varepsilon;y)dy(\widetilde{L}_t^\varepsilon-\widetilde{L}_t)\Big]+\widetilde{E}^{Q}\big[\partial_\mu \sigma(\mu_t^{\varepsilon,\lambda},u_t^\varepsilon;\widetilde{U}_t^{\varepsilon,\lambda})\widetilde{L}_t^{\varepsilon,\lambda}(\widetilde{U}_t^\varepsilon-\widetilde{U}_t)\big]\Big\}d\lambda.\\
\end{split}
\end{equation*}
Then we have
\begin{equation}\label{sqr*}
\begin{split}
&X_t^\varepsilon-(X_t+Y_t^{1,\varepsilon})\\
=&\int_0^t\int_0^1\widetilde{E}^Q\Big[\Big(\int_0^{\widetilde{U}_s^{\varepsilon,\lambda}}\partial_\mu\sigma(\mu_s^{\varepsilon,\lambda},u_s^\varepsilon;y)dy-\int_0^{\widetilde{U}_s}
\partial_\mu\sigma(\mu_s,u_s^\varepsilon;y)dy\Big)(\widetilde{L}_s^\varepsilon-\widetilde{L}_s)\Big]d\lambda dB_s^1\\
&+\int_0^t\int_0^1\widetilde{E}^Q\Big[\int_0^{\widetilde{U}_s}\partial_\mu\sigma(\mu_s,u_s^\varepsilon;y)dy\big(\widetilde{L}_s^\varepsilon-(\widetilde{L}_s+\widetilde{K}_s^{1,
\varepsilon})\big)\Big]d\lambda dB_s^1\\
&+\int_0^t\int_0^1\widetilde{E}^Q\big[\big(\partial_\mu\sigma(\mu_s^{\varepsilon,\lambda},u_s^\varepsilon;\widetilde{U}_s^{\varepsilon,\lambda})\widetilde{L}_s^{\varepsilon,\lambda}-
\partial_\mu\sigma(\mu_s,u_s^\varepsilon;\widetilde{U}_s)\widetilde{L}_s\big)(\widetilde{U}_s^\varepsilon-\widetilde{U}_s)\big]d\lambda dB_s^1\\
&+\int_0^t\int_0^1\widetilde{E}^Q\big[\partial_\mu\sigma(\mu_s,u_s^\varepsilon;\widetilde{U}_s)\widetilde{L}_s\big(\widetilde{U}_s^\varepsilon-(\widetilde{U}_s+\widetilde{V}_s^{1,
\varepsilon})\big)\big]d\lambda dB_s^1\\
&+\int_0^t\Big(\widetilde{E}^Q\Big[\int_0^{\widetilde{U}_s}\delta(\partial_\mu\sigma)(s,y)dy\widetilde{K}_s^{1,\varepsilon}\Big]+\widetilde{E}^Q\big[\delta(\partial_\mu\sigma)(s,
\widetilde{U}_s)\widetilde{L}_s\widetilde{V}_s^{1,\varepsilon}\big]\Big)\mathbf{1}_{E_\varepsilon}(s)dB_s^1,
\end{split}
\end{equation}
and we denote by $I_t^{1,\varepsilon}$, $I_t^{2,\varepsilon}$ and $I_t^{3,\varepsilon}$ the first, the third and the last term of the right side in above equation, respectively. To compute $I_t^{1,\varepsilon}$ and  $I_t^{2,\varepsilon}$, note that \\
(1) $W_1(\mu_s^{\varepsilon,\rho},\mu_s)\leq C\sqrt{\varepsilon},\ s\in[0,T],\ \varepsilon>0,\ \rho\in[0,1]$. Indeed, for all $\varphi\in\mbox{Lip}_1(\mathbb{R})$ with $\varphi(0)=0$,
\begin{equation*}
\begin{split}
&\int_\mathbb{R}\varphi d\mu_s^{\varepsilon,\rho}-\int_\mathbb{R}\varphi d\mu_s=E^Q[L_s^{\varepsilon,\rho}\varphi(U_s^{\varepsilon,\rho})]-E^Q[L_s\varphi(U_s)]\\
=&E^Q\big[(L_s^{\varepsilon,\rho}-L_s)\varphi(U_s^{\varepsilon,\rho})\big]+E^Q\big[L_s\big(\varphi(U_s^{\varepsilon,\rho})-\varphi(U_s)\big)\big]\leq E^Q\big[|L_s^{\varepsilon}-L_s||U_s^{\varepsilon,\rho}|\big]
+E^Q\big[L_s|U_s^{\varepsilon,\rho}-U_s|\big]\\
\leq& E^Q\big[|L_s^{\varepsilon}-L_s|\big(|U_s|+|U_s^{\varepsilon}|\big)\big]+E^Q\big[L_s|U_s^{\varepsilon}-U_s|\big]\leq C\sqrt{\varepsilon}.
\end{split}
\end{equation*}
The latter inequality comes from Proposition \ref{EstofXL} and Corollary \ref{EstofUV}. Thus, from the Kantorovich-Rubinstein duality it follows that $W_1(\mu_s^{\varepsilon,\rho},\mu_s)\leq C\sqrt{\varepsilon}$.\\
(2) Now we estimate $I_t^{1,\varepsilon}$. It is easy to see that
\begin{equation*}
\begin{split}
I_t^{1,\varepsilon}=&\int_0^t\int_0^1\widetilde{E}^Q\Big[\Big(\int_0^{\widetilde{U}_s^{\varepsilon,\lambda}}\partial_\mu\sigma(\mu_s^{\varepsilon,\lambda},u_s^\varepsilon;y)dy-\int_0^{
\widetilde{U}_s}\partial_\mu\sigma(\mu_s,u_s^\varepsilon;y)dy\Big)(\widetilde{L}_s^\varepsilon-\widetilde{L}_s)\Big]d\lambda dB_s^1\\
=&\int_0^t\int_0^1\int_0^\lambda\partial_\rho\Big\{\widetilde{E}^Q\Big[\Big(\int_0^{\widetilde{U}_s^{\varepsilon,\rho}}\partial_\mu\sigma(\mu_s^{\varepsilon,\rho},u_s^\varepsilon;y)dy\Big)
(\widetilde{L}_s^\varepsilon-\widetilde{L}_s)\Big]\Big\}d\rho d\lambda dB_s^1\\
=&\int_0^t\int_0^1\int_0^\lambda\Big(\widetilde{E}^Q\big[\partial_\mu\sigma(\mu_s^{\varepsilon,\rho},u_s^\varepsilon;\widetilde{U}_s^{\varepsilon,\rho})(\widetilde{U}_s^\varepsilon-
\widetilde{U}_s)(\widetilde{L}_s^\varepsilon-\widetilde{L}_s)\big]\Big)d\rho d\lambda dB_s^1+I_t^{1,1,\varepsilon}+I_t^{1,2,\varepsilon},
\end{split}
\end{equation*}
with
\small
\begin{equation*}
\begin{split}
&I_t^{1,1,\varepsilon}=\int_0^t\int_0^1\int_0^\lambda\widetilde{E}^Q\Big[\int_0^{\widetilde{U}_s^{\varepsilon,\rho}}dy\widehat{E}^Q\Big[\Big(\int_0^{\widehat{U}_s^{\varepsilon,\rho}}
\partial_\mu(\partial_\mu\sigma)(\mu_s^{\varepsilon,\rho},u_s^\varepsilon;y,y')dy'\Big)(\widehat{L}_s^\varepsilon-\widehat{L}_s)\Big](\widetilde{L}_s^\varepsilon-\widetilde{L}_s)\Big]d\rho d\lambda dB_s^1,\\
&I_t^{1,2,\varepsilon}=\int_0^t\int_0^1\int_0^\lambda\widetilde{E}^Q\Big[\int_0^{\widetilde{U}_s^{\varepsilon,\rho}}dy\widehat{E}^Q\big[\partial_\mu(\partial_\mu\sigma)(\mu_s^{
\varepsilon,\rho},u_s^\varepsilon;y,\widehat{U}_s^{\varepsilon,\rho})\widehat{L}_s^{\varepsilon,\rho}(\widehat{U}_s^{\varepsilon}-\widehat{U}_s)\big](\widetilde{L}_s^\varepsilon
-\widetilde{L}_s)\Big]d\rho d\lambda dB_s^1.
\end{split}
\end{equation*}
\normalsize
Combining the Lipschitz continuity of $(\gamma,y')\rightarrow\partial_\mu(\partial_\mu\sigma)(\gamma,u_s^\varepsilon;y,y')$ (in $\gamma\in\mathcal{P}_2(\mathbb{R})$ with respect to the 1-Wasserstein distance) with (1) and Proposition \ref{EstofXL} as well as Corollary \ref{EstofUV} we see that
\begin{equation*}
\begin{split}
I_t^{1,2,\varepsilon}&=\int_0^t\frac{1}{2}\widetilde{E}^Q\Big[\int_0^{\widetilde{U}_s}dy\widehat{E}^Q\big[\partial_\mu(\partial_\mu\sigma)(\mu_s,u_s^\varepsilon;y,\widehat{U}_s)
\widehat{L}_s(\widehat{U}_s^{\varepsilon}-\widehat{U}_s)\big](\widetilde{L}_s^\varepsilon-\widetilde{L}_s)\Big]dB_s^1+R_t^{1,2,\varepsilon}\\
&=\int_0^t\frac{1}{2}\widetilde{E}^Q\Big[\int_0^{\widetilde{U}_s}dy\widehat{E}^Q\big[\partial_\mu(\partial_\mu\sigma)(\mu_s,u_s;y,\widehat{U}_s)\widehat{L}_s(\widehat{U}_s^{\varepsilon}
-\widehat{U}_s)\big](\widetilde{L}_s^\varepsilon-\widetilde{L}_s)\Big]dB_s^1+R_t^{1,3,\varepsilon},
\end{split}
\end{equation*}
where $\displaystyle \Big(E^Q\big[\sup_{t\leq T}|R_t^{1,2,\varepsilon}|^p\big]\Big)^{\frac{1}{p}}\leq C_p\varepsilon^{\frac{3}{2}},\ \varepsilon>0,\ p\geq 2,$\
and
$$ R_t^{1,3,\varepsilon}=R_t^{1,2,\varepsilon}+\int_0^t\frac{1}{2}\widetilde{E}^Q\Big[\int_0^{\widetilde{U}_s}dy\widehat{E}^Q\big[\delta\big(\partial_\mu(\partial_\mu\sigma)\big)(s,y,
\widehat{U}_s)\widehat{L}_s(\widehat{U}_s^{\varepsilon}-\widehat{U}_s)\big](\widetilde{L}_s^\varepsilon-\widetilde{L}_s)\Big]\mathbf{1}_{E_\varepsilon}(s)dB_s^1.
$$
Obviously, also $\displaystyle  \Big(E^Q\big[\sup_{t\leq T}|R_t^{1,3,\varepsilon}|^p\big]\Big)^{\frac{1}{p}}\leq C_p\varepsilon^{\frac{3}{2}},\ \varepsilon>0,\ p\geq 2. $ \\
On the other hand, again from Proposition \ref{EstofXL} and Corollary \ref{EstofUV} we obtain
\begin{equation}\label{4.16-1}
\begin{split}
I_t^{1,2,\varepsilon}&=\int_0^t\frac{1}{2}\widetilde{E}^Q\Big[\int_0^{\widetilde{U}_s}dy\widehat{E}^Q\big[\partial_\mu(\partial_\mu\sigma)(\mu_s,u_s^\varepsilon;y,\widehat{U}_s)
\widehat{L}_s(\widehat{U}_s^{\varepsilon}-\widehat{U}_s)\big](\widetilde{L}_s^\varepsilon-\widetilde{L}_s)\Big]dB_s^1+R_t^{1,2,\varepsilon}\\
&=:\int_0^t\frac{1}{2}\widetilde{E}^Q\Big[\int_0^{\widetilde{U}_s}dy\widehat{E}^Q\big[\partial_\mu(\partial_\mu\sigma)(\mu_s,u_s^\varepsilon;y,\widehat{U}_s)\widehat{L}_s\widehat{V}_s^{1,
\varepsilon}\big]\widetilde{K}_s^{1,\varepsilon}\Big]dB_s^1+R_t^{1,4,\varepsilon}\\
&=:\int_0^t\frac{1}{2}\widetilde{E}^Q\Big[\int_0^{\widetilde{U}_s}dy\widehat{E}^Q\big[\partial_\mu(\partial_\mu\sigma)(\mu_s,u_s;y,\widehat{U}_s)\widehat{L}_s\widehat{V}_s^{1,
\varepsilon}\big]\widetilde{K}_s^{1,\varepsilon}\Big]dB_s^1+R_t^{1,4,\varepsilon}+R_t^{1,5,\varepsilon},
\end{split}
\end{equation}
where $\displaystyle  \Big(E^Q\big[\sup_{t\leq T}|R_t^{1,4,\varepsilon}|^p\big]\Big)^{\frac{1}{p}}\leq C_p\varepsilon^{\frac{3}{2}},\ \varepsilon>0,\ p\geq 2,$\ and
$$ \Big(E^Q\big[\sup_{t\leq T}|R_t^{1,5,\varepsilon}|^p\big]\Big)^{\frac{1}{p}}\leq C_p\varepsilon\Big(\int_0^T \mathbf{1}_{E_\varepsilon}(s)ds\Big)^{\frac{1}{2}}
\leq C_p\varepsilon^{\frac{3}{2}},\ \varepsilon>0,\ p\geq 2. $$
Put $\sigma_{\mu\mu}(s,y,y'):=\partial_\mu(\partial_\mu\sigma)(\mu_s,u_s;y,y')$, and observe that
\begin{equation*}
\begin{split}
&\widehat{E}^Q\big[\sigma_{\mu\mu}(s,y,\widehat{U}_s)\widehat{L}_s\widehat{V}_s^{1,\varepsilon}\big]\\
=&\widehat{E}^Q\Big[\sigma_{\mu\mu}(s,y,\widehat{U}_s)\widehat{L}_s \Big(\frac{\widehat{E}^Q[\widehat{L}_s \widehat{Y}_s^{1,\varepsilon}+\widehat{X}_s \widehat{K}_s^{1,\varepsilon}\,|\,\mathcal{F}_s^{\widehat{Y}}]}{\widehat{E}^Q[\widehat{L}_s\,|\,\mathcal{F}_s^{\widehat{Y}}]}-\frac{\widehat{E}^{Q}[\widehat{L}_s \widehat{X}_s\,|\,\mathcal{F}_s^{\widehat{Y}}] }{(\widehat{E}^{Q}[\widehat{L}_s\,|\,\mathcal{F}_s^{\widehat{Y}}])^2}\widehat{E}^{Q}[\widehat{K}_s^{1,\varepsilon}\,|\,\mathcal{F}_s^{\widehat{Y}}] \Big)\Big]\\
=&\widehat{E}^Q\Big[\Big(\frac{\widehat{E}^{Q}[\sigma_{\mu\mu}(s,y,\widehat{U}_s)\widehat{L}_s\,|\,\mathcal{F}_s^{\widehat{Y}}]}{\widehat{E}^{Q}[\widehat{L}_s\,|\,\mathcal{F}_s^{
\widehat{Y}}]}\widehat{L}_s\Big)\widehat{Y}_s^{1,\varepsilon}\Big]\\
&\qquad+\widehat{E}^Q\Big[\bigg(\frac{\widehat{E}^{Q}[\sigma_{\mu\mu}(s,y,\widehat{U}_s)\widehat{L}_s\,|\,\mathcal{F}_s^{\widehat{Y}}]}{\widehat{E}^{Q}[\widehat{L}_s\,|\,\mathcal{F}_s^{
\widehat{Y}}]}\Big(\widehat{X}_s-\frac{\widehat{E}^{Q}[\widehat{L}_s \widehat{X}_s\,|\,\mathcal{F}_s^{\widehat{Y}}] }{\widehat{E}^{Q}[\widehat{L}_s\,|\,\mathcal{F}_s^{\widehat{Y}}]}\Big)\bigg)\widehat{K}_s^{1,\varepsilon}\Big].
\end{split}
\end{equation*}
For simplicity we put $ \displaystyle \widehat{\theta}_s^1(y):=\frac{\widehat{E}^{Q}[\sigma_{\mu\mu}(s,y,\widehat{U}_s)\widehat{L}_s\,|\,\mathcal{F}_s^{\widehat{Y}}]}{\widehat{E}^{Q}[\widehat{L}_s\,|\,\mathcal{F}_s^{
\widehat{Y}}]}\widehat{L}_s, $ and
$$ \widehat{\theta}_s^2(y):=\frac{\widehat{E}^{Q}[\sigma_{\mu\mu}(s,y,\widehat{U}_s)\widehat{L}_s\,|\,\mathcal{F}_s^{\widehat{Y}}]}{\widehat{E}^{Q}[\widehat{L}_s\,|\,\mathcal{F}_s^{
\widehat{Y}}]}\Big(\widehat{X}_s-\frac{\widehat{E}^{Q}[\widehat{L}_s \widehat{X}_s\,|\,\mathcal{F}_s^{\widehat{Y}}] }{\widehat{E}^{Q}[\widehat{L}_s\,|\,\mathcal{F}_s^{\widehat{Y}}]}\Big). $$
With these notations, we have $\displaystyle  \widehat{E}^{Q}[\sigma_{\mu\mu}(s,y,\widehat{U}_s)\widehat{L}_s\widehat{V}_s^{1,\varepsilon}]=\widehat{E}^{Q}[\widehat{\theta}_s^1(y)\widehat{Y}_s^{1,\varepsilon}]-\widehat{E}^{Q}[
\widehat{\theta}_s^2(y)\widehat{K}_s^{1,\varepsilon}]. $\\
Note that for $p\geq2$, $\displaystyle \big(\widehat{E}^Q[|\widehat{\theta}_s^i(y)|^p]\big)^{\frac{1}{p}}\leq C_p,\ s\in[0,T],\ y\in\mathbb{R},\ i=1,2. $\
Then, from Proposition \ref{tech},
$$ \big|\widehat{E}^Q[\widehat{\theta}_s^1(y)\widehat{Y}_s^{1,\varepsilon}]\big|+\big|\widehat{E}^Q[\widehat{\theta}_s^2(y)\widehat{K}_s^{1,\varepsilon}]\big|\leq \rho_s(\varepsilon,y)\sqrt{\varepsilon},
\ \varepsilon>0,\ s\in[0,T], $$
with $\rho_s(\varepsilon,y)\rightarrow0\ (\varepsilon\searrow0)$, and $\displaystyle \rho_s(\varepsilon,y)\leq C\widehat{E}^Q\big[|\widehat{\theta}_s^1(y)|^2+|\widehat{\theta}_s^2(y)\widehat{L}_s|^2\big]\leq C,\ s\in[0,T],
\ \varepsilon>0,\ y\in\mathbb{R}. $\\

\noindent Consequently, from \eqref{4.16-1} we have
\begin{equation*}
\begin{split}
&\big(E^Q[\sup_{t\in[0,T]}|I_t^{1,2,\varepsilon}|^p]\big)^{\frac{1}{p}}\leq C_p\varepsilon^{\frac{3}{2}}+C_p\bigg(\!E^Q\!\Big[\Big(\!\int_0^T\!\Big|\widetilde{E}^Q\Big[
    \int_0^{\widetilde{U}_s}dy\widehat{E}^Q\big[\sigma_{\mu\mu}(s,y,\widehat{U}_s)\widehat{L}_s
\widehat{V}_s^{1,\varepsilon}\big]\widetilde{K}_s^{1,\varepsilon}\Big]\Big|^2ds\Big)^{\frac{p}{2}}\!\Big]\!\bigg)^{\frac{1}{p}}\\
\leq & C_p\varepsilon^{\frac{3}{2}}+C_p\sqrt{\varepsilon}\bigg(E^Q\Big[\Big(\int_0^T\Big|\widetilde{E}^Q\big[\big|\int_0^{\widetilde{U}_s}\rho_s(\varepsilon,y)dy\big|\,\big|
\widetilde{K}_s^{1,\varepsilon}\big|\big]
\Big|^2ds\Big)^{\frac{p}{2}}\Big]\bigg)^{\frac{1}{p}}\\
\leq & C_p\varepsilon^{\frac{3}{2}}+C_p\varepsilon\bigg(E^Q\Big[\Big(\int_0^T\widetilde{E}^Q\big[\big|\int_0^{\widetilde{U}_s}\rho_s(\varepsilon,y)dy\big|^2\big]ds\Big)^{\frac{p}{2}}
\Big]\bigg)^{\frac{1}{p}}
=: C_p\varepsilon^{\frac{3}{2}}+C_p\varepsilon\rho_p(\varepsilon),
\end{split}
\end{equation*}
where, thanks to the dominated convergence theorem, $\rho_p(\varepsilon)\rightarrow 0$, as $\varepsilon \searrow 0$.\ This proves that
\begin{equation}\label{xing}
\big(E^Q[\sup_{t\in[0,T]}|I_t^{1,2,\varepsilon}|^p]\big)^{\frac{1}{p}}\leq C_p\varepsilon\rho_p(\varepsilon),\ \varepsilon>0,
\end{equation}
with $\rho_p(\varepsilon)\rightarrow0\ (\varepsilon\searrow0)$, $p\geq2$. Similarly,
\begin{equation}\label{xingxing}
\big(E^Q[\sup_{t\in[0,T]}|I_t^{1,1,\varepsilon}|^p]\big)^{\frac{1}{p}}\leq C_p\varepsilon\rho_p(\varepsilon),\ \varepsilon>0,
\end{equation}
with $\rho_p(\varepsilon)\rightarrow0\ (\varepsilon\searrow0)$, $p\geq2$. Indeed, comparing $I_t^{1,1,\varepsilon}$ with $I_t^{1,2,\varepsilon}$ and taking into account the approach to estimate
$I_t^{1,2,\varepsilon}$, we realize that the hard kernel concerns the estimate
\begin{equation*}
\begin{split}
&\widehat{E}^Q\Big[\Big(\int_0^{\widehat{U}_s}\sigma_{\mu\mu}(s,y,y')dy'\Big)(\widehat{L}_s^\varepsilon-\widehat{L}_s)\Big]\\
=&\widehat{E}^Q\Big[\Big(\int_0^{\widehat{U}_s}\sigma_{\mu\mu}(s,y,y')dy'\Big)\big(\widehat{L}_s^\varepsilon-(\widehat{L}_s+\widehat{K}_s^{1,\varepsilon})\big)\Big]+\widehat{E}^Q\Big[
\Big(\int_0^{\widehat{U}_s}\sigma_{\mu\mu}(s,y,y')dy'\Big)\widehat{K}_s^{1,\varepsilon}\Big].
\end{split}
\end{equation*}
Putting $\displaystyle \widehat{\theta}_s^3(y):=\int_0^{\widehat{U}_s}\sigma_{\mu\mu}(s,y,y')dy'$, $s\in[0,T]$, $y\in\mathbb{R}$, we have $|\widehat{\theta}_s^3(y)|\leq C|\widehat{U}_s|$, $s\in[0,T]$,
$y\in\mathbb{R}$, and from Proposition \ref{EstofXL} and Corollary \ref{EstofUV} we get
$\displaystyle \big|\widehat{E}^Q\big[\widehat{\theta}_s^3(y)\big(\widehat{L}_s^\varepsilon-(\widehat{L}_s+\widehat{K}_s^{1,\varepsilon})\big)\big]\big|\leq C\varepsilon. $\\
On the other hand, to $\displaystyle \widehat{E}^Q\Big[\int_0^{\widehat{U}_s}\sigma_{\mu\mu}(s,y,y')dy'\cdot\widehat{K}_s^{1,\varepsilon}\Big]$ we apply Proposition \ref{tech}. This yields
$$ \big|\widehat{E}^Q\big[\widehat{\theta}_s^3(y)\widehat{K}_s^{1,\varepsilon}\big]\big|\leq \rho_s(\varepsilon,y)\sqrt{\varepsilon},\ \varepsilon>0,\ s\in[0,T], $$
with $\rho_s(\varepsilon,y)\rightarrow0\ (\varepsilon\searrow0)$ and $\displaystyle \rho_s(\varepsilon,y)\leq C\widehat{E}^Q\big[|\widehat{\theta}_s^3(y)\widehat{L}_s|^2\big]\leq C,\ s\in[0,T],\ y\in\mathbb{R},
\ \varepsilon>0. $\ We obtain $\displaystyle \Big(E^Q\big[\sup_{t\in[0,T]}|I_t^{1,1,\varepsilon}|^p\big]\Big)^{\frac{1}{p}}\leq C_p\varepsilon\rho_p(\varepsilon),\ \varepsilon>0,\ p\geq 2.
$\ Consequently, having \eqref{xing} and \eqref{xingxing}, again using the Proposition \ref{EstofXL} and Corollary \ref{EstofUV} yields
\begin{equation*}
\begin{split}
I_t^{1,\varepsilon}&=\int_0^t\int_0^1\int_0^\lambda\Big(\widetilde{E}^Q\big[\partial_\mu\sigma(\mu_s^{\varepsilon,\rho},u_s^\varepsilon;\widetilde{U}_s^{\varepsilon,\rho})(
\widetilde{U}_s^\varepsilon-\widetilde{U}_s)(\widetilde{L}_s^\varepsilon-\widetilde{L}_s)\big]\Big)d\rho d\lambda dB_s^1+I_t^{1,1,\varepsilon}+I_t^{1,2,\varepsilon}\\
&=:\frac{1}{2}\int_0^t\widetilde{E}^Q\big[\partial_\mu\sigma(\mu_s,u_s^\varepsilon;\widetilde{U}_s)\widetilde{V}_s^{1,\varepsilon}\widetilde{K}_s^{1,\varepsilon}\big]dB_s^1+R_t^{1,1,\varepsilon}\\
&=\frac{1}{2}\int_0^t\widetilde{E}^Q\big[\widetilde{\sigma}_\mu(s)\widetilde{V}_s^{1,\varepsilon}\widetilde{K}_s^{1,\varepsilon}\big]dB_s^1+R_t^{1,2,\varepsilon},
\end{split}
\end{equation*}
where $\displaystyle  R_t^{1,2,\varepsilon}:=R_t^{1,1,\varepsilon}+\frac{1}{2}\int_0^t \widetilde{E}^Q\big[\delta\widetilde{\sigma}_\mu(s)\widetilde{V}_s^{1,\varepsilon}\widetilde{K}_s^{1,\varepsilon}\big]
\mathbf{1}_{E_\varepsilon}(s)dB_s^1$, and
$$ \Big(E^Q\big[\sup_{t\in[0,T]}|R_t^{1,2,\varepsilon}|^p\big]\Big)^{\frac{1}{p}}\leq C_p\varepsilon\rho_p(\varepsilon),\ \varepsilon>0,\ p\geq 2.$$
(3) In analogy to $I_t^{1,\varepsilon}$ we study $I_t^{2,\varepsilon}$ (refer to (\ref{sqr*})):
\begin{equation*}
\begin{split}
I_t^{2,\varepsilon}=&\int_0^t\int_0^1\widetilde{E}^Q\big[\big(\partial_\mu\sigma(\mu_s^{\varepsilon,\lambda},u_s^\varepsilon;\widetilde{U}_s^{\varepsilon,\lambda})\widetilde{L}_s^{
\varepsilon,\lambda}-\partial_\mu\sigma(\mu_s,u_s^\varepsilon;\widetilde{U}_s)\widetilde{L}_s\big)(\widetilde{U}_s^\varepsilon-\widetilde{U}_s)\big]d\lambda dB_s^1\\
=&\int_0^t\int_0^1\int_0^\lambda\widetilde{E}^Q\big[\partial_\rho\big(\partial_\mu\sigma(\mu_s^{\varepsilon,\rho},u_s^\varepsilon;\widetilde{U}_s^{\varepsilon,\rho})\widetilde{L}_s^{
\varepsilon,\rho}\big)(\widetilde{U}_s^\varepsilon-\widetilde{U}_s)\big]d\rho d\lambda dB_s^1\\
=&\int_0^t\int_0^1\int_0^\lambda\Big(\widetilde{E}^Q\big[\partial_\mu\sigma(\mu_s^{\varepsilon,\rho},u_s^\varepsilon;\widetilde{U}_s^{\varepsilon,\rho})(\widetilde{L}_s^{\varepsilon
}-\widetilde{L}_s)(\widetilde{U}_s^\varepsilon-\widetilde{U}_s)\big]\\
&\qquad\qquad+\widetilde{E}^Q\big[\partial_z(\partial_\mu\sigma)(\mu_s^{\varepsilon,\rho},u_s^\varepsilon;\widetilde{U}_s^{\varepsilon,\rho})\widetilde{L}_s^{\varepsilon,\rho}
(\widetilde{U}_s^\varepsilon-\widetilde{U}_s)^2\big]\Big)d\rho d\lambda dB_s^1\\
&+\int_0^t\int_0^1\int_0^\lambda\bigg(\widetilde{E}^Q\Big[\widehat{E}^Q\Big[\int_0^{\widehat{U}_s^{\varepsilon,\rho}}\partial_\mu(\partial_\mu\sigma)(\mu_s^{\varepsilon,\rho},
u_s^\varepsilon;\widetilde{U}_s^{\varepsilon,\rho},y)dy\widetilde{L}_s^{\varepsilon,\rho}(\widehat{L}_s^{\varepsilon}-\widehat{L}_s)\Big](\widetilde{U}_s^{\varepsilon}-\widetilde{U}_s)\Big]\\
&\qquad\qquad+\widetilde{E}^Q\Big[\widehat{E}^Q\big[\partial_\mu(\partial_\mu\sigma)(\mu_s^{\varepsilon,\rho},u_s^\varepsilon;\widetilde{U}_s^{\varepsilon,\rho},\widehat{U}_s^{
\varepsilon,\rho})\widehat{L}_s^{\varepsilon,\rho}(\widehat{U}_s^{\varepsilon}-\widehat{U}_s)\big]\widetilde{L}_s^{\varepsilon,\rho}(\widetilde{U}_s^{\varepsilon}-\widetilde{U}_s)\Big]\bigg)d\rho d\lambda dB_s^1\\
=:&I_t^{2,1,\varepsilon}+I_t^{2,2,\varepsilon}.
\end{split}
\end{equation*}
Similar to (2), we obtain $\displaystyle \big(E^Q[\sup_{t\in[0,T]}|I_t^{2,2,\varepsilon}|^p]\big)^{\frac{1}{p}}\leq C_p\varepsilon\rho_p(\varepsilon),\ \varepsilon>0, $\
with $\rho_p(\varepsilon)\rightarrow0\ (\varepsilon\searrow0)$, $p\geq2$. Moreover, with arguments already developed in (2), we also see that
\begin{equation*}
\begin{split}
I_t^{2,1,\varepsilon}=&\frac{1}{2}\int_0^t\Big(\widetilde{E}^Q\big[\partial_\mu\sigma(\mu_s,u_s;\widetilde{U}_s)\widetilde{K}_s^{1,\varepsilon}\widetilde{V}_s^{1,\varepsilon}\big]+
\widetilde{E}^Q\big[\partial_z(\partial_\mu\sigma)(\mu_s,u_s;\widetilde{U}_s)\widetilde{L}_s(\widetilde{V}_s^{1,\varepsilon})^2\big]\Big)dB_s^1+R_t^{2,1,\varepsilon},
\end{split}
\end{equation*}
with $\displaystyle \big(E^Q[\sup_{t\in[0,T]}|R_t^{2,1,\varepsilon}|^p]\big)^{\frac{1}{p}}\leq C_p\varepsilon\rho_p(\varepsilon),\ \varepsilon>0, $\
where $\rho_p(\varepsilon)\rightarrow0\ (\varepsilon\searrow0)$, $p\geq2$.
Hence,
$$ I_t^{2,\varepsilon}=\frac{1}{2}\int_0^t\Big(\widetilde{E}^Q\big[\widetilde{\sigma}_\mu(s)\widetilde{K}_s^{1,\varepsilon}\widetilde{V}_s^{1,\varepsilon}\big]+\widetilde{E}^Q\big[
\widetilde{\sigma}_{z\mu}(s)\widetilde{L}_s(\widetilde{V}_s^{1,\varepsilon})^2\big]\Big)dB_s^1+R_t^{2,\varepsilon}, $$
with $\displaystyle  \big(E^Q[\sup_{t\in[0,T]}|R_t^{2,\varepsilon}|^p]\big)^{\frac{1}{p}}\leq C_p\varepsilon\rho_p(\varepsilon),\ \varepsilon>0. $\
Consequently, coming back to \eqref{sqr*}, this together with Proposition \ref{EstofXL} and Corollary \ref{EstofUV} yields
\begin{equation}\label{2ndVE1}
\begin{split}
&X_t^\varepsilon-(X_t+Y_t^{1,\varepsilon})\\
=&\int_0^t\Big\{\widetilde{E}^Q\Big[\int_0^{\widetilde{U}_s}\sigma_\mu(s,y)dy\big(\widetilde{L}_s^\varepsilon-(\widetilde{L}_s+\widetilde{K}_s^{1,\varepsilon})\big)\Big]+\widetilde{E}^Q
\big[\widetilde{\sigma}_\mu(s)\widetilde{L}_s\big(\widetilde{U}_s^\varepsilon-(\widetilde{U}_s+\widetilde{V}_s^{1,\varepsilon})\big)\big]\Big\}dB_s^1\\
&+\int_0^t\Big(\widetilde{E}^Q\big[\widetilde{\sigma}_\mu(s)\widetilde{V}_s^{1,\varepsilon}\widetilde{K}_s^{1,\varepsilon}\big]+\frac{1}{2}\widetilde{E}^Q\big[\widetilde{\sigma}_{z\mu}
(s)\widetilde{L}_s(\widetilde{V}_s^{1,\varepsilon})^2\big]\Big)dB_s^1\\
&+\int_0^t\Big(\widetilde{E}^Q\Big[\int_0^{\widetilde{U}_s}\delta\big(\sigma_\mu(s,y)\big)dy\cdot\widetilde{K}_s^{1,\varepsilon}\Big]+\widetilde{E}^Q\big[\delta\big(
\widetilde{\sigma}_{\mu}(s)\big)\widetilde{L}_s\widetilde{V}_s^{1,\varepsilon}\big]\Big)\mathbf{1}_{E_\varepsilon}(s)dB_s^1+R_t^{X,\varepsilon},
\end{split}
\end{equation}
where $\displaystyle  \big(E^Q[\sup_{t\in[0,T]}|R_t^{X,\varepsilon}|^p]\big)^{\frac{1}{p}}\leq C_p\varepsilon\rho_p(\varepsilon),\ \varepsilon>0, $\
and $\rho_p(\varepsilon)\rightarrow0\ (\varepsilon\searrow0)$, $p\geq2$. On the other hand,
\begin{equation*}
\begin{split}
&L_t^\varepsilon-L_t=\int_0^t\int_0^1\partial_\lambda\big[h(X_s^{\varepsilon,\lambda},u_s^\varepsilon)L_s^{\varepsilon,\lambda}\big]d\lambda dY_s+\int_0^t\delta h(s)L_s dY_s\\
&=\int_0^t\int_0^1\big\{h(X_s^{\varepsilon,\lambda},u_s^\varepsilon)(L_s^{\varepsilon}-L_s)+h_x(X_s^{\varepsilon,\lambda},u_s^\varepsilon)L_s^{\varepsilon,\lambda}(X_s^{\varepsilon}-X_s)\big\}
d\lambda dY_s+\int_0^t\delta h(s)L_s\mathbf{1}_{E_\varepsilon}(s)dY_s.
\end{split}
\end{equation*}
Then we have
\begin{equation*}
\begin{split}
&L_t^\varepsilon-(L_t+K_t^{1,\varepsilon})=\int_0^t\int_0^1\int_0^\lambda\partial_\rho\big\{h(X_s^{\varepsilon,\rho},u_s^\varepsilon)(L_s^{\varepsilon}-L_s)+h_x(X_s^{\varepsilon,\rho},
u_s^\varepsilon)L_s^{\varepsilon,\rho}(X_s^{\varepsilon}-X_s)\big\}d\rho d\lambda dY_s\\
&\ \hskip3cm+\int_0^t\Big\{h(X_s,u_s^\varepsilon)\big(L_s^{\varepsilon}-(L_s+K_s^{1,\varepsilon})\big)+h_x(X_s,u_s^\varepsilon)L_s\big(X_s^{\varepsilon}-(X_s+Y_s^{1,\varepsilon})\big)\Big\}dY_s\\
&\hskip3cm+\int_0^t\big(\delta h(s)K_s^{1,\varepsilon}+\delta h_x(s)L_s Y_s^{1,\varepsilon}\big)\mathbf{1}_{E_\varepsilon}(s)dY_s\\
=&\int_0^t\int_0^1\int_0^\lambda\Big\{2h_x(X_s^{\varepsilon,\rho},u_s^\varepsilon)(X_s^\varepsilon-X_s)(L_s^\varepsilon-L_s)+h_{xx}(X_s^{\varepsilon,\rho},u_s^\varepsilon)L_s^{
\varepsilon,\rho}(X_s^\varepsilon-X_s)^2\Big\}d\rho d\lambda dY_s\\
&+\int_0^t\Big\{h(X_s,u_s^\varepsilon)\big(L_s^{\varepsilon}-(L_s+K_s^{1,\varepsilon})\big)+h_x(X_s,u_s^\varepsilon)L_s\big(X_s^{\varepsilon}-(X_s+Y_s^{1,\varepsilon})\big)\Big\}dY_s\\
&+\int_0^t\big(\delta h(s)K_s^{1,\varepsilon}+\delta h_x(s)L_s Y_s^{1,\varepsilon}\big)\mathbf{1}_{E_\varepsilon}(s)dY_s,\ t\in[0,T],
\end{split}
\end{equation*}
and with arguments similar to those for $X_t^\varepsilon-(X_t+Y_t^{1,\varepsilon})$ we obtain:
\begin{equation}\label{2ndVE2}
\begin{split}
&L_t^\varepsilon-(L_t+K_t^{1,\varepsilon})=\int_0^t\Big\{h(s)\big(L_s^{\varepsilon}-(L_s+K_s^{1,\varepsilon})\big)+h_x(s)L_s\big(X_s^{\varepsilon}-(X_s+Y_s^{1,\varepsilon})\big)\Big\}dY_s\\
&\ \hskip3cm+\int_0^t\Big\{h_x(s)Y_s^{1,\varepsilon}K_s^{1,\varepsilon}+\frac{1}{2}h_{xx}(s)L_s(Y_s^{1,\varepsilon})^2\Big\}dY_s\\
&\ \hskip3cm+\int_0^t\big(\delta h(s)K_s^{1,\varepsilon}+\delta h_x(s)L_s Y_s^{1,\varepsilon}\big)\mathbf{1}_{E_\varepsilon}(s)dY_s+R_t^{L,\varepsilon},\ t\in[0,T],
\end{split}
\end{equation}
with, for all $p\geq2$, $\displaystyle \big(E^Q[\sup_{t\in[0,T]}|R_t^{L,\varepsilon}|^p]\big)^{\frac{1}{p}}\leq C_p\varepsilon\rho_p(\varepsilon),\ \varepsilon>0,\ (\rho_p(\varepsilon)\rightarrow0\ \mbox{ as } \varepsilon\searrow0 ). $\\
\eqref{2ndVE1} and \eqref{2ndVE2} suggest to have the following second order variational equation:\\
\begin{equation}\label{2ndVE3}
\left\{
\begin{split}
&dY_t^{2,\varepsilon}=\Big\{\widetilde{E}^Q\Big[\int_0^{\widetilde{U}_t}\sigma_\mu(t,y)dy\cdot\widetilde{K}_t^{2,\varepsilon}\Big]+\widetilde{E}^Q\big[\widetilde{\sigma}_\mu(t)
\widetilde{L}_t\widetilde{V}_t^{2,\varepsilon}\big]\Big\}dB_t^1\\
&\qquad\qquad+\Big\{\widetilde{E}^Q\big[\widetilde{\sigma}_\mu(t)\widetilde{V}_t^{1,\varepsilon}\widetilde{K}_t^{1,\varepsilon}\big]+\frac{1}{2}\widetilde{E}^Q\big[
\widetilde{\sigma}_{z\mu}(t)\widetilde{L}_t(\widetilde{V}_t^{1,\varepsilon})^2\big]\Big\}dB_t^1\\
&\qquad+\Big\{\widetilde{E}^Q\Big[\int_0^{\widetilde{U}_t}\delta\big(\sigma_\mu(t,y)\big)dy\cdot\widetilde{K}_t^{1,\varepsilon}\Big]+\widetilde{E}^Q\big[\delta\big(
\widetilde{\sigma}_\mu(t)\big)\widetilde{L}_t\widetilde{V}_t^{1,\varepsilon}\big]\Big\}\mathbf{1}_{E_\varepsilon}(t)dB_t^1,\ t\in[0,T],\\
&Y_0^{2,\varepsilon}=0; \\
&dK_t^{2,\varepsilon}=\Big\{h(t)K_t^{2,\varepsilon}+h_x(t)L_t Y_t^{2,\varepsilon}\Big\}dY_t+\Big\{h_x(t)Y_t^{1,\varepsilon}K_t^{1,\varepsilon}+\frac{1}{2}h_{xx}(t)L_t (Y_t^{1,\varepsilon})^2\Big\}dY_t\\
&\qquad\qquad+\big(\delta h(t)K_t^{1,\varepsilon}+\delta h_x(t)L_t Y_t^{1,\varepsilon}\big)\mathbf{1}_{E_\varepsilon}(t)dY_t,\ t\in[0,T],\\
&K_0^{2,\varepsilon}=0,
\end{split}\right.
\end{equation}
and
\begin{equation}\label{xingxingxing}
\begin{split}
V_t^{2,\varepsilon}=&\frac{E^Q[L_t Y_t^{2,\varepsilon}+X_t K_t^{2,\varepsilon}\,|\,\mathcal{F}_t^Y]}{E^Q[L_t\,|\,\mathcal{F}_t^Y]}-\frac{E^{Q}[L_t X_t\,|\,\mathcal{F}_t^Y] E^{Q}[K_t^{2,\varepsilon}\,|\,
\mathcal{F}_t^Y]}{(E^{Q}[L_t\,|\,\mathcal{F}_t^Y])^2}+\frac{E^Q[K_t^{1,\varepsilon}Y_t^{1,\varepsilon}\,|\,\mathcal{F}_t^Y]}{E^Q[L_t\,|\,\mathcal{F}_t^Y]}\\
&-\frac{E^Q[K_t^{1,\varepsilon}\,|\,\mathcal{F}_t^Y]}{E^Q[L_t\,|\,\mathcal{F}_t^Y]}\bigg(\frac{E^Q[L_t Y_t^{1,\varepsilon}+X_t K_t^{1,\varepsilon}\,|\,\mathcal{F}_t^Y]}{E^Q[L_t\,|\,\mathcal{F}_t^Y]}-
\frac{E^{Q}[L_t X_t\,|\,\mathcal{F}_t^Y] E^{Q}[K_t^{1,\varepsilon}\,|\,\mathcal{F}_t^Y]}{(E^{Q}[L_t\,|\,\mathcal{F}_t^Y])^2}\bigg)\\
\bigg(=&\theta_t(Y_t^{2,\varepsilon},K_t^{2,\varepsilon})+\frac{E^Q[K_t^{1,\varepsilon}Y_t^{1,\varepsilon}\,|\,\mathcal{F}_t^Y]}{E^Q[L_t\,|\,\mathcal{F}_t^Y]}-\frac{E^Q[K_t^{1,
\varepsilon}\,|\,\mathcal{F}_t^Y]}{E^Q[L_t\,|\,\mathcal{F}_t^Y]}\theta_t(Y_t^{1,\varepsilon},K_t^{1,\varepsilon})\bigg).
\end{split}
\end{equation}
With standard arguments we get the following statement:
\begin{lemma}
Under Assumption (H2), \eqref{2ndVE3} has a unique solution $(Y^{2,\varepsilon},K^{2,\varepsilon})\in S^2_{\mathbb{F}}([0,T],Q)\times S^2_{\mathbb{F}}([0,T],Q)$. Moreover,
$Y^{2,\varepsilon},\,K^{2,\varepsilon}\in S^{\infty -}_{\mathbb{F}}([0,T],Q)$ with $S^p_{\mathbb{F}}([0,T],Q)$-bounds independent of $\varepsilon>0$, for all $p\geq2$.
\end{lemma}
Let us now give an important and very technical but somehow also beautiful result.
\begin{proposition}\label{tech2}
For all $p\geq1$, there is a constant $C_p>0$ such that for $t\in[0,T],\ \varepsilon>0$,
$$ \Big(\!E^Q\big[\big|\big(U_t^\varepsilon-(U_t+V_t^{1,\varepsilon}+V_t^{2,\varepsilon})\big)-\theta_t\big(X_t^\varepsilon-(X_t+Y_t^{1,\varepsilon}+Y_t^{2,\varepsilon}),L_t^\varepsilon-
(L_t+K_t^{1,\varepsilon}+K_t^{2,\varepsilon})\big)\big|^p\big]\!\Big)^{\frac{1}{p}}\!\leq\! C_p\varepsilon^{\frac{3}{2}}. $$
\end{proposition}
\begin{proof}
Recall the notation $(X^{\varepsilon,\lambda},L^{\varepsilon,\lambda},U^{\varepsilon,\lambda}):=(1-\lambda)(X,L,U)+\lambda(X^{\varepsilon},L^{\varepsilon},U^{\varepsilon}),\ \lambda\in[0,1]$. Then
\begin{equation*}
\begin{split}
K_t^\varepsilon:=&U_t^\varepsilon-(U_t+V_t^{1,\varepsilon}+V_t^{2,\varepsilon})=\int_0^1\partial_\lambda\Big[\frac{E^{Q}[L_t^{\varepsilon,\lambda}X_t^{\varepsilon,\lambda}\,|\,
\mathcal{F}_t^Y]}{E^{Q}[L_t^{\varepsilon,\lambda}\,|\,\mathcal{F}_t^Y]}\Big]d\lambda-V_t^{1,\varepsilon}-V_t^{2,\varepsilon}\\
=&\int_0^1\Big\{\frac{E^{Q}[L_t^{\varepsilon,\lambda}(X_t^\varepsilon-X_t)+X_t^{\varepsilon,\lambda}(L_t^\varepsilon-L_t)\,|\,\mathcal{F}_t^Y]}{E^{Q}[L_t^{\varepsilon,\lambda}\,|\,
\mathcal{F}_t^Y]}-\frac{E^{Q}[L_t^{\varepsilon,\lambda}X_t^{\varepsilon,\lambda}\,|\,\mathcal{F}_t^Y]}{(E^{Q}[L_t^{\varepsilon,\lambda}\,|\,\mathcal{F}_t^Y])^2}E^Q[L_t^\varepsilon-
L_t\,|\,\mathcal{F}_t^Y]\Big\}d\lambda\\
&-V_t^{1,\varepsilon}-V_t^{2,\varepsilon}.
\end{split}
\end{equation*}
We also recall that
$$ V_t^{1,\varepsilon}=\theta_t(Y_t^{1,\varepsilon},K_t^{1,\varepsilon})=\frac{E^Q[L_t Y_t^{1,\varepsilon}+X_t K_t^{1,\varepsilon}\,|\,\mathcal{F}_t^Y]}{E^Q[L_t\,|\,\mathcal{F}_t^Y]}-\frac{E^{Q}[L_t X_t\,|\,\mathcal{F}_t^Y]}{(E^{Q}[L_t\,|\,\mathcal{F}_t^Y])^2}E^{Q}[K_t^{1,\varepsilon}\,|\,\mathcal{F}_t^Y]. $$
Hence,
\begin{equation*}
\begin{split}
K_t^\varepsilon=\!&\int_0^1\!\int_0^\lambda\!\partial_\rho\Big\{\frac{E^{Q}[L_t^{\varepsilon,\rho}(X_t^\varepsilon-X_t)\!+\!X_t^{\varepsilon,\rho}(L_t^\varepsilon-L_t)\,|\,\mathcal{F}_t^Y]}
{E^{Q}[L_t^{\varepsilon,\rho}\,|\,\mathcal{F}_t^Y]}\!-\!\frac{E^{Q}[L_t^{\varepsilon,\rho}X_t^{\varepsilon,\rho}\,|\,\mathcal{F}_t^Y]}{(E^{Q}[L_t^{\varepsilon,\rho}\,|\,\mathcal{F}_t^Y])^2}
E^Q[L_t^\varepsilon\!-\! L_t\,|\,\mathcal{F}_t^Y]\Big\}d\rho d\lambda\\
&\qquad+\theta_t(X_t^\varepsilon-X_t-Y_t^{1,\varepsilon},L_t^\varepsilon-L_t-K_t^{1,\varepsilon})-V_t^{2,\varepsilon}\\
=&\int_0^1\int_0^\lambda\Big\{\frac{2E^{Q}[(L_t^\varepsilon-L_t)(X_t^\varepsilon-X_t)\,|\,\mathcal{F}_t^Y]}{E^{Q}[L_t^{\varepsilon,\rho}\,|\,\mathcal{F}_t^Y]}-\frac{2E^{Q}
[L_t^{\varepsilon,\rho}
(X_t^\varepsilon-X_t)+X_t^{\varepsilon,\rho}(L_t^\varepsilon-L_t)\,|\,\mathcal{F}_t^Y]}{(E^{Q}[L_t^{\varepsilon,\rho}\,|\,\mathcal{F}_t^Y])^2}\\
&\displaystyle\qquad\qquad \times E^{Q}[L_t^\varepsilon-L_t\,|\,
\mathcal{F}_t^Y]+\frac{2E^{Q}[L_t^{\varepsilon,\rho}X_t^{\varepsilon,\rho}\,|\,\mathcal{F}_t^Y]}{(E^{Q}[L_t^{\varepsilon,\rho}\,|\,\mathcal{F}_t^Y])^3}(E^{Q}[L_t^\varepsilon-L_t\,|\,
\mathcal{F}_t^Y])^2\Big\}d\rho d\lambda \\
&\displaystyle\qquad+\theta_t(X_t^\varepsilon-X_t-Y_t^{1,\varepsilon},L_t^\varepsilon-L_t-K_t^{1,\varepsilon})-V_t^{2,\varepsilon}.
\end{split}
\end{equation*}
As $\theta_t\big(X_t^\varepsilon-(X_t+Y_t^{1,\varepsilon}+Y_t^{2,\varepsilon}),L_t^\varepsilon-(L_t+K_t^{1,\varepsilon}+K_t^{2,\varepsilon})\big)=\theta_t\big(X_t^\varepsilon-(X_t+Y_t^{1,
\varepsilon}),L_t^\varepsilon-(L_t+K_t^{1,\varepsilon})\big)-\theta_t(Y_t^{2,\varepsilon},K_t^{2,\varepsilon})$ (recall that $(\zeta,\eta)\rightarrow\theta_t(\zeta,\eta)$ is linear, see the definition (\ref{theta})) and
\begin{equation*}
\begin{split}
V_t^{2,\varepsilon}=&\theta_t(Y_t^{2,\varepsilon},K_t^{2,\varepsilon})+\frac{E^Q[K_t^{1,\varepsilon}Y_t^{1,\varepsilon}\,|\,\mathcal{F}_t^Y]}{E^Q[L_t\,|\,\mathcal{F}_t^Y]}\\
&-\frac{E^Q[L_t Y_t^{1,\varepsilon}+X_t K_t^{1,\varepsilon}\,|\,\mathcal{F}_t^Y]}{(E^Q[L_t\,|\,\mathcal{F}_t^Y])^2}E^Q[K_t^{1,\varepsilon}\,|\,\mathcal{F}_t^Y]
+\frac{E^{Q}[L_t X_t\,|\,\mathcal{F}_t^Y]}{(E^{Q}[L_t\,|\,\mathcal{F}_t^Y])^3}(E^{Q}[K_t^{1,\varepsilon}\,|\,\mathcal{F}_t^Y])^2,
\end{split}
\end{equation*}
(see, \eqref{xingxingxing}), we get
\begin{equation}\label{sqrtri}
\begin{split}
K_t^\varepsilon=&\theta_t\big(X_t^\varepsilon-(X_t+Y_t^{1,\varepsilon}+Y_t^{2,\varepsilon}),L_t^\varepsilon-(L_t+K_t^{1,\varepsilon}+K_t^{2,\varepsilon})\big)+I_t^1(\varepsilon)+I_t^2
(\varepsilon)+I_t^3(\varepsilon),
\end{split}
\end{equation}
where
\begin{equation*}
\begin{split}
I_t^1(\varepsilon)=&\int_0^1\int_0^\lambda2\Big(\frac{E^{Q}[(L_t^\varepsilon-L_t)(X_t^\varepsilon-X_t)\,|\,\mathcal{F}_t^Y]}{E^{Q}[L_t^{\varepsilon,\rho}\,|\,\mathcal{F}_t^Y]}-
\frac{E^Q[Y_t^{1,\varepsilon}K_t^{1,\varepsilon}\,|\,\mathcal{F}_t^Y]}{E^Q[L_t\,|\,\mathcal{F}_t^Y]}\Big)d\rho d\lambda,\\
I_t^2(\varepsilon)=&-\int_0^1\int_0^\lambda2\Big(\frac{E^{Q}[L_t^{\varepsilon,\rho}(X_t^\varepsilon-X_t)+X_t^{\varepsilon,\rho}(L_t^\varepsilon-L_t)\,|\,\mathcal{F}_t^Y]}{(E^{Q}
[L_t^{\varepsilon,\rho}\,|\,\mathcal{F}_t^Y])^2}E^{Q}[L_t^\varepsilon-L_t\,|\,\mathcal{F}_t^Y]\\
&\qquad\qquad\quad-\frac{E^Q[L_t Y_t^{1,\varepsilon}+X_t K_t^{1,\varepsilon}\,|\,\mathcal{F}_t^Y]}{(E^Q[L_t\,|\,\mathcal{F}_t^Y])^2}E^Q[K_t^{1,\varepsilon}\,|\,\mathcal{F}_t^Y]\Big)d\rho d\lambda,\\
I_t^3(\varepsilon)=&\int_0^1\int_0^\lambda2\Big(\frac{E^{Q}[L_t^{\varepsilon,\rho}X_t^{\varepsilon,\rho}\,|\,\mathcal{F}_t^Y]}{(E^{Q}[L_t^{\varepsilon,\rho}\,|\,\mathcal{F}_t^Y])^3}
(E^{Q}[L_t^\varepsilon-L_t\,|\,\mathcal{F}_t^Y])^2-\frac{E^{Q}[L_t X_t\,|\,\mathcal{F}_t^Y]}{(E^{Q}[L_t\,|\,\mathcal{F}_t^Y])^3}(E^{Q}[K_t^{1,\varepsilon}\,|\,\mathcal{F}_t^Y])^2 \Big)d\rho d\lambda.
\end{split}
\end{equation*}
Let us estimate them, and we begin with $I_t^1(\varepsilon)$. From Proposition \ref{EstofXL} and the fact that $L^{-1},\,(L^{\varepsilon,\rho})^{-1}$, $\varepsilon>0$, $\rho\in[0,1]$, are bounded in $S^p_{\mathbb{F}}([0,T],Q)$, uniformly w.r.t. $\varepsilon$ and $\rho$, we have
$$ \Big(E^Q\Big[\Big|\frac{E^{Q}[(L_t^\varepsilon-L_t)(X_t^\varepsilon-X_t-Y_t^{1,\varepsilon})\,|\,\mathcal{F}_t^Y]}{E^{Q}[L_t^{\varepsilon,\rho}\,|\,\mathcal{F}_t^Y]}\Big|^p\Big]\Big)^{\frac{1}{p}}\leq C_p\varepsilon^{\frac{3}{2}}, $$
$$ \Big(E^Q\Big[\Big|\frac{E^{Q}[(L_t^\varepsilon-L_t-K_t^{1,\varepsilon})Y_t^{1,\varepsilon}\,|\,\mathcal{F}_t^Y]}{E^{Q}[L_t^{\varepsilon,\rho}\,|\,\mathcal{F}_t^Y]}\Big|^p\Big]\Big)^{\frac{1}{p}}\leq C_p\varepsilon^{\frac{3}{2}}, $$
and
\begin{equation*}
\begin{split}
&\Big(E^Q\Big[\Big|\Big(\frac{1}{E^{Q}[L_t^{\varepsilon,\rho}\,|\,\mathcal{F}_t^Y]}-\frac{1}{E^{Q}[L_t\,|\,\mathcal{F}_t^Y]}\Big)E^Q[K_t^{1,\varepsilon}Y_t^{1,\varepsilon}\,|\,
\mathcal{F}_t^Y]\Big|^p\Big]\Big)^{\frac{1}{p}}\\
&\leq \Big(E^Q\Big[\Big|\frac{E^{Q}[L_t^\varepsilon-L_t\,|\,\mathcal{F}_t^Y]}{E^{Q}[L_t^{\varepsilon,\rho}\,|\,\mathcal{F}_t^Y]E^{Q}[L_t\,|\,\mathcal{F}_t^Y]}\Big|^{2p}\Big]\Big)^{\frac{1}{2p}} \big(E^Q[|K_t^{1,\varepsilon}Y_t^{1,\varepsilon}|^{2p}]\big)^{\frac{1}{2p}} \leq C_p\varepsilon^{\frac{3}{2}},
\end{split}
\end{equation*}
and, hence, $\displaystyle \big(E^Q[|I_t^1(\varepsilon)|^p]\big)^{\frac{1}{p}}\leq C_p\varepsilon^{\frac{3}{2}},\ t\in[0,T],\ \varepsilon>0,\ p\geq1. $\
Similarly, but easier, also for $I_t^2(\varepsilon)$ and $I_t^3(\varepsilon)$, we obtain $\displaystyle \big(E^Q[|I_t^i(\varepsilon)|^p]\big)^{\frac{1}{p}}\leq C_p\varepsilon^{\frac{3}{2}},\ t\in[0,T],\ \varepsilon>0,\ p\geq1,\ i=2,3. $

Let us now define
\begin{equation*}
\begin{split}
&\eta_t^\varepsilon:=X_t^\varepsilon-(X_t+Y_t^{1,\varepsilon}+Y_t^{2,\varepsilon}),\\
&\vartheta_t^\varepsilon:=L_t^\varepsilon-(L_t+K_t^{1,\varepsilon}+K_t^{2,\varepsilon}),\ t\in[0,T],\ \varepsilon>0.
\end{split}
\end{equation*}
Then, due to \eqref{sqrtri} combined with the above estimates,
\begin{equation}\label{2ndVE4}
K_t^\varepsilon=U_t^\varepsilon-(U_t+V_t^{1,\varepsilon}+V_t^{2,\varepsilon})=\theta_t(\eta_t^\varepsilon,\vartheta_t^\varepsilon)+\rho_t^\varepsilon,\ t\in[0,T],\ \varepsilon>0,
\end{equation}
with $\big(E^Q[|\rho_t^\varepsilon|^p]\big)^{\frac{1}{p}}\leq C_p\varepsilon^{\frac{3}{2}},\ t\in[0,T],\ \varepsilon>0,\ p\geq2$. The proof is complete now.
\end{proof}

As a consequence of above computations, we can state the following proposition.
\begin{proposition}\label{EstofYK}
For all $p\geq2$, there exists a constant $C_p\in\mathbb{R}_+$\ such that,
 \begin{equation*}
\begin{split}
&\qquad\qquad\qquad\quad{\rm (i)}\ E^Q\big[\sup_{t\in[0,T]}\big|X_t^\varepsilon-(X_t+Y_t^{1,\varepsilon}+Y_t^{2,\varepsilon})\big|^{p}\big]\leq C_p\varepsilon^{p}\rho_p(\varepsilon);\\
&\qquad\qquad\qquad\quad{\rm (ii)}\ E^Q\big[\sup_{t\in[0,T]}\big|L_t^\varepsilon-(L_t+K_t^{1,\varepsilon}+K_t^{2,\varepsilon})\big|^{p}\big]\leq C_p\varepsilon^{p}\rho_p(\varepsilon);\\
&\qquad\qquad\qquad\quad{\rm (iii)}\ E^Q\big[\sup_{t\in[0,T]}\big|U_t^\varepsilon-(U_t+V_t^{1,\varepsilon}+V_t^{2,\varepsilon})\big|^{p}\big]\leq C_p\varepsilon^{p}\rho_p(\varepsilon),\hspace{8cm}
\end{split}
\end{equation*}
with $\rho_p(\varepsilon)\rightarrow0$, as $\varepsilon\searrow0$. Moreover,
$${\rm (iv)}\ E^Q\big[\sup_{t\in[0,T]}|Y_t^{2,\varepsilon}|^p+|K_t^{2,\varepsilon}|^p\big]\leq C_p\varepsilon^{p};\quad  E^Q\big[|V_t^{2,\varepsilon}|^p\big]\leq C_p\varepsilon^{p},\ \varepsilon>0,\ t\in[0,T]. $$
\end{proposition}
\begin{proof}
The estimates (i)-(iii) are a consequence of Proposition \ref{tech2} combined with Burkholder-Davis-Gundy's inequality. Indeed, following the proof of Proposition \ref{tech2}, from \eqref{2ndVE1}, \eqref{2ndVE2}, \eqref{2ndVE4} and the second order variational equation, we have
\begin{equation*}
\begin{split}
&\eta_t^\varepsilon=\int_0^t\Big\{\widetilde{E}^Q\Big[\int_0^{\widetilde{U}_s}\sigma_\mu(s,y)dy\cdot\widetilde{\vartheta}_s^\varepsilon\Big]+\widetilde{E}^Q\big[\widetilde{\sigma}_\mu(s)
\widetilde{L}_s\widetilde{\theta_s(\eta_s^\varepsilon,\vartheta_s^\varepsilon)}\big]\Big\}dB_s^1+R_t^1(\varepsilon),\\
&\vartheta_t^\varepsilon=\int_0^t\big(h(s)\vartheta_s^\varepsilon+h_x(s)L_s\eta_s^\varepsilon\big)dY_s+R_t^2(\varepsilon),\ t\in[0,T],
\end{split}
\end{equation*}
where $\displaystyle \big(E^Q[\sup_{t\leq T}|R_t^i(\varepsilon)|^p]\big)^{\frac{1}{p}}\leq C_p\varepsilon\rho_p(\varepsilon),\ \varepsilon>0,\ p\geq2,\ i=1,2$.
Morever, recalling the estimate \eqref{Estoftheta}, we have
\begin{equation*}
\begin{split}
&\Big(E^Q\big[L_s^q\big|\theta_s(\eta_s^\varepsilon,\vartheta_s^\varepsilon)\big|^q\big]\Big)^{\frac{1}{q}}=\Big(E^Q\big[L_s^q\big|\theta_s(\eta_s^\varepsilon,\vartheta_s^\varepsilon)
-\theta_s(0,0)\big|^q\big]\Big)^{\frac{1}{q}}\\
&\leq C\big(E^Q\big[|\eta_s^\varepsilon|^{2q}+|\vartheta_s^\varepsilon|^{2q}\big]\big)^{\frac{1}{2q}},\ s\in[0,T],\ \varepsilon>0,\ q\geq1.
\end{split}
\end{equation*}
Thus, for all $p\geq2$,
\begin{equation*}
\begin{split}
E^Q[\sup_{s\leq t}|\eta_s^\varepsilon|^{2p}]&\leq C_p\int_0^t\big\{\big(E^Q[|\vartheta_s^\varepsilon|^{p}]\big)^2+E^Q[|\eta_s^\varepsilon|^{2p}]\big\}ds+C_p\varepsilon^{2p}\rho_p(\varepsilon),\\
\big(E^Q[\sup_{s\leq t}|\vartheta_s^\varepsilon|^{p}]\big)^2&\leq C_p\int_0^t\big\{\big(E^Q[|\vartheta_s^\varepsilon|^{p}]\big)^2+E^Q[|\eta_s^\varepsilon|^{2p}]\big\}ds+C_p\varepsilon^{2p}\rho_p(\varepsilon),
\end{split}
\end{equation*}
and from Gronwall's Lemma we get $\displaystyle \big(E^Q[\sup_{s\leq T}|\eta_s^\varepsilon|^{2p}]\big)^{\frac{1}{2}}+E^Q[\sup_{s\leq T}|\vartheta_s^\varepsilon|^{p}]\leq C_p\varepsilon^p\rho_p(\varepsilon),\
\varepsilon>0, $\
for $p\geq2$, and $\rho_p(\varepsilon)\rightarrow0$, as $\varepsilon\searrow0$. Finally, from \eqref{2ndVE4}, it follows that
$$\displaystyle\begin{array}{rcl}
& \big(E^Q[\sup_{t\in[0, T]}|K_t^\varepsilon|^p]\big)^{\frac{1}{p}}\leq\big(E^Q[\sup_{t\in[0, T]}|\theta_t(\eta_t^\varepsilon,\vartheta_t^\varepsilon)|^p]\big)^{\frac{1}{p}}+C_p \varepsilon^{\frac{3}{2}}\\
&\leq C_p\big(E^Q\big[\sup_{t\in[0, T]}|
\eta_t^\varepsilon|^{2p}+\sup_{t\in[0, T]}|\vartheta_t^\varepsilon|^{2p}\big]\big)^{\frac{1}{2p}}+C_p \varepsilon^{\frac{3}{2}}\leq C\varepsilon\rho_p(\varepsilon).
\end{array}$$
Finally, (iv) is an immediate consequence of (i), Proposition \ref{EstofXL} and
Corollary \ref{EstofUV}.
\end{proof}
This proposition also holds for the general case, and the proof is based on a similar approach using the ideas developed above, hence we
can omit it here. In particular, we obtain for the general case the following second order variational equation:
\begin{equation}\label{2ndVE5}
\left\{
\begin{split}
&dY_t^{2,\varepsilon}=\Big\{\sigma_x(t)Y_t^{2,\varepsilon}+\frac{1}{2}\sigma_{xx}(t)(Y_t^{1,\varepsilon})^2+\widetilde{E}^Q\Big[\int_0^{\widetilde{U}_t}\sigma_\mu(t,y)dy\cdot
\widetilde{K}_t^{2,\varepsilon}\Big]+\widetilde{E}^Q\big[\widetilde{\sigma}_\mu(t)\widetilde{L}_t\widetilde{V}_t^{2,\varepsilon}\big]\\
&\qquad\qquad+\widetilde{E}^Q\big[\widetilde{\sigma}_\mu(t)\widetilde{V}_t^{1,\varepsilon}\widetilde{K}_t^{1,\varepsilon}\big]+\frac{1}{2}\widetilde{E}^Q\big[\widetilde{\sigma}_{z\mu}
(t)\widetilde{L}_t(\widetilde{V}_t^{1,\varepsilon})^2\big]\\
&\qquad\qquad+\Big(\delta\sigma_x(t)Y_t^{1,\varepsilon}+\widetilde{E}^Q\Big[\int_0^{\widetilde{U}_t}\delta\sigma_\mu(t,y)dy\cdot\widetilde{K}_t^{1,\varepsilon}\Big]+\widetilde{E}^Q
\big[\delta\widetilde{\sigma}_\mu(t)\widetilde{L}_t\widetilde{V}_t^{1,\varepsilon}\big]\Big)\mathbf{1}_{E_\varepsilon}(t)\Big\}dB_t^1,\ t\in[0,T],\\
&Y_0^{2,\varepsilon}=0; \\
&dK_t^{2,\varepsilon}=\Big\{h(t)K_t^{2,\varepsilon}+h_x(t)L_t Y_t^{2,\varepsilon}+h_x(t)Y_t^{1,\varepsilon}K_t^{1,\varepsilon}+\frac{1}{2}h_{xx}(t)L_t (Y_t^{1,\varepsilon})^2\\
&\qquad\qquad+L_t\widetilde{E}^Q\Big[\int_0^{\widetilde{U}_t}h_\mu(t,y)dy\cdot\widetilde{K}_t^{2,\varepsilon}\Big]+L_t\widetilde{E}^Q\big[\widetilde{h}_\mu(t)\widetilde{L}_t
\widetilde{V}_t^{2,\varepsilon}\big]\\
&\qquad\qquad+L_t\widetilde{E}^Q\big[\widetilde{h}_\mu(t)\widetilde{V}_t^{1,\varepsilon}\widetilde{K}_t^{1,\varepsilon}\big]+\frac{1}{2}L_t\widetilde{E}^Q\big[\widetilde{h}_{z\mu}(t)
\widetilde{L}_t(\widetilde{V}_t^{1,\varepsilon})^2\big]\\
&\qquad\qquad+\Big(\delta h(t)K_t^{1,\varepsilon}+\delta h_x(t)L_t Y_t^{1,\varepsilon}+L_t\widetilde{E}^Q\Big[\int_0^{\widetilde{U}_t}\delta h_\mu(t,y)dy\cdot\widetilde{K}_t^{1,\varepsilon}\Big]\\
&\qquad\qquad\quad\ \  +L_t\widetilde{E}^Q\big[\delta\widetilde{h}_\mu(t)\widetilde{L}_t\widetilde{V}_t^{1,\varepsilon}\big]\Big)\mathbf{1}_{E_\varepsilon}(t)\Big\}dY_t,\ t\in[0,T],\\
&K_0^{2,\varepsilon}=0.
\end{split}\right.
\end{equation}
(Recall that ${V}_t^{1,\varepsilon}$ and ${V}_t^{2,\varepsilon}$\ are defined in (\ref{SMP2}) and (\ref{xingxingxing}), respectively.)
\begin{proposition}\label{Wellof2ndVar}
Under Assumption (H2), equation \eqref{2ndVE5} has a unique solution $(Y^{2,\varepsilon},K^{2,\varepsilon})\in S^2_{\mathbb{F}}([0,T],Q)\times S^2_{\mathbb{F}}([0,T],Q)$. Moreover, $Y^{2,\varepsilon},\,K^{2,\varepsilon}$, $\varepsilon>0$, are bounded in $S^p_{\mathbb{F}}([0,T],Q)$, for all $p\geq2$.
\end{proposition}

\subsection{Duality}
It is well known that for deducing the stochastic maximum principle we shall make use of the duality relations between the solutions of the variational equations and the adjoint backward SDEs. The first order adjoint BSDE over our reference probability space $(\Omega,\mathcal{F},Q)$ is the following:
%\begin{equation*}\left\{
%\begin{split}
%&p_t^1=p_T^1+\int_t^T\alpha_s ds-\int_t^T q_s^1dB_s^1,\ t\in[0,T];\\
%&p_t^2=p_T^2+\int_t^T\beta_s ds-\int_t^T q_s^2dY_s,\ t\in[0,T].
%\end{split}\right.
%\end{equation*}
%where processes $\alpha$, $\beta$ and terminal conditions $p_T^1$, $p_T^2$ are to be confirmed.
\begin{equation}\label{SMP4}\left\{
\begin{split}
&dp_t^1=-\alpha_t(q_t^1,q_t^2)dt+q_t^1dB_t^1+\check{q}_t^1dY_t,\ t\in[0,T],\\
&p_T^1=-\Phi_x(T)-L_T\widetilde{E}^Q\big[E^P[\widetilde{\Phi}^*_\mu(T)\,|\,\mathcal{F}_T^Y]\big];\\
&dp_t^2=-\beta_t(q_t^1,q_t^2)dt+\check{q}_t^2dB_t^1+q_t^2dY_t,\ t\in[0,T],\\
&p_T^2=-(X_T-U_T)\widetilde{E}^Q\big[E^P[\widetilde{\Phi}^*_\mu(T)\,|\,\mathcal{F}_T^Y]\big]-\widetilde{E}^Q\Big[\int_0^{U_T}\Phi^*_\mu(T,y)
dy\Big],
\end{split}\right.
\end{equation}
where
\begin{align*}
\alpha_t^0(q_t^1,q_t^2):=&\sigma_x(t) q_t^1+L_t\widetilde{E}^Q\big[\widetilde{q}_t^1 E^P[\widetilde{\sigma}^*_\mu(t)\,|\,\mathcal{F}_t^Y]\big]+h_x(t) L_t q_t^2+L_t\widetilde{E}^Q\big[\widetilde{q}_t^2 \widetilde{L}_t E^P[\widetilde{h}^*_\mu(t) \,|\,\mathcal{F}_t^Y]\big];\\
\alpha_t(q_t^1,q_t^2):=&\alpha_t^0(q_t^1,q_t^2)-f_x(t)-L_t\widetilde{E}^Q\big[E^P[\widetilde{f}^*_\mu(t)\,|\,\mathcal{F}_t^Y]\big];\\
\beta_t^0(q_t^1,q_t^2):=&(X_t-U_t)\widetilde{E}^Q\big[\widetilde{q}_t^1 E^P[\widetilde{\sigma}^*_\mu(t)\,|\,\mathcal{F}_t^Y]\big]+\widetilde{E}^Q\big[\widetilde{q}_t^1 \int_0^{U_t}\sigma^*_\mu(t,y)dy\big]\\
&+h(t)q_t^2+(X_t-U_t)\widetilde{E}^Q\big[\widetilde{q}_t^2 \widetilde{L}_t E^P[\widetilde{h}^*_\mu(t) \,|\,\mathcal{F}_t^Y]\big]+\widetilde{E}^Q\big[\widetilde{q}_t^2 \widetilde{L}_t \int_0^{U_t}h^*_\mu(t,y)dy\big];\\
\beta_t(q_t^1,q_t^2):=&\beta_t^0(q_t^1,q_t^2)-(X_t-U_t)\widetilde{E}^Q\big[E^P[\widetilde{f}^*_\mu(t)\,|\,\mathcal{F}_t^Y]\big]-\widetilde{E}^Q\big[\int_0^{U_t}f^*_\mu(t,y)dy\big],\ t\in[0,T].
\end{align*}
Before showing the well-posedness of the first order adjoint equation, assuming the existence and the uniqueness of the solution for \eqref{SMP4} we study the duality. From It\^{o}'s formula we get (for simplicity we omit the $(q_t^1,q_t^2)$ in the generators)
\begin{equation}\label{SMP3}
\begin{split}
&dE^Q[p_t^1 Y_t^{1,\varepsilon}]=E^Q\Big[-\alpha_t Y_t^{1,\varepsilon}+q_t^1\Big\{\sigma_x(t)Y_t^{1,\varepsilon}+\widetilde{E}^Q\Big[\int_0^{\widetilde{U}_t}\sigma_\mu(t,y)dy\cdot\widetilde{K}_t^{1,\varepsilon}\Big]\\
&\qquad\qquad\qquad\qquad\qquad\qquad+\widetilde{E}^Q\big[\widetilde{\sigma}_\mu(t)\widetilde{L}_t\widetilde{V}_t^{1,\varepsilon}\big]+\delta\sigma(t)\mathbf{1}_{E_\varepsilon}(t)\Big\}\Big]dt,\\
&dE^Q[p_t^2 K_t^{1,\varepsilon}]=E^Q\Big[-\beta_t K_t^{1,\varepsilon}+q_t^2 h(t)K_t^{1,\varepsilon}+q_t^2 L_t\Big\{h_x(t)Y_t^{1,\varepsilon}+\widetilde{E}^Q\Big[\int_0^{\widetilde{U}_t}h_\mu(t,y)dy\cdot\widetilde{K}_t^{1,\varepsilon}\Big]\\
&\qquad\qquad\qquad\qquad\qquad\qquad+\widetilde{E}^Q\big[\widetilde{h}_\mu(t)\widetilde{L}_t\widetilde{V}_t^{1,\varepsilon}\big]+\delta h(t)\mathbf{1}_{E_\varepsilon}(t)\Big\}\Big]dt.\\
\end{split}
\end{equation}
Note that
\begin{equation*}
\begin{split}
&\widetilde{E}^Q\big[\widetilde{\sigma}_\mu(t)\widetilde{L}_t\widetilde{V}_t^{1,\varepsilon}\big]=\widetilde{E}^Q\Big[\widetilde{\sigma}_\mu(t)\widetilde{L}_t\Big\{\frac{\widetilde{E}^Q[\widetilde{L}_t \widetilde{Y}_t^{1,\varepsilon}+\widetilde{X}_t \widetilde{K}_t^{1,\varepsilon}\,|\,\mathcal{F}_t^{\widetilde{Y}}]}{\widetilde{E}^Q[\widetilde{L}_t\,|\,\mathcal{F}_t^{\widetilde{Y}}]}-\frac{\widetilde{E}^{Q}[\widetilde{L}_t \widetilde{X}_t\,|\,\mathcal{F}_t^{\widetilde{Y}}] }{(\widetilde{E}^{Q}[\widetilde{L}_t\,|\,\mathcal{F}_t^{\widetilde{Y}}])^2}\widetilde{E}^{Q}[\widetilde{K}_t^{1,\varepsilon}\,|\,\mathcal{F}_t^{\widetilde{Y}}]\Big\}\Big]\\
=&\widetilde{E}^Q\Big[\frac{\widetilde{E}^Q[\widetilde{\sigma}_\mu(t)\widetilde{L}_t\,|\,\mathcal{F}_t^{\widetilde{Y}}]}{\widetilde{E}^Q[\widetilde{L}_t\,|\,\mathcal{F}_t^{
\widetilde{Y}}]}\Big(\big(\widetilde{L}_t \widetilde{Y}_t^{1,\varepsilon}+\widetilde{X}_t \widetilde{K}_t^{1,\varepsilon}\big)-\frac{\widetilde{E}^{Q}[\widetilde{L}_t \widetilde{X}_t\,|\,\mathcal{F}_t^{\widetilde{Y}}] }{\widetilde{E}^{Q}[\widetilde{L}_t\,|\,\mathcal{F}_t^{\widetilde{Y}}]}\widetilde{K}_t^{1,\varepsilon}\Big)\Big]\\
=&\widetilde{E}^Q\Big[\widetilde{E}^P[\widetilde{\sigma}_\mu(t)\,|\,\mathcal{F}_t^{\widetilde{Y}}]\Big(\big(\widetilde{L}_t \widetilde{Y}_t^{1,\varepsilon}+\widetilde{X}_t \widetilde{K}_t^{1,\varepsilon}\big)-\widetilde{E}^{P}[\widetilde{X}_t\,|\,\mathcal{F}_t^{\widetilde{Y}}] \widetilde{K}_t^{1,\varepsilon}\Big)\Big].
\end{split}
\end{equation*}
Thus, from Fubini's theorem we obtain that
\begin{equation*}
\begin{split}
&E^Q\big[q_t^1 \widetilde{E}^Q[\widetilde{\sigma}_\mu(t)\widetilde{L}_t\widetilde{V}_t^{1,\varepsilon}]\big]=E^Q\Big[\Big\{\big(L_t Y_t^{1,\varepsilon}+X_t K_t^{1,\varepsilon}\big)-E^{P}[X_t\,|\,\mathcal{F}_t^Y] K_t^{1,\varepsilon}\Big\} \widetilde{E}^Q\big[\widetilde{q}_t^1 E^P[\widetilde{\sigma}^*_\mu(t) \,|\,\mathcal{F}_t^Y]\big]\Big],\\
&E^Q\big[q_t^2 L_t \widetilde{E}^Q[\widetilde{h}_\mu(t)\widetilde{L}_t\widetilde{V}_t^{1,\varepsilon}]\big]=E^Q\Big[\Big\{\big(L_t Y_t^{1,\varepsilon}+X_t K_t^{1,\varepsilon}\big)-E^{P}[X_t\,|\,\mathcal{F}_t^Y] K_t^{1,\varepsilon}\Big\} \widetilde{E}^Q\big[\widetilde{q}_t^2 \widetilde{L}_t E^P[\widetilde{h}^*_\mu(t) \,|\,\mathcal{F}_t^Y]\big]\Big],\\
\end{split}
\end{equation*}
and similarly we also have
$$ E^Q\Big[q_t^2 L_t \widetilde{E}^Q\big[\int_0^{\widetilde{U}_t}h_\mu(t,y)dy\cdot\widetilde{K}_t^{1,\varepsilon}\big] \Big]=E^Q\Big[K_t^{1,\varepsilon}\widetilde{E}^Q\big[\widetilde{q}_t^2 \widetilde{L}_t \int_0^{U_t}h^*_\mu(t,y)dy\big]\Big]. $$
Hence,
\begin{equation*}
\begin{split}
&d\big(E^Q[p_t^1 Y_t^{1,\varepsilon}]+E^Q[p_t^2 K_t^{1,\varepsilon}]\big)\\
=&E^Q\Big[Y_t^{1,\varepsilon}\Big\{-\alpha_t+\sigma_x(t) q_t^1+L_t \widetilde{E}^Q\big[\widetilde{q}_t^1 E^P[\widetilde{\sigma}^*_\mu(t) \,|\,\mathcal{F}_t^Y]\big]+h_x(t)q_t^2 L_t\\
&\qquad\qquad+L_t\widetilde{E}^Q\big[\widetilde{q}_t^2 \widetilde{L}_t E^P[\widetilde{h}^*_\mu(t) \,|\,\mathcal{F}_t^Y]\big]\Big\}+K_t^{1,\varepsilon}\Big\{-\beta_t+\widetilde{E}^Q\big[\widetilde{q}_t^1 \int_0^{U_t}\sigma^*_\mu(t,y)dy\big]\\
&\qquad\qquad+\big(X_t-E^P[X_t\,|\,\mathcal{F}_t^Y]\big)\widetilde{E}^Q\big[\widetilde{q}_t^1 E^P[\widetilde{\sigma}^*_\mu(t) \,|\,\mathcal{F}_t^Y]\big]+h(t)q_t^2+\widetilde{E}^Q\big[\widetilde{q}_t^2 \widetilde{L}_t \int_0^{U_t}h^*_\mu(t,y)dy\big]\\
&\qquad\qquad+\big(X_t-E^P[X_t\,|\,\mathcal{F}_t^Y]\big)\widetilde{E}^Q\big[\widetilde{q}_t^2 \widetilde{L}_t E^P[\widetilde{h}^*_\mu(t) \,|\,\mathcal{F}_t^Y]\big]\Big\}+\big(q_t^1\delta\sigma(t)+q_t^2 L_t\delta h(t)\big)\mathbf{1}_{E_\varepsilon}(t)\Big]dt.
\end{split}
\end{equation*}
Recall that $U_t=E^P[X_t\,|\,\mathcal{F}_t^Y]$. Using the definition of $\alpha_t$ and  $\beta_t$ we get the following duality relation
\begin{equation}\label{eq4.29star}
\begin{split}
&E^Q[p_T^1 Y_T^{1,\varepsilon}+p_T^2 K_T^{1,\varepsilon}]\\
=&E^Q\Big[\int_0^T\Big\{Y_t^{1,\varepsilon}\big(f_x(t)+L_t\widetilde{E}^Q\big[E^P[\widetilde{f}^*_\mu(t)\,|\,\mathcal{F}_t^Y]\big]\big)\\
&\qquad\qquad\qquad+K_t^{1,\varepsilon}\Big((X_t-U_t)\widetilde{E}^Q\big[E^P[\widetilde{f}^*_\mu(t)\,|\,\mathcal{F}_t^Y]\big]+\widetilde{E}^Q\big[\int_0^{U_t}f^*_\mu(t,y)dy\big]\Big)\\
&\qquad\qquad\qquad+\big(q_t^1\delta\sigma(t)+q_t^2 L_t\delta h(t)\big)\mathbf{1}_{E_\varepsilon}(t)\Big\}dt\Big].
\end{split}
\end{equation}

As the mean-field BSDE \eqref{SMP4} does not have Lipschitz coefficients, to the best of our knowledge, it is beyond the existing frameworks of BSDEs, so we need the following result.
\begin{proposition}{\label{4.777}}
Under Assumption (H2), BSDE \eqref{SMP4} has a unique strong solution $\big((p^1,(q^1,\check{q}^1)),\\ (p^2,(\check{q}^2,q^2))\big)$. Furthermore, for any $p\geq 2$, it holds that $\big((p^1,(q^1,\check{q}^1)),\, (p^2,(\check{q}^2,q^2))\big)\in \big(S_{\mathbb{F}}^p([0,T],Q) \times (L_{\mathbb{F}}^p([0,T],Q))^2\big)\times\big(S_{\mathbb{F}}^{2p}([0,T],Q) \times (L_{\mathbb{F}}^{2p}([0,T],Q))^2\big)$.
\end{proposition}
\begin{proof}
Note that the functionals $\alpha_t^0,\,\beta_t^0:\, L^1(\mathcal{F}_t,Q)\times L^2(\mathcal{F}_t,Q)\rightarrow L^1(\mathcal{F}_t,Q)\times L^2(\mathcal{F}_t,Q)$ are linear, and
\begin{equation*}
\begin{split}
&|\alpha_t^0(q_t^1,q_t^2)|\leq C\Big(|q_t^1|+L_t E^Q[|q_t^1|]+L_t\big\{|q_t^2|+\big(E^Q[|q_t^2|^2]\big)^{\frac{1}{2}} \big\}\Big);\\
&|\beta_t^0(q_t^1,q_t^2)|\leq C\Big(|q_t^2|+K(t)\big\{ E^Q[|q_t^1|]+\big(E^Q[|q_t^2|^2]\big)^{\frac{1}{2}}\big\}\Big),\ (q_t^1,q_t^2)\in L^1(\mathcal{F}_t,Q)\times L^2(\mathcal{F}_t,Q),
\end{split}
\end{equation*}
where $K(t)=|X_t|+|U_t|$. Note that $K(\cdot),\, L,\, U\in S_{\mathbb{F}}^{\infty -}([0,T],Q)\big(\!:=\bigcap\limits_{p\geq2}S_{\mathbb{F}}^p([0,T],Q)\big)$. Given $p\geq2$, and arbitrary $(q^1,q^2)\in
 L_{\mathbb{F}}^p([0,T],Q)\times L_{\mathbb{F}}^{2p}([0,T],Q)$, we remark that, $t\in [0,T]$,
\begin{equation*}
\begin{split}
&E^Q\Big[\Big(\int_t^T\big|\alpha_s^0(q_s^1,q_s^2)\big|^2ds\Big)^{\frac{p}{2}}\Big]\leq C_p \Big(E^Q\Big[\Big(\int_t^T|q_s^1|^2ds\Big)^{\frac{p}{2}}\Big]+\Big(E^Q\Big[\Big(\int_t^T|q_s^2|^2ds\Big)^p\Big]\Big)^{\frac{1}{2}}\Big),\\
&E^Q\Big[\Big(\int_t^T\big|\beta_s^0(q_s^1,q_s^2)\big|^2ds\Big)^{p}\Big]\leq C_p \Big(\Big(E^Q\Big[\Big(\int_t^T|q_s^1|^2ds\Big)^{\frac{p}{2}}\Big]\Big)^2+E^Q\Big[\Big(\int_t^T|q_s^2|^2ds\Big)^p\Big]\Big).
\end{split}
\end{equation*}
Let us denote by $\big((\overline{p}^1,(\overline{q}^1,\overline{\overline{q}}^1)),\, (\overline{p}^2,(\overline{\overline{q}}^2,\overline{q}^2))\big)\in \big(S_{\mathbb{F}}^p([0,T],Q)\times (L_{\mathbb{F}}^p([0,T],Q))^2\big)\times\big( S_{\mathbb{F}}^{2p}([0,T],Q)\times (L_{\mathbb{F}}^{2p}([0,T],Q))^2 \big)$ the unique solution of the 2-dimensional BSDE
\begin{equation*}\left\{
\begin{split}
&d\overline{p}_t^1=-\alpha_t(q_t^1,q_t^2)dt+\overline{q}_t^1dB_t^1+\overline{\overline{q}}_t^1dY_t,\ t\in[0,T],\ \overline{p}_T^1=p_T^1,\\
&d\overline{p}_t^2=-\beta_t(q_t^1,q_t^2)dt+\overline{\overline{q}}_t^2dB_t^1+\overline{q}_t^2dY_t,\ t\in[0,T],\ \overline{p}_T^2=p_T^2,
\end{split}\right.
\end{equation*}
where $p_T^1$ and $p_T^2$ are defined by \eqref{SMP4}. The existence and the uniqueness are guaranteed by the fact that $p_T^1,\, p_T^2\in L^2(\mathcal{F}_T,Q)$ and
\begin{equation*}
\begin{split}
|\alpha_t(q_t^1,q_t^2)|&\leq|\alpha_t^0(q_t^1,q_t^2)|+C(1+L_t),\\
|\beta_t(q_t^1,q_t^2)|&\leq|\beta_t^0(q_t^1,q_t^2)|+C(|X_t|+|U_t|),\ t\in[0,T],
\end{split}
\end{equation*}
proving that the driving coefficients of the BSDE are $\mathbb{F}$-adapted and square integrable. This allows to define the mapping $\Phi(q^1,q^2):=\big((\overline{p}^1,(\overline{q}^1,\overline{\overline{q}}^1)),\, (\overline{p}^2,(\overline{\overline{q}}^2,\overline{q}^2))\big)$ as a function over $L_{\mathbb{F}}^{p}([0,T],Q)\times L_{\mathbb{F}}^{2p}([0,T],Q)$. Remark that, putting
$\Phi\big((p^1,(q^1,\check{q}^1)),\, (p^2,(\check{q}^2,q^2))\big):=\Phi(q^1,q^2)$, for $\big((p^1,(q^1,\check{q}^1)),\\ (p^2,(\check{q}^2,q^2))\big)\in \big(S_{\mathbb{F}}^{p}([0,T],Q)\times (L_{\mathbb{F}}^{p}([0,T],Q))^2\big) \times \big(S_{\mathbb{F}}^{2p}([0,T],Q)\times (L_{\mathbb{F}}^{2p}([0,T],Q))^2\big)$, this mapping can be interpreted as one on $\big(S_{\mathbb{F}}^{p}([0,T],Q)\times (L_{\mathbb{F}}^{p}([0,T],Q))^2\big) \times \big(S_{\mathbb{F}}^{2p}([0,T],Q)\times (L_{\mathbb{F}}^{2p}([0,T],Q))^2\big)$ to itself.
Given now any $(q^{i,1},q^{i,2})\in L_{\mathbb{F}}^{p}([0,T],Q)\times L_{\mathbb{F}}^{2p}([0,T],Q),\, i=1,2$, we denote $\Big(\big(\overline{p}^{i,1},(\overline{q}^{i,1},\overline{\overline{q}}^{i,1})\big),\\ \big(\overline{p}^{i,2},(\overline{\overline{q}}^{i,2},\overline{q}^{i,2})\big)\Big):=\Phi(q^{i,1},q^{i,2})$, and we put
$$ \widehat{p}^j:=\overline{p}^{1,j}-\overline{p}^{2,j},\ \widehat{q}^j:=\overline{q}^{1,j}-\overline{q}^{2,j},\ \widehat{\overline{q}}^j:=\overline{\overline{q}}^{1,j}-\overline{\overline{q}}^{2,j},\ j=1,2;\ \ \ \ \widetilde{q}^j:=q^{1,j}-q^{2,j},\ j=1,2. $$
Then, as $\alpha_s^0,\, \beta_s^0$ are linear, we get
\begin{equation*}\left\{
\begin{split}
&d\widehat{p}_t^1=-\alpha^0_t(\widetilde{q}_t^1,\widetilde{q}_t^2)dt+\widehat{q}_t^1dB_t^1+\widehat{\overline{q}}_t^1dY_t,\ t\in[0,T],\ \widehat{p}_T^1=0,\\
&d\widehat{p}_t^2=-\beta^0_t(\widetilde{q}_t^1,\widetilde{q}_t^2)dt+\widehat{\overline{q}}_t^2dB_t^1+\widehat{q}_t^2dY_t,\ t\in[0,T],\ \widehat{p}_T^2=0,
\end{split}\right.
\end{equation*}
and standard BSDE estimates yield
\begin{equation*}
\begin{split}
&E^Q\Big[\big|\widehat{p}_t^1\big|^p+\Big(\int_t^T\big(\big|\widehat{q}_s^1\big|^2+\big|\widehat{\overline{q}}_s^1\big|^2  \big)ds\Big)^{\frac{p}{2}}\Big]\leq C_p E^Q\Big[\Big(\int_t^T\big|\alpha_s^0(\widetilde{q}_s^1,\widetilde{q}_s^2)\big|ds\Big)^{p}\Big]\\
&\leq C_p(T-t)^{\frac{p}{2}}E^Q\Big[\Big(\int_t^T\big|\alpha_s^0(\widetilde{q}_s^1,\widetilde{q}_s^2)\big|^2ds\Big)^{\frac{p}{2}}\Big]\\
&\leq C_p(T-t)^{\frac{p}{2}}\Big(E^Q\Big[\Big(\int_t^T|\widetilde{q}_s^1|^2ds\Big)^{\frac{p}{2}}\Big]+\Big(E^Q\Big[\Big(\int_t^T|\widetilde{q}_s^2|^2ds\Big)^p\Big]\Big)^{\frac{1}{2}}\Big),\ t\in[0,T],\\
\end{split}
\end{equation*}
and
\begin{equation*}
\begin{split}
&E^Q\Big[\big|\widehat{p}_t^2\big|^{2p}+\Big(\int_t^T\big(\big|\widehat{q}_s^2\big|^2+\big|\widehat{\overline{q}}_s^2\big|^2\big) ds\Big)^p\Big]\leq C_p E^Q\Big[\Big(\int_t^T\big|\beta_s^0(\widetilde{q}_s^1,\widetilde{q}_s^2)\big|ds\Big)^{2p}\Big]\\
&\leq C_p(T-t)^{p}E^Q\Big[\Big(\int_t^T\big|\beta_s^0(\widetilde{q}_s^1,\widetilde{q}_s^2)\big|^2ds\Big)^{p}\Big]\\
&\leq C_p(T-t)^{p}\Big(\Big(E^Q\Big[\Big(\int_t^T|\widetilde{q}_s^1|^2ds\Big)^{\frac{p}{2}}\Big]\Big)^2+E^Q\Big[\Big(\int_t^T|\widetilde{q}_s^2|^2ds\Big)^p\Big]\Big),\ t\in[0,T].
\end{split}
\end{equation*}
This allows to conclude that the BSDE has a unique solution, first on a small interval $[T-\delta,T]$ ($\delta>0$ small enough) by using the contraction mapping arguments, and after by iteration, on $[T-2\delta,T-\delta]$,
$[T-3\delta,T-2\delta]$, ..., and thus on $[0,T]$.
\end{proof}

\subsection{Computations for the main results}
With the same notations as those introduced in the proof of Proposition \ref{EstofXL}, and so in particular $(X^{\varepsilon,\lambda},L^{\varepsilon,\lambda},U^{\varepsilon,\lambda}):=(1-\lambda)(X,L,U)+\lambda(X^{\varepsilon},L^{\varepsilon},U^{\varepsilon}),\ \lambda\in[0,1]$, we have that
\begin{equation*}
\begin{split}
&\Phi(X_T^\varepsilon,\mu_T^\varepsilon)-\Phi(X_T,\mu_T)=\int_0^1\partial_\lambda\big[\Phi(X_T^{\varepsilon,\lambda},\mu_T^{\varepsilon,\lambda})\big]d\lambda\\
&=\int_0^1\Big\{\Phi_x(X_T^{\varepsilon,\lambda},\mu_T^{\varepsilon,\lambda})(X_T^\varepsilon-X_T)+\widetilde{E}^Q\Big[\int_0^{\widetilde{U}_T^{\varepsilon,\lambda}}
\Phi_\mu(X_T^{\varepsilon,\lambda},\mu_T^{\varepsilon,\lambda};y)dy(\widetilde{L}_T^\varepsilon-\widetilde{L}_T)\Big]\\
&\qquad\qquad+\widetilde{E}^{Q}\big[\Phi_\mu(X_T^{\varepsilon,\lambda},\mu_T^{\varepsilon,\lambda};\widetilde{U}_T^{\varepsilon,\lambda})\widetilde{L}_T^{\varepsilon,\lambda}
(\widetilde{U}_T^\varepsilon-\widetilde{U}_T)\big]\Big\}d\lambda\\
&=I_\varepsilon^1+I_\varepsilon^2,
\end{split}
\end{equation*}
where
\begin{equation*}
\begin{split}
I_\varepsilon^1&=\Phi_x(T)(X_T^\varepsilon-X_T)+\widetilde{E}^Q\Big[\int_0^{\widetilde{U}_T}\Phi_\mu(T,y)dy(\widetilde{L}_T^\varepsilon-\widetilde{L}_T)\Big]+\widetilde{E}^{Q}
\big[\widetilde{\Phi}_\mu(T)\widetilde{L}_T(\widetilde{U}_T^\varepsilon-\widetilde{U}_T)\big],\\
I_\varepsilon^2&=\int_0^1\int_0^\lambda\partial_\rho\Big\{\Phi_x(X_T^{\varepsilon,\rho},\mu_T^{\varepsilon,\rho})(X_T^\varepsilon-X_T)+\widetilde{E}^Q\Big[\int_0^{\widetilde{U}_T^{\varepsilon,
\rho}}\Phi_\mu(X_T^{\varepsilon,\rho},\mu_T^{\varepsilon,\rho};y)dy(\widetilde{L}_T^\varepsilon-\widetilde{L}_T)\Big]\\
&\qquad\qquad+\widetilde{E}^{Q}\big[\Phi_\mu(X_T^{\varepsilon,\rho},\mu_T^{\varepsilon,\rho};\widetilde{U}_T^{\varepsilon,\rho})\widetilde{L}_T^{\varepsilon,\rho}(\widetilde{U}_T^\varepsilon
-\widetilde{U}_T)\big]\Big\}d\rho d\lambda.
\end{split}
\end{equation*}
From Proposition \ref{EstofYK}, we get
\begin{equation*}
\begin{split}
I_\varepsilon^1=&\Phi_x(T)(Y_T^{1,\varepsilon}+Y_T^{2,\varepsilon})+\widetilde{E}^Q\Big[\int_0^{\widetilde{U}_T}\Phi_\mu(T,y)dy(\widetilde{K}_T^{1,\varepsilon}+\widetilde{K}_T^{2,\varepsilon})\Big]\\
&+\widetilde{E}^{Q}\big[\widetilde{\Phi}_\mu(T)\widetilde{L}_T(\widetilde{V}_T^{1,\varepsilon}+\widetilde{V}_T^{2,\varepsilon})\big]+R_\varepsilon^1,
\end{split}
\end{equation*}
%where, for $\mu^v=\mu^{X^v,Y}$ we have $\overline{V}^{1,\varepsilon}+\overline{V}^{2,\varepsilon}$ instead of $V^{1,\varepsilon}+V^{2,\varepsilon}$, and, in both cases,
where $\big(E^Q[|R_\varepsilon^1|^p]\big)^{\frac{1}{p}}\leq C_p\varepsilon\rho_p(\varepsilon)$, $\varepsilon>0$, with $\rho_p(\varepsilon)\rightarrow0\ (\varepsilon\searrow0)$, for all $p\geq2$.
%Again with $(\overline{V}^{1,\varepsilon},\overline{V}^{2,\varepsilon})$ instead of $(V^{1,\varepsilon},V^{2,\varepsilon})$, if $\mu^v=\mu^{X^v,Y}$.
Considering again our arguments of the computations for the second order variational equations, we get from the Propositions \ref{EstofXL}, \ref{tech} and \ref{EstofYK} that
\begin{equation*}
\begin{split}
I_\varepsilon^2=&\frac{1}{2}\Phi_{xx}(T)(X_T^\varepsilon-X_T)^2+\widetilde{E}^{Q}\big[\widetilde{\Phi}_\mu(T)(\widetilde{U}_T^\varepsilon-\widetilde{U}_T)(\widetilde{L}_T^\varepsilon-\widetilde{L}_T)\big]
+\frac{1}{2}\widetilde{E}^{Q}\big[\widetilde{\Phi}_{z\mu}(T)\widetilde{L}_T(\widetilde{U}_T^\varepsilon-\widetilde{U}_T)^2\big]+R_\varepsilon^2\\
=&\frac{1}{2}\Phi_{xx}(T)(Y_T^{1,\varepsilon})^2+\widetilde{E}^{Q}\big[\widetilde{\Phi}_\mu(T)\widetilde{V}_T^{1,\varepsilon}\widetilde{K}_T^{1,\varepsilon}\big]+\frac{1}{2}\widetilde{E}^{Q}
\big[\widetilde{\Phi}_{z\mu}(T)\widetilde{L}_T(\widetilde{V}_T^{1,\varepsilon})^2\big]+R_\varepsilon^3,\\
\end{split}
\end{equation*}
again with $\big(E^Q[|R_\varepsilon^i|^p]\big)^{\frac{1}{p}}\leq C_p\varepsilon\rho_p(\varepsilon)$, $\varepsilon>0$, $i=2,3$, with $\rho_p(\varepsilon)\rightarrow0\ (\varepsilon\searrow0)$, for all $p\geq2$.
Consequently,
\begin{equation}\label{xingsqr}
\begin{split}
 \Phi(X_T^\varepsilon,\mu_T^\varepsilon)-\Phi(X_T,\mu_T)
=&\Phi_x(T)(Y_T^{1,\varepsilon}+Y_T^{2,\varepsilon})+\widetilde{E}^Q\Big[\int_0^{\widetilde{U}_T}\Phi_\mu(T,y)dy(\widetilde{K}_T^{1,\varepsilon}+\widetilde{K}_T^{2,\varepsilon})\Big]\\
&+\widetilde{E}^{Q}\big[\widetilde{\Phi}_\mu(T)\widetilde{L}_T(\widetilde{V}_T^{1,\varepsilon}+\widetilde{V}_T^{2,\varepsilon})\big]+\frac{1}{2}\Phi_{xx}(T)(Y_T^{1,\varepsilon})^2\\
&+\widetilde{E}^{Q}\big[\widetilde{\Phi}_\mu(T)\widetilde{V}_T^{1,\varepsilon}\widetilde{K}_T^{1,\varepsilon}\big]+\frac{1}{2}\widetilde{E}^{Q}\big[\widetilde{\Phi}_{z\mu}(T)
\widetilde{L}_T(\widetilde{V}_T^{1,\varepsilon})^2\big]+R_\varepsilon^4,
\end{split}
\end{equation}
where $\big(E^Q[|R_\varepsilon^4|^p]\big)^{\frac{1}{p}}\leq C_p\varepsilon\rho_p(\varepsilon)$, $\varepsilon>0$, with $\rho_p(\varepsilon)\rightarrow0\ (\varepsilon\searrow0)$, for all $p\geq2$. Now, observing that
\begin{equation}\label{2xing}
\begin{split}
&f(t,X_t^\varepsilon,\mu_t^\varepsilon,u_t^\varepsilon)-f(t,X_t,\mu_t,u_t)=f(t,X_t^\varepsilon,\mu_t^\varepsilon,u_t^\varepsilon)-f(t,X_t,\mu_t,u_t^\varepsilon)+\delta f(t)\mathbf{1}_{E_\varepsilon}(t)\\
=&\big\{f(t,X_t^\varepsilon,\mu_t^\varepsilon,v_t)-f(t,X_t,\mu_t,v_t)\big\}\mathbf{1}_{E_\varepsilon}(t)\\
&+\big\{f(t,X_t^\varepsilon,\mu_t^\varepsilon,u_t)-f(t,X_t,\mu_t,u_t)\big\}\mathbf{1}_{E_\varepsilon^c}(t)+\delta f(t)\mathbf{1}_{E_\varepsilon}(t),\ t\in[0,T].
\end{split}
\end{equation}
By using the same arguments as for \eqref{xingsqr}, but now for $f(t,\cdot,\cdot,v_t)$ and $f(t,\cdot,\cdot,u_t)$, respectively, instead of $\Phi$, we obtain for $\phi\in\{u,v\}$:
\begin{equation}\label{3xing}
\begin{split}
&f(t,X_t^\varepsilon,\mu_t^\varepsilon,\phi_t)-f(t,X_t,\mu_t,\phi_t)\\
=&\partial_x f(t,X_t,\mu_t,\phi_t)(Y_t^{1,\varepsilon}+Y_t^{2,\varepsilon})+\widetilde{E}^Q\Big[\int_0^{\widetilde{U}_t}\partial_\mu f(t,X_t,\mu_t,\phi_t;y)dy(\widetilde{K}_t^{1,\varepsilon}+\widetilde{K}_t^{2,\varepsilon})\Big]\\
&+\widetilde{E}^{Q}\big[\partial_\mu f(t,X_t,\mu_t,\phi_t;\widetilde{U}_t)\widetilde{L}_t(\widetilde{V}_t^{1,\varepsilon}+\widetilde{V}_t^{2,\varepsilon})\big]\\
&+\frac{1}{2}\partial_{xx}f(t,X_t,\mu_t,\phi_t)(Y_t^{1,\varepsilon})^2+\widetilde{E}^{Q}\big[\partial_\mu f(t,X_t,\mu_t,\phi_t;\widetilde{U}_t)\widetilde{V}_t^{1,\varepsilon}\widetilde{K}_t^{1,\varepsilon}\big]\\
&+\frac{1}{2}\widetilde{E}^{Q}\big[\partial_z(\partial_\mu f)(t,X_t,\mu_t,\phi_t;\widetilde{U}_t)\widetilde{L}_t(\widetilde{V}_t^{1,\varepsilon})^2\big]+R_\varepsilon(t,\phi_t),\ t\in[0,T],
\end{split}
\end{equation}
with $\big(E^Q[|R_\varepsilon(t,\phi_t)|^p]\big)^{\frac{1}{p}}\leq C_p\varepsilon\rho_p(\varepsilon)$, $\varepsilon>0$, for $\rho_p(\varepsilon)\rightarrow0\ (\varepsilon\searrow0)$, $p\geq2$, where $C_p$ does not depend on $t$ nor on  $\varepsilon>0$. Consequently, from \eqref{2xing} and \eqref{3xing} we conclude that, with the notations introduced in Subsection 4.2,
\begin{equation*}
\begin{split}
&f(t,X_t^\varepsilon,\mu_t^\varepsilon,u_t^\varepsilon)-f(t,X_t,\mu_t,u_t)=f_x(t)(Y_t^{1,\varepsilon}+Y_t^{2,\varepsilon})+\widetilde{E}^Q\Big[\int_0^{\widetilde{U}_t}f_\mu(t,y)dy(\widetilde{K}_t^{1,\varepsilon}+\widetilde{K}_t^{2,\varepsilon})\Big]\\
&+\widetilde{E}^{Q}\big[\widetilde{f}_\mu(t)\widetilde{L}_t(\widetilde{V}_t^{1,\varepsilon}+\widetilde{V}_t^{2,\varepsilon})\big]+\frac{1}{2}f_{xx}(t)(Y_t^{1,\varepsilon})^2
+\widetilde{E}^{Q}\big[\widetilde{f}_\mu(t)\widetilde{V}_t^{1,\varepsilon}\widetilde{K}_t^{1,\varepsilon}\big]\\
&+\frac{1}{2}\widetilde{E}^{Q}\big[\widetilde{f}_{z\mu}(t)\widetilde{L}_t(\widetilde{V}_t^{1,\varepsilon})^2\big]+ \delta f(t)\mathbf{1}_{E_\varepsilon}(t)+R_\varepsilon^5(t),\ t\in[0,T],
\end{split}
\end{equation*}
where
\begin{equation*}
\begin{split}
R_\varepsilon^5(t)=&\Big\{\delta f_x(t)(Y_t^{1,\varepsilon}+Y_t^{2,\varepsilon})+\widetilde{E}^Q\Big[\int_0^{\widetilde{U}_t}\delta f_\mu(t,y)dy(\widetilde{K}_t^{1,\varepsilon}+\widetilde{K}_t^{2,\varepsilon})\Big]\\
&+\widetilde{E}^{Q}\big[\delta \widetilde{f}_\mu(t)\widetilde{L}_t(\widetilde{V}_t^{1,\varepsilon}+\widetilde{V}_t^{2,\varepsilon})\big]+\frac{1}{2}\delta f_{xx}(t)(Y_t^{1,\varepsilon})^2+\widetilde{E}^{Q}\big[\delta\widetilde{f}_\mu(t)\widetilde{V}_t^{1,\varepsilon}\widetilde{K}_t^{1,\varepsilon}\big]\\
&+\frac{1}{2}\widetilde{E}^{Q}\big[\delta\widetilde{f}_{z\mu}(t)\widetilde{L}_t(\widetilde{V}_t^{1,\varepsilon})^2\big]\Big\}\mathbf{1}_{E_\varepsilon}(t)
+\big(R_\varepsilon(t,v_t)-R_\varepsilon(t,u_t)\big)\mathbf{1}_{E_\varepsilon}(t)+R_\varepsilon(t,u_t).
\end{split}
\end{equation*}
With by now standard arguments based on Corollary \ref{EstofUV}, Propositions \ref{EstofXL} and \ref{EstofYK} we get that $\displaystyle R_\varepsilon^5:=\int_0^T|R_\varepsilon^5(t)|dt\ \ \mbox{satisfies}\ \big(E^Q[|R_\varepsilon^5|^p]\big)^{\frac{1}{p}}\leq C_p\varepsilon\rho_p(\varepsilon), \ \varepsilon>0,$ with $\rho_p(\varepsilon)\rightarrow0\ (\varepsilon\searrow0)$,\ $p\geq2$.

Consequently, from the definition of the cost functional and the optimality of $u$, it holds that
\begin{equation}\label{MR1}
\begin{split}
0&\leq J(u^\varepsilon)-J(u)=E^Q\big[\Phi(X_T^\varepsilon,\mu_T^\varepsilon)-\Phi(X_T,\mu_T)\big]\!+\!E^Q\Big[\!\int_0^T\!\big(f(t,X_t^\varepsilon,\mu_t^\varepsilon,u_t^\varepsilon)-f(t,X_t,\mu_t,u_t)\big)dt\Big]\\
&=E^Q\Big[\Phi_x(T)(Y_T^{1,\varepsilon}+Y_T^{2,\varepsilon})\!+\!\widetilde{E}^Q\Big[\int_0^{\widetilde{U}_T}\Phi_\mu(T,y)dy(\widetilde{K}_T^{1,\varepsilon}+\widetilde{K}_T^{2,\varepsilon})\Big]
+\widetilde{E}^{Q}\big[\widetilde{\Phi}_\mu(T)\widetilde{L}_T(\widetilde{V}_T^{1,\varepsilon}+\widetilde{V}_T^{2,\varepsilon})\big]\Big]\\
&\ +E^Q\Big[\!\int_0^T\!\!\Big(f_x(t)(Y_t^{1,\varepsilon}+Y_t^{2,\varepsilon})\!+\!\widetilde{E}^Q\Big[\!\int_0^{\widetilde{U}_t}\!f_\mu(t,y)dy(\widetilde{K}_t^{1,\varepsilon}+\widetilde{K}_t^{2,\varepsilon})\Big]
\! +\!\widetilde{E}^{Q}\big[\widetilde{f}_\mu(t)\widetilde{L}_t(\widetilde{V}_t^{1,\varepsilon}+\widetilde{V}_t^{2,\varepsilon})\!\big]\!\Big)dt\!\Big]\\
&\ +E^Q\Big[\frac{1}{2}\Phi_{xx}(T)(Y_T^{1,\varepsilon})^2+\widetilde{E}^{Q}\big[\widetilde{\Phi}_\mu(T)\widetilde{V}_T^{1,\varepsilon}\widetilde{K}_T^{1,\varepsilon}\big]
+\frac{1}{2}\widetilde{E}^{Q}\big[\widetilde{\Phi}_{z\mu}(T)\widetilde{L}_T(\widetilde{V}_T^{1,\varepsilon})^2\big]\Big]\\
&\ +E^Q\Big[\int_0^T\Big(\frac{1}{2}f_{xx}(t)(Y_t^{1,\varepsilon})^2+\widetilde{E}^{Q}\big[\widetilde{f}_\mu(t)\widetilde{V}_t^{1,\varepsilon}\widetilde{K}_t^{1,\varepsilon}\big]
+\frac{1}{2}\widetilde{E}^{Q}\big[\widetilde{f}_{z\mu}(t)\widetilde{L}_t(\widetilde{V}_t^{1,\varepsilon})^2\big]\Big)dt\Big]\\
&\ +E^Q\Big[\int_0^T\delta f(t)\mathbf{1}_{E_\varepsilon}(t)dt\Big]+o(\varepsilon),\ \mbox{ as } \varepsilon\rightarrow0.
\end{split}
\end{equation}

\subsection{Peng's stochastic maximum principle}
Our objective is to rewrite the right-hand side of \eqref{MR1}. For simplicity we continue for our computations to use the setting of the proof of Proposition \ref{EstofXL}, i.e., we assume that $\sigma(t,x,\gamma,u)=\sigma(\gamma,u)$, $h(t,x,\gamma,u)=h(x,u)$, and we put:
$$\gamma_t^{1,\varepsilon}:= \widetilde{E}^{Q}\Big[\Big(\int_0^{\widetilde{U}_t}(\partial_\mu \sigma)(\mu_t,u_t;y)dy\Big)\widetilde{K}^{1,\varepsilon}_t \Big]+ \widetilde{E}^{Q}\Big[(\partial_\mu \sigma)(\mu_t,u_t;\widetilde{U}_t)\widetilde{L}_t \widetilde{V}_t^{1,\varepsilon}\Big].$$
From Proposition \ref{tech} we have
$$\Big|\widetilde{E}^{Q}\Big[\Big(\int_0^{\widetilde{U}_t}(\partial_\mu \sigma)(\mu_t,u_t;y)dy\Big)\widetilde{K}^{1,\varepsilon}_t \Big]\Big|\leq \rho_t (\varepsilon) \sqrt\varepsilon, $$
where $\rho_t (\varepsilon)\leq C,\ t\in [0,T],\ \varepsilon> 0,$ and $ \rho_t (\varepsilon)\rightarrow 0\ (\varepsilon\searrow 0)$. Moreover, we observe that
\begin{equation*}
\begin{split}
E^{Q}\Big[(\partial_\mu \sigma)(\mu_t,\widetilde{u}_t;U_t))L_t& V_t^{1,\varepsilon}\Big]=E^Q\Big[\Big(\frac{E^Q[L_t(\partial_\mu \sigma)(\mu_t,\widetilde{u}_t;U_t)\,|\,\mathcal{F}_t^Y]L_t}{E^Q[L_t\,|\,\mathcal{F}_t^Y]}\Big)Y_t^{1,\varepsilon}\Big]\\
&\qquad\quad +E^Q\Big[\Big(\frac{E^Q[L_t(\partial_\mu \sigma)(\mu_t,\widetilde{u}_t;U_t)\,|\,\mathcal{F}_t^Y]X_t}{E^Q[L_t\,|\,\mathcal{F}_t^Y]}\Big)K_t^{1,\varepsilon}\Big]\\
& -E^Q\Big[\Big(\frac{E^Q[L_t(\partial_\mu \sigma)(\mu_t,\widetilde{u}_t;U_t)\,|\,\mathcal{F}_t^Y]E^Q[L_t X_t|\mathcal{F}_t^Y]}{(E^Q[L_t\,|\,\mathcal{F}_t^Y])^2}\Big)K_t^{1,\varepsilon}\Big].
\end{split}
\end{equation*}
Thus, again, from Proposition \ref{tech} and standard estimates,
$$\Big| E^{Q}\big[(\partial_\mu \sigma)(\mu_t,\widetilde{u}_t;U_t))L_t V_t^{1,\varepsilon}\big]\Big|\leq \rho_t (\varepsilon) \sqrt\varepsilon,$$
where $\rho_t (\varepsilon)\leq C,\ t\in [0,T],\ \varepsilon> 0,$ and $ \rho_t (\varepsilon)\rightarrow 0\ (\varepsilon\searrow 0)$. Hence, also for $\gamma_t^{1,\varepsilon}: $
$$|\gamma_t^{1,\varepsilon}|\leq \rho_t (\varepsilon)\sqrt\varepsilon,\ \text{where}\ \ \rho_t (\varepsilon)\leq C,\ \rho_t (\varepsilon)\rightarrow 0\ (\varepsilon\searrow 0).$$
Recall, for the simplified setting in the proof of Proposition \ref{EstofXL} we have (the general case uses the same arguments):
$$\displaystyle dY_t^{1,\varepsilon}=\big(\gamma_t^{1,\varepsilon} + \delta\sigma(t)\mathbf{1}_{E_\varepsilon}(t)\big)dB_t^1,\ t\in [0,T],\ Y_0^{1,\varepsilon} = 0, \ \ $$
and now we introduce
$$\displaystyle dY_t^{1,\varepsilon,1}=  \delta\sigma(t)\mathbf{1}_{E_\varepsilon}(t)dB_t^1,\ t\in [0,T],\ \ Y_0^{1,\varepsilon,1} = 0.$$
Obviously, for all $p\geq 2$,
\begin{equation}\label{3starstar}
\begin{split}
&\Big(E^{Q}\big[\sup_{s\in[0,t]} |Y_s^{1,\varepsilon}- Y_s^{1,\varepsilon,1}|^p \,\big|\, \mathcal{F}_t^Y\big]\Big)^{\frac{1}{p}}= \Big(E^{Q}\big[\sup_{s\in[0,t]} |\int_0^s \gamma_r^{1,\varepsilon}
dB_r^1|^p \big| \mathcal{F}_t^Y\big]\Big)^{\frac{1}{p}}\\
 &\leq\ C_p \Big( E^Q[\int_0^t |\gamma_r^{1,\varepsilon}|^{p} dr \,\big|\, \mathcal{F}_t^Y]\Big)^{\frac{1}{p}}\leq \rho_p (\varepsilon) \sqrt\varepsilon,\ \ \varepsilon> 0,
\end{split}
\end{equation}
for $\rho_p(\varepsilon)\rightarrow0\ (\varepsilon \searrow 0)  $ and $\rho_p (\varepsilon)\leq C_p,\ \varepsilon > 0,\ p\geq 2.$

Given an arbitrary process $\theta \in S^{\infty -}_{\mathbb{F}}([0,T],Q)$, we denote by $\theta_{\cdot, \cdot} =(\theta_{t,s})_{0\leq s\leq t \leq T}$ the jointly measurable process such that,
\begin{equation*}
\begin{array}{lll}
&\qquad\qquad\qquad{\rm(i)}\ \theta_{t,s}\ \text{is}\ \mathcal{F}_s^{B^1} \vee \mathcal{F}_t^{Y}\text{-measurable},\ 0\leq s \leq t \leq T;\\
&\qquad\qquad\qquad{\rm(ii)}\ E^Q\Big[\int^T_0\int^t_0 |\theta_{t,s}|^2 dsdt\Big]< +\infty;\\
&\qquad\qquad\qquad{\rm(iii)}\ \theta_t = E^Q\big[\theta_t \,|\,\mathcal{F}_t^{Y}\big] + \int^t_0 \theta_{t,s} dB_s^1.\hspace{10cm}\\
\end{array}
\end{equation*}
Note that, for all $p\geq 2$, from Burkholder-Davis-Gundy inequality and Doob's martingale inequality,
\begin{equation*}
\begin{split}
&\Big(\!E^Q\Big[\Big(\int^t_0\big| \theta_{t,s}\big|^2 ds\Big)^{\frac{p}{2}} \,\big|\, \mathcal{F}_t^Y\Big]\!\Big)^{\frac{1}{p}}\leq C_p \Big(\!E^{Q}\big[\sup_{s\in[0,t]} \big|\int^s_0 \theta_{t,r}
dB_r^1\big|^p \,\big|\, \mathcal{F}_t^Y\big]\!\Big)^{\frac{1}{p}}\\
=&C_p \Big(\!E^{Q}\big[\sup_{s\in[0,t]}\big|\theta_s -E^Q\big[\theta_s \,\big|\,\mathcal{F}_s^{Y}\big]\big|^p \big|\, \mathcal{F}_t^Y\big] \!\Big)^{\frac{1}{p}}\leq C_p\Big(E^{Q}
\big[\sup_{s\in[0,t]}\big|\theta_s\big|^p  \,\big|\, \mathcal{F}_t^Y\big]\Big)^{\frac{1}{p}}.
\end{split}
\end{equation*}
Consequently, as
$$ E^Q\Big[ \theta_t Y_t^{1,\varepsilon,1}\,\Big|\,\mathcal{F}_t^{Y}\Big]=E^Q\Big[ \theta_t \,\Big|\,\mathcal{F}_t^{Y}\Big] E^Q\Big[Y_t^{1,\varepsilon,1}\,\Big|\,\mathcal{F}_t^{Y}\Big]
+E^Q\Big[\int^t_0 \theta_{t,s}\delta\sigma(s)\mathbf{1}_{E_\varepsilon}(s) ds\,\Big|\,\mathcal{F}_t^{Y} \Big],$$
and using that $ E^Q\big[Y_t^{1,\varepsilon,1}\big|\mathcal{F}_t^{Y}\big]=0$, we get\\
\begin{equation}\label{4starstar}
\Big|E^Q\big[ \theta_t Y_t^{1,\varepsilon,1}\,\Big|\,\mathcal{F}_t^{Y}\big]\Big| \leq C E^Q\Big[\Big(\int^t_0 \big|\theta_{t,s}\big|^2\mathbf{1}_{E_\varepsilon}(s) ds\Big)^{\frac{1}{2}}\,\Big|\,\mathcal{F}_t^{Y}
 \Big]\sqrt\varepsilon \leq \varphi_t (\varepsilon)\sqrt\varepsilon,
\end{equation}
for $\displaystyle\varphi_t (\varepsilon):=C E^Q\Big[\Big(\int^t_0 \big|\theta_{t,s}\big|^2\mathbf{1}_{E_\varepsilon}(s) ds\Big)^{\frac{1}{2}}\,\big|\,\mathcal{F}_t^{Y} \Big]\leq C\Big(E^{Q}\Big[\sup_{s\in[0,t]}|
\theta_s|^2 \,\Big|\, \mathcal{F}_t^{Y}\Big]\Big)^{\frac{1}{2}}.$\\
Then, obviously, for all $p\geq 2$,
$$\Big(E^{Q}\big[\sup_{t\in[0,T]}\big|\varphi_t (\varepsilon)\big|^p\big]\Big)^{\frac{1}{p}}\leq C_p\Big(E^{Q}\big[\sup_{t\in[0,T]}\big|\theta_t\big|^p \big]\Big)^{\frac{1}{p}}\leq C_p \in \mathbb{R}_{+} ,$$
and from the dominated convergence theorem, it follows that $\varphi_t (\varepsilon)\rightarrow 0\ (\varepsilon\searrow 0)$ in $L^p (Q),\ t\in [0,T].$\\
On the other hand, from \eqref{3starstar},
\begin{equation}\label{pound}
    \Big|E^Q\big[ \theta_t (Y_t^{1,\varepsilon}-Y_t^{1,\varepsilon,1})\,\big|\,\mathcal{F}_t^{Y}\big]\Big|\leq C\, \Big(E^Q\big[\,\big| \theta_t\big|^2 \,\big| \,\mathcal{F}_t^Y\big]\Big)^{\frac{1}{2}}\cdot
    \rho_2(\varepsilon)\sqrt{\varepsilon},
\end{equation}
where
\begin{equation*}
    \begin{split}
         &\Big(E^Q\big[\sup_{t\in[0,T]}\big(E^Q\big[\,\big| \theta_t\big|^2 \,\big| \,\mathcal{F}_t^Y\big]\big)^{\frac{p}{2}}\big]\Big)^{\frac{1}{p}}
       =\Big(E^Q\big[\sup_{t\in[0,T]}\big(E^Q\big[\,\big| \theta_t\big|^2 \,\big| \,\mathcal{F}_T^Y\big]\big)^{\frac{p}{2}}\big]\Big)^{\frac{1}{p}}\\
        \leq&\ \Big(E^Q\big[\sup_{t\in[0,T]} \big| \theta_t\big|^p\big]\Big)^{\frac{1}{p}} \leq C_p < +\infty,\ \ p\geq 2.
    \end{split}
\end{equation*}
Combining \eqref{4starstar} and \eqref{pound} we obtain that, for all $\theta\in S^{\infty -}_{\mathbb{F}}([0,T],Q)$,
\begin{equation*}
\begin{split}
E^Q\big[ \theta_t Y_t^{1,\varepsilon}\,\big|\,\mathcal{F}_t^{Y}\big]&=E^Q\big[ \theta_t (Y_t^{1,\varepsilon}-Y_t^{1,\varepsilon,1})\,\big|\,\mathcal{F}_t^{Y}\big]
+E^Q\big[ \theta_t Y_t^{1,\varepsilon,1}\,\big|\,\mathcal{F}_t^{Y}\big]
\end{split}
\end{equation*}
satisfies
\begin{equation}\label{5star}
\begin{split}
\big|E^Q\big[ \theta_t Y_t^{1,\varepsilon}\,\big|\,\mathcal{F}_t^{Y}\big]\big|\leq \varphi_t (\varepsilon)\sqrt\varepsilon,\ \ t\in [0,T],\ \ \varepsilon>0,
\end{split}
\end{equation}
where $\displaystyle\Big(E^{Q}\big[\sup_{t\in[0,T]}\big|\varphi_t (\varepsilon)\big|^p\big]\Big)^{\frac{1}{p}}\leq C_p<+\infty,\ p\geq 2,$ \text{and} $\varphi_t (\varepsilon)\rightarrow 0\ (\varepsilon\searrow 0)\ $\text{in probability Q}, \text{for all} $t\in [0,T].$\ Recall
\begin{equation*}
\left\{
\begin{split}
&dK_t^{1,\varepsilon}=\big\{h(X_t,u_t)K_t^{1,\varepsilon}+L_t h_x(X_t,u_t)Y_t^{1,\varepsilon}+L_t \delta h(t)\mathbf{1}_{E_\varepsilon}(t)\big\}dY_t,\ t\in [0,T],\\
&K_0^{1,\varepsilon}=0.
\end{split}\right.
\end{equation*}
Let $K_t^{1,\varepsilon, 1}\in S^{\infty-}_{\mathbb{F}}([0,T],Q)$ be defined by the following SDE:
\begin{equation*}
\left\{
\begin{split}
&dK_t^{1,\varepsilon,1}=\big\{h(X_t,u_t)K_t^{1,\varepsilon,1}+L_t h_x(X_t,u_t)Y_t^{1,\varepsilon,1}+L_t \delta h(t)\mathbf{1}_{E_\varepsilon}(t)\big\}dY_t,\ t\in [0,T],\\
&K_0^{1,\varepsilon,1}=0.
\end{split}\right.
\end{equation*}
From \eqref{3starstar} we have
\begin{equation}\label{6star}
\begin{split}
\Big(E^Q[\sup_{t\in[0,T]} \big|K_t^{1,\varepsilon}- K_t^{1,\varepsilon,1}\big|^p \,\big|\, \mathcal{F}_t^Y\big]\Big)^{\frac{1}{p}}\leq \rho_p (\varepsilon) \sqrt\varepsilon,\ \ \varepsilon> 0.
\end{split}
\end{equation}
Let us estimate $\displaystyle E^Q\big[ \theta_t K_t^{1,\varepsilon}\,\big|\,\mathcal{F}_t^{Y}\big]$.
For $\theta \in S^{\infty -}_{\mathbb{F}}([0,T],Q)$,
\begin{equation*}
\begin{split}
&d\big(L^{-1}_t K_t^{1,\varepsilon,1}\big)=L_t^{-1}\big(h(X_t,u_t) K_t^{1,\varepsilon,1} + L_t h_x (X_t, u_t) Y_t^{1,\varepsilon,1} +L_t \delta h(t) \mathbf{1}_{E_\varepsilon}(t)\big)dY_t\\
&\qquad\qquad\quad\quad\quad - L_t^{-1}h(X_t,u_t) K_t^{1,\varepsilon,1} dY_t +L_t^{-1}\big(h(X_t,u_t)\big)^2 K_t^{1,\varepsilon,1} dt\\
&\qquad\qquad\quad\quad\quad - L_t^{-1}h(X_t,u_t)\big(h(X_t,u_t) K_t^{1,\varepsilon,1}+L_t h_x(X_t,u_t)  Y_t^{1,\varepsilon,1}+L_t \delta h(t) \mathbf{1}_{E_\varepsilon}(t)\big)dt\\
=&\big(h_x(X_t,u_t)Y_t^{1,\varepsilon,1}+\delta h(t)\mathbf{1}_{E_\varepsilon}(t)\big)dY_t-\big((h\cdot h_x)(X_t,u_t)Y_t^{1,\varepsilon,1}+h(X_t,u_t)\delta h(t)\mathbf{1}_{E_\varepsilon}(t)\big)dt,
\end{split}
\end{equation*}
which leads to that
\begin{equation*}
\begin{split}
&E^Q\big[ \theta_t K_t^{1,\varepsilon,1}\big|\mathcal{F}_t^{Y}\big]=E^Q\big[(\theta_t L_t)(K_t^{1,\varepsilon,1} L_t^{-1})\,\big|\,\mathcal{F}_t^{Y}\big]\\
=&E^Q\big[(\theta_t L_t)\big(\int^t_0 h_x(X_s,u_s)Y_s^{1,\varepsilon,1} dY_s\big)\,\big|\,\mathcal{F}_t^{Y}\big]+E^Q\big[(\theta_t L_t)\big(\int^t_0 \delta h(s)\mathbf{1}_{E_\varepsilon}(s)dY_s\big)\,\big|\,\mathcal{F}_t^{Y}\big]\\
&-\!E^Q\big[(\theta_t L_t)\big(\!\int^t_0\!(h\cdot h_x)(X_s,u_s)Y_s^{1,\varepsilon,1} ds\big)\,\big|\,\mathcal{F}_t^{Y}\big]\!-\!E^Q\big[(\theta_t L_t)\big(\!\int^t_0\! h(X_s,u_s)\delta h(s)\mathbf{1}_{E_\varepsilon}(s)ds\big)\,\big|\,\mathcal{F}_t^{Y}\big]\\
=&:I_t^{1,\varepsilon}+I_t^{2,\varepsilon}+I_t^{3,\varepsilon}+I_t^{4,\varepsilon},
\end{split}
\end{equation*}
and we denote these terms by $I_t^{i,\varepsilon},\ i=1,2,3,4,\ $ respectively. Let us estimate these terms.\\
\noindent$\bullet$ For $I_t^{4,\varepsilon}$ we easily see
 $$ \big|I_t^{4,\varepsilon}\big|=\Big|E^Q\big[(\theta_t L_t)\big(\int^t_0 h(X_s,u_s)\delta h(s)\mathbf{1}_{E_\varepsilon}(s)ds\big)\big|\mathcal{F}_t^{Y}\big]\Big|\leq C E^Q\big[|\theta_t L_t|\big|\mathcal{F}_t^{Y}\big]\cdot\varepsilon,\ \varepsilon >0.$$
\noindent$\bullet$ For $I_t^{3,\varepsilon}$ we notice that
\begin{equation*}
\begin{split}
I_t^{3,\varepsilon}&=-E^Q\big[(\theta_t L_t)\big(\int^t_0(h\cdot h_x)(X_s,u_s)Y_s^{1,\varepsilon,1} ds\big)\big|\mathcal{F}_t^{Y}\big]\\
&=-\int^t_0 E^Q\big[(\theta_t L_t)(h\cdot h_x)(X_s,u_s)\int_0^s\delta \sigma(r)\mathbf{1}_{E_\varepsilon}(r) dB_r^1\,\big|\,\mathcal{F}_t^{Y}\big]ds.
\end{split}
\end{equation*}
Let $\theta^s_{t,\cdot} \in L^2_{\mathbb{F}^{B^1} \vee \mathcal{F}_t^Y} ([0,t])$ be jointly measurable w.r.t. $(s,t,r,\omega)$ and such that:
$$\theta_t L_t (h\cdot h_x)(X_s,u_s) = E^Q \big[\theta_t L_t (h\cdot h_x)(X_s,u_s)\,\big|\,\mathcal{F}_t^{Y}\big]+\int_0^t \theta^s_{t,r} dB_r^1, \ \ 0\leq s\leq t,$$
which yields that
$$I_t^{3,\varepsilon} = -\int_0^t E^Q\big[\int^s_0\theta^s_{t,r} \delta \sigma(r)\mathbf{1}_{E_\varepsilon}(r)dr\,\big|\,\mathcal{F}_t^{Y}\big]ds$$
and
$$\big|I_t^{3,\varepsilon}\big|^p\leq C_p\int_0^t E^Q\big[(\int^s_0|\theta^s_{t,r}|^2 \mathbf{1}_{E_\varepsilon}(r)dr)^{\frac{p}{2}}\,\big|\,\mathcal{F}_t^{Y}\big]ds\cdot \varepsilon^{\frac{p}{2}},\ \ p>1.$$
We observe that, for $p>1$, thanks to the Burkholder-Davis-Gundy inequality and Doob's martingale inequality and the fact that $h\cdot h_x$ is bounded, we have
\begin{equation*}
\begin{split}
&\int_0^t E^Q\big[(\int^s_0|\theta^s_{t,r}|^2 dr)^{\frac{p}{2}}\,\big|\,\mathcal{F}_t^{Y}\big]ds
\leq C_p \int_0^t E^Q\big[\big|\theta_t L_t (h\cdot h_x)(X_s,u_s)\big|^p\,\big|\,\mathcal{F}_t^{Y}\big]ds\\
\leq&\ C_p E^Q\big[\;\big|\theta_t L_t\big|^p\,\big|\,\mathcal{F}_t^{Y}\big] = C_p E^Q\big[\big|\theta_t L_t\big|^p\,\big|\,\mathcal{F}_T^{Y}\big]
\leq C_p E^Q[\;\sup_{0\leq t\leq T} |\theta_t L_t|^p\,\big|\,\mathcal{F}_T^{Y}\big]=:\eta^{(p)},\ \
\end{split}
\end{equation*}
where $\eta^{(p)}\in L^{\infty -} (\Omega, \mathcal{F}, Q)$ \big(Recall $ L^{\infty -} (\Omega, \mathcal{F}, Q) := \bigcap_{p\geq 1}L^{p} (\Omega, \mathcal{F}, Q)$\big). Thus, we have
 $$\displaystyle\big|I_t^{3,\varepsilon}\big| \leq (\eta_t(\varepsilon))^{\frac{1}{p}}\cdot\sqrt\varepsilon,\ \ \varepsilon>0,\
\mbox{where}\ \ \displaystyle\eta_t(\varepsilon):=C_p \int_0^t E^Q\big[(\int^s_0|\theta^s_{t,r}|^2 \mathbf{1}_{E_\varepsilon}(r) dr)^{\frac{p}{2}}\big|\mathcal{F}_t^{Y}\big]ds,$$
and from the dominated convergence theorem, $\displaystyle\eta_t(\varepsilon)\rightarrow 0\,(\varepsilon\searrow 0)$, and $\displaystyle\big|\eta_t(\varepsilon)\big|\leq \eta^{(p)},
\ t\in[0,T].$\\
In the following we use the notation $\eta_t(\varepsilon)$ for random variables depending on $(t,\varepsilon)$, such that $\eta_t(\varepsilon)\rightarrow 0 \,
(\varepsilon\searrow 0),\ {Q}$-a.s., and $|\eta_t(\varepsilon)|\leq \eta,\ t\in[0,T],\ \varepsilon >0,$ for some $\eta \in  L^{\infty -} (\Omega,\mathcal{F}, Q).$\\
\noindent$\bullet$ For $I_t^{1,\varepsilon}$ we observe that
\begin{equation*}
\begin{split}
&\int^t_0 h_x(X_s,u_s)Y_s^{1,\varepsilon,1} dY_s=\int_0^t h_x(X_s,u_s)\big(\int_0^s \delta \sigma(r)\mathbf{1}_{E_\varepsilon}(r) dB_r^1\big)dY_s\\
=&\int_0^t h_x(X_s,u_s)dY_s\cdot Y_t^{1,\varepsilon,1} - \int^t_0\big(\int^s_0h_x(X_r,u_r)dY_r\big)\delta \sigma(s)\mathbf{1}_{E_\varepsilon}(s) dB_s^1.
\end{split}
\end{equation*}
Then,
\begin{equation*}
\begin{split}
I_t^{1,\varepsilon}=&E^Q\big[\big(\theta_t L_t \int^t_0h_x(X_s,u_s)dY_s\big)\cdot Y_t^{1,\varepsilon,1} \,\big|\,\mathcal{F}_t^{Y}\big]\\
&-E^Q\big[(\theta_t L_t)\int^t_0\big(\int^s_0h_x(X_r,u_r)dY_r\big)\delta \sigma(s)\mathbf{1}_{E_\varepsilon}(s) dB_s^1\,\big|\,\mathcal{F}_t^{Y}\big]=:I_t^{1,1,\varepsilon} - I_t^{1,2,\varepsilon}.\\
\end{split}
\end{equation*}
For $I_t^{1,1,\varepsilon}$, let $\displaystyle\zeta_t := \theta_t L_t \int^t_0h_x(X_s,u_s)dY_s,\ t\in [0,T]$, and denote by $\zeta_{t,\cdot} \in L^2_{\mathbb{F}^{B^1} \vee \mathcal{F}_t^Y} ([0,t])$ the jointly measurable two-parameter process, such that
$$\zeta_t = E^Q \big[\zeta_t\,\big|\,\mathcal{F}_t^{Y}\big]+\int_0^t \zeta_{t,s} dB_s^1,\ Q\text{-a.s.},\ t\in [0,T].$$
Then, obviously, $\displaystyle I_t^{1,1,\varepsilon} = E^Q\Big[\int_0^t \zeta_{t,s} \delta \sigma(s) \mathbf{1}_{E_\varepsilon}(s)ds\big|\mathcal{F}_t^{Y}\Big].$\ Thus, we have
$$\big|I_t^{1,1,\varepsilon} \big|\leq C E^Q\big[\big(\int_0^T |\zeta_{t,s}|^2\mathbf{1}_{E_\varepsilon}(s)ds\big)^{\frac{1}{2}}\,\big|\,\mathcal{F}_t^{Y}\big]\cdot \sqrt \varepsilon=: \eta_t(\varepsilon)\sqrt \varepsilon,$$
where
\begin{equation*}
\begin{split}
&\big|\eta_t(\varepsilon)\big|^p \leq C_p E^Q\Big[\Big(\int_0^T |\zeta_{t,s}|^2\mathbf{1}_{E_\varepsilon}(s)ds\Big)^{\frac{p}{2}}\,\big|\,\mathcal{F}_t^{Y}\Big]\leq C_p
E^Q\Big[\Big(\int_0^T |\zeta_{t,s}|^2ds\Big)^{\frac{p}{2}}\,\Big|\,\mathcal{F}_t^{Y}\Big]\\
\leq&\ C_p E^Q\Big[\,\big |\zeta_t\big|^p\big|\mathcal{F}_t^{Y}\Big]=C_p E^Q\Big[\,\big|\zeta_t\big|^p \,\big|\,\mathcal{F}_T^{Y}\Big]\leq C_p E^Q\Big[\;\sup_{t\leq T}
\big|\zeta_t\big|^p \,\big|\,\mathcal{F}_T^{Y}\Big]=:\eta ^p, \ \eta\!\in\! L^{\infty -} (\Omega,\mathcal{F}_T, Q).
\end{split}
\end{equation*}
It follows from the dominated convergence theorem that $\eta_t(\varepsilon)\rightarrow 0\ (\varepsilon\searrow0)$ in $L^p (Q)\ ( p\geq 2),$ and $\big|\eta_t(\varepsilon)\big| \leq \eta,\ t\in [0,T],\ \varepsilon > 0.$\\
To estimated $I_t^{1,2,\varepsilon}$, we introduce the jointly measurable process $\gamma_{t,\cdot} \in L^2_{\mathbb{F}^{B^1} \vee \mathcal{F}_t^Y} ([0,t])$ such that
$$\theta_t L_t= E^Q\big[\theta_t L_t\,\big|\,\mathcal{F}_t^{Y}\big]+\int_0^t \gamma_{t,s} dB_s^1,\ \text{Q-a.s.}, \ t\in [0,T].$$
Then, $\displaystyle I_t^{1,2,\varepsilon} = E^Q\Big[\int^t_0 \gamma_{t,s}\big(\int^s_0 h_x (X_r,u_r)dY_r\big)\delta \sigma(s)
\mathbf{1}_{E_\varepsilon}(s)ds\,\big|\,\mathcal{F}_t^{Y}\Big]$, and
\begin{equation*}
\begin{split}
\big| I_t^{1,2,\varepsilon}\big| &\leq C E^Q\Big[\big(\int^t_0 \Big|\gamma_{t,s}\cdot\int^s_0 h_x (X_r,u_r)dY_r\Big|^2 \mathbf{1}_{E_\varepsilon}(s)ds\big)^{\frac{1}{2}}\,\Big|\,\mathcal{F}_t^{Y}\Big]\cdot\sqrt\varepsilon\\
&\leq C E^Q\Big[\Big(\sup_{s\in[0,t]}\Big| \int_0^s h_x (X_r,u_r)dY_r\Big|\Big)\Big(\int_0^t \big|\gamma_{t,s}\big|^2\mathbf{1}_{E_\varepsilon}(s)ds\Big)^{\frac{1}{2}}\,\Big|\,\mathcal{F}_t^{Y}\Big]\cdot\sqrt\varepsilon\\
&=:\eta_t(\varepsilon)\sqrt\varepsilon,
\end{split}
\end{equation*}
where we have
\begin{align*}
\big|\eta_t(\varepsilon)\big|&\leq \Big(E^Q[\sup_{s\in[0,T]} \big|\int^s_0h_x (X_r,u_r)dY_r|^2\,\big|\,\mathcal{F}_T^{Y}\big]\Big)^{\frac{1}{2}} \Big(E^Q\big[\big(\int_0^t |\gamma_{t,s}|^2
\mathbf{1}_{E_\varepsilon}(s)ds\big)\,\big|\,\mathcal{F}_t^{Y}\big]\Big)^{\frac{1}{2}}\\
&\leq  \Big(E^Q[\sup_{s\in[0,T]} \big|\int^s_0h_x (X_r,u_r)dY_r|^2\,\big|\,\mathcal{F}_T^{Y}\big]\Big)^{\frac{1}{2}} \Big(E^Q\big[\big(\int_0^t |\gamma_{t,s}|^2ds\big)\,\big|\,
\mathcal{F}_t^{Y}\big]\Big)^{\frac{1}{2}}\\
&\leq C  \Big(E^Q[\sup_{s\in[0,T]} \big|\int^s_0h_x (X_r,u_r)dY_r|^2\,\big|\,\mathcal{F}_T^{Y}\big]\Big)^{\frac{1}{2}} \Big(E^Q\big[\big|\theta_tL_t\big|^2\,\big|\,\mathcal{F}_t^{Y}\big]
\Big)^{\frac{1}{2}}\\
&=C  \Big(E^Q[\sup_{s\in[0,T]} \big|\int^s_0h_x (X_r,u_r)dY_r|^2\,\big|\,\mathcal{F}_T^{Y}\big]\Big)^{\frac{1}{2}} \Big(E^Q\big[\big|\theta_tL_t\big|^2\,\big|\,\mathcal{F}_T^{Y}\big]
\Big)^{\frac{1}{2}}\\
&\leq C  \Big(E^Q\big[\sup_{s\in[0,T]} \big|\int^s_0h_x (X_r,u_r)dY_r|^2\,\big|\,\mathcal{F}_T^{Y}\big]\Big)^{\frac{1}{2}} \Big(E^Q\big[\sup_{t\in[0,T]}|\theta_tL_t|^2\,\big|\,
\mathcal{F}_T^{Y}\big]\Big)^{\frac{1}{2}}\\
&=:\eta,\ \ \ \eta \in  L^{\infty -} (\Omega,\mathcal{F}_T, Q).
\end{align*}
Hence, $|\eta_t(\varepsilon)|\leq \eta,\ t\in[0,T],\ \varepsilon > 0,$ and from the dominated convergence theorem, we deduce
 $$\displaystyle \eta_t(\varepsilon)\rightarrow 0,\ t\in [0,T]\ (\varepsilon\searrow 0),\ \ \mbox{in}\  L^p(Q),\  p\geq 2.$$
Consequently, from the above estimates of $I^{i,\varepsilon}_t,\ i=1,2,3,4,$ it follows that, for all $\theta \in S^{\infty -}_{\mathbb{F}} ([0,T]; Q)$,
\begin{equation}\label{theta K 11}
\begin{split}
E^Q\big[\theta_t K_t^{1,\varepsilon,1}\,\big|\,\mathcal{F}_t^{Y}\big] = E^Q\Big[(\theta_t L_t)\cdot\int^t_0 \delta h(s)\mathbf{1}_{E_\varepsilon}(s)dY_s\,\big|\,\mathcal{F}_t^{Y}\Big]
+\eta_t(\varepsilon)\sqrt\varepsilon,
\end{split}
\end{equation}
where $\displaystyle \big|\eta_t(\varepsilon)\big|\leq \eta \in L^{\infty -} (\Omega,\mathcal{F}_T, Q),\ \eta_t(\varepsilon)\rightarrow0\ (\varepsilon\searrow 0),\ t\in [0,T],
\ \mbox{in}\  L^p(Q),\  p\geq 2.$\ Since
\begin{equation*}
\begin{split}
&\big|E^Q\big[ \theta_t (K_t^{1,\varepsilon}-K_t^{1,\varepsilon,1})\,\big|\,\mathcal{F}_t^{Y}\big]\big|=\big|E^Q\big[ \theta_t (K_t^{1,\varepsilon}-K_t^{1,\varepsilon,1})\,
\big|\,\mathcal{F}_T^{Y}\big]\big|\leq E^Q \big[\sup_{t\in[0,T]} |
\theta_t| |K_t^{1,\varepsilon}-K_t^{1,\varepsilon,1}|\,\big|\,\mathcal{F}_T^{Y}\big]\\
\leq& \Big(E^Q\big[\sup_{t\in[0,T]} |\theta_t|^2\,
\big|\,\mathcal{F}_T^{Y}\big]\Big)^{\frac{1}{2}}\Big(E^Q [\sup_{t\in[0,T]}|K_t^{1,\varepsilon}-K_t^{1,\varepsilon,1}|^2\,\big|\,
\mathcal{F}_T^{Y}\big]\Big)^{\frac{1}{2}}=:\sqrt\varepsilon\eta(\varepsilon),
\end{split}
\end{equation*}
and from \eqref{6star}, we obtain $\displaystyle \Big(E^Q\big[\,|\eta(\varepsilon)|^p\big]\Big)^{\frac{1}{p}}\rightarrow 0\ (\varepsilon\searrow0),\ \ p\geq 2$.\ Combined with (\ref{theta K 11}), this yields
\begin{equation}\label{7star}
\begin{split}
E^Q\big[\theta_t K_t^{1,\varepsilon}\big|\mathcal{F}_t^{Y}\big] = E^Q\big[(\theta_t L_t) \int_0^t\delta h(s) \mathbf{1}_{E_\varepsilon}(s)dY_s\,\big|\,\mathcal{F}_t^{Y}\big] + \sqrt \varepsilon \eta_t (\varepsilon),
\end{split}
\end{equation}
where $\displaystyle \big|\eta_t(\varepsilon)\big|\leq \eta \in L^{\infty-} (\Omega,\mathcal{F}_T, Q),\ \eta_t(\varepsilon)\rightarrow0\ (\varepsilon\searrow 0),\ t\in[0,T],
\ \mbox{in}\  L^p(Q),\  p\geq 2.$\\

Let us now estimate $\displaystyle E^Q\big[K_t^{1,\varepsilon}Y_t^{1,\varepsilon}\,\big|\,\mathcal{F}_t^{Y}\big]$. For this we recall that
$$\Big(E^Q \big[\sup_{t\in[0,T]}(\big|K_t^{1,\varepsilon}\big|^p + \big|K_t^{1,\varepsilon,1}\big|^p +\big|Y_t^{1,\varepsilon}\big|^p +\big|Y_t^{1,\varepsilon,1}\big|^p)\big]\Big)^{\frac{1}{p}}\leq C_p \sqrt \varepsilon,\ \varepsilon>0,\ p\geq2,$$
and
$$\Big(E^Q \big[\sup_{t\in[0,T]}(\big|K_t^{1,\varepsilon}-K_t^{1,\varepsilon,1}\big|^p +\big|Y_t^{1,\varepsilon}-Y_t^{1,\varepsilon,1}\big|^p)\big]\Big)^{\frac{1}{p}}\leq \rho_p(\varepsilon) \sqrt \varepsilon,\ \varepsilon>0,\ p\geq2,\ \rho_p(0+)=0.$$
Then,
\begin{equation*}
\begin{split}
E^Q\big[K_t^{1,\varepsilon} Y_t^{1,\varepsilon}\,\big|\,\mathcal{F}_t^{Y}\big]=&E^Q\big[K_t^{1,\varepsilon} Y_t^{1,\varepsilon}\,\big|\,\mathcal{F}_T^{Y}\big]
=E^Q\big[K_t^{1,\varepsilon,1} Y_t^{1,\varepsilon,1}\,\big|\,\mathcal{F}_t^{Y}\big]+\varepsilon \eta_t(\varepsilon),
\end{split}
\end{equation*}
where $\displaystyle\Big(E^Q\big[\sup_{t\in[0,T]} \big|\eta_t(\varepsilon)\big|^p\big]\Big)^{\frac{1}{p}}\rightarrow 0,\  \varepsilon\searrow0.$\ Recall that $\displaystyle L_t = 1\!+\! \int_0^t\! h(X_s,u_s)L_sdY_s,\ t\in[0,T]$,\ and
\begin{equation*}
\begin{split}
L^{-1}_t K_t^{1,\varepsilon,1}=&\int_0^t \big(h_x(X_s,u_s)Y_s^{1,\varepsilon,1}+\delta h(s)\mathbf{1}_{E_\varepsilon}(s)\big)dY_s\\
&-\int^t_0\big((h\cdot h_x)(X_s,u_s)Y_s^{1,\varepsilon,1}+h(X_s,u_s)\delta h(s)\mathbf{1}_{E_\varepsilon}(s)\big)ds,\ t\in[0,T].
\end{split}
\end{equation*}

Now let us introduce the following assumption for $\sigma$.\\
\noindent\textbf{Assumption (H3)}:\\
(i)\ ${\sigma(t,x,\gamma,u)=\sigma(t,\gamma,u),\ (t,x,\gamma,u)\in [0,T]\times\mathbb{R}\times\mathcal{P}_{2}(\mathbb{R})\times U}.$\\

\noindent Let $D^1 = (D^1_s)_{s\in[0,T]}$ denote the Malliavin derivative w.r.t. the Brownian motion $B^1 = (B_s^1)_{s\in[0,T]}.$
We use the well-known basic results and commonly used notions, for which the reader is, for example, referred to Nualart \cite{N06}.
Noting that $\displaystyle Y_r^{1,\varepsilon,1}=\int_0^r \delta \sigma(s) \mathbf{1} _{E_\varepsilon}(s)dB_s^1, $ and $\delta \sigma (s)$ is
independent of $B^1$, we see that, for $s\leq r,\ dsdQ\text{-a.e.}$, $D_s^1[Y_r^{1,\varepsilon,1}] = \delta \sigma (s)\mathbf{1}_{E_\varepsilon}(s).$ Thus,
\begin{equation*}
\begin{split}
D_s^1[L_t^{-1} K_t^{1,\varepsilon,1}]&=\int_s^t\big\{D_s^1 [h_x(X_r,u_r)]Y_r^{1,\varepsilon,1} + D_s^1[\delta h(r)]\mathbf{1}_{E_\varepsilon}(r)\big\}dY_r\\
&\quad-\int_s^t\big\{D_s^1\big[(h\cdot h_x)(X_r,u_r)\big]Y_r^{1,\varepsilon,1} +D_s^1\big[h(X_r,u_r)\delta h(r)\big]\mathbf{1}_{E_\varepsilon}(r)\big\}dr\\
&\quad+\Big(\int_s^t h_x(X_r,u_r)dY_r\Big)\delta \sigma(s)\mathbf{1}_{E_\varepsilon}(s)-\int_s^t(h\cdot h_x)(X_r,u_r)dr\cdot\delta\sigma(s)\mathbf{1}_{E_\varepsilon}(s).
\end{split}
\end{equation*}
We remark that all the coefficients
\begin{equation*}
\begin{split}
&D_s^1\big[h_x(X_r,u_r)\big] = h_{xx} (X_r,u_r)\sigma(\mu_s,u_s),\\
&D_s^1\big[\delta h(r)\big] = D_s\big[h(X_r,v_r)-h(X_r,u_r)\big] = \big(h_x(X_r,v_r)-h_x(X_r,u_r)\big)\sigma(\mu_s,u_s)=\delta h_x(r)\sigma(\mu_s,u_s),\\
&{D_s^1\big[h(X_r,u_r)\delta h(r)\big]=h_x(X_r,u_r)\delta h(r)\sigma(\mu_s,u_s)+h(X_r,u_r)\delta h_x(r)\sigma(\mu_s,u_s)},\\
&{D_s^1\big[(h\cdot h_x)(X_r,u_r)\big]=|h_x(X_r,u_r)|^2\sigma(\mu_s,u_s)+(h\cdot h_{xx})(X_r,u_r)\sigma(\mu_s,u_s)},
\end{split}
\end{equation*}
are bounded. Then obviously,
$$D_s^1\big[L_t^{-1} K_t^{1,\varepsilon,1}\big]=\Big(\int_s^t h_x(X_r,u_r)dY_r-\int_s^t(h\cdot h_x)(X_r,u_r)dr\Big)\delta \sigma(s)\mathbf{1}_{E_\varepsilon}(s)+{\varphi_{s,t}(\varepsilon)}
\sqrt\varepsilon,$$
{with $\displaystyle \sup_{s\in[0,T]}\big(E^Q[\sup_{s\leq t\leq T} |\varphi_{s,t}(\varepsilon)|^p]\big)^{\frac{1}{p}}\leq C_p.$}  On the other hand, putting
\begin{equation*}
\begin{split}
&h(s,t):=\int_s^t h_x(X_r,u_r)dY_r-\int_s^t(h\cdot h_x)(X_r, u_r)dr,\
H_{\varepsilon}(t):=\int_0^t h(s,t)|\delta\sigma(s)|^2 \mathbf{1}_{E_\varepsilon}(s)ds,
\end{split}
\end{equation*}
we observe that
\begin{equation*}
\begin{split}
&E^Q\big[K_t^{1,\varepsilon,1} Y_t^{1,\varepsilon,1}\,\big|\,\mathcal{F}_t^{Y}\big]
=E^Q\big[L_t(L_t^{-1} K_t^{1,\varepsilon,1})\int_0^t \delta\sigma(s)\mathbf{1}_{E_\varepsilon}(s)dB_s^1\,\big|\,\mathcal{F}_T^{Y}\big]\\
=&E^Q\big[(L_t^{-1} K_t^{1,\varepsilon,1})\int^t_0 D_s^1[L_t] \delta\sigma(s)\mathbf{1}_{E_\varepsilon}(s)ds\,\big|\,\mathcal{F}_T^{Y}\big]+E^Q\big[L_t \int_0^t D_s^1[L_t^{-1} K_t^{1,\varepsilon,1}]\delta\sigma(s)\mathbf{1}_{E_\varepsilon}(s)ds\,\big|\,\mathcal{F}_T^{Y}\big]\\
=&:J_t^{1,\varepsilon} + J_t^{2,\varepsilon}.
\end{split}
\end{equation*}
As $\displaystyle L_t=\exp\{\int_0^t h(X_r,u_r)dY_r - \frac{1}{2} \int_0^t\big|h(X_r,u_r)\big|^2dr\}$, and, for $s\leq t$,
$$D_s^1[L_t]=L_t\{\int_s^t h_x(X_r,u_r)dY_r - \int_s^t (h\cdot h_x)(X_r,u_r)dr\}\sigma (\mu_s,u_s),$$
we have $\displaystyle\big( E^Q[\sup_{0\leq s\leq t\leq T}|D_s^1[L_t]|^p]\big)^{1/p}\leq C_p,\ p\geq 2.$ This yields the following estimates
for $J_t^{1,\varepsilon}$:
$$\Big(E^Q\big[\sup_{0\leq t\leq T}\big|J_t^{1,\varepsilon}\big|^p\big]\Big)^{\frac{1}{p}}\leq C_p\Big(E^Q\big[\sup_{0\leq t\leq T}\big|L_t^{-1} K_t^{1,\varepsilon,1}\big|^{2p}\big]\Big)^{\frac{1}{2p}}\cdot \varepsilon\leq C_p\cdot \varepsilon^{\frac{3}{2}},\ \varepsilon>0,\ p\geq 2.$$
On the other hand, recalling the formula for $D_s^1[L_t^{-1}K_t^{1,\varepsilon,1}]$, for $J_t^{2,\varepsilon}$ we have:
\begin{equation*}
\begin{split}
J_t^{2,\varepsilon}&=E^Q\big[L_t\int^t_0\big\{h(s,t)\delta\sigma(s)\mathbf{1}_{E_\varepsilon}(s)+{\varphi_{s,t}(\varepsilon)}\sqrt\varepsilon\big\}\delta\sigma(s)
\mathbf{1}_{E_\varepsilon}(s)ds\,\big|\,\mathcal{F}_T^{Y}\big]\\
&=\int_0^tE^Q\big[L_t{\varphi_{s,t}(\varepsilon)}\delta\sigma(s)\,\big|\,\mathcal{F}_T^{Y}\big]\sqrt\varepsilon\mathbf{1}_{E_\varepsilon}(s)ds +E^Q\big[L_t\int_0^th(s,t)|\delta \sigma(s)|^2\mathbf{1}_{E_\varepsilon}(s)ds\,\big|\,\mathcal{F}_T^{Y}\big].\\
\end{split}
\end{equation*}
We observe that $\displaystyle \int_0^tE^Q\big[L_t{\varphi_{s,t}(\varepsilon)}\delta\sigma(s)\,\big|\,\mathcal{F}_T^{Y}\big]\sqrt\varepsilon \mathbf{1}_{E_\varepsilon}(s)ds=:
{\varphi_t(\varepsilon)\varepsilon^{\frac{3}{2}}},$ where from our estimate for $\varphi_{s,t}(\varepsilon)$ it follows that
 $\displaystyle\Big(E^Q\big[\sup_{t\in[0,T]} |\varphi_t(\varepsilon)|^p\big]\Big)^{\frac{1}{p}}\leq C_p.$\ Consequently, combining the above estimates we see that
\begin{equation}\label{8star}
\begin{split}
&E^Q[K_t^{1,\varepsilon}Y_t^{1,\varepsilon}\,\big|\,\mathcal{F}_t^{Y}\big]=E^Q[K_t^{1,\varepsilon}Y_t^{1,\varepsilon,1}\,\big|\,
\mathcal{F}_t^{Y}\big]+\varepsilon \eta_t(\varepsilon) = J_t^{1,\varepsilon}+J_t^{2,\varepsilon} +\varepsilon\eta_t(\varepsilon)\\
=&E^Q\big[L_t H_{\varepsilon} (t)\,\big|\,\mathcal{F}_t^{Y}\big]+\varepsilon\eta_t(\varepsilon)=E^Q\big[L_t\int_0^t h(s,t)|\delta\sigma(s)|^2\mathbf{1}_{E_\varepsilon}(s)ds\,
\big|\,\mathcal{F}_t^{Y}\big]+\varepsilon\eta_t(\varepsilon), \ \varepsilon>0,
\end{split}
\end{equation}
where $\displaystyle E^Q[\sup_{t\in[0,T]} |\eta_t(\varepsilon)|^p]\rightarrow 0,\ \varepsilon \searrow 0.$\\

\noindent Notice that
\begin{equation*}
\begin{split}
H_{\varepsilon}(t) &= \int_0^t h(s,t) |\delta\sigma(s)|^2\mathbf{1}_{E_\varepsilon}(s)ds\\
&=\int_0^t\big(\int_s^t h_x(X_r,u_r)dY_r - \int^t_s(h\cdot h_x)(X_r,u_r)dr\big)|\delta\sigma(s)|^2\mathbf{1}_{E_\varepsilon}(s)ds\\
&=\int_0^t\big(\int_0^r|\delta\sigma(s)|^2\mathbf{1}_{E_\varepsilon}(s)ds\big)\big(h(X_r,u_r)dY_r-(h\cdot h_x)(X_r,u_r)dr\big).
\end{split}
\end{equation*}
Now we need the following assumption for $h$.\\
\noindent\textbf{Assumption (H3)}:\\
(ii)\ $h(t,x,\gamma,u)=h_0(t,x,\gamma)+\phi(x)h_1(t,\gamma,u),\ (t,x,\gamma,u)\in [0,T]\times\mathbb{R}\times\mathcal{P}_{2}(\mathbb{R})\times U,$\\
where $\phi:\ \mathbb{R}\rightarrow \mathbb{R}$ is bounded Lipschitz function.\\

\noindent Note that now
\begin{equation*}
    \begin{split}
\delta h(t) = \phi(X_t)\delta h_1(t) = \phi(X_t)\big(h_1(t,\mu_t,v_t)-h_1(t,\mu_t,u_t)\big),
\end{split}
\end{equation*}
and $\delta h_1$ is $\mathbb{F}^Y $-adapted and bounded. Now let us introduce the following notations:\\
$\displaystyle\Gamma^1_t :=\frac{L_t}{E^Q[L_t\,|\,\mathcal{F}_t^Y]},\ \ \Gamma_t := \frac{1}{E^Q[L_t\,|\,\mathcal{F}_t^Y]}
\Big(X_t - \frac{E^Q[X_tL_t\,|\,\mathcal{F}_t^Y]}{E^Q[L_t\,|\,\mathcal{F}_t^Y]}\Big)=\frac{1}{E^Q[L_t\,|\,\mathcal{F}_t^Y]}
\Big(X_t -E^P[X_t\,|\,\mathcal{F}_t^Y]\Big).$\\
Recall:
\begin{equation*}
\begin{split}
V_t^{1,\varepsilon}=&\theta_t(Y_t^{1,\varepsilon},K_t^{1,\varepsilon})=E^Q\big[\Gamma_t^1 Y_t^{1,\varepsilon}\,|\,\mathcal{F}_t^Y\big]+E^Q\big[\Gamma_t K_t^{1,\varepsilon}\,|\,\mathcal{F}_t^Y\big]\\
=&E^P\big[Y_t^{1,\varepsilon}|\mathcal{F}_t^Y\big]+\frac{E^Q[(X_t -E^P[X_t\,|\,\mathcal{F}_t^Y])K_t^{1,\varepsilon}\,|\,\mathcal{F}_t^Y]}{E^Q[L_t\,|\,\mathcal{F}_t^Y]},
\end{split}
\end{equation*}
and
$$
V_t^{2,\varepsilon}=\theta_t(Y_t^{2,\varepsilon},K_t^{2,\varepsilon})+\frac{E^Q[K_t^{1,\varepsilon}Y_t^{1,\varepsilon}\,|\,\mathcal{F}_t^Y]}{E^Q[L_t\,|\,\mathcal{F}_t^Y]}
-\frac{E^Q[K_t^{1,\varepsilon}|\,\mathcal{F}_t^Y]}{E^Q[L_t\,|\,\mathcal{F}_t^Y]}V_t^{1,\varepsilon}.
$$
From \eqref{8star} we immediately have:\ $\displaystyle \frac{ E^Q[K_t^{1,\varepsilon}Y_t^{1,\varepsilon}\,|\,\mathcal{F}_t^Y]}{E^Q[L_t\,|
\,\mathcal{F}_t^Y]}={\varepsilon \eta_t (\varepsilon)}+{E^P\big[H_\varepsilon (t)\,\big|\,\mathcal{F}_t^Y\big]}$, \
where $\eta_t(\varepsilon)$ stands here and in the following for $\displaystyle\Big(E^Q\big[\sup_{t\in[0,T]} \big|\eta_t(\varepsilon)\big|^p
\big]\Big)^{\frac{1}{p}}\rightarrow 0\ (\varepsilon\searrow0),$\ and
$\displaystyle \Big(E^Q\big[\big|\eta_t(\varepsilon)\big|^p\big]\Big)^{\frac{1}{p}}
\leq C_p,\ t\in[0,T],\ \varepsilon >0,\ p\geq 2.$\\

From \eqref{5star} and Proposition \ref{EstofXL} we obtain
$$\frac{E^Q[K_t^{1,\varepsilon}\,|\,\mathcal{F}_t^Y] E^Q\big[\Gamma_t^1 Y_t^{1,\varepsilon}\,|\,\mathcal{F}_t^Y\big]}
{E^Q[L_t\,|\,\mathcal{F}_t^Y]} = \varepsilon\eta_t(\varepsilon).$$
Moreover, from \eqref{7star} we get
$$E^Q\big[\Gamma_t K_t^{1,\varepsilon}\,|\,\mathcal{F}_t^Y\big]=E^Q\big[\Gamma_t L_t \int^t_0\delta h(s)\mathbf{1}_{E_\varepsilon}(s) dY_s\,\big|\,\mathcal{F}_t^{Y}\big]+\sqrt\varepsilon \eta_t(\varepsilon),\ t\in[0,T],\ \varepsilon>0.$$
Assuming $\displaystyle E_\varepsilon =  \big[t_0,t_0 +\varepsilon](\subset[0,T])$, we have
$$\int^t_0\delta h(s)\mathbf{1}_{E_\varepsilon}(s) dY_s = \int^{(t_0 +\varepsilon) \wedge t}_{t_0 \wedge t} \phi(X_s)\delta h_1(s) dY_s = \phi(X_{t_0 \wedge t})\int^{(t_0 +\varepsilon) \wedge t}_{t_0 \wedge t}\delta h_1(s) dY_s + r_t^{1,\varepsilon},$$
where $\displaystyle r_t^{1,\varepsilon} = \int^{(t_0 +\varepsilon) \wedge t}_{t_0 \wedge t}\big(  \phi(X_s)- \phi(X_{t_0 \wedge t} )\big)
\delta h_1(s)dY_s,$ and
\begin{equation*}
\begin{split}
&\Big(E^Q[\sup_{0\leq t\leq T}|r_t^{1,\varepsilon}|^p]\Big)^{\frac{1}{p}}\leq C_p\Big(E^Q\Big[\big(\int^{(t_0 +\varepsilon) \wedge t}_{t_0 \wedge t}|\phi(X_s)-\phi(X_{t_0 \wedge t})|^2ds\big)^{\frac{p}{2}}\Big]\Big)^{\frac{1}{p}}\\
    \leq&\ C_p \sqrt\varepsilon\big(E^Q\big[\sup_{s\in[t_0,t_0 +\varepsilon]}\big|\phi(X_s)-\phi(X_{t_0 })\big|^p\big]\big)^{\frac{1}{p}}\leq C_p \varepsilon.
\end{split}
\end{equation*}
Consequently,
\begin{equation*}
\begin{split}
E^Q\big[\Gamma_t K_t^{1,\varepsilon}\,|\,\mathcal{F}_t^Y\big]&=E^Q\big[\Gamma_t L_t\phi(X_{t_0 \wedge t} )\int^{t}_{0} \delta h_1(s) \mathbf{1}_{E_\varepsilon}(s) dY_s\,\big|\,\mathcal{F}_t^{Y}\big]+\sqrt\varepsilon\eta_t (\varepsilon)\\
&=E^Q\big[\Gamma_t L_t\phi(X_{t_0 \wedge t} )\,\big|\,\mathcal{F}_t^{Y}\big]\int^{t}_{0} \delta h_1(s) \mathbf{1}_{E_\varepsilon}(s) dY_s
+\sqrt\varepsilon\eta_t (\varepsilon),
\end{split}
\end{equation*}
and
\begin{equation}\label{4.39-1}
  E^Q\big[ K_t^{1,\varepsilon}\,|\,\mathcal{F}_t^Y\big]=E^Q\big[ L_t\phi(X_{t_0 \wedge t} )\,\big|\,\mathcal{F}_t^{Y}\big]\int^{t}_{0}
\delta h_1(s) \mathbf{1}_{E_\varepsilon}(s) dY_s+\sqrt\varepsilon\eta_t (\varepsilon).
\end{equation}
Combined with \eqref{5star} we have:
\begin{equation}\label{tri1}
\begin{split}
V_t^{1,\varepsilon}=&\theta_t(Y_t^{1,\varepsilon},K_t^{1,\varepsilon})=E^Q\big[\Gamma_t K_t^{1,\varepsilon}|\mathcal{F}_t^Y\big]
+\sqrt\varepsilon\eta_t (\varepsilon)\\
=&E^Q\big[\Gamma_t L_t\phi(X_{t_0 \wedge t} )\,\big|\,\mathcal{F}_t^{Y}\big]\int^{t}_{0} \delta h_1(s) \mathbf{1}_{E_\varepsilon}(s) dY_s
+\sqrt\varepsilon\eta_t (\varepsilon),
\end{split}
\end{equation}
and from \eqref{4.39-1} and \eqref{tri1}
$$\frac{E^Q[K_t^{1,\varepsilon}|\,\mathcal{F}_t^Y]}{E^Q[L_t\,|\,\mathcal{F}_t^Y]}V_t^{1,\varepsilon} =
E^Q\big[\Gamma_t^1 \phi(X_{t_0 \wedge t} )\big|\mathcal{F}_t^{Y}\big]E^Q\big[\Gamma_t L_t\;\phi(X_{t_0 \wedge t} )\big|\mathcal{F}_t^{Y}\big]\Big(\int^{t}_{0}
\delta h_1(s)\mathbf{1}_{E_\varepsilon}(s) dY_s\Big)^2 +\varepsilon\eta_t(\varepsilon).$$
Furthermore, from \eqref{8star} we get
\begin{equation*}
\begin{split}
V_t^{2,\varepsilon}=&\theta_t(Y_t^{2,\varepsilon},K_t^{2,\varepsilon})+{{E^P\big[H_\varepsilon (t)\,\big|\,\mathcal{F}_t^Y\big]}}\\
&-E^Q\big[\Gamma^1_t\phi(X_{t_0 \wedge t} )\,\big|\,\mathcal{F}_t^{Y}\big]E^Q\big[\Gamma_t L_t\phi(X_{t_0 \wedge t} )\,\big|\,\mathcal{F}_t^{Y}\big]\Big(\int^{t}_{0}\delta h_1(s)\mathbf{1}_{E_\varepsilon}(s) dY_s\Big)^2+\varepsilon\eta_t(\varepsilon).
\end{split}
\end{equation*}
Finally, we obtain
\begin{equation}\label{1jing}
\begin{split}
V_t^{1,\varepsilon}+&V_t^{2,\varepsilon}=\theta_t(Y_t^{1,\varepsilon}+Y_t^{2,\varepsilon} , K_t^{1,\varepsilon}+K_t^{2,\varepsilon})
+{{E^P\big[H_\varepsilon (t)\,\big|\,\mathcal{F}_t^Y\big]}}\\
&-E^Q\big[ \Gamma^1_t\phi(X_{t_0 \wedge t} )\,\big|\,\mathcal{F}_t^{Y}\big]E^Q\big[\Gamma_t L_t\phi(X_{t_0 \wedge t} )\,\big|\,\mathcal{F}_t^{Y}\big]\cdot\Big(\int^{t}_{0} \delta h_1(s) \mathbf{1}_{E_\varepsilon}(s) dY_s\Big)^2+\varepsilon\eta_t(\varepsilon).
\end{split}
\end{equation}
Putting $\displaystyle N_t^{1,\varepsilon}:=\int_0^t\delta h_1(s) \mathbf{1}_{E_\varepsilon}(s) dY_s,$ we have
$$\Big(\int^{t}_{0}\delta h_1(s)\mathbf{1}_{E_\varepsilon}(s) dY_s\Big)^2=\big(N_t^{1,\varepsilon}\big)^2=2\int^t_0 N_t^{1,\varepsilon}\delta h_1(s)\mathbf{1}_{E_\varepsilon}(s) dY_s +\int^t_0 \big(\delta h_1(s)\big)^2 \mathbf{1}_{E_\varepsilon}(s)ds.$$
Now consider $\displaystyle\theta\in L^{\infty-}_{\mathbb{F}}\big([0,T],L^p(Q)\big)(:=\bigcap_{q\geq 2}L^{q}_{\mathbb{F}}\big([0,T],L^p(Q)\big)),\ p\geq 2$, where\ $\displaystyle L^{q}_{\mathbb{F}}\big([0,T],L^p(Q)\big):=\big\{\theta\in L^{2}_{\mathbb{F}}\big([0,T],Q\big)\big|\int_0^T(E^Q[|\theta_t|^p])^{\frac{q}{p}}dt<+\infty \big\},\ q\geq 2$, then there exists a jointly measurable random
field $\{\theta_{t,s},0\leq s\leq t\}$ such that $\theta_{t,\cdot}=(\theta_{t,\cdot}^1,\theta_{t,\cdot}^2) \in L^{2}_{\mathbb{F}}\big([0,T],Q;\mathbb{R}^2\big),$ with
$$\theta_t = E^Q[\theta_t]+\int_0^t \theta^1_{t,s} dY_s +\int_0^t \theta_{t,s}^2 dB_s^1,\ Q\text{-a.s.}$$
Note that, thanks to the Burkholder-Davis-Gundy inequality,
$$\Big(E^Q\big[\big(\int_0^t (|\theta^1_{t,s}|^2+| \theta^2_{t,s}|^2) ds\big)^{\frac{p}{2}}\big]\Big)^{\frac{1}{p}} \leq C_p\Big(E^Q\big[|\theta_t|^p\big]\Big)^{\frac{1}{p}}\leq C_p,\ t\in[0,T],\ p\geq 2.$$
Thus, we can deduce that
\begin{equation*}
\begin{split}
&\Big|E^Q\big[\theta_t(\int^t_0 N_t^{1,\varepsilon}\delta h_1(s)\mathbf{1}_{E_\varepsilon}(s) dY_s )\big]\Big|=\Big|E^Q\big[\int^t_0 \theta^1_{t,s} N_t^{1,\varepsilon}\delta h_1(s)\mathbf{1}_{E_\varepsilon}(s) ds\big]\Big|\\
\leq &C E^Q\Big[\sup_{s\leq t} \big|N_s^{1,\varepsilon}\big|\int^t_0|\theta^1_{t,s}| \mathbf{1}_{E_\varepsilon}(s)ds\Big]
\leq C E^Q\Big[\sup_{s\leq t} \big|N_s^{1,\varepsilon}\big|\big(\int^t_0|\theta^1_{t,s}|^2 \mathbf{1}_{E_\varepsilon}(s)ds\big)^{\frac{1}{2}}\Big]\sqrt\varepsilon\\
\leq &\varepsilon\cdot C \Big(E^Q[\int^t_0|\theta^1_{t,s}|^2 \mathbf{1}_{E_\varepsilon}(s)ds]\Big)^{\frac{1}{2}} = \varepsilon \rho_t(\varepsilon),\ \ \mbox{where}
\end{split}
\end{equation*}
\small$$\rho_t(\varepsilon):=C\Big(E^Q\big[\int^t_0|\theta^1_{t,s}|^2 \mathbf{1}_{E_\varepsilon}(s)ds\big]\Big)^{\frac{1}{2}} \leq C  \Big(E^Q\big[\int^t_0\big|\theta^1_{t,s}\big|^2 ds\big]\Big)^{\frac{1}{2}}\leq C
\Big(E^Q\big[\big|\theta_t\big|^2\big]\Big)^{\frac{1}{2}}\leq C,\ t\in [0,T],\ \varepsilon>0,$$
\normalsize
and, from the dominated convergence theorem we see that $\rho_t(\varepsilon)\rightarrow 0\ (\varepsilon\searrow0)$. Consequently,
\begin{equation}\label{2jing}
    \begin{split}
\Big(\int_0^t \delta h_1(s)\mathbf{1}_{E_\varepsilon}(s)dY_s\Big)^2 = \int_0^t\big(\delta h_1(s)\big)^2
\mathbf{1}_{E_\varepsilon}(s)ds + \varepsilon \varphi_t(\varepsilon),
\end{split}
\end{equation}
where $E^Q\big[\theta_t\varphi_t(\varepsilon)\big] \rightarrow 0\ (\varepsilon \searrow 0)$ and $\big|E^Q\big[\theta_t\varphi_t(\varepsilon)\big]\big|\leq C,\ t\in [0,T],\ \varepsilon > 0.$
We now rewrite the terms of the right-hand side of (\ref{MR1}). For this we notice that combining \eqref{1jing} and \eqref{2jing} yields
\small\begin{equation*}
\begin{split}
{\rm i)}\ \ \ \ &E^Q\Big[\widetilde{f}_\mu^* (t) L_t(V_t^{1,\varepsilon}+V_t^{2,\varepsilon})\Big]\\
=&E^Q\Big[\widetilde{f}_\mu^* (t) L_t \theta_t(Y_t^{1,\varepsilon}+Y_t^{2,\varepsilon} , K_t^{1,\varepsilon}+K_t^{2,\varepsilon})\Big]
+{E^Q\Big[\widetilde{f}_\mu^* (t) L_tE^P\big[H_\varepsilon (t)\,\big|\,\mathcal{F}_t^Y\big]\Big]}\\
&-E^Q\Big[\widetilde{f}_\mu^* (t) L_t E^Q\big[\Gamma_t^1 \phi(X_{t_0 \wedge t} )\,\big|\,\mathcal{F}_t^{Y}\big]E^Q\big[\Gamma_t L_t
\phi(X_{t_0 \wedge t} )\,\big|\,\mathcal{F}_t^{Y}\big]\Big(\int_0^t \delta h_1(s)\mathbf{1}_{E_\varepsilon}(s)dY_s\Big)^2\Big]+
\varepsilon \rho_t(\varepsilon)\\
=&E^Q\Big[E^Q\big[\widetilde{f}_\mu^* (t) L_t \,\big|\,\mathcal{F}_t^{Y}\big]\big(\Gamma^1_t(Y_t^{1,\varepsilon}+Y_t^{2,\varepsilon})
+\Gamma_t(K_t^{1,\varepsilon}+K_t^{2,\varepsilon})\big)\Big]+{E^Q\Big[\widetilde{f}_\mu^* (t) L_tE^P\big[H_\varepsilon (t)\,\big|\,\mathcal{F}_t^Y\big]\Big]}\\
&-E^Q\Big[\widetilde{f}_\mu^* (t) L_t E^Q\big[\Gamma_t^1 \phi(X_{t_0 \wedge t} )\,\big|\,\mathcal{F}_t^{Y}\big]E^Q\big[\Gamma_t L_t\phi(X_{t_0 \wedge t} )\,\big|\,\mathcal{F}_t^{Y}\big]
\Big(\int_0^t \delta h_1(s)\mathbf{1}_{E_\varepsilon}(s)dY_s\Big)^2\Big]+\varepsilon \rho_t(\varepsilon)\\
=&E^Q\Big[E^Q\big[\widetilde{f}_\mu^* (t) L_t \,\big|\,\mathcal{F}_t^{Y}\big]\big(\Gamma^1_t(Y_t^{1,\varepsilon}+Y_t^{2,\varepsilon})
+\Gamma_t(K_t^{1,\varepsilon}+K_t^{2,\varepsilon})\big)\Big]+{E^Q\Big[\widetilde{f}_\mu^* (t) L_tE^P\big[H_\varepsilon (t)\,\big|\,\mathcal{F}_t^Y\big]\Big]}\\
&-E^Q\Big[\widetilde{f}_\mu^* (t) L_t E^Q\big[\Gamma_t^1 \phi(X_{t_0 \wedge t} )\,\big|\,\mathcal{F}_t^{Y}\big]E^Q\big[\Gamma_t L_t\phi(X_{t_0 \wedge t} )\,\big|\,\mathcal{F}_t^{Y}\big]
\Big(\int_0^t (\delta h_1(s))^2\mathbf{1}_{E_\varepsilon}(s)ds\Big)\Big]+\varepsilon \rho_t(\varepsilon)\\
=&E^Q\Big[E^Q\big[\widetilde{f}_\mu^* (t) L_t \,\big|\,\mathcal{F}_t^{Y}\big]\big(\Gamma^1_t(Y_t^{1,\varepsilon}+Y_t^{2,\varepsilon})
+\Gamma_t(K_t^{1,\varepsilon}+K_t^{2,\varepsilon})\big)\Big]+{E^Q\Big[\widetilde{f}_\mu^* (t) L_tE^P\big[H_\varepsilon (t)\,\big|\,\mathcal{F}_t^Y\big]\Big]}\\
&-E^Q\Big[\!\widetilde{f}_\mu^* (t) L_t\! \int^{(t_0 +\varepsilon) \wedge t}_{t_0 \wedge t}\!\!E^Q\big[\Gamma^1_t\phi(X_{t_0\wedge t})\delta h_1(s)\,\big|\,\mathcal{F}_t^{Y}\!\big]E^Q
\big[\!\Gamma_t L_t\phi(X_{t_0 \wedge t} )\delta h_1(s)\,\big|\,\mathcal{F}_t^{Y}\!\big]\mathbf{1}_{E_\varepsilon}(s)ds\!\Big]\!+\!\varepsilon \rho_t(\varepsilon)\\
=&E^Q\Big[E^Q\big[\widetilde{f}_\mu^* (t) L_t \,\big|\,\mathcal{F}_t^{Y}\big]\big(\Gamma^1_t(Y_t^{1,\varepsilon}+Y_t^{2,\varepsilon})
+\Gamma_t(K_t^{1,\varepsilon}+K_t^{2,\varepsilon})\big)\Big]+{E^Q\Big[\widetilde{f}_\mu^* (t) L_tE^P\big[H_\varepsilon (t)\,\big|\,\mathcal{F}_t^Y\big]\Big]}\\
&-E^Q\Big[\widetilde{f}_\mu^* (t) L_t \int^{t}_{0}E^Q\big[\Gamma^1_t\delta h(s)\,\big|\,\mathcal{F}_t^{Y}\big]E^Q\big[\Gamma_t L_t \delta h(s)\,\big|\,\mathcal{F}_t^{Y}\big]\mathbf{1}
_{E_\varepsilon}(s)ds\Big]+\varepsilon \rho_t(\varepsilon),
\end{split}
\end{equation*}
\normalsize
for $\rho_t(\varepsilon)\rightarrow 0\ (\varepsilon\searrow0),\ \big|\rho_t(\varepsilon)\big|\leq C,\ \varepsilon>0,\ t\in [0,T].$\ This shows that\\
\begin{equation*}
\begin{split}
E^Q\big[\widetilde{f}_\mu^* (t) L_t(V_t^{1,\varepsilon}+V_t^{2,\varepsilon})\big]=&E^Q\Big[E^Q\big[\widetilde{f}_\mu^* (t) L_t \,\big|\,\mathcal{F}_t^{Y}\big]\big\{\Gamma^1_t(Y_t^{1,
\varepsilon}+Y_t^{2,\varepsilon})+\Gamma_t(K_t^{1,\varepsilon}+K_t^{2,\varepsilon})\\
&\qquad-\int_0^t E^Q\big[\Gamma^1_t\delta h(s)\,\big|\,\mathcal{F}_t^{Y}\big]E^Q\big[\Gamma_t L_t \delta h(s)\,\big|\,\mathcal{F}_t^{Y}\big]\mathbf{1}_{E_\varepsilon}(s)ds  \big\}\Big]\\
&+{E^Q\Big[\widetilde{f}_\mu^* (t) L_tE^P\big[H_\varepsilon (t)\,\big|\,\mathcal{F}_t^Y\big]\Big]}+\varepsilon \rho_t(\varepsilon).
\end{split}
\end{equation*}
{\rm ii)} We now combine \eqref{tri1}, \eqref{7star} and \eqref{2jing}. This allows to get
\begin{align}
\notag&E^Q\Big[\widetilde{f}_\mu^* (t) V_t^{1,\varepsilon}K_t^{1,\varepsilon}\Big]\\
\notag=&E^Q\Big[\widetilde{f}_\mu^* (t)E^Q\big[\Gamma_t L_t\phi(X_{t_0 \wedge t} )\,\big|\,\mathcal{F}_t^{Y}\big]\big(\int_0^t \delta h_1(s)\mathbf{1}_{E_\varepsilon}(s)dY_s\big)
K_t^{1,\varepsilon}\Big]+\varepsilon \rho_t(\varepsilon)\\
\label{P10,ii}=& E^Q\Big[ E^Q[\Gamma_t L_t\phi(X_{t_0 \wedge t} )\,\big|\,\mathcal{F}_t^{Y}]\big(\int_0^t \delta h_1(s)\mathbf{1}_{E_\varepsilon}(s)dY_s\big)E^Q\big[\widetilde{f}_\mu^* (t)
K_t^{1,\varepsilon}\,\big|\,\mathcal{F}_t^{Y}\big]\Big] +\varepsilon \rho_t(\varepsilon)\\
\notag=& E^Q\Big[ E^Q[\Gamma_t L_t\phi(X_{t_0 \wedge t} )\,\big|\,\mathcal{F}_t^{Y}]E^Q\big[\widetilde{f}_\mu^* (t)L_t\! \int_0^t\! \delta h(s)\mathbf{1}_{E_\varepsilon}(s)dY_s
\,\big|\,\mathcal{F}_t^{Y}\!\big]\big(\!\int_0^t\! \delta h_1(s)\mathbf{1}_{E_\varepsilon}(s)dY_s\big)\!\Big]+\varepsilon \rho_t(\varepsilon)\\
\notag=& E^Q\Big[ E^Q[\Gamma_t L_t\phi(X_{t_0 \wedge t} )\,\big|\,\mathcal{F}_t^{Y}]E^Q\big[\widetilde{f}_\mu^* (t)L_t\phi(X_{t_0 \wedge t} )\,\big|\,\mathcal{F}_t^{Y}\big]\big(
    \int_0^t \delta h_1(s)\mathbf{1}_{E_\varepsilon}(s)dY_s\big)^2\Big]+\varepsilon \rho_t(\varepsilon)\\
\notag=& E^Q\Big[ E^Q[\Gamma_t L_t\phi(X_{t_0 \wedge t} )\,\big|\,\mathcal{F}_t^{Y}]E^Q\big[\widetilde{f}_\mu^* (t)L_t\phi(X_{t_0 \wedge t} )\,\big|\,\mathcal{F}_t^{Y}\big]\int_0^t
(\delta h_1(s))^2\mathbf{1}_{E_\varepsilon}(s)ds\Big]+\varepsilon \rho_t(\varepsilon)\\
\notag=& E^Q\Big[\int_0^t E^Q[\Gamma_t L_t \delta h(s)\,\big|\,\mathcal{F}_t^{Y}]E^Q\big[\widetilde{f}_\mu^* (t)L_t\delta h(s)\big|\mathcal{F}_t^{Y}\big]\mathbf{1}_{E_\varepsilon}(s)ds
\Big]+\varepsilon \rho_t(\varepsilon).
\end{align}
{\rm iii)} In the following we use \eqref{tri1} and \eqref{2jing}. Then,
\begin{equation}\label{P10,iii}
\begin{split}
&E^Q\Big[\widetilde{f}_{z\mu}^* (t) L_t (V_t^{1,\varepsilon})^2\Big]\\
=& E^Q\Big[\widetilde{f}_{z\mu}^* (t) L_t\big( E^Q[\Gamma_t L_t\phi(X_{t_0 \wedge t} )\,\big|\,\mathcal{F}_t^{Y}\big]\big)^2 \big(\int_0^t \delta h_1(s)\mathbf{1}_{E_\varepsilon}(s)dY_s\big)^2\Big]+\varepsilon \rho_t(\varepsilon)\\
=& E^Q\Big[\widetilde{f}_{z\mu}^* (t) L_t\big( E^Q[\Gamma_t L_t\phi(X_{t_0 \wedge t} )\,\big|\,\mathcal{F}_t^{Y}\big]\big)^2 \int_0^t (\delta h_1(s))^2\mathbf{1}_{E_\varepsilon}(s)ds\Big]+\varepsilon \rho_t(\varepsilon)\\
=& E^Q\Big[\widetilde{f}_{z\mu}^* (t) L_t \int_0^t \big( E^Q[\Gamma_t L_t\delta h(s)\,\big|\,\mathcal{F}_t^{Y}\big]\big)^2\mathbf{1}_{E_\varepsilon}(s)ds\Big]+\varepsilon \rho_t(\varepsilon).
\end{split}
\end{equation}
{\rm iv)}\ In analogy to {\rm i)} we have:
\begin{equation*}
\begin{split}
E^Q\Big[\widetilde{\Phi}_{\mu}^{*} (T) L_T (V_T^{1,\varepsilon} + V_T^{2,\varepsilon} )\Big]
=&E^Q \Big[E^Q[\widetilde{\Phi}_{\mu}^{*} (T) L_T \,\big|\,\mathcal{F}_T^{Y}\big]\big\{\Gamma^1_T(Y_T^{1,\varepsilon} +Y_T^{2,\varepsilon})+\Gamma_T(K_T^{1,\varepsilon} +K_T^{2,\varepsilon})\\
&\quad-\int^T_0 E^Q \big[\Gamma^1_T \delta h(s)\,\big|\,\mathcal{F}_T^{Y}\big] E^Q\big[\Gamma_T L_T \delta h(s)\,\big|\,\mathcal{F}_T^{Y}\big]\mathbf{1}_{E_\varepsilon}(s)ds\big\}\Big]\\
&+{E^Q\Big[\widetilde{\Phi}_{\mu}^{*} (T) L_T E^P\big[H_\varepsilon (T)\,\big|\,\mathcal{F}_T^Y\big]\Big]}+\varepsilon \rho_T(\varepsilon).\qquad\qquad\qquad\qquad\qquad\qquad\qquad\qquad\qquad\qquad\qquad\qquad\qquad\qquad\qquad
\end{split}
\end{equation*}
{\rm v)} Moreover, similar to {\rm ii) we get
\begin{equation*}
\begin{split}
E^Q\Big[\widetilde{\Phi}_{\mu}^{*} (T)V_T^{1,\varepsilon}K_T^{1,\varepsilon}\Big]
=E^Q\Big[\int^T_0 E^Q\big[\Gamma_T L_T \delta h(s)\,\big|\,\mathcal{F}_T^{Y}\big]E^Q\big[\widetilde{\Phi}_{\mu}^{*} (T)L_T \delta h(s)\big|\mathcal{F}_T^{Y}\big]\mathbf{1}_{E_\varepsilon}(s)ds\Big]+\varepsilon \rho_T(\varepsilon).\qquad
\end{split}
\end{equation*}
{\rm vi)} And finally, in analogy to {\rm iii)} we see
\begin{equation*}
\begin{split}
E^Q\Big[\widetilde{\Phi}_{z\mu}^* (T) L_T (V_T^{1,\varepsilon})^2\Big]= E^Q\Big[\widetilde{\Phi}_{z\mu}^* (T) L_T \int_0^T \big( E^Q[\Gamma_T L_T \delta h(s)\big|\mathcal{F}_T^{Y}\big]\big)^2 \mathbf{1}_{E_\varepsilon}(s)ds\Big]+\varepsilon \rho_T(\varepsilon).\qquad\qquad
\end{split}
\end{equation*}

We substitute ${\rm i)}-{\rm vi)}$ in (\ref{MR1}). This yields
\begin{align}\label{new4.33}
\notag 0\leq J(u^\varepsilon)-J&(u)=E^Q\Big[\Phi_x (T)(Y_T^{1,\varepsilon}+Y_T^{2,\varepsilon})+\widetilde{E}^Q\big[\int_0^{U_T} \Phi_{\mu}^{*} (T,y)dy(K_T^{1,\varepsilon} +K_T^{2,\varepsilon})\big]\Big]\\
\notag&+\Big\{E^Q\Big[\widetilde{E}^Q\big[E^Q[\widetilde{\Phi}_{\mu}^{*} (T)L_T\,\big|\,\mathcal{F}_T^{Y}\big]\big\{\Gamma^1_T(Y_T^{1,\varepsilon} +Y_T^{2,\varepsilon})+\Gamma_T(K_T^{1,\varepsilon} +K_T^{2,\varepsilon})\big\}\big]\\
\notag&\ \ \qquad -\widetilde{E}^Q\big[E^Q[\widetilde{\Phi}_{\mu}^{*} (T)L_T\,\big|\,\mathcal{F}_T^{Y}\big]\int_0^{T}E^Q[\Gamma_T^1\delta h(s) \,\big|\,\mathcal{F}_T^{Y}\big]\\
\notag&\qquad\qquad\qquad\qquad\qquad\qquad\cdot E^Q\big[\Gamma_T L_T \delta h(s) \,\big|\,\mathcal{F}_T^{Y}\big]\mathbf{1}_{E_\varepsilon}(s)ds\big]\Big]\Big\}\\
\notag&+{E^Q\Big[\widetilde{E}^Q\big[\widetilde{\Phi}_{\mu}^{*} (T)\big]L_TE^P\big[H_\varepsilon (T)\,\big|\,\mathcal{F}_T^Y\big]\Big]}\\
\notag&+E^Q\Big[\int^T_0\Big(f_x(t)(Y_t^{1,\varepsilon}+Y_t^{2,\varepsilon})+\widetilde{E}^Q\big[\int^{U_t}_0 f_{\mu}^* (t,y)dy(K_t^{1,\varepsilon} +K_t^{2,\varepsilon})\big]\Big)dt\Big]\\
\notag&+\Big\{E^Q\Big[\widetilde{E}^Q\big[\int^T_0 E^Q[\widetilde{f}^*_{\mu}(t)L_t\,\big|\,\mathcal{F}_t^{Y}\big]\big(\Gamma^1_t(Y_t^{1,\varepsilon} +Y_t^{2,\varepsilon}) + \Gamma_t(K_t^{1,\varepsilon} + K_t^{2,\varepsilon})\big)dt\big]\\
&\ \ \qquad -\widetilde{E}^Q\big[\int^T_0 E^Q[\widetilde{f}^*_{\mu}(t)L_t\,\big|\,\mathcal{F}_t^{Y}\big](\int^t_0E^Q[\Gamma_t^1\delta h(s)\,\big|\,\mathcal{F}_t^{Y}\big]\\
\notag&\qquad\qquad\qquad\qquad\qquad\qquad\cdot E^Q[\Gamma_t L_t\delta h(s)\,\big|\,\mathcal{F}_t^{Y}\big]\mathbf{1}_{E_\varepsilon}(s)ds\big)dt\big]\Big]\Big\}\\
\notag&+{E^Q\Big[\int^T_0\widetilde{E}^Q\big[ \widetilde{f}^*_{\mu}(t)\big]L_t E^P\big[H_\varepsilon (t)\,\big|\,\mathcal{F}_t^{Y}\big]dt\Big]}\\
\notag&+\frac{1}{2}E^Q\Big[\Phi_{xx} (T) (Y_T^{1,\varepsilon})^2\Big]+E^Q\Big[\widetilde{E}^Q\big[\int^T_0E^Q[\Gamma_T L_T \delta h(s)\,\big|\,\mathcal{F}_T^{Y}\big]\\
\notag&\qquad\qquad\qquad\qquad\qquad\qquad\cdot E^Q[\widetilde{\Phi}_{\mu}^* (T) L_T \delta h(s)\,\big|\,\mathcal{F}_T^{Y}\big]\mathbf{1}_{E_\varepsilon}(s)ds\big]\Big]\\
\notag&+\frac{1}{2}E^Q\Big[\widetilde{E}^Q\big[\widetilde{\Phi}_{z\mu}^* (T)L_T \int^T_0(E^Q[\Gamma_T L_T\delta h(s)\,\big|\,\mathcal{F}_T^{Y}\big])^2\mathbf{1}_{E_\varepsilon}(s)ds\big]\Big]\\
\notag&+E^Q\Big[\widetilde{E}^Q\big[\int^T_0\int^t_0E^Q[\Gamma_t L_t \delta h(s)\,\big|\,\mathcal{F}_t^{Y}\big]E^Q[\widetilde{f}^*_{\mu} (t)L_t \delta h(s) \,\big|\,\mathcal{F}_t^{Y}\big]\mathbf{1}_{E_\varepsilon}(s)dsdt\big]\Big]\\
\notag&+\frac{1}{2}E^Q\Big[\widetilde{E}^Q\big[\int^T_0 \int^t_0 \widetilde{f} ^*_{z\mu} (t) L_t\big(E^Q[\Gamma_t L_t \delta h(s)\,\big|\,\mathcal{F}_t^{Y}\big]\big)^2\mathbf{1}_{E_\varepsilon}(s)dsdt\big]\Big]\\
\notag&+\frac{1}{2}E^Q\Big[\int^T_0 f_{xx}(t)(Y_t^{1,\varepsilon})^2dt\Big]+E^Q\Big[\int_0^T\delta f(t)\mathbf{1}_{E_\varepsilon}(t)dt\Big] + o(\varepsilon),\ \ \mbox{ as } \varepsilon\searrow0.
\end{align}
Restructuring the above inequality we get
\begin{align}\label{new4.34}
\notag 0\leq J(u^\varepsilon)-&J(u)=E^Q\Big[\big(\Phi_x (T)+E^Q\big[\widetilde{E}^Q[\widetilde{\Phi}_{\mu}^*(T)]L_T\,\big|\,\mathcal{F}_T^{Y}]\Gamma_T^1\big)\big(Y_T^{1,\varepsilon}+Y_T^{2,\varepsilon}\big)\Big]\\
\notag&+E^Q\Big[\big(\widetilde{E}^Q\big[\int_0^{U_T} \Phi^*_{\mu} (T,y) dy \big] + E^Q[\widetilde{E}^Q[\widetilde{\Phi}_{\mu}^*(T)]L_T\,\big|\,\mathcal{F}_T^{Y}]\Gamma_T\big)\big(K_T^{1,\varepsilon}+K_T^{2,\varepsilon}\big)\Big]\\
\notag&+\frac{1}{2} E^Q\big[\Phi_{xx} (T)(Y_T^{1,\varepsilon})^2\big]+{E^Q\Big[\widetilde{E}^Q\big[\widetilde{\Phi}_{\mu}^{*} (T)\big]L_TE^P\big[H_\varepsilon (T)\,\big|\,\mathcal{F}_T^Y\big]\Big]}\\
\notag&+E^Q\Big[\int^T_0\Big(E^Q[\Gamma_T L_T \delta h(t)\,\big|\,\mathcal{F}_T^{Y}\big]\big\{E^Q[\widetilde{E}^Q[\widetilde{\Phi}^*_{\mu} (T)]L_T\delta h(t) \,\big|\,\mathcal{F}_T^{Y}\big]\\
\notag&\ \ \qquad-E^Q\big[\widetilde{E}^Q[\widetilde{\Phi}^*_{\mu} (T)]L_T\,\big|\,\mathcal{F}_T^{Y}\big]E^Q\big[\Gamma_T^1\delta h(t)\,\big|\,\mathcal{F}_T^{Y}\big]\big\}\\
\notag&\ \ +\frac{1}{2}\big(E^Q[\Gamma_T L_T\delta h(t) \,\big|\,\mathcal{F}_T^{Y}\big]\big)^2\widetilde{E}^Q[\Phi^*_{z\mu} (T)]L_T\Big)\mathbf{1}_{E_\varepsilon}(t)dt\Big]\\
&+E^Q\Big[\int^T_0 \big(f_x(t) + E^Q[\widetilde{E}^Q[\widetilde{f}_{\mu}^*(t)]L_t\,\big|\,\mathcal{F}_t^{Y}\big]\Gamma^1_t\big)\big(Y_t^{1,\varepsilon}+Y_t^{2,\varepsilon}\big)dt\Big]\\
\notag&+E^Q\Big[\int_0^T\big(\widetilde{E}^Q\big[\int_0^{U_t} f^*_{\mu}(t,y)dy\big] + E^Q[\widetilde{E}^Q[\widetilde{f}^*_{\mu}(t)]L_t\,\big|\,\mathcal{F}_t^{Y}]\Gamma_t\big)(K_t^{1,\varepsilon}+K_t^{1,\varepsilon})dt\Big]\\
\notag&+\frac{1}{2}E^Q\Big[\int^T_0 f_{xx} (t)(Y_t^{1,\varepsilon})^2dt\Big]+{E^Q\Big[\int^T_0\widetilde{E}^Q\big[ \widetilde{f}^*_{\mu}(t)\big]L_t E^P\big[H_\varepsilon (t)\,\big|\,\mathcal{F}_t^{Y}\big]dt\Big]}\\
\notag&+ E^Q\Big[\int_0^T \int_0^t\big(E^Q[\Gamma_t L_t \delta h(s)\,\big|\,\mathcal{F}_t^{Y}]\big\{E^Q\big[\widetilde{E}^Q[\widetilde{f}^*_{\mu} (t)]L_t \delta h(s)\,\big|\,\mathcal{F}_t^{Y}\big]\\
\notag&\ \ \qquad-E^Q\big[\widetilde{E}^Q[\widetilde{f}^*_{\mu} (t)]L_t\,\big|\,\mathcal{F}_t^{Y}\big]E^Q\big[\Gamma^1_t \delta h(s)\,\big|\,\mathcal{F}_t^{Y}\big]\big\}\\
\notag&\ \ +\frac{1}{2} \big(E^Q[\Gamma_t L_t \delta h(s)\,\big|\,\mathcal{F}_t^{Y}]\big)^2\widetilde{E}^Q[\widetilde{f}^*_{z\mu} (t)]L_t\big)\mathbf{1}_{E_\varepsilon}(s)dsdt\Big]\\
\notag&+E^Q\Big[\int^T_0\delta f(t)\mathbf{1}_{E_\varepsilon}(t)dt\Big] + o(\varepsilon),\ \ \mbox{ as } \varepsilon\searrow0.
\end{align}
Put, for $\phi = \Phi,\ f:$
\small
\begin{align}\label{gamma_s^phi}
\notag\gamma_t^{\phi} \big(\delta h(s)\big) &= E^Q\big[\Gamma_t L_t\delta h(s)\,\big|\,\mathcal{F}_t^{Y}\big]\big\{E^Q\big[\widetilde{E}^Q[
    \widetilde{\phi}^*_{\mu} (t)] L_t\delta h(s)\,\big|\,\mathcal{F}_t^{Y}\big]-E^Q\big[\widetilde{E}^Q[\widetilde{\phi}^*_{\mu} (t)]
    L_t\,\big|\,\mathcal{F}_t^{Y}\big]E^Q\big[\Gamma_t^1\delta h(s)\,\big|\,\mathcal{F}_t^{Y}\big]\big\}\\
    &\quad +\frac{1}{2}\Big(E^Q[\Gamma_t L_t \delta h(s)\,\big|\,
    \mathcal{F}_t^{Y}]\Big)^2\widetilde{E}^Q[\widetilde{\phi}^*_{z\mu} (t)] L_t,\ \ \ 0\leq s\leq t\leq T.
\end{align}
\normalsize
Then, from above we get
\begin{align}\label{reorganized4.34}
0\leq J(u^\varepsilon)-J(u)=&E^Q\Big[\big(\Phi_x (T)+E^Q\big[\widetilde{E}^Q[\widetilde{\Phi}_{\mu}^*(T)]L_T\,\big|\,\mathcal{F}_T^{Y}]\Gamma_T^1\big)\big(Y_T^{1,\varepsilon}+Y_T^{2,\varepsilon}\big)\Big]\notag\\
&+E^Q\Big[\big(\widetilde{E}^Q[\int_0^{U_T} \Phi^*_{\mu} (T,y) dy ] + E^Q[\widetilde{E}^Q[\widetilde{\Phi}_{\mu}^*(T)]L_T\,\big|\,\mathcal{F}_T^{Y}]\Gamma_T\big)\big(K_T^{1,\varepsilon}+K_T^{2,\varepsilon}\big)\Big]\notag\\
&+\frac{1}{2} E^Q\Big[\Phi_{xx} (T)(Y_T^{1,\varepsilon})^2\Big]+\frac{1}{2}E^Q\Big[\int^T_0 f_{xx} (t)(Y_t^{1,\varepsilon})^2dt\Big]\notag\\
&+{E^Q\Big[\widetilde{E}^Q\big[\widetilde{\Phi}_{\mu}^{*} (T)\big]L_TE^P\big[H_\varepsilon (T)\,\big|\,\mathcal{F}_T^Y\big]\Big]}\\
&+{E^Q\Big[\int^T_0\widetilde{E}^Q\big[ \widetilde{f}^*_{\mu}(t)\big]L_t E^P\big[H_\varepsilon (t)\,\big|\,\mathcal{F}_t^{Y}\big]dt\Big]}\notag\\
&+E^Q\Big[\int^T_0 \big(f_x(t) + E^Q[\widetilde{E}^Q[\widetilde{f}_{\mu}^*(t)]L_t\,\big|\,\mathcal{F}_t^{Y}\big]\Gamma^1_t\big)(Y_t^{1,\varepsilon}+Y_t^{2,\varepsilon}\big)dt\Big]\notag\\
&+E^Q\Big[\int_0^T\big(\widetilde{E}^Q\big[\int_0^{U_t} f^*_{\mu}(t,y)dy\big]+E^Q[\widetilde{E}^Q[\widetilde{f}^*_{\mu}(t)]L_t\big|\mathcal{F}_t^{Y}]\Gamma_t\big)(K_t^{1,\varepsilon}+K_t^{1,\varepsilon})dt\Big]\notag\\
&+E^Q\Big[\int^T_0\big\{E^Q\big[\gamma^{\Phi}_T(\delta h(t))+\int^T_t \gamma^f_s (\delta h(t))ds\,\big|\,\mathcal{F}_t^{Y}\big]+\delta f(t)\big\}\mathbf{1}_{E_\varepsilon}(t)dt\Big]\notag\\
&+o(\varepsilon),\ \ \mbox{ as } \varepsilon\searrow0.\notag
\end{align}
In addition to the duality (\ref{eq4.29star}), we have to calculate the duality described by $\displaystyle E^Q\big[p^1_T Y_T^{2,\varepsilon}
 +p^2_T K_T^{2,\varepsilon}\big]$.\\ Recall the second order variational equation \eqref{2ndVE5} under the Assumption (H3)-(i) (i.e., $\sigma(t,x,\gamma,u) = \sigma(t,\gamma, u)$):
\begin{equation*}
\left\{
\begin{split}
dY_t^{2,\varepsilon} =& \Big\{\widetilde{E}^Q\big[\int^{\widetilde{U}_t}_0 \sigma_{\mu} (t,y)dy\widetilde{K}^{2,\varepsilon}_t\big]+\widetilde{E}^Q\big[\widetilde{\sigma}_{\mu} (t) \widetilde{L}_t \widetilde{V}_t^{2,\varepsilon}\big]\\
&\quad+\widetilde{E}^Q\big[\widetilde{\sigma}_{\mu}(t)\widetilde{V}_t^{1,\varepsilon} \widetilde{K}_t^{1,\varepsilon}\big]+\frac{1}{2}\widetilde{E}^Q\big[\widetilde{\sigma}_{z\mu}(t)\widetilde{L}_t (\widetilde{V}_t^{1,\varepsilon})^2\big]\\
&\quad+\Big(\widetilde{E}^Q\Big[\int_0^{\widetilde{U}_t} \delta \sigma_{\mu}(t,y)dy\cdot\widetilde{K}_t^{1,\varepsilon}\Big]+\widetilde{E}^Q\big[\delta\widetilde{\sigma}_{\mu} (t)\widetilde{L}_t\widetilde{V}_t^{1,\varepsilon}\big]\Big)\mathbf{1}_{E_\varepsilon}(t)\Big\}dB_t^1;\\
dK_t^{2,\varepsilon}=&\Big\{h(t)K_t^{2,\varepsilon}+h_x(t)L_t Y_t^{2,\varepsilon}+h_x(t)Y_t^{1,\varepsilon}K_t^{1,\varepsilon}+\frac{1}{2}h_{xx}(t)L_t (Y_t^{1,\varepsilon})^2\\
&\ +L_t\widetilde{E}^Q\Big[\int_0^{\widetilde{U}_t}h_\mu(t,y)dy\cdot\widetilde{K}_t^{2,\varepsilon}\Big]+L_t\widetilde{E}^Q\big[\widetilde{h}_\mu(t)\widetilde{L}_t
\widetilde{V}_t^{2,\varepsilon}\big]\\
&\ +L_t\widetilde{E}^Q\big[\widetilde{h}_\mu(t)\widetilde{V}_t^{1,\varepsilon}\widetilde{K}_t^{1,\varepsilon}\big]+\frac{1}{2}L_t\widetilde{E}^Q\big[\widetilde{h}_{z\mu}(t)
\widetilde{L}_t(\widetilde{V}_t^{1,\varepsilon})^2\big]\\
&\ +\Big(\delta h(t)K_t^{1,\varepsilon}+\delta h_x(t)L_t Y_t^{1,\varepsilon}+L_t\widetilde{E}^Q\Big[\int_0^{\widetilde{U}_t}\delta h_\mu(t,y)dy\cdot\widetilde{K}_t^{1,\varepsilon}\Big]\\
&\qquad   +L_t\widetilde{E}^Q\big[\delta\widetilde{h}_\mu(t)\widetilde{L}_t\widetilde{V}_t^{1,\varepsilon}\big]\Big)\mathbf{1}_{E_\varepsilon}(t)\Big\}dY_t;\\
Y_0^{2,\varepsilon}=& K_0^{2,\varepsilon}=0,\ t\in[0,T].
\end{split}\right.
\end{equation*}
Also recall the first adjoint BSDE (\ref{SMP4}) reduced to our simplified setting
\begin{equation*}
\begin{split}
&dp_t^1=-\alpha_t(q_t^1,q_t^2)dt+q_t^1dB_t^1+\check{q}_t^1dY_t,\ t\in[0,T];\\
&dp_t^2=-\beta_t(q_t^1,q_t^2)dt+\check{q}_t^2dB_t^1+q_t^2dY_t,\ t\in[0,T],
\end{split}
\end{equation*}
where
\begin{equation*}
\begin{split}
\alpha_t(q^1_t,q^2_t) :=&L_t\widetilde{E}^Q\big[\widetilde{q}^1_t E^P[\widetilde{\sigma}^*_{\mu} (t) \,\big|\,\mathcal{F}_t^{Y}]\big]+h_x(t)L_t q^2_t +L_t\widetilde{E}^Q\big[\widetilde{q}^2_t \widetilde{L}_t E^P[\widetilde{h}^*_{\mu}(t)\,\big|\,\mathcal{F}_t^{Y}]\big]\\
&- f_x(t)-L_t\widetilde{E}^Q\big[E^P[\widetilde{f}^*_{\mu}(t)\,\big|\,\mathcal{F}_t^{Y}]\big],\ t\in[0,T];\\
\beta_t(q^1_t,q^2_t):=& (X_t-U_t)\widetilde{E}^Q\Big[\widetilde{q}_t^1 E^P[\widetilde{\sigma}^*_{\mu}(t)\,\big|\,\mathcal{F}_t^{Y}]\Big]+\widetilde{E}^Q\Big[\widetilde{q}_t^1 \int_0^{U_t} \sigma^*_{\mu}(t,y)dy\Big]+h(t)q^2_t\\
&+(X_t -U_t)\widetilde{E}^Q\Big[\widetilde{q}_t^2 \widetilde{L}_t E^P[\widetilde{h}^*_{\mu}(t)\,\big|\,\mathcal{F}_t^{Y}]\Big]+\widetilde{E}^Q\Big[\widetilde{q}_t^2 \widetilde{L}_t \int^{U_t}_0 h_{\mu}^* (t,y)dy\Big]\\
&-(X_t -U_t)\widetilde{E}^Q\Big[E^P[\widetilde{f}_{\mu}^*(t)\,\big|\,\mathcal{F}_t^{Y}]\Big]-\widetilde{E}^Q\Big[\int^{U_t}_0 f_{\mu}^* (t,y)dy\Big],\ t\in[0,T].
\end{split}
\end{equation*}
Then, using It\^{o}'s formula, we obtain
\begin{equation*}
\begin{split}
&dE^Q\big[p_t^1 Y_t^{2,\varepsilon}\big]=E^Q\Big[q^1_t\big\{\widetilde{E}^Q \Big[\int^{\widetilde{U}_t}_0 \sigma_{\mu} (t,y) dy\widetilde{K}_t^{2,\varepsilon}\Big] +
\widetilde{E}^Q\big[\widetilde{\sigma}_{\mu} (t) \widetilde{L}_t \widetilde{V}_t^{2,\varepsilon}\big] +\widetilde{E}^Q\big[\widetilde{\sigma}_{\mu}(t)\widetilde{V}_t^{1,\varepsilon}
\widetilde{K}_t^{1,\varepsilon}\big]\\
&\qquad+\frac{1}{2}\widetilde{E}^Q\big[\widetilde{\sigma}_{z\mu} (t)\widetilde{L}_t(\widetilde{V}_t^{1,\varepsilon})^2\big]+\big(\widetilde{E}^Q\big[\int_0^{\widetilde{U}_t}
\delta \sigma_{\mu}(t,y)dy\cdot\widetilde{K}_t^{1,\varepsilon}\big]+\widetilde{E}^Q\big[\delta\widetilde{\sigma}_{\mu} (t)\widetilde{L}_t\widetilde{V}_t^{1,\varepsilon}\big]\big)
\mathbf{1}_{E_\varepsilon}(t)\big\}\Big]dt\\
&\qquad-E^Q\Big[Y_t^{2,\varepsilon}\big(L_t\widetilde{E}^Q\big[\widetilde{q}_t^1 E^P\big[\widetilde{\sigma}^*_{\mu}(t)\,\big|\,\mathcal{F}_t^{Y}\big]\big]+h_x(t)L_t q_t^2 +
 L_t\widetilde{E}^Q\big[\widetilde{q}_t^2 \widetilde{L}_t E^P\big[\widetilde{h}_{\mu}^* (t)\,\big|\,\mathcal{F}_t^{Y}\big]\big]\\
&\qquad\qquad\quad-f_x(t)-L_t \widetilde{E}^Q\big[E^P\big[\widetilde{f}^*_{\mu} (t)\,\big|\,\mathcal{F}_t^{Y}\big]\big]\big)\Big]dt,
\end{split}
\end{equation*}
and
\begin{equation*}
    \begin{split}
&dE^Q\big[p_t^2 K_t^{2,\varepsilon}\big]= E^Q\Big[q^2_t\big\{h(t)K   _t^{2,\varepsilon} +h_x(t) L_t Y_t^{2,\varepsilon}
 + h_x(t) Y_t^{1,\varepsilon}K_t^{1,\varepsilon} +\frac{1}{2} h_{xx}(t) L_t \big(Y_t^{1,\varepsilon}\big)^2\\
 &\quad +L_t \widetilde{E}^Q\big[\int_0^{\widetilde{U}_t} h_{\mu} (t,y) dy \cdot \widetilde{K}_t^{2,\varepsilon}\big]
+L_t\widetilde{E}^Q\big[\widetilde{h}_{\mu} (t)\widetilde{V}_t^{1,\varepsilon}\widetilde{K}_t^{1,\varepsilon}\big]\\
&\quad +L_t\widetilde{E}^Q\big[\widetilde{h}_{\mu} (t)\widetilde{L}_t\widetilde{V}_t^{2,\varepsilon}\big]+\frac{1}{2}L_t
 \widetilde{E}^Q\big[\widetilde{h}_{z\mu} (t)\widetilde{L}_t\big(\widetilde{V}_t^{1,\varepsilon}\big)^2\big]\\
&\quad +\big(\delta h(t)K_t^{1,\varepsilon}\! + \!\delta h_x(t)L_tY_t^{1,\varepsilon}\! +\! L_t\widetilde{E}^Q\big[\!\int_0^{\widetilde{U}_t}
 \delta h_{\mu} (t,y)dy\cdot \widetilde{K}_t^{1,\varepsilon}\!\big]
 \!+\!L_t\widetilde{E}^Q\big[\delta \widetilde{h}_{\mu}(t)\widetilde{L}_t\widetilde{V}_t^{1,\varepsilon}\big]\big)\mathbf{1}_{E_\varepsilon}
 (t)\}\Big]dt\\
 &-E^Q\Big[K_t^{2,\varepsilon}\{(X_t-U_t)\widetilde{E}^Q\big[\widetilde{q}_t^1 E^P\big[\widetilde{\sigma}^*_{\mu} (t)\,\big|\,
 \mathcal{F}_t^{Y}\big]\big]+\widetilde{E}^Q\big[\widetilde{q}_t^1 \int_0^{U_t} \sigma_{\mu}^* (t,y)dy\big]\\
& \quad +h(t)q^2_t +(X_t-U_t)\widetilde{E}^Q\big[\widetilde{q}_t^2 \widetilde{L}_t E^P\big[\widetilde{h}_{\mu}^* (t)\,\big|\,\mathcal{F}_t^{Y}\big]
\big]+\widetilde{E}^Q\big[\widetilde{q}_t^2 \widetilde{L}_t \int_0^{U_t} h_{\mu}^* (t,y)dy\big]\\
&\quad  -(X_t - U_t)\widetilde{E}^Q\big[E^P\big[\widetilde{f}_{\mu}^* (t)\,\big|\,\mathcal{F}_t^{Y}\big]\big]-\widetilde{E}^Q
\big[\int_0^{U_t}f_{\mu}^* (t,y)dy\big]\big\}\Big]dt.
\end{split}
\end{equation*}
The following estimates will be needed to simplify the above relations.\\
{\rm{i)} First, due to Proposition \ref{EstofXL} and Corollary \ref{EstofUV}, we have\\
\begin{equation*}
\begin{split}
&\Big|E^Q\Big[\int_0^T q_t^1\big(\widetilde{E}^Q\big[\int^{\widetilde{U}_t}_0 \delta\sigma_{\mu} (t,y)dy\cdot
\widetilde{K}_t^{1,\varepsilon}\big]
+\widetilde{E}^Q\big[\delta \widetilde{\sigma}_{\mu} (t)\widetilde{L}_t\widetilde{V}_t^{1,\varepsilon}\big]\big)
\mathbf{1}_{E_\varepsilon}(t)dt\Big]\Big|\\
&\leq C\sqrt\varepsilon E^Q\Big[\int^T_0 \big|q^1_t\big|\mathbf{1}_{E_\varepsilon}(t)dt\Big]\leq C\varepsilon
\Big(E^Q\Big[\int^T_0 |q^1_t|^2\mathbf{1}_{E_\varepsilon}(t)dt\Big]\Big)^{\frac{1}{2}} = o(\varepsilon),\ \ \mbox{ as }
\varepsilon\rightarrow 0.
\end{split}
\end{equation*}
The same argument also shows that
\begin{equation*}
\begin{split}
\Big|E^Q\Big[&\int_0^Tq_t^2\Big(\delta h(t)K_t^{1,\varepsilon} +\delta h_x(t) L_t Y_t^{1,\varepsilon} +L_t \widetilde{E}^Q
\big[\int^{\widetilde{U}_t}_0 \delta h_{\mu} (t,y)dy\cdot\widetilde{K}_t^{1,\varepsilon}\big]\\
&+L_t\widetilde{E}^Q\big[\delta \widetilde{h}_{\mu} (t) \widetilde{L}_t\widetilde{V}_t^{1,\varepsilon}\big]\Big)\mathbf{1}_{E_\varepsilon}(t)dt\Big]\Big|
\leq C\varepsilon\Big(E^Q\big[\int_0^T |q_t^2|^2\mathbf{1}_{E_\varepsilon}(t)dt\big]\Big)^{\frac{1}{2}} = o(\varepsilon),\ \ \mbox{ as } \varepsilon\rightarrow 0.
\end{split}
\end{equation*}
{\rm{ii)} We now make the computations for $\displaystyle E^Q\big[q^1_t \widetilde{E}^Q[\widetilde{\sigma}_{\mu} (t) \widetilde{L}_t
\widetilde{V}_t^{2,\varepsilon}]\big]$. Recall that,
\begin{equation*}
    \begin{split}
V_t^{2,\varepsilon} = &\theta_t \big(Y_t^{2,\varepsilon},K_t^{2,\varepsilon}\big)+{E^P\big[H_\varepsilon (t)\,\big|\,
\mathcal{F}_t^{Y}\big]}\\
&-\int_0^t E^Q\big[\Gamma_t^1 \delta h(s) \,\big|\,\mathcal{F}_t^{Y}\big] E^Q\big[\Gamma_t L_t \delta h(s)\,
\big|\,\mathcal{F}_t^{Y}\big] \mathbf{1}_{E_\varepsilon}(s)ds+\varepsilon\eta_t(\varepsilon),
\end{split}
\end{equation*}
where\ \ $\displaystyle \theta_t\big(Y_t^{2,\varepsilon},K_t^{2,\varepsilon}\big)=E^Q\big[\Gamma_t^1 Y_t^{2,\varepsilon} \,\big|\,\mathcal{F}_t^{Y}\big]
+E^Q\big[\Gamma_t K_t^{2,\varepsilon} \,\big|\,\mathcal{F}_t^{Y}\big],$\ with
$$\displaystyle\Gamma^1_t =\frac{L_t}{E^Q\big[L_t\,|\,\mathcal{F}_t^Y\big]},\ \ \Gamma_t =\frac{1}{E^Q\big[L_t\,|\,\mathcal{F}_t^Y\big]}
 \big(X_t -U_t\big),$$
and, $\eta_t(\varepsilon)$ stands for a stochastic process such that
$$E^Q\big[\sup_{t\in[0,T]}|\eta_t(\varepsilon)|^p\big]\rightarrow 0\ (\varepsilon\searrow 0),\ \ \big(E^Q\big[\big|\eta_t(\varepsilon)\big|^p\big]
\big)^{\frac{1}{p}}\leq C_p,\ t\in[0,T],\ \varepsilon>0,\ p\geq 2.$$
Then,
\begin{equation*}
\begin{split}
&E^Q\Big[q^1_t\widetilde{E}^Q\big[\widetilde{\sigma}_{\mu} (t)\widetilde{L}_t\widetilde{V}_t^{2,\varepsilon}\big]\Big]=E^Q\Big[\widetilde{E}^Q\big[\widetilde{q}_t^1 \widetilde{\sigma}_{\mu}^* (t)\big] L_t V_t^{2,\varepsilon}\Big]\\
=&\ E^Q\Big[\widetilde{E}^Q\big[\widetilde{q}_t^1 \widetilde{\sigma}_{\mu}^* (t)\big] L_t E^Q\big[\Gamma_t^1 Y_t^{2,\varepsilon} +\Gamma_t K_t^{2,\varepsilon}\,\big|\,\mathcal{F}_t^{Y}\big]\Big]\\
&\ - E^Q\big[\widetilde{E}^Q\big[\widetilde{q}_t^1 \widetilde{\sigma}_{\mu}^* (t)\big] L_t \int_0^t E^Q\big[\Gamma_t^1\delta h(s)\,\big|\,\mathcal{F}_t^{Y}\big]E^Q\big[\Gamma_t L_t \delta h(s)\,\big|\,\mathcal{F}_t^{Y}]\mathbf{1}_{E_\varepsilon}(s)ds\big]\\
&\ +{E^Q\Big[\widetilde{E}^Q\big[\widetilde{q}_t^1 \widetilde{\sigma}_{\mu}^* (t)\big] L_tE^P\big[H_\varepsilon(t)\,\big|\,\mathcal{F}_t^{Y}\big]\Big]} + \varepsilon \rho_t(\varepsilon)\\
=&\ I_t^{1,\varepsilon} - I_t^{2,\varepsilon}+{E^Q\Big[\widetilde{E}^Q\big[\widetilde{q}_t^1 \widetilde{\sigma}_{\mu}^* (t)\big] L_tE^P\big[H_\varepsilon(t)\,\big|\,\mathcal{F}_t^{Y}\big]\Big]} + \varepsilon \rho_t (\varepsilon),\ \
\end{split}
\end{equation*}
\mbox{where}
\begin{equation*}
\begin{split}
I_t^{1,\varepsilon}:&=E^Q\Big[E^Q\big[L_t\widetilde{E}^Q\big[\widetilde{q}_t^1\widetilde{\sigma}^*_{\mu} (t)\big]\,\big|\,\mathcal{F}_t^{Y}\big](\Gamma_t^1 Y_t^{2,\varepsilon} + \Gamma_t K_t^{2,\varepsilon})\Big]\\
&=E^Q\Big[\widetilde{E}^Q\big[E^Q\big[\frac{L_t}{E^Q[L_t\,|\,\mathcal{F}_t^Y]}\widetilde{\sigma}_{\mu}^* (t)\,\big|\,\mathcal{F}_t^{Y}\big]\widetilde{q}_t^1\big](L_t Y_t^{2,\varepsilon} + (X_t - U_t)K_t^{2,\varepsilon})\Big]\\
&=E^Q\Big[\widetilde{E}^Q\big[\widetilde{q}_t^1 E^P\big[\widetilde{\sigma}_{\mu}^* (t)\,\big|\,\mathcal{F}_t^{Y}\big]\big](L_t Y_t^{2,\varepsilon})\Big] +E^Q\Big[\widetilde{E}^Q\big[\widetilde{q}_t^1  E^P\big[\widetilde{\sigma}_{\mu}^* (t)\,\big|\,\mathcal{F}_t^{Y}\big]\big](X_t - U_t)K_t^{2,\varepsilon}\Big].
\end{split}
\end{equation*}
Thus,
\begin{equation*}
\begin{split}
E^Q&\Big[q^1_t\widetilde{E}^Q\big[\widetilde{\sigma}_{\mu}(t) \widetilde{L}_t\widetilde{V}_t^{2,\varepsilon}\big]\Big]
=E^Q\Big[\widetilde{E}^Q\big[\widetilde{q}_t^1 E^P \big[\widetilde{\sigma}_{\mu}^* (t)\,\big|\,\mathcal{F}_t^{Y}\big]\big](L_t
Y_t^{2,\varepsilon})\Big]\\
&+E^Q\Big[\widetilde{E}^Q\big[\widetilde{q}_t^1 E^P \big[\widetilde{\sigma}_{\mu}^* (t)\,\big|\,\mathcal{F}_t^{Y}\big]\big](X_t - U_t)
K_t^{2,\varepsilon}\Big]
+{E^Q\Big[\widetilde{E}^Q\big[\widetilde{q}_t^1 \widetilde{\sigma}_{\mu}^* (t)\big] L_tE^P\big[H_\varepsilon(t)\,\big|\,
\mathcal{F}_t^{Y}\big]\Big]}\\
& -E^Q\Big[\widetilde{E}^Q\big[\widetilde{q}_t^1 \widetilde{\sigma}_{\mu}^* (t)\big]L_t \int^t_0 E^P \big[\delta h(s) \,\big|\,
\mathcal{F}_t^{Y}\big] E^P\big[(X_t - U_t)\delta h(s)\,\big|\,\mathcal{F}_t^{Y}\big]\mathbf{1}_{E_\varepsilon}(s)ds\Big]\\
& +\varepsilon\rho_t(\varepsilon),\ \ \ \ \rho_t(\varepsilon)\rightarrow 0\ (\varepsilon\searrow 0).
\end{split}
\end{equation*}
{\rm{iii)} Moreover, in analogy to (\ref{P10,ii}) we deduce that
\begin{equation}\label{P16,iii}
\begin{split}
&E^Q\Big[{q}_t^1\widetilde{E}^Q[\widetilde{\sigma}_{\mu}(t)\widetilde{V}_t^{1,\varepsilon}\widetilde{K}_t^{1,\varepsilon}]\Big] = E^Q\Big[\widetilde{E}^Q\big[\widetilde{q}_t^1\widetilde{\sigma}_{\mu}^*(t)\big]V_t^{1,\varepsilon} K_t^{1,\varepsilon}\Big]\\
=&\ E^Q\Big[\int_0^t E^Q\big[\Gamma_t L_t\delta h(s) \,\big|\,\mathcal{F}_t^{Y}\big]E^Q\big[\widetilde{E}^Q\big[\widetilde{q}^1_t \widetilde{\sigma}_{\mu}^*(t)\big]L_t \delta h(s) \big|\mathcal{F}_t^{Y}\big]\mathbf{1}_{E_\varepsilon}(s)ds\Big]+\varepsilon\rho_t(\varepsilon)\\
=&\ E^Q\Big[\int_0^t E^P\big[(X_t - U_t)\delta h(s)\,\big|\,\mathcal{F}_t^{Y}\big]E^Q\big[\widetilde{E}^Q\big[\widetilde{q}^1_t \widetilde{\sigma}_{\mu}^*(t)\big]L_t \delta h(s) \,\big|\,\mathcal{F}_t^{Y}\big]\mathbf{1}_{E_\varepsilon}(s)ds\Big]+\varepsilon\rho_t(\varepsilon),
\end{split}
\end{equation}
{\rm{iv)} Similarly to (\ref{P10,iii}) we have
\begin{equation}\label{P16,iv}
\begin{split}
&E^Q\Big[{q}_t^1\widetilde{E}^Q[\widetilde{\sigma}_{z\mu}(t)\widetilde{L}_t(\widetilde{V}_t^{1,\varepsilon})^2]\Big]\\
=& E^Q\Big[\widetilde{E}^Q\big[\widetilde{q}_t^1\widetilde{\sigma}_{z\mu}^*(t)\big]
L_t(V_t^{1,\varepsilon})^2\Big]\\
=&\ E^Q\Big[\widetilde{E}^Q\big[\widetilde{q}_t^1\widetilde{\sigma}_{z\mu}^*(t)\big]L_t\int_0^t\big(E^Q\big[\Gamma_t L_t \delta h(s)\,\big|\,\mathcal{F}_t^{Y}\big]\big)^2\mathbf{1}_{E_\varepsilon}(s)ds\Big] +\varepsilon\rho_t(\varepsilon)\\
=& \ E^Q\Big[\widetilde{E}^Q\big[\widetilde{q}_t^1\widetilde{\sigma}_{z\mu}^*(t)\big]L_t\int_0^t\big(E^P\big[(X_t -U_t) \delta h(s)\,\big|\,\mathcal{F}_t^{Y}\big]\big)^2\mathbf{1}_{E_\varepsilon}(s)ds\Big] +\varepsilon\rho_t(\varepsilon).
\end{split}
\end{equation}
{\rm{v)} In analogy to {\rm{ii)} we have
\begin{equation}\label{P16,v}
\begin{split}
E^Q\Big[q^2_t L_t &\widetilde{E}^Q\big[\widetilde{h}_{\mu} (t)\widetilde{L}_t \widetilde{V}_t^{2,\varepsilon}\big]\Big]
=E^Q\Big[L_tY_t^{2,\varepsilon} \widetilde{E}^Q\big[\widetilde{q}^2_t \widetilde{L}_t E^P \big[\widetilde{h}_{\mu}^* (t)\,\big|\,\mathcal{F}_t^{Y}\big]\big]\Big]\\
& +E^Q\Big[(X_t-U_t)K_t^{2,\varepsilon} \widetilde{E}^Q\big[\widetilde{q}_t^2 \widetilde{L}_t E^P\big[\widetilde{h}_{\mu}^*(t)\,\big|\,\mathcal{F}_t^{Y}\big]\big]\Big]\\
&+{E^Q\Big[\widetilde{E}^Q\big[\widetilde{q}^2_t \widetilde{L}_t\widetilde{h}_{\mu}^* (t)\big]L_tE^P\big[H_\varepsilon(t)\,\big|\,\mathcal{F}_t^{Y}\big]\Big]}\\
& -E^Q\Big[\widetilde{E}^Q\big[\widetilde{q}_t^2 \widetilde{L}_t \widetilde{h}_{\mu}^* (t)\big]L_t \int_0^t E^P \big[\delta h(s)\,\big|\,\mathcal{F}_t^{Y}\big] E^P\big[(X_t- U_t)\delta
h(s)\,\big|\,\mathcal{F}_t^{Y}\big]\mathbf{1}_{E_\varepsilon}(s)ds\Big]\\
& +\varepsilon\rho_t (\varepsilon).
\end{split}
\end{equation}
{\rm{vi)} From the arguments for \rm{iii)} we see
\begin{equation*}
\begin{split}
&E^Q\Big[q^2_t L_t \widetilde{E}^Q\big[\widetilde{h}_{\mu}(t)\widetilde{V}_t^{1,\varepsilon} \widetilde{K}_t^{1,\varepsilon}\big]\Big]\\
=&\ E^Q\Big[\int_0^t E^P\big[(X_t -U_t)\delta h(s)\,\big|\,\mathcal{F}_t^{Y}\big]E^Q\big[\widetilde{E}^Q\big[\widetilde{q}^2_t\widetilde{L}_t\widetilde{h}^*_{\mu} (t)\big]L_t \delta h(s)
 \big|\mathcal{F}_t^{Y}\big]\mathbf{1}_{E_\varepsilon}(s)ds\Big]+\varepsilon\rho_t (\varepsilon).
\end{split}
\end{equation*}
{\rm{vii)} Following \rm{iv)} we obtain
\begin{equation*}
\begin{split}
&E^Q\Big[q^2_t L_t \widetilde{E}^Q\big[\widetilde{h}_{z\mu}(t) \widetilde{L}_t (\widetilde{V}_t^{1,\varepsilon})^2\big]\Big]\\
=&\ E^Q\Big[\widetilde{E}^Q\big[\widetilde{q}_t^2 \widetilde{L}_t \widetilde{h}_{z\mu}^* (t) \big]L_t \int_0^t \big(E^P \big[ (X_t - U_t) \delta h(s) \big|\mathcal{F}_t^{Y}\big]\big)^2
\mathbf{1}_{E_\varepsilon}(s)ds\Big]+\varepsilon\rho_t (\varepsilon).
\end{split}
\end{equation*}
{\rm{viii)} Similarly, we have
\begin{equation*}
    \begin{split}
        &E^Q\Big[q^2_t h_x(t) Y^{1,\varepsilon}_t K^{1,\varepsilon}_t\Big]
        =\ E^Q\Big[\int^t_0 E^Q\big[q_t^2h_x(t)\,\big|\,\mathcal{F}_t^{Y}\big]L_t h(s,t)|\delta\sigma(s)|^2\mathbf{1}_{E_\varepsilon}(s)ds\Big] + \varepsilon \rho_t (\varepsilon)\\
        =&E^Q\Big[E^Q\big[q_t^2h_x(t)\,\big|\,\mathcal{F}_t^{Y}\big]L_t E^P\big[H_\varepsilon (t)\,\big|\,\mathcal{F}_t^Y\big]\Big]+ \varepsilon \rho_t (\varepsilon).
    \end{split}
\end{equation*}
Indeed, if we define $\zeta_t:=q_t^2h_x(t),\ t\in [0,T]$, then we have $\displaystyle \zeta_t = E^Q [\zeta_t\,|\,\mathcal{F}_t^Y]+\int_0^t \zeta_{s,t} dB_s^1,$\ where $(\zeta_{s,t})$ is a jointly measurable two-parameter process such that,
\begin{equation*}
    \begin{array}{lll}
\quad\quad\quad(\mbox{i})\  \zeta_{s, t}\ \mbox{is}\ \mathcal{F}^{B^1}_s \vee  \mathcal{F}_t^Y\mbox{-adapted},\ 0\leq s\leq t\leq T;\\
\quad\quad\quad(\mbox{ii})\ E^Q\Big[\int_0^T \int_0^t \big|\zeta_{s,t}\big|^2dsdt \Big] \leq E^Q \Big[\int_0^T \big|\zeta_t\big|^2 dt \Big]<+\infty.\\
    \end{array}
\end{equation*}
Then, we have that
\begin{equation*}
    \begin{split}
        &E^Q\Big[q^2_t h_x(t) Y^{1,\varepsilon}_t K^{1,\varepsilon}_t\Big]= E^Q\Big[\zeta_t Y^{1,\varepsilon}_t K^{1,\varepsilon}_t\Big]\\
        =& E^Q\Big[ E^Q\big[\zeta_t \,\big|\,\mathcal{F}_t^Y\big]Y^{1,\varepsilon}_t K^{1,\varepsilon}_t\Big] + E^Q\Big[ \big( \int_0^t \zeta_{s,t}
         dB_s^1\big)Y^{1,\varepsilon}_t K^{1,\varepsilon}_t\Big]=:I_t^{1,\varepsilon}+ I_t^{2,\varepsilon}.
    \end{split}
\end{equation*}
\noindent$\bullet$ For $I_t^{2,\varepsilon}$ we observe that
\begin{equation}\label{p55-last}
    I_t^{2,\varepsilon}=E^Q\Big[ \big( \int_0^t \zeta_{s,t}dB_s^1\big)Y^{1,\varepsilon,1}_t K^{1,\varepsilon}_t\Big]+
E^Q\Big[ \big( \int_0^t \zeta_{s,t}dB_s^1\big)\big(Y^{1,\varepsilon}_t-Y^{1,\varepsilon,1}_t \big)K^{1,\varepsilon}_t\Big].
\end{equation}
It follows from from Proposition \ref{EstofXL} and \eqref{3starstar} that
\begin{equation*}
    \begin{split}
      &\Big|E^Q\Big[ \big( \int_0^t \zeta_{s,t}dB_s^1\big)\big(Y^{1,\varepsilon}_t-Y^{1,\varepsilon,1}_t \big)K^{1,\varepsilon}_t\Big] \Big|\\
      \leq& \Big(E^Q\big[ \int_0^t |\zeta_{s,t}|^2 ds\big] \Big)^{\frac{1}{2}} \Big(E^Q\big[ \big|Y^{1,\varepsilon}_t-Y^{1,\varepsilon,1}_t\big|^4\big] \Big)^{\frac{1}{4}}
      \Big(E^Q\big[ \big|K^{1,\varepsilon}_t\big|^4\big] \Big)^{\frac{1}{4}}\\
     \leq & C\, \Big(E^Q\big[ |\zeta_{t}|^2 \big] \Big)^{\frac{1}{2}} \rho_4 (\varepsilon)\sqrt{\varepsilon} \cdot \sqrt{\varepsilon} =  \rho_t (\varepsilon)\,{\varepsilon},
    \end{split}
\end{equation*}
where $\rho_t (\varepsilon) = C\,\Big(E^Q\big[ |\zeta_{t}|^2 \big] \Big)^{\frac{1}{2}}\rho_4 (\varepsilon)$, and
$$ \int_0^T |\rho_t(\varepsilon)|^2dt \leq |\rho_4 (\varepsilon)|^2\!\cdot\!E^Q \big[\int_0^T |\zeta_t|^2 dt \big] \leq C\,
|\rho_4 (\varepsilon)|^2\!\cdot\!E^Q \big[\int_0^T |q_t^2|^2 dt \big] \rightarrow 0\ (\varepsilon \searrow 0).$$
On the other hand, notice that
\begin{equation*}
    \begin{split}
     & \big( \int_0^t \zeta_{s,t} dB_s^1 \big)Y^{1,\varepsilon,1}_t = \int_0^t \Big( \zeta_{s,t} Y^{1,\varepsilon,1}_t+\big(\int_0^s \zeta_{r,t}dB_r^1\big)\delta\sigma(s)
      \mathbf{1}_{E_\varepsilon}(s) \Big)dB_s^1 + \int_0^t \zeta_{s,t}\delta\sigma(s) \mathbf{1}_{E_\varepsilon}(s)ds\\
      =&\int_0^t \eta^\varepsilon _{s,t} dB_s^1+  \int_0^t \zeta_{s,t}\delta\sigma(s) \mathbf{1}_{E_\varepsilon}(s)ds,
    \end{split}
\end{equation*}
where $\displaystyle \eta^\varepsilon _{s,t}:=\zeta_{s,t} Y^{1,\varepsilon,1}_t+\big(\int_0^s \zeta_{r,t}dB_r^1\big)\delta\sigma(s)
      \mathbf{1}_{E_\varepsilon}(s)$, $0\leq s\leq t\leq T$, and
$\displaystyle E^Q \Big[\big(\int_0^t \eta^\varepsilon _{s,t} dB_s^1\big)K_t^{1,\varepsilon}\Big]=
E^Q \Big[\big\langle \int_0^\cdot \eta^\varepsilon _{s,t} dB_s^1, K_\cdot^{1,\varepsilon}\big\rangle_t\Big]=0.$\\
Hence,
\begin{equation*}
    \begin{split}
         &\Big|E^Q\Big[ \big( \int_0^t \zeta_{s,t}dB_s^1\big)Y^{1,\varepsilon,1}_t K^{1,\varepsilon}_t\Big]\Big| = \Big|E^Q\Big[ \big( \int_0^t \zeta_{s,t}\delta\sigma(s)
         \mathbf{1}_{E_\varepsilon}(s)ds\big) K^{1,\varepsilon}_t\Big]\Big|\\
         \leq &C\, \Big(E^Q\Big[ \int_0^t\big| \zeta_{s,t}\big| ^2\mathbf{1}_{E_\varepsilon}(s)ds\Big]\Big)^{\frac{1}{2}}\sqrt \varepsilon \cdot
         \Big(E^Q \Big[ \big| K_t^{1,\varepsilon}\big|^2\Big] \Big)^{\frac{1}{2}}\leq C\, \rho_t^1(\varepsilon)\cdot \varepsilon,
    \end{split}
\end{equation*}
where $\displaystyle \rho_t^1(\varepsilon) := \Big(E^Q\Big[ \int_0^t\big| \zeta_{s,t}\big| ^2\mathbf{1}_{E_\varepsilon}(s)ds\Big]\Big)^{\frac{1}{2}}\rightarrow 0\ (\varepsilon\searrow 0),$
and
$$\big|\rho_t^1(\varepsilon)\big|^2 \leq E^Q\Big[ \int_0^t\big| \zeta_{s,t}\big|^2 ds\Big]\leq
E^Q\big[ \big| \zeta_{t}\big|^2\big] \leq E^Q\big[ \big| q_{t}^2\big|^2\big],\ \varepsilon>0.$$
Thus, we get $$|I_t^{2,\varepsilon}|\leq \rho_t(\varepsilon)\varepsilon,\ \mbox{where}\ \displaystyle \rho_t(\varepsilon)\rightarrow 0\ (\varepsilon\searrow 0),\ \mbox{and} \int_0^T\! |\rho_t(\varepsilon)|^2 dt \leq CE^Q\big[\!\int_0^T\! |q_t^2|^2dt \big],\ \varepsilon>0.$$
\noindent$\bullet$ For $I_t^{1,\varepsilon}$, from \eqref{8star} it yields that
\begin{align*}
    I_t^{1,\varepsilon} & = E^Q\Big[E^Q\big[ \zeta_t \,\big|\,\mathcal{F}_t^Y\big]Y_t^{1,\varepsilon}K_t^{1,\varepsilon} \Big] =
    E^Q\Big[\zeta_t E^Q\big[ Y_t^{1,\varepsilon}K_t^{1,\varepsilon} \,\big|\,\mathcal{F}_t^Y\big] \Big]                          \\
                        & = E^Q\Big[\zeta_t E^Q\big[L_t\int_0^t h(s,t)|\delta\sigma(s)|^2\mathbf{1}_{E_\varepsilon}(s)ds\,
    \big|\,\mathcal{F}_t^{Y}\big]\Big]+\varepsilon E^Q\big[\zeta_t \eta_t(\varepsilon)\big].
\end{align*}
Finally, from above results we obtain
\begin{equation*}
    \begin{split}
        E^Q\Big[q^2_t h_x(t) Y^{1,\varepsilon}_t K^{1,\varepsilon}_t\Big]=E^Q\Big[\int_0^t E^Q\big[q_t^2h_x(t)\,\big|\,\mathcal{F}_t^{Y}\big]L_t
        h(s,t)|\delta \sigma(s)|^2\mathbf{1}_{E_\varepsilon}(s)ds\Big]+ \varepsilon \rho_t (\varepsilon),
    \end{split}
\end{equation*}
where $\rho_t(\varepsilon)\rightarrow 0\ (\varepsilon\searrow 0),\ \int_0^T |\rho_t(\varepsilon)|^2 dt \leq C, \ \varepsilon>0.$\\

\noindent Consequently, as
\begin{equation*}
\begin{split}
& dE^Q\Big[p_t^1 Y_t^{2,\varepsilon} + p^2_t K_t^{2,\varepsilon}\Big]\\
=&\ \Big\{\Big( E^Q \Big[q^1_t\widetilde{E}^Q\big[\widetilde{\sigma}_{\mu}(t) \widetilde{L}_t\widetilde{V}_t^{2,\varepsilon}\big]\Big]
- E^Q\Big[Y_t^{2,\varepsilon} L_t \widetilde{E}^Q \big[ \widetilde{q_t^1} E^P \big[\widetilde{\sigma}_{\mu}^* (t) \,\big|\,\mathcal{F}_t^{Y}\big]\big]\Big]\\
&\qquad\qquad -E^Q\Big[(X_t -U_t) K_t^{2,\varepsilon} \widetilde{E}^Q\big[\widetilde{q}_t^1 E^P \big[\widetilde{\sigma}_{\mu}^* (t)
 \,\big|\,\mathcal{F}_t^{Y}\big]\big]\Big]\Big)\\
&\qquad +E^Q\Big[q_t^1 \widetilde{E}^Q\big[\widetilde{\sigma}_{\mu} (t) \widetilde{V}_t^{1,\varepsilon} \widetilde{K}_t^ {1,\varepsilon}
\big]\Big]+\frac{1}{2} E^Q \Big[q_t^1 \widetilde{E}^Q\big[\widetilde{\sigma}_{z\mu} (t) \widetilde{L}_t (\widetilde{V}_t^{1,\varepsilon}
)^2\big]\Big]\\
&\qquad +\Big(E^Q \Big[q_t^2 L_t \widetilde{E}^Q\big[\widetilde{h}_{\mu} (t) \widetilde{L}_t \widetilde{V}_t^{2,\varepsilon}\big]\Big]-E^Q\Big[L_t Y_t^{2,\varepsilon} \widetilde{E}^Q\big[
    \widetilde{q}_t^2 \widetilde{L}_t E^P\big[\widetilde{h}_{\mu}^* (t)\,\big|\,\mathcal{F}_t^{Y}\big]\big]\Big]\\
&\qquad\qquad -E^Q\big[(X_t -U_t) K_t^{2,\varepsilon}\widetilde{E}^Q\big[\widetilde{q}_t^2 \widetilde{L}_tE^P\big[\widetilde{h}^*_{\mu} (t) \,\big|\,\mathcal{F}_t^{Y}\big]\big]\big]\Big)\\
&\qquad +E^Q\Big[q_t^2 L_t \widetilde{E}^Q \big[\widetilde{h}_{\mu}(t) \widetilde{V}_t^{1,\varepsilon} \widetilde{K}_t^{1,\varepsilon}\big]\Big] + \frac{1}{2}E^Q \Big[q_t^2 L_t
\widetilde{E}^Q\big[\widetilde{h}_{z\mu} (t)\widetilde{L}_t (\widetilde{V}_t^{1,\varepsilon})^2\big]\Big]\\
&\qquad +E^Q\Big[Y_t^{2,\varepsilon}\big(f_x(t) + L_t \widetilde{E}^Q\big[E^P\big[\widetilde{f}_{\mu}^*(t)\,\big|\,\mathcal{F}_t^{Y}\big]\big]\big)\Big] +\frac{1}{2}E^Q\Big[q_t^2 h_{xx}
(t)L_t\big(Y_t^{1,\varepsilon}\big)^2\Big]\\
&\qquad + E^Q\Big[K_t^{2,\varepsilon}\big((X_t-U_t)\widetilde{E}^Q\big[E^P\big[\widetilde{f}_{\mu}^*(t)\,\big|\,\mathcal{F}_t^{Y}\big]\big]+\widetilde{E}^Q\big[\int^{U_t}_0 f_{\mu} ^* (t,y)
 dy\big]\big)\Big]\\
&\qquad + E^Q\Big[q_t^2 h_x (t) Y_t^{1,\varepsilon} K_t^{1,\varepsilon}\Big] + \varepsilon \rho_t(\varepsilon)\Big\}dt,
\end{split}
\end{equation*}

\noindent combining our previous estimates \rm{i)}--\rm{viii)} yields
\begin{align*}
&dE^Q\Big[p_t^1 Y_t^{2,\varepsilon} + p^2_t K_t^{2,\varepsilon}\Big]\\
=&\ \Big\{E^Q\Big[Y_t^{2,\varepsilon}\big(f_x(t) +L_t \widetilde{E}^Q\big[E^P\big[\widetilde{f}^*_{\mu} (t) \,\big|\,\mathcal{F}_t^{Y}\big]
\big]\big)\Big]+E^Q\Big[E^Q\big[q_t^2h_x(t)\,\big|\,\mathcal{F}_t^{Y}\big]L_t E^P\big[H_\varepsilon (t)\,\big|\,\mathcal{F}_t^Y\big]\Big]\\
&\quad +E^Q\Big[K_t^{2,\varepsilon}\big((X_t-U_t)\widetilde{E}^Q\big[E^P\big[\widetilde{f}^*_{\mu} (t) \,\big|\,\mathcal{F}_t^{Y}\big]\big]+\widetilde{E}^Q\big[\int_0^{U_t} f_{\mu}^* (t,y) dy\big]\big)\Big]\\
&\quad +{E^Q\Big[\widetilde{E}^Q\big[\widetilde{q}_t^1 \widetilde{\sigma}^*_\mu (t) + \widetilde{q}_t^2\widetilde{L}_t \widetilde{h}^*_\mu (t)\big] L_t E^P\big[H_\varepsilon (t)
\,\big|\,\mathcal{F}_t^{Y}\big]\Big]}+\frac{1}{2} E^Q\big[q_t^2h_{xx}(t) L_t \big(Y_t^{1,\varepsilon}\big)^2\big]\\
&\quad +E^Q\Big[\int^t_0\big\{E^P\big[(X_t-U_t)\delta h(s)\,\big|\,\mathcal{F}_t^{Y}\big]\Big(E^Q\big[\!\widetilde{E}^Q\big[\widetilde{q}_t^1\widetilde{\sigma}_{\mu}^*(t) + \widetilde{q}_t^2 \widetilde{L}_t
\widetilde{h}_{\mu}^*(t)\big]L_t \delta h(s)\,\big|\,\mathcal{F}_t^{Y}\big]\!\\
&\quad\quad\quad -\!\widetilde{E}^Q\big[\widetilde{q}^1_t \widetilde{\sigma}^*_{\mu}
(t) +\widetilde{q}_t^2 \widetilde{L}_t \widetilde{h}_{\mu}^*(t)\big]L_t E^P\big[\delta h(s)\,\big|\,\mathcal{F}_t^{Y}\big]\Big)\\
&\quad\quad +\frac{1}{2}\widetilde{E}^Q\big[\widetilde{q}_t^1\widetilde{\sigma}_{z\mu}^*(t) + \widetilde{q}_t^2\widetilde{L}_t \widetilde{h}^*_{z\mu} (t)\big] L_t \big(E^P \big[(X_t-U_t)\delta h(s)\,\big|\,\mathcal{F}_t^{Y}\big]\big)^2\big\}\mathbf{1}_{E_\varepsilon}(s)ds\Big]\\
&\quad +\varepsilon \rho_t(\varepsilon)\Big\}dt.
\end{align*}
This shows that
\begin{equation*}
\begin{split}
&E^Q\Big[p_T^1 Y_T^{2,\varepsilon} + p^2_T K_T^{2,\varepsilon}\Big]\\
=&\ E^Q\Big[\int^T_0 Y_t^{2,\varepsilon} \big(f_x(t) + L_t \widetilde{E}^Q\big[E^P[\widetilde{f}_{\mu}^*(t) \,\big|\,\mathcal{F}_t^{Y}]\,\big]\big)dt\Big]
+E^Q\Big[\frac{1}{2}\int^T_0 h_{xx} (t) q_t^2 (Y_t^{1,\varepsilon})^2 L_tdt\Big]\\
\end{split}
\end{equation*}
\begin{equation*}
\begin{split}
&\quad +E^Q\Big[\int^T_0 K_t^{2,\varepsilon} \Big\{(X_t-U_t) \widetilde{E}^Q\big[E^P[\widetilde{f}_{\mu}^*(t)\,\big|\,\mathcal{F}_t^{Y}]\,\big]
+\widetilde{E}^Q\big[\int_0^{U_t} f_{\mu}^* (t,y)dy\big]\Big\}dt\Big]\\
&\quad +{E^Q\Big[\int_0^T\big( \widetilde{E}^Q\big[\widetilde{q}_t^1 \widetilde{\sigma}^*_\mu (t) + \widetilde{q}_t^2\widetilde{L}_t \widetilde{h}^*_\mu (t)\big]
+E^Q\big[q_t^2h_x(t)\,\big|\,\mathcal{F}_t^{Y}\big]\big) L_t E^P\big[H_\varepsilon (t)\,\big|\,\mathcal{F}_t^{Y}\big]dt\Big]}\\
&\quad +E^Q\Big[\int_0^T E^Q \Big[\int_t^T\Big(E^P\big[(X_s -U_s)\delta h(t) \big|\mathcal{F}_s^{Y}\big]\Big\{ E^Q\Big[\widetilde{E}^Q\big[\widetilde{q}_s^1 \widetilde{\sigma}_{\mu}^* (s)
+\widetilde{q}_s^2\widetilde{L}_s \widetilde{h}_{\mu}^*(s) \big]L_s \delta h(t)\big|\mathcal{F}_s^{Y}\Big]\\
&\quad\quad \ - \widetilde{E}^Q\big[\widetilde{q}_s^1 \widetilde{\sigma}_{\mu}^* (s) +\widetilde{q}_s^2 \widetilde{L}_s\widetilde{h}_{\mu}^* (s)\big]L_s E^P\big[\delta h(t)\big|
\mathcal{F}_s^{Y}\big]\Big\}\\
&\quad\ +\frac{1}{2} \widetilde{E}^Q\big[\widetilde{q}_s^1 \widetilde{\sigma}^*_{z\mu}(s) +\widetilde{q}_s^2\widetilde{L}_s \widetilde{h}^*_{z\mu} (s)\big]L_s \big(E^P
\big[(X_s -U_s) \delta h(t)\big|\mathcal{F}_s^{Y}\big]\big)^2\Big)ds\,\big|\,\mathcal{F}_t^{Y}\Big]\mathbf{1}_{E_\varepsilon}(t)dt\Big]\\
&\quad + o(\varepsilon)\ \ \ (\varepsilon\searrow 0).
\end{split}
\end{equation*}
Recalling \eqref{eq4.29star} for $\displaystyle E^Q\big[p_T^1 Y_T^{1,\varepsilon} + p_T^2 K_T^{1,\varepsilon} \big]$, it follows that
\begin{align}
&E^Q\Big[p_T^1 \big(Y_T^{1,\varepsilon}+Y_T^{2,\varepsilon}\big) + p_T^2 \big(K_T^{1,\varepsilon}+K_T^{2,\varepsilon}\big)\Big]\notag \\
=&\ E^Q\Big[\int^T_0\Big\{\big(Y_t^{1,\varepsilon}+Y_t^{2,\varepsilon}\big)\Big(f_x (t) +L_t \widetilde{E}^Q \big[E^P\big[\widetilde{f}_{\mu}^*(t) \,\big|\,\mathcal{F}_t^{Y}\big]\big]\Big)\notag\\
\notag&\ \quad  +\big(K_t^{1,\varepsilon} +K_t^{2,\varepsilon}\big)\Big( (X_t-U_t)\widetilde{E}^Q\big[E^P\big[\widetilde{f}_{\mu}^*(t) \,\big|
\,\mathcal{F}_t^{Y}\big]\big]+\widetilde{E}^Q\big[\int_0^{U_t} f_{\mu}^* (t,y)dy\big]\Big)\Big\}dt\Big]\\
&\  +E^Q\Big[\int^T_0 \frac{1}{2} h_{xx} (t) q_t^2 L_t ( Y_t^{1,\varepsilon})^2 dt\Big]\notag\\
&\ +{E^Q\Big[\int_0^T\big(\widetilde{E}^Q \big[\widetilde{q}_t^1 \widetilde{\sigma}^*_\mu (t) + \widetilde{q}_t^2\widetilde{L}_t \widetilde{h}^*_\mu (t)\big]+
 E^Q\big[q_t^2h_x(t)\,\big|\,\mathcal{F}_t^{Y}\big]\big)L_t E^P\big[ H_\varepsilon (t)\,\big|\,\mathcal{F}_t^{Y}\big]dt\Big]}\notag\\
&\  +E^Q\Big[\!\int^T_0\! \Big\{\big(q_t^1 \delta \sigma (t) + q^2_t L_t \delta h(t) \big)\!+\! E^Q \big[\!\int_t^T\! \Big(E^P\big[(X_s - U_s)
\delta h(t)\,\big|\,\mathcal{F}_s^{Y}\big]\notag \\
&\quad\quad \cdot\Big[E^Q \Big[\widetilde{E}^Q\big[\widetilde{q}_s^1 \widetilde{\sigma}_{\mu}^* (s) + \widetilde{q}_s^2 \widetilde{L}_s
\widetilde{h}^*_{\mu} (s)\big]L_s \delta h(t)\,\big|\,\mathcal{F}_s^{Y}\Big]-\widetilde{E}^Q\big[\widetilde{q}_s^1 \widetilde{\sigma}^*_{\mu}
 (s) + \widetilde{q}_s^2 \widetilde{L}_s\widetilde{h}_{\mu}^* (s)\big]L_s E^P \big[\delta h(t) \big|\mathcal{F}_s^{Y}\big]\Big]\notag\\
&\ \quad +\frac{1}{2} \widetilde{E}^Q\Big[\widetilde{q}_s^1 \widetilde{\sigma}^*_{z\mu} (s) +\widetilde{q}^2_s\widetilde{L}_s
\widetilde{h}^*_{z\mu}(s)\Big]L_s \big(E^P\big[(X_s -U_s)\delta h(t)\,\big|\,\mathcal{F}_s^{Y}\big]\big)^2\Big)ds\,\big|\,
\mathcal{F}_t^{Y}\big]\Big\}\mathbf{1}_{E_\varepsilon}(t)dt\Big]\notag\\
 \label{tau}&\ + o(\varepsilon),\ \ \ \mbox{as}\ \varepsilon \searrow 0.
\end{align}
Recall from (\ref{gamma_s^phi}): For $\phi = \Phi,\ f,\ 0\leq t\leq s\leq T$,
\begin{equation*}
\begin{split}
&\gamma_s^{\phi} \big(\delta h(t)\big) \!=\! E^Q\big[\Gamma_s L_s\delta h(t)\,\big|\,\mathcal{F}_s^{Y}\big]\big\{E^Q\big[\widetilde{E}^Q\big[
    \widetilde{\phi}^*_{\mu} (s)\big] L_s\delta h(t)\,\big|\,\mathcal{F}_s^{Y}\big]\!\\
&\ \ \qquad\quad\!-\!E^Q\big[\widetilde{E}^Q\big[\widetilde{\phi}^*_{\mu} (s)\big] L_s\,\big|\,\mathcal{F}_s^{Y} \big]
E^Q\big[\Gamma_s^1\delta h(t)\,\big|\,\mathcal{F}_s^{Y}\big]\big\}
 +\frac{1}{2}\Big(E^Q\big[\Gamma_s L_s \delta h(t)\,\big|\,\mathcal{F}_s^{Y}\big]\Big)^2\widetilde{E}^Q\big[\widetilde{\phi}^*_{z\mu} (s)\big] L_s\\
\!=\!&\ E^P\big[(X_s \!-\!U_s)\delta h(t)\,\big|\,\mathcal{F}_s^{Y}\big]\big\{E^Q\big[\widetilde{E}^Q\big[\widetilde{\phi}^*_{\mu} (s)\big]
L_s\delta h(t)\,\big|\,\mathcal{F}_s^{Y}\big]\!-\!E^Q\big[\widetilde{E}^Q\big[\widetilde{\phi}^*_{\mu} (s)\big] L_s E^P\big[\delta h(t)\,\big|\,\mathcal{F}_s^{Y}\big]\,\big|\,
\mathcal{F}_s^{Y}\big]\big\}\\
&\  +\frac{1}{2}\Big(E^P\big[(X_s-U_s) \delta h(t)\,\big|\,\mathcal{F}_s^{Y}\big]\Big)^2\widetilde{E}^Q\big[\widetilde{\phi}^*_{z\mu} (s)\big] L_s,\ \ \ 0\leq t\leq s\leq T.
\end{split}
\end{equation*}

As $\ \ \displaystyle p_T^1 = -\Phi_x (T) - L_T\widetilde{E}^Q\Big[E^P \big[\widetilde{\Phi}_{\mu}^* (T) \,\big|\,\mathcal{F}_T^{Y}\big]\Big],\ \ \ $ and\ \ \
 $ \displaystyle p_T^2 = -(X_T -U_T) \widetilde{E}^Q\Big[E^P\big[\widetilde{\Phi}_{\mu}^* (T) \,\big|\,\mathcal{F}_T^{Y}\big]\Big]$ $\displaystyle-\widetilde{E}^Q\Big[\int_0^{U_T}
  \Phi_{\mu}^*(T,y) dy\Big],$\ we deduce from (\ref{reorganized4.34}}) that
\begin{equation*}\label{again reo 4.34}
\begin{split}
0\leq &\ \ J(u^\varepsilon)-J(u)=-E^Q\Big[p_T^1\big(Y_T^{1,\varepsilon}+Y_T^{2,\varepsilon}\big)+p_T^2 \big(K_T^{1,\varepsilon}+K_T^{2,\varepsilon}\big)\Big]\hskip5cm\\
\end{split}
\end{equation*}
\begin{equation}\label{again reo 4.34}
\begin{split}
&+E^Q\Big[\int^T_0\big\{ \big(f_x(t)+L_t \widetilde{E}^Q\big[E^P\big[\widetilde{f}_{\mu}^* (t)\,\big|\,\mathcal{F}_t^{Y}\big]\big]\big)
\big(Y_t^{1,\varepsilon}+Y_t^{2,\varepsilon}\big)\\
&\qquad+\big( (X_t- U_t)\widetilde{E}^Q\big[E^P\big[\widetilde{f}_{\mu}^* (t)\,\big|\,\mathcal{F}_t^{Y}\big]\big]+\widetilde{E}^Q\big[
    \int_0^{U_t} f_{\mu}^* (t,y)dy \big]\big)\big(K_t^{1,\varepsilon}+K_t^{2,\varepsilon}\big)\big\}dt\Big]\\
&+\frac{1}{2} E^Q\Big[\Phi_{xx} (T)\big(Y_T^{1,\varepsilon}\big)^2\Big]+\frac{1}{2}E^Q\Big[\int^T_0 f_{xx} (t)\big(Y_t^{1,\varepsilon}\big)^2dt\Big]\\
&+{E^Q\Big[\widetilde{E}^Q\big[\widetilde{\Phi}_{\mu}^{*} (T)\big]L_TE^P\big[H_\varepsilon (T)\,\big|\,\mathcal{F}_T^Y
\big]\Big]}+{E^Q\Big[\int^T_0\widetilde{E}^Q\big[ \widetilde{f}^*_{\mu}(t)\big]L_t E^P\big[H_\varepsilon (t)\,\big|\,\mathcal{F}_t^{Y}\big]dt\Big]}\\
&+E^Q\Big[\int^T_0\big\{E^Q\big[\gamma^{\Phi}_T\big(\delta h(t)\big)+\int^T_t \gamma^f_s \big(\delta h(t)\big)ds\,\big|\,\mathcal{F}_t^{Y}
\big]+\delta f(t)\big\}\mathbf{1}_{E_\varepsilon}(t)dt\Big]+o(\varepsilon).
\end{split}
\end{equation}
Thus, substituting the above formula \eqref{tau} in \eqref{again reo 4.34} we obtain:
\begin{align}
&0\leq E^Q\Big[\frac{1}{2}\Phi_{xx} (T)\big(Y_T^{1,\varepsilon}\big)^2\Big]\notag\\
&\!-\!{E^Q\Big[\int^T_0\Big(\widetilde{E}^Q\big[\widetilde{q}_t^1 \widetilde{\sigma}^*_\mu (t) + \widetilde{q}_t^2
\widetilde{L}_t \widetilde{h}^*_\mu (t)\big]+E^Q\big[q_t^2h_x(t)\,\big|\,\mathcal{F}_t^{Y}\big]\Big) L_t E^P\big[H_\varepsilon (t)\,\big|\,\mathcal{F}_t^{Y}\big]dt\Big]}\notag\\
&+{E^Q\Big[\widetilde{E}^Q\big[\widetilde{\Phi}_{\mu}^{*} (T)\big]L_TE^P\big[H_\varepsilon (T)\,\big|\,\mathcal{F}_T^Y
\big]\Big]}+{E^Q\Big[\int^T_0\widetilde{E}^Q\big[ \widetilde{f}^*_{\mu}(t)\big]L_t E^P\big[H_\varepsilon (t)\,
\big|\,\mathcal{F}_t^{Y}\big]dt\Big]}\notag\\
&-E^Q\Big[\int_0^T \frac{1}{2}\big(h_{xx}(t) L_tq^2_t - f_{xx} (t) \big)\big(Y_t^{1,\varepsilon}\big)^2dt\Big]-E^Q\Big[\int^T_0
\big(q_t^1\delta \sigma (t) + q_t^2 L_t \delta h(t) -\delta f(t)\big)\mathbf{1}_{E_\varepsilon}(t)dt\Big]\notag\\
\label{tau1}&+E^Q\Big[\int_0^T E^Q\big[ E^P\big[(X_T-U_T) \delta h(t) \big|\mathcal{F}_T^{Y}\big]\big\{E^Q\big[\widetilde{E}^Q\big[
    \widetilde{\Phi}_{\mu}^* (T)\big]L_T \delta h(t)  \,\big|\,\mathcal{F}_T^{Y}\big]\\
&\quad\!-\! \widetilde{E}^Q\big[\widetilde{\Phi}_{\mu}^* (T)\big]L_T E^P \big[\delta h(t) \,\big|\,\mathcal{F}_T^{Y}\big]\big\}
 \!+\!\frac{1}{2} \big(E^P\big[(X_T-U_T) \delta h(t) \,\big|\,\mathcal{F}_T^{Y}\big]\big)^2 \widetilde{E}^Q\big[\widetilde{\Phi}_{z\mu}^*
 (T)\big]L_T \big|\mathcal{F}_t^{Y}\big]\mathbf{1}_{E_\varepsilon}(t)dt\Big]\notag\\
&\!-\!E^Q\Big[\int^T_0 E^Q \Big[\int^T_t \Big( E^P \big[(X_s\!-\!U_s)\delta h(t) \,\big|\,\mathcal{F}_s^{Y}\big]\Big\{E^Q\Big[\widetilde{E}^Q\big[
    \widetilde{q}_s^1 \widetilde{\sigma}^*_{\mu} (s)\! +\! \widetilde{q}_s^2\widetilde{L}_s
\widetilde{h}_{\mu}^*(s)\! -\!\widetilde{f}_{\mu}^* (s)\big] L_s \delta h(t) \,\big|\,\mathcal{F}_s^{Y}\Big]\notag\\
&\quad\quad\quad -\widetilde{E}^Q\Big[\widetilde{q}^1_s \widetilde{\sigma}_{\mu}^* (s) + \widetilde{q}^2_s\widetilde{L}_s \widetilde{h}_{\mu}^*(s)
 -\widetilde{f}_{\mu}^* (s)\Big] L_sE^P\big[\delta h(t)\,\big|\,\mathcal{F}_s^{Y}\big]\Big\}\notag\\
&\quad+\frac{1}{2}\big(E^P \big[(X_s -U_s)\delta h(t) \,\big|\,\mathcal{F}_s^{Y}\big]\big)^2\widetilde{E}^Q\big[\widetilde{q}_s^1
\widetilde{\sigma}^*_{z\mu} (s) +\widetilde{q}^2_s\widetilde{L}_s \widetilde{h}_{z\mu}^*(s)-\widetilde{f}_{z\mu}^* (s)\big]L_s\Big)ds
\,\big|\,\mathcal{F}_t^{Y}\Big]\mathbf{1}_{E_\varepsilon}(t)dt\Big]\notag\\
&+o(\varepsilon),\ \ \ \mbox{as}\ \varepsilon \searrow 0.\notag
\end{align}
With the Hamiltonian
\begin{equation}\label{Hamilton 1}
H(t,x,l,\gamma,v,q_1,q_2):=\sigma(t,\gamma,v)q_1+h(t,x,\gamma,v)lq_2-f(t,x,\gamma,v),
\end{equation}
for $(t,x,l,\gamma,v,q_1,q_2)\in[0,T]\times\mathbb{R}\times\mathbb{R}_+\times\mathcal{P}_2(\mathbb{R})\times U\times\mathbb{R}
\times\mathbb{R}$, and the notations
\begin{equation*}
\begin{split}
\delta H(t)&:=\delta\sigma(t)q_t^1+\delta h(t)L_t q_t^2-\delta f(t),\ \ H_{xx}(t):=h_{xx}(t)L_t q_t^2-f_{xx}(t),\\
H_{x}(t)&:=h_x(t)L_tq_t^2-f_{x}(t),\ \ \qquad\quad\quad\ \ \widetilde{H}^*_\mu(t):=\widetilde{\sigma}^*_\mu(t)\widetilde{q}_t^1+\widetilde{h}^*_\mu(t)\widetilde{L}_t \widetilde{q}_t^2-\widetilde{f}^*_\mu(t),\\
\widetilde{H}^*_{z\mu}(t)&:=\widetilde{\sigma}^*_{z\mu}(t)\widetilde{q}_t^1+\widetilde{h}^*_{z\mu}(t)\widetilde{L}_t \widetilde{q}_t^2-\widetilde{f}^*_{z\mu}(t),
\end{split}
\end{equation*}
where $\big((p^1,q^1),\,(p^2,q^2)\big)$ is the solution of the first adjoint BSDE.\\
\noindent Then we rewrite \eqref{tau1} as follows:

\begin{equation}\label{conclude re 4.34}
\begin{split}
0\leq &E^Q\Big[\frac{1}{2}\Phi_{xx} (T)(Y_T^{1,\varepsilon})^2\Big]-E\Big[\int_0^T \frac{1}{2}H_{xx}(t)\big(Y_t^{1,\varepsilon}\big)^2
dt\Big]-E^Q\Big[\int_0^T \delta H(t) \mathbf{1}_{E_\varepsilon}(t)dt\Big]\\
&+{E^Q\Big[\widetilde{E}^Q\big[\widetilde{\Phi}_{\mu}^{*} (T)\big]L_TE^P\big[H_\varepsilon (T)\,\big|\,\mathcal{F}_T^Y\big]\Big]}\\
&-E^Q\Big[\int^T_0 \Big( \widetilde{E}^Q\big[\widetilde{H}_\mu^*(t)\big]+E^Q \big[ \big(H_x(t)+f_x(t)\big)L_t^{-1}\,\big|\,\mathcal{F}_t^Y\big]\Big)L_t E^P\big[H_\varepsilon (t)\,
\big|\,\mathcal{F}_t^{Y}\big]dt\Big]\\
&+E^Q\Big[\int^T_0 E^Q\big[E^P \big[(X_T -U_T) \delta h(t) \,\big|\,\mathcal{F}_T^{Y}\big]\big\{E^Q\big[\widetilde{E}^Q\big[
    \widetilde{\Phi}_{\mu}^* (T)\big]L_T \delta h(t)\big|\mathcal{F}_T^{Y}\big]\\
&\quad\quad -\widetilde{E}^Q\big[\widetilde{\Phi}^*_{\mu} (T)\big]L_T E^P\big[\delta h(t) \big|\mathcal{F}_T^{Y}\big]\big\}\\
&\quad +\frac{1}{2} \big(E^P \big[(X_T-U_T) \delta h(t)\,\big|\,\mathcal{F}_T^{Y}\big]\big)^2\widetilde{E}^Q\big[\widetilde{\Phi}^*_{z\mu}
 (T)\big] L_T\,\big|\,\mathcal{F}_t^{Y}\big]\mathbf{1}_{E_\varepsilon}(t)dt\Big]\\
&-E^Q\Big[\int_0^T E^Q \Big[\int_t^T\Big(E^P\big[(X_s -U_s) \delta h(t)  \,\big|\,\mathcal{F}_s^{Y}\big]\big\{E^Q\big[\widetilde{E}^Q\big[\widetilde{H}_{\mu}^* (s)\big]
L_s \delta h(t) \,\big|\,\mathcal{F}_s^{Y}\big]\\
&\quad\quad - \widetilde{E}^Q\big[\widetilde{H}_{\mu}^* (s)\big]L_s E^P\big[\delta h(t)\,\big|\,\mathcal{F}_s^{Y}\big]\big\}\\
&\quad +\frac{1}{2} \big(E^P \big[(X_s -U_s)\delta h(t)\,\big|\,\mathcal{F}_s^{Y}\big]\big)^2\widetilde{E}^Q\big[\widetilde{H}_{z\mu}^*
(s)\big]L_s\Big) ds\,\big|\,\mathcal{F}_t^{Y}\Big]\mathbf{1}_{E_\varepsilon}(t)dt\Big]\\
&+o(\varepsilon),\ \ \mbox{as}\ \varepsilon \searrow 0.
\end{split}
\end{equation}
Now we make computations for the second adjoint BSDE, where $\alpha_t^1$ is to be detemined,
\begin{equation}\label{4.56starstar}
 dP_t^1 = -\alpha_t^1 dt + Q_t^{1,1} dB_t^1 + Q_t^{1,2} dY_t,\ \ \ P_T^1 = -\Phi_{xx} (T).
\end{equation}
Recall for $\sigma(t,x,\gamma,u) = \sigma(t,\gamma,u)$\\
\begin{equation*}
\left\{
\begin{split}
&dY_t^{1,\varepsilon}=\Big\{\widetilde{E}^Q\Big[\int_0^{\widetilde{U}_t} \sigma_{\mu}(t,y) dy \widetilde{K}_t^{1,\varepsilon}\Big] + \widetilde{E}^Q\big[\widetilde{\sigma}_{\mu} (t) \widetilde{L}_t \widetilde{V}_t^{1,\varepsilon}\big]+\delta\sigma (t)\mathbf{1}_{E_\varepsilon}(t)\Big\}dB_t^1,\\
&dK_t^{1,\varepsilon} =\Big\{h(t)K_t^{1,\varepsilon} + L_t h_x (t) Y_t^{1,\varepsilon} + L_t \widetilde{E}^Q\big[\int^{\widetilde{U}_t}_0 h_{\mu} (t,y) dy\cdot \widetilde{K}_t^{1,\varepsilon}\big]\\
&\qquad\qquad\qquad\qquad\qquad\qquad + L_t \widetilde{E}^Q\big[\widetilde{h}_{\mu} (t) \widetilde{L}_t \widetilde{V}_t^{1,\varepsilon}\big] + \delta h(t) \mathbf{1}_{E_\varepsilon}(t) L_t\Big\}dY_t,\\
&Y_0^{1,\varepsilon} = K_0^{1,\varepsilon} = 0.
\end{split}\right.
\end{equation*}
Then,
\begin{equation*}
\begin{split}
& dE^Q\Big[P_t^1\big(Y_t^{1,\varepsilon}\big)^2\Big]\\
 =&\ E^Q\Big[\!-\alpha_t^1\big(Y_t^{1,\varepsilon}\big)^2 \!+\!P_t^1\big\{\widetilde{E}^Q\big[\!\int_0^{\widetilde{U}_t} \!\sigma_{\mu}(t,y) dy \widetilde{K}_t^{1,\varepsilon}\big]\!+\!\widetilde{E}^Q\big[\widetilde{\sigma}_{\mu} (t) \widetilde{L}_t \widetilde{V}_t^{1,\varepsilon}\big]+\delta\sigma (t)\mathbf{1}_{E_\varepsilon}(t)\big\}^2\\
&\qquad + 2 Q_t^{1,1} Y_t^{1,\varepsilon} \big\{\widetilde{E}^Q\big[\int^{\widetilde{U}_t}_0 \sigma_{\mu} (t,y) dy \widetilde{K}_t^{1,\varepsilon} \big]+\widetilde{E}^Q\big[\widetilde{\sigma}_{\mu} (t) \widetilde{L}_t\widetilde{V}_t^{1,\varepsilon}\big]+\delta \sigma(t) \mathbf{1}_{E_\varepsilon}(t)\big\}\Big]dt.
\end{split}
\end{equation*}
Also recalling that $\displaystyle \gamma_t^{1,\varepsilon}:= \widetilde{E}^{Q}\Big[\Big(\int_0^{\widetilde{U}_t}\sigma_\mu(t,y)dy\Big)\widetilde{K}^{1,\varepsilon}_t \Big]+ \widetilde{E}^{Q}\Big[\widetilde{\sigma}_\mu(t)\widetilde{L}_t \widetilde{V}_t^{1,\varepsilon}\Big],\ |\gamma_t^{1,\varepsilon}|\leq \rho_t (\varepsilon)\sqrt\varepsilon,$\\
where $\rho_t (\varepsilon)\leq C,\ \rho_t (\varepsilon)\rightarrow 0\ \big(\varepsilon\searrow 0\big),\ t\in [0,T].$\ From Proposition \ref{EstofXL} we get
$$E^Q\Big[P_T^1\big(Y_t^{1,\varepsilon}\big)^2\Big] =  E^Q\Big[\int_0^T\big\{(-\alpha_t^1)\big(Y_t^{1,\varepsilon}\big)^2+
P_t^1\big(\delta\sigma(t)\big)^2\mathbf{1}_{E_\varepsilon}(t)\big\}dt\Big] + o(\varepsilon).$$
Hence, as $\displaystyle \frac{1}{2} E^Q\Big[\Phi_{xx} (T)\big(Y_T^{1,\varepsilon}\big)^2\Big] = -\frac{1}{2} E^Q\Big[P_T^1\big(Y_T^{1,
\varepsilon}\big)^2\Big]$, and recall that
\begin{equation*}
    \begin{split}
\displaystyle h(s,t)=\int_s^t\! h_x(X_r,u_r)dY_r\!-\!\int_s^t (h\cdot h_x)(X_r,u_r)dr,\ \
H_{\varepsilon}(t)=\int_0^t h(s,t)|\delta\sigma(s)|^2 \mathbf{1}_{E_\varepsilon}(s)ds,\ \qquad
\end{split}
\end{equation*}
we have
\begin{equation*}
\begin{split}
&E^Q\Big[\int^T_0\widetilde{E}^Q\big[\widetilde{H}_\mu^*(t)\big] L_t E^P\big[H_\varepsilon (t)\,\big|\,\mathcal{F}_t^{Y}\big]dt\Big]\\
=&E^Q\Big[\int^T_0\widetilde{E}^Q\big[\widetilde{H}_\mu^*(t)\big] L_t E^P\big[\int_0^t h(s,t)|\delta\sigma(s)|^2 \mathbf{1}_{E_\varepsilon}(s)ds\,\big|\,\mathcal{F}_t^{Y}\big]dt\Big]\\
=&E^Q\Big[\int^T_0\int_0^t\widetilde{E}^Q\big[\widetilde{H}_\mu^*(t)\big] L_t E^P\big[ h(s,t)|\delta\sigma(s)|^2 \,\big|\,\mathcal{F}_t^{Y}\big]\mathbf{1}_{E_\varepsilon}(s)ds dt\Big]\\
=&E^Q\Big[\int^T_0\int_s^T\widetilde{E}^Q\big[\widetilde{H}_\mu^*(t)\big] L_t E^P\big[ h(s,t)|\delta\sigma(s)|^2 \,\big|\,\mathcal{F}_t^{Y}\big]dt\mathbf{1}_{E_\varepsilon}(s) ds\Big]\\
=&E^Q\Big[\int^T_0\int_t^T\widetilde{E}^Q\big[\widetilde{H}_\mu^*(s)\big] L_sE^P\big[ h(t,s)|\delta\sigma(t)|^2 \,\big|\,\mathcal{F}_s^{Y}\big]ds\mathbf{1}_{E_\varepsilon}(t) dt\Big],\\
\end{split}
\end{equation*}
and similarly, we can deduce that
\begin{equation*}
\begin{split}
E&^Q\Big[\widetilde{E}^Q\big[\widetilde{\Phi}_{\mu}^{*} (T)\big]L_TE^p\big[H_\varepsilon (T)\,\big|\,\mathcal{F}_T^Y\big]\Big]\\
=&E^Q\Big[\int_0^T\widetilde{E}^Q\big[ \widetilde{\Phi}_{\mu}^{*} (T)L_T E^P[h(t,T)|\delta\sigma(t)|^2\,\big|\,\mathcal{F}_T^Y]\big]\mathbf{1}_{E_\varepsilon}(t)dt\Big].
\end{split}
\end{equation*}
Combined these relations with \eqref{conclude re 4.34}, we have
\begin{align*}
&0\leq-E^Q\Big[\int_0^T \frac{1}{2}\big\{(-\alpha_t^1)\big(Y_t^{1,\varepsilon}\big)^2+P_t^1\big(\delta\sigma(t)\big)^2
\mathbf{1}_{E_\varepsilon}(t)\big\}dt\Big]\\
&\!+\!{E^Q\Big[\int_0^T\widetilde{E}^Q\big[ \widetilde{\Phi}_{\mu}^{*} (T)L_T E^P[h(t,T)|\delta\sigma(t)|^2\,\Big|\,
\mathcal{F}_T^Y]\big]\mathbf{1}_{E_\varepsilon}(t)dt\Big]}\\
&\!-\!{E^Q\Big[\int^T_0\int_t^T\Big(\widetilde{E}^Q\big[\widetilde{H}_\mu^*(s)\big]\!+\!E^Q\big[\big( H_x(s)+f_x(s)\big)L_s^{-1}\,\big|\,\mathcal{F}_s^Y\big] \Big)
L_s E^P\big[h(t,s)|\delta\sigma(t)|^2\, \Big|\,\mathcal{F}_s^{Y}\big]ds\mathbf{1}_{E_\varepsilon}(t)dt\Big]}\\
&\!-\!E^Q\Big[\int_0^T \frac{1}{2} H_{xx} (t)\cdot\big(Y_t^{1,\varepsilon}\big)^2 dt\Big]-E^Q\Big[\int^T_0 \delta H(t) \cdot
\mathbf{1}_{E_\varepsilon}(t)dt\Big]\\
&\!+\!E^Q\Big[\int^T_0E^Q\Big[E^P\big[(X_T-U_T) \delta h(t)\,\big|\,\mathcal{F}_T^{Y}\big]\big\{E^Q\big[\widetilde{E}^Q\big[
\widetilde{\Phi}_{\mu}^* (T)\big]L_T \delta h(t)\,\big|\,\mathcal{F}_T^{Y}\big]\\
&\quad\!-\!\widetilde{E}^Q\big[\widetilde{\Phi}_{\mu}^* (T)\big]L_T E^P\big[\delta h(t) \,\big|\,\mathcal{F}_T^{Y}\big]\big\}
\!+\!\frac{1}{2}\big(E^P\big[(X_T-U_T) \delta h(t)\,\big|\,\mathcal{F}_T^{Y}\big]\big)^2 \widetilde{E}^Q\big[\widetilde{\Phi}^*_{z\mu}
 (T)\big]L_T\,\Big|\,\mathcal{F}_t^{Y}\Big]\mathbf{1}_{E_\varepsilon}(t)dt\Big]\\
&\!-\!E^Q\Big[\int^T_0 E^Q\Big[\int^T_t\Big(E^P\big[(X_s -U_s)\delta h(t)\big|\mathcal{F}_s^{Y}\big]\big\{E^Q\big[\widetilde{E}^Q\big[\widetilde{H}^*_{\mu}
 (s)\big]L_s \delta h(t) \,\big|\,\mathcal{F}_s^{Y}\big]\\
&\quad\!-\!\widetilde{E}^Q\big[ \widetilde{H}_{\mu}^* (s)\big]L_s E^P\big[\delta h(t)\,\big|\,\mathcal{F}_s^{Y}\big]\big\}
 \!+\! \frac{1}{2}\big(E^P\big[(X_s - U_s) \delta h(t)\, \big|\,\mathcal{F}_s^{Y}\big]\big)^2 \widetilde{E}^Q\big[\widetilde{H}^*_{z\mu}
(s)\big]L_s\Big)ds\,\Big|\,\mathcal{F}_t^{Y}\Big]\mathbf{1}_{E_\varepsilon}(t)dt\Big]\\
&+o(\varepsilon),\ \ \ \mbox{as}\ \varepsilon\searrow 0.
\end{align*}
Recall that $\displaystyle\delta h(t) = \phi (X_t)\delta h_1(t),$ and $\displaystyle\delta h_1(t)$ is $\mathcal{F}_t^{Y}$-measurable.
Thus, choosing $\displaystyle\alpha^1_t = H_{xx}(t)$, we obtain that
\begin{equation*}
\begin{split}
0\leq& -E^Q\Big[\int^T_0\big(\delta H(t) +\frac{1}{2} P_t^1\big(\delta \sigma (t)\big)^2\big)\mathbf{1}_{E_\varepsilon}(t)dt\Big]-E^Q\Big[\int_0^T M_t\big(
    \delta\sigma(t)\big)^2\mathbf{1}_{E_\varepsilon}(t)dt\Big]\\
&-E^Q\Big[\int_0^T R_t \big(\delta h_1(t)\big)^2\mathbf{1}_{E_\varepsilon}(t)dt\Big] + o(\varepsilon),\ \ \ \mbox{as}\ \varepsilon\searrow 0,
\end{split}
\end{equation*}
where, for all $t\in[0,T],$
\begin{equation*}
\begin{split}
M_t:=&-\widetilde{E}^Q\big[ \widetilde{\Phi}_{\mu}^{*} (T)L_T E^P[h(t,T)\,\big|\,\mathcal{F}_T^Y]\big]\\
&+\int_t^T\Big(\widetilde{E}^Q\big[\widetilde{H}_\mu^*(s)\big]+E^Q\big[(H_x(s) +f_x(s)) L_s^{-1}\,\big|\,\mathcal{F}_s^Y\big] \Big)L_s E^P\big[h(t,s)\,\big|\,\mathcal{F}_s^{Y}\big]ds,\\
\end{split}
\end{equation*}
\begin{equation*}
\begin{split}
R_t: = &-E^Q\Big[E^P\big[\big(X_T -U_T\big)\phi(X_t)\,\big|\,\mathcal{F}_T^{Y}\big]\Big\{E^Q\big[\widetilde{E}^Q\big[\Phi^*_{\mu}
 (T)\big]L_T \phi(X_t)\,\big|\,\mathcal{F}_T^{Y}\big]\\
&\ \ \!-\!\widetilde{E}^Q\big[\widetilde{\Phi}_{\mu}^* (T)\big]L_T E^P\big[\phi(X_t)\,\big|\,\mathcal{F}_T^{Y}\big]\Big\}\!+\!\frac{1}{2} \big(E^P\big[(X_T -U_T)
\phi(X_t)\,\big|\,\mathcal{F}_T^{Y}\big]\big)^2\widetilde{E}^Q\big[\widetilde{\Phi}^*_{z \mu}(T)\big]L_T\,\big|\,\mathcal{F}_t^{Y}\Big]\\
&+ E^Q\Big[\int^T_t\Big(E^P\big[(X_s -U_s)\phi(X_t) \,\big|\,\mathcal{F}_s^{Y}\big]\Big\{E^Q\big[\widetilde{E}^Q\big[\widetilde{H}_{\mu}^*(s)\big]L_s\phi(X_t)\,\big|\,\mathcal{F}_s^{Y}\big]\\
&\ \  \!-\!\widetilde{E}^Q\big[\widetilde{H}_{\mu}^* (s)\big]L_sE^P\big[\phi(X_t)\,\big|\,\mathcal{F}_s^{Y}\big]\Big\}
\!+\!\frac{1}{2} \big(E^P\big[(X_s -U_s)\phi(X_t)\,\big|\,\mathcal{F}_s^{Y}\big]\big)^2\widetilde{E}^Q\big[\widetilde{H}^*_{z\mu}
 (s)\big]L_s\Big)ds\,\big|\,\mathcal{F}_t^{Y}\Big] .
\end{split}
\end{equation*}
Notice now from \eqref{4.56starstar} the second adjoint BSDE is the following
\begin{equation}\label{4.56starstarstar}
    dP_t^1 = -H_{xx}(t) dt +Q_t^{1,1} dB_t^1 + Q^{1,2}_t dY_t,\ \ P_T^1 = -\varPhi_{xx}(T).
\end{equation}
Thus we have:
$$0\leq -E^Q\Big[\int_0^T\big(\delta H(t) + \frac{1}{2} P_t^1\big(\delta \sigma(t)\big)^2 + R_t \big(\delta h_1(t)\big)^2
+{M_t\big(\delta\sigma(t)\big)^2}\big)\mathbf{1}_{E_\varepsilon}(t)dt\Big]+o(\varepsilon),$$
and, as $\displaystyle v\in \mathcal{U}_{ad}$ has been fixed arbitrarily, Lebesgue's differentiation theorem implies
\begin{equation*}
\begin{split}
E^Q\Big[H\big(t,X_t,&L_t,v_t,q_t^1,q_t^2\big)-H\big(t,X_t,L_t,u_t,q_t^1,q_t^2\big)+\frac{1}{2} P_t^1\big|\sigma(t,\mu_t,v_t)-
\sigma(t,\mu_t,u_t)\big|^2\\
+&{M_t\big|\sigma(t,\mu_t,v_t)-\sigma(t,\mu_t,u_t)\big|^2}+R_t\big|h_1(t,\mu_t,v_t)-h_1(t,\mu_t,u_t)\big|^2
\,\big|\,\mathcal{F}_t^{Y}\Big]\leq 0,\ \ dtdQ\text{-a.s.,}
\end{split}
\end{equation*}
for all $\displaystyle v\in \mathcal{U}_{ad}$. (The fact that we have to take in this formula $E^Q\big[\ \cdot\ \big|\mathcal{F}_t^{Y}
\big]$ stems from the fact the control processes are $\mathbb{F}^Y$-adapted).

Now, finally, we can obtain our Peng's stochastic maximum principle.
\begin{theorem}\label{th4.2}
Under the Assumptions (H2) and (H3), let $u\in\mathcal{U}_{ad}$ be optimal and $(X,L)$ be the associated solution of system
\eqref{system}. Then, for all $v\in U$, it holds that for $dtdQ$-a.e. $(t,\omega)\in[0,T]\times\Omega$,
\begin{equation}\label{finalSMP}
\begin{split}
E^Q\Big[H\big(t,&X_t,L_t,v,q_t^1,q_t^2\big)-H\big(t,X_t,L_t,u_t,q_t^1,q_t^2\big)+\frac{1}{2} P_t^1\big(\sigma(t,\mu_t,v)-
\sigma(t,\mu_t,u_t)\big)^2\\
+&{M_t\big(\sigma(t,\mu_t,v)-\sigma(t,\mu_t,u_t)\big)^2}+R_t\big(h_1(t,\mu_t,v)-h_1(t,\mu_t,u_t)\big)^2\,
\big|\,\mathcal{F}_t^{Y}\Big]\leq 0,\
\end{split}
\end{equation}
where $\big((p^1,(q^1,\check{q}^1)),\, (p^2,(\check{q}^2,q^2))\big)$ and $(P^1,(Q^{1,1},Q^{1,2}))$
are the unique solutions to \eqref{SMP4} and \eqref{4.56starstarstar}, respectively.
\end{theorem}

\begin{remark}
Comparing the formula with the SMP got in previous works by different authors, namely in the classical case (no mean field, and no
conditional expectation), the terms with $R=(R_t)_{t\in[0,T]}$ and $M=(M_t)_{t\in [0,T]}$ are new here. Note that $R=(R_t)_{t\in[0,T]}$ depends in a nonlocal
way on $(X,L,U)$. This comes from the fact that we have a mean-field control problem involving the law of the conditional expectation
of the controlled state process.
\end{remark}

\section{Appendix}

For a better appreciation of the optimal condition, let us shortly discuss the case, where $\mu_t = P_{E^P[X_t\,|\,\mathcal{F}_t^Y]}$
is replaced now by mathematically easier $\mu_t = P_{\varphi(X_t,Y_{\cdot\wedge t})}$, for some Borel function $\varphi$.

More precisely, let us consider the case
$\displaystyle U_t=\varphi(X_t,Y_{\cdot\wedge t})$ instead of $\displaystyle U_t=E^P[X_t\,|\,\mathcal{F}_t^Y]$, where $\varphi:\mathbb{R}\times C_T\rightarrow
\mathbb{R}$ is a Borel measurable function differentiable w.r.t. $x\in\mathbb{R}$, and with bounded derivative $\varphi_x:\mathbb{R}\times C_T\rightarrow\mathbb{R}$. That is, we consider the full observation as a more classical case to have a comparison with Peng's SMP for $\displaystyle U_t=E^P[X_t\,|\,\mathcal{F}_t^Y]$. We can deduce the second variational equations with an approach similar to that for \eqref{2ndVE5}. The form of the first and second order variational equations, respectively, is the following:
{\small
\begin{equation}\label{2ndFofVE}
\left\{
\begin{split}
&dY_t^{1,\varepsilon}=\Big\{\sigma_x(t)Y_t^{1,\varepsilon}\!+\!\widetilde{E}^Q\Big[\!\int_0^{\widetilde{U}_t}\!\sigma_\mu(t,y)dy\cdot\widetilde{K}_t^{1,\varepsilon}\Big]
\!+\!\widetilde{E}^Q\big[\widetilde{\sigma}_\mu(t)\widetilde{L}_t\widetilde{\overline{V}}_t^{1,\varepsilon}\big]\!+\!\delta\sigma(t)\mathbf{1}_{E_\varepsilon}(t)\Big\}dB_t^1,\ \\
 &dK_t^{1,\varepsilon}=\Big\{h(t)K_t^{1,\varepsilon}+\Big(h_x(t)Y_t^{1,\varepsilon}+\widetilde{E}^Q\Big[\int_0^{\widetilde{U}_t}h_\mu(t,y)dy\cdot\widetilde{K}_t^{1,\varepsilon}\Big]\\
&\qquad\qquad+\widetilde{E}^Q\big[\widetilde{h}_\mu(t)\widetilde{L}_t\widetilde{\overline{V}}_t^{1,\varepsilon}\big]+\delta h(t)\mathbf{1}_{E_\varepsilon}(t)\Big)L_t\Big\}dY_t,\\
&Y_0^{1,\varepsilon}=0;\ K_0^{1,\varepsilon}=0,  \\
&\overline{V}_t^{1,\varepsilon}=\varphi_x(X_t,Y_{\cdot\wedge t})Y_t^{1,\varepsilon},\ \ \ t\in[0,T].
\end{split}\right.
\end{equation}
\begin{equation}\label{2ndVE6}
\left\{
\begin{split}
&dY_t^{2,\varepsilon}=\Big\{\sigma_x(t)Y_t^{2,\varepsilon}+\frac{1}{2}\sigma_{xx}(t)(Y_t^{1,\varepsilon})^2+\widetilde{E}^Q\Big[\int_0^{\widetilde{U}_t}\sigma_\mu(s,y)dy\cdot
\widetilde{K}_t^{2,\varepsilon}\Big]+\widetilde{E}^Q\big[\widetilde{\sigma}_\mu(t)\widetilde{L}_t\widetilde{\overline{V}}_t^{2,\varepsilon}\big]\\
&\qquad\quad+\widetilde{E}^Q\big[\widetilde{\sigma}_\mu(t)\widetilde{\overline{V}}_t^{1,\varepsilon}\widetilde{K}_t^{1,\varepsilon}\big]+\frac{1}{2}\widetilde{E}^Q
\big[\widetilde{\sigma}_{z\mu}(t)\widetilde{L}_t(\widetilde{\overline{V}}_t^{1,\varepsilon})^2\big]\\
&\qquad\quad+\Big(\delta\sigma_x(t)Y_t^{1,\varepsilon}+\widetilde{E}^Q\Big[\int_0^{\widetilde{U}_t}\delta\sigma_\mu(t,y)dy\cdot\widetilde{K}_t^{1,\varepsilon}\Big]+\widetilde{E}^Q
\big[\delta\widetilde{\sigma}_\mu(t)\widetilde{L}_t\widetilde{\overline{V}}_t^{1,\varepsilon}\big]\Big)\mathbf{1}_{E_\varepsilon}(t)\Big\}dB_t^1,\\
&dK_t^{2,\varepsilon}=\Big\{h(t)K_t^{2,\varepsilon}+h_x(t)L_t Y_t^{2,\varepsilon}+h_x(t)Y_t^{1,\varepsilon}K_t^{1,\varepsilon}+\frac{1}{2}h_{xx}(t)L_t (Y_t^{1,\varepsilon})^2\\
&\qquad\quad+L_t\widetilde{E}^Q\Big[\int_0^{\widetilde{U}_t}h_\mu(t,y)dy\cdot\widetilde{K}_t^{2,\varepsilon}\Big]+L_t\widetilde{E}^Q\big[\widetilde{h}_\mu(t)\widetilde{L}_t
\widetilde{\overline{V}}_t^{2,\varepsilon}\big]\\
&\qquad\quad+L_t\widetilde{E}^Q\big[\widetilde{h}_\mu(t)\widetilde{\overline{V}}_t^{1,\varepsilon}\widetilde{K}_t^{1,\varepsilon}\big]+\frac{1}{2}L_t\widetilde{E}^Q\big[
\widetilde{h}_{z\mu}(t)\widetilde{L}_t(\widetilde{\overline{V}}_t^{1,\varepsilon})^2\big]\\
&\qquad\quad+\Big(\delta h(t)K_t^{1,\varepsilon}+\delta h_x(t)L_t Y_t^{1,\varepsilon}+L_t\widetilde{E}^Q\Big[\int_0^{\widetilde{U}_t}\delta h_\mu(t,y)dy\cdot\widetilde{K}_t^{1,\varepsilon}\Big]\\
&\qquad\quad\quad\ \ +L_t\widetilde{E}^Q\big[\delta\widetilde{h}_\mu(t)\widetilde{L}_t\widetilde{\overline{V}}_t^{1,\varepsilon}\big]\Big)\mathbf{1}_{E_\varepsilon}(t)\Big\}dY_t,\\
&Y_0^{2,\varepsilon}=0;\ K_0^{2,\varepsilon}=0,\\
&\overline{V}_t^{2,\varepsilon}=\varphi_x(X_t,Y_{\cdot\wedge t})Y_t^{2,\varepsilon},\ \ \ t\in[0,T].
\end{split}\right.
\end{equation}}

Also here we get the same results as those in the Propositions \ref{Wellof1stVar}--\ref{Wellof2ndVar}, and for Proposition \ref{EstofYK}
we have (iii) for  $U_t^\varepsilon=\varphi(X_t^\varepsilon,Y_{\cdot\wedge t})$, $U_t=\varphi(X_t,Y_{\cdot\wedge t})$,\
$\overline{V}_t^{1,\varepsilon}=\varphi_x(X_t,Y_{\cdot\wedge t})Y_t^{1,\varepsilon}$ and $\overline{V}_t^{2,\varepsilon}=\varphi_x
(X_t,Y_{\cdot\wedge t})Y_t^{2,\varepsilon}$ instead of $V_t^{1,\varepsilon}$ and $V_t^{2,\varepsilon}$, respectively.

The first order adjoint BSDEs and their well-posedness can be got similarly to \eqref{SMP4} and Proposition \ref{4.777}.
We give here only the equations without the proof:
\begin{equation}\label{Firadj}\left\{
\begin{split}
&dp_t^1=-\alpha_t(q_t^1,q_t^2)dt+q_t^1dB_t^1+\check{q}_t^1dY_t,\ t\in[0,T],\\
&p_T^1=-\Phi_x(T)-L_T\varphi_x(X_T,Y_{\cdot\wedge T})\widetilde{E}^Q[\widetilde{\Phi}^*_\mu(T)];\\
&dp_t^2=-\beta_t(q_t^1,q_t^2)dt+\check{q}_t^2dB_t^1+q_t^2dY_t,\ t\in[0,T],\\
&p_T^2=-\widetilde{E}^Q\Big[\int_0^{U_T}\Phi^*_\mu(T,y)dy\Big],
\end{split}\right.
\end{equation}
where
\begin{equation*}
\begin{split}
\alpha_t(q_t^1,q_t^2):=&\sigma_x(t) q_t^1+L_t\varphi_x(X_t,Y_{\cdot\wedge t})\widetilde{E}^Q[\widetilde{q}_t^1 \widetilde{\sigma}^*_\mu(t)]+h_x(t) L_t q_t^2+L_t\varphi_x(X_t,Y_{\cdot\wedge t})\widetilde{E}^Q[\widetilde{q}_t^2 \widetilde{L}_t\widetilde{h}^*_\mu(t)]\\
&-f_x(t)-L_t\varphi_x(X_t,Y_{\cdot\wedge t})\widetilde{E}^Q[\widetilde{f}^*_\mu(t)];\\
\beta_t(q_t^1,q_t^2):=&h(t)q_t^2+\widetilde{E}^Q\big[\widetilde{q}_t^1 \int_0^{U_t}\sigma^*_\mu(t,y)dy\big]+\widetilde{E}^Q\big[\widetilde{q}_t^2 \widetilde{L}_t \int_0^{U_t}h^*_\mu(t,y)dy\big]-\widetilde{E}^Q\big[\int_0^{U_t}f^*_\mu(t,y)dy\big],\
\end{split}
\end{equation*}
$t\in[0,T].$\ A straight-forward computations give also here the duality result:
\begin{equation}\label{Dua}
\begin{split}
E&^Q\big[p_T^1 Y_T^{1,\varepsilon}+p_T^2 K_T^{1,\varepsilon}\big]\\
=&E^Q\Big[\int_0^T\Big(Y_t^{1,\varepsilon}\big(f_x(t)+\widetilde{E}^Q[\widetilde{f}^*_\mu(t)]L_t\varphi_x(X_t,Y_{\cdot \wedge t})\big)+K_t^{1,\varepsilon}\widetilde{E}^Q\Big[\int_0^{U_t}f^*_\mu(t,y)dy\Big]\Big)dt\Big]\\
&+E^Q\Big[\int_0^T\big(q_t^1\delta\sigma(t)+q_t^2 L_t\delta h(t)\big)\mathbf{1}_{E_\varepsilon}(t)dt\Big].
\end{split}
\end{equation}
Then, similar to our computations for the case $U_t = E^P[X_t\,|\mathcal{F}_t^Y]$, we see that, for
$U_t=\varphi(X_t,Y_{\cdot\wedge t})$,
\begin{align}\label{MR2}
0\leq& J(u^\varepsilon)-J(u) \notag\\
=&E^Q\Big[(Y_T^{1,\varepsilon}+Y_T^{2,\varepsilon})\Big(\Phi_x(T)+\widetilde{E}^{Q}\big[\widetilde{\Phi}^*_\mu(T)\big]L_T\varphi_x(X_T,Y_{\cdot \wedge T})\Big)\Big]\notag\\
&+E^Q\Big[(K_T^{1,\varepsilon}+K_T^{2,\varepsilon})\widetilde{E}^Q\Big[\int_0^{U_T}\Phi^*_\mu(T,y)dy\Big]\Big]\notag\\
&+E^Q\Big[\int_0^T(Y_t^{1,\varepsilon}+Y_t^{2,\varepsilon})\Big(f_x(t)+\widetilde{E}^{Q}\big[\widetilde{f}^*_\mu(t)\big]L_t\varphi_x(X_t,Y_{\cdot \wedge t})\Big)dt\Big]\notag\\
&+E^Q\Big[\int_0^T(K_t^{1,\varepsilon}+K_t^{2,\varepsilon})\widetilde{E}^Q\Big[\int_0^{U_t}f^*_\mu(t,y)dy\Big]\Big]\\
&+E^Q\Big[\frac{1}{2}(Y_T^{1,\varepsilon})^2\Big(\Phi_{xx}(T)+\widetilde{E}^{Q}\big[\widetilde{\Phi}^*_{z\mu}(T)\big]L_T\varphi_x(X_T,Y_{\cdot \wedge T})^2\Big)\Big]\notag\\
&+E^Q\Big[Y_T^{1,\varepsilon}K_T^{1,\varepsilon}\widetilde{E}^{Q}\big[\widetilde{\Phi}^*_\mu(T)\big]\varphi_x(X_T,Y_{\cdot \wedge T})\Big]\notag\\
&+E^Q\Big[\frac{1}{2}\int_0^T(Y_t^{1,\varepsilon})^2\Big(f_{xx}(t)+\widetilde{E}^{Q}\big[\widetilde{f}^*_{z\mu}(t)\big]L_t\varphi_x(X_t,Y_{\cdot \wedge t})^2\Big)dt\Big]\notag\\
&+E^Q\Big[\int_0^T Y_t^{1,\varepsilon}K_t^{1,\varepsilon}\widetilde{E}^{Q}\big[\widetilde{f}^*_\mu(t)\big]\varphi_x(X_t,Y_{\cdot \wedge t})dt\Big]+E^Q\Big[\int_0^T\delta f(t)\mathbf{1}_{E_\varepsilon}(t)dt\Big]+o(\varepsilon).\notag
\end{align}
On the other hand, using \eqref{2ndVE6} and \eqref{Firadj} and applying It\^{o}'s formula to $p_t^1 Y_t^{2,\varepsilon}+p_t^2 K_t^{2,\varepsilon}$ we compute $E^Q\big[p_T^1 Y_T^{2,\varepsilon}+p_T^2 K_T^{2,\varepsilon}\big]$ by using Fubini's Theorem:
\begin{align}\label{MRII}
\notag&E^Q\big[p_T^1 Y_T^{2,\varepsilon}+p_T^2 K_T^{2,\varepsilon}\big]\\
\notag=&E^Q\Big[\int_0^T \Big(Y_t^{2,\varepsilon}\big(f_x(t)+\widetilde{E}^Q[\widetilde{f}^*_\mu(t)]L_t\varphi_x(X_t,Y_{\cdot \wedge t})\big)+K_t^{2,\varepsilon}\widetilde{E}^Q\Big[\int_0^{U_t}f^*_\mu(t,y)dy\Big]\Big)dt\Big]\\
\notag&+E^Q\Big[\frac{1}{2}\int_0^T(Y_t^{1,\varepsilon})^2\Big\{\sigma_{xx}(t)q_t^1+h_{xx}(t)q_t^2 L_t+L_t\varphi_x(X_t,Y_{\cdot\wedge t})^2\widetilde{E}^Q\big[\widetilde{q}_t^1\widetilde{\sigma}^*_{z\mu}(t)+\widetilde{L}_t\widetilde{q}_t^2\widetilde{h}^*_{z\mu}(t)\big]\Big\}dt\Big]\\
&+E^Q\Big[\int_0^T Y_t^{1,\varepsilon}K_t^{1,\varepsilon}\Big(h_x(t)q_t^2+\varphi_x(X_t,Y_{\cdot\wedge t})\widetilde{E}^Q\big[\widetilde{q}_t^1\widetilde{\sigma}^*_{\mu}(t)+\widetilde{L}_t\widetilde{q}_t^2\widetilde{h}^*_{\mu}(t)\big]\Big)dt\Big]\\
\notag&+E^Q\Big[\int_0^T Y_t^{1,\varepsilon}\Big(q_t^1\delta\sigma_x(t)+L_t q_t^2\delta h_x(t)+L_t\varphi_x(X_t,Y_{\cdot\wedge t})\widetilde{E}^Q\big[\widetilde{q}_t^1\delta\widetilde{\sigma}^*_{\mu}(t)+\widetilde{L}_t\widetilde{q}_t^2\delta\widetilde{h}^*_{\mu}(t)\big]\Big)\mathbf{1}_{E_\varepsilon}(t)dt\Big]\\
\notag&+E^Q\Big[\int_0^T K_t^{1,\varepsilon}\Big(q_t^2\delta h(t)+\widetilde{E}^Q\Big[\widetilde{q}_t^1\int_0^{U_t}\delta\sigma^*_\mu(t,y)dy+\widetilde{L}_t\widetilde{q}_t^2\int_0^{U_t}\delta h^*_\mu(t,y)dy\Big]\Big)\mathbf{1}_{E_\varepsilon}(t)dt\Big].
\end{align}
Thus, combining \eqref{Dua} and \eqref{MRII}, we obtain the following duality
\begin{align}\label{MRIII}
&E^Q\big[p_T^1(Y_T^{1,\varepsilon}+Y_T^{2,\varepsilon})+p_T^2(K_T^{1,\varepsilon}+ K_T^{2,\varepsilon})\big] \notag\\
=&E^Q\Big[\int_0^T (Y_t^{1,\varepsilon}+Y_t^{2,\varepsilon})\big(f_x(t)+\widetilde{E}^Q[\widetilde{f}^*_\mu(t)]L_t\varphi_x(X_t,Y_{\cdot \wedge t})\big)dt\Big]\notag\\
&+E^Q\Big[\int_0^T (K_t^{1,\varepsilon}+K_t^{2,\varepsilon})\widetilde{E}^Q\Big[\int_0^{U_t}f^*_\mu(t,y)dy\Big]dt\Big]\\
&+E^Q\Big[\frac{1}{2}\int_0^T|Y_t^{1,\varepsilon}|^2\Big\{\sigma_{xx}(t)q_t^1+h_{xx}(t)q_t^2 L_t+L_t\varphi_x(X_t,Y_{\cdot\wedge t})^2\widetilde{E}^Q\big[\widetilde{q}_t^1
\widetilde{\sigma}^*_{z\mu}(t)+\widetilde{L}_t\widetilde{q}_t^2\widetilde{h}^*_{z\mu}(t)\big]\Big\}dt\Big]\notag\\
&+E^Q\Big[\int_0^T Y_t^{1,\varepsilon}K_t^{1,\varepsilon}\Big(h_x(t)q_t^2+\varphi_x(X_t,Y_{\cdot\wedge t})\widetilde{E}^Q\big[\widetilde{q}_t^1\widetilde{\sigma}^*_{\mu}(t)
+\widetilde{L}_t\widetilde{q}_t^2\widetilde{h}^*_{\mu}(t)\big]\Big)dt\Big]\notag\\
&+E^Q\Big[\int_0^T\big(q_t^1\delta\sigma(t)+q_t^2 L_t\delta h(t)\big)\mathbf{1}_{E_\varepsilon}(t)dt\Big]+R_\varepsilon,\notag
\end{align}
where
\begin{equation*}
\begin{split}
R_\varepsilon=&E^Q\Big[\int_0^T Y_t^{1,\varepsilon}\Big(q_t^1\delta\sigma_x(t)+L_t q_t^2\delta h_x(t)+L_t\varphi_x(X_t,Y_{\cdot\wedge t})\widetilde{E}^Q\big[\widetilde{q}_t^1\delta\widetilde{\sigma}^*_{\mu}(t)+\widetilde{L}_t\widetilde{q}_t^2\delta\widetilde{h}^*_{\mu}(t)\big]\Big)\mathbf{1}_{E_\varepsilon}(t)dt\Big]\\
&+E^Q\Big[\int_0^T K_t^{1,\varepsilon}\Big(q_t^2\delta h(t)+\widetilde{E}^Q\Big[\widetilde{q}_t^1\int_0^{U_t}\delta\sigma^*_\mu(t,y)dy+\widetilde{L}_t\widetilde{q}_t^2\int_0^{U_t}\delta h^*_\mu(t,y)dy\Big]\Big)\mathbf{1}_{E_\varepsilon}(t)dt\Big].
\end{split}
\end{equation*}
Notice that $|R_\varepsilon|\leq C\varepsilon\rho(\varepsilon),\ \varepsilon>0$, with $\rho(\varepsilon)\rightarrow0$, as $\varepsilon\searrow0$. Indeed, putting
$\theta_t:=q_t^1\delta\sigma_x(t)+L_t q_t^2\delta h_x(t)+L_t\varphi_x(X_t,Y_{\cdot\wedge t})\widetilde{E}^Q\big[\widetilde{q}_t^1\delta\widetilde{\sigma}^*_{\mu}(t)
+\widetilde{L}_t\widetilde{q}_t^2\delta\widetilde{h}^*_{\mu}(t)\big]$, $t\in[0,T]$,
from Proposition \ref{EstofXL} we have
$$ \Big|E^Q\Big[\int_0^T Y_t^{1,\varepsilon}\theta_t\mathbf{1}_{E_\varepsilon}(t)dt\Big]\Big|\leq E^Q\Big[\sup_{t\in[0,T]}|Y_t^{1,\varepsilon}|\Big(\int_0^T |\theta_t|^2
\mathbf{1}_{E_\varepsilon}(t)dt\Big)^{\frac{1}{2}}\Big(\int_0^T \mathbf{1}_{E_\varepsilon}(t)dt\Big)^{\frac{1}{2}}\Big]\leq\varepsilon\rho_1(\varepsilon), $$
with $\displaystyle \rho_1(\varepsilon):=\Big(E^Q\Big[\int_0^T |\theta_t|^2\mathbf{1}_{E_\varepsilon}(t)dt\Big]\Big)^{\frac{1}{2}}\rightarrow0\ (\varepsilon\searrow0)$, as
$\theta\in L^2_{\mathbb{F}}([0,T],Q)$. Similarly we estimate the second term of $R_\varepsilon$, and we get the wished estimate for $R_\varepsilon$. Using the duality \eqref{MRIII}
and recalling the terminal conditions of \eqref{Firadj}:
$$
p_T^1=-\Phi_x(T)-L_T\varphi_x(X_T,Y_{\cdot\wedge T})\widetilde{E}^Q[\widetilde{\Phi}^*_\mu(T)],\ \ \ p_T^2=-\widetilde{E}^Q\Big[\int_0^{U_T}\Phi^*_\mu(T,y)dy\Big],
$$
we deduce from \eqref{MR2} that
\begin{align}\label{xingjing}
0\leq& J(u^\varepsilon)-J(u) \notag\\
=&E^Q\Big[\frac{1}{2}(Y_T^{1,\varepsilon})^2\Big(\Phi_{xx}(T)+L_T\varphi_x(X_T,Y_{\cdot \wedge T})^2\widetilde{E}^{Q}\big[\widetilde{\Phi}^*_{z\mu}(T)\big]\Big)\Big]\notag\\
&+E^Q\Big[Y_T^{1,\varepsilon}K_T^{1,\varepsilon}\varphi_x(X_T,Y_{\cdot \wedge T})\widetilde{E}^{Q}\big[\widetilde{\Phi}^*_\mu(T)\big]\Big]\notag\\
&-E^Q\Big[\frac{1}{2}\int_0^T(Y_t^{1,\varepsilon})^2\Big(\big(\sigma_{xx}(t)q_t^1+h_{xx}(t)L_t q_t^2-f_{xx}(t)\big)\\
&\qquad\qquad\qquad\qquad+L_t\varphi_x(X_t,Y_{\cdot \wedge t})^2\widetilde{E}^{Q}\big[\widetilde{\sigma}^*_{z\mu}(t)\widetilde{q}_t^1+\widetilde{h}^*_{z\mu}(t)\widetilde{L}_t\widetilde{q}_t^2-\widetilde{f}^*_{z\mu}(t)\big]\Big)dt\Big]\notag\\
&-E^Q\Big[\int_0^T Y_t^{1,\varepsilon}K_t^{1,\varepsilon}\Big(h_x(t)q_t^2+\varphi_x(X_t,Y_{\cdot \wedge t})\widetilde{E}^{Q}\big[\widetilde{\sigma}^*_{\mu}(t)\widetilde{q}_t^1+\widetilde{h}^*_{\mu}(t)\widetilde{L}_t\widetilde{q}_t^2-\widetilde{f}^*_\mu(t)\big]\Big)dt\Big]\notag\\
&-E^Q\Big[\int_0^T\big(\delta\sigma(t)q_t^1+\delta h(t)L_t q_t^2-\delta f(t)\big)\mathbf{1}_{E_\varepsilon}(t)dt\Big]+o(\varepsilon),\ \mbox{ as } \varepsilon\searrow0.\notag
\end{align}
With the Hamiltonian
\begin{equation}\label{Hamilton}
H(t,x,l,\gamma,v,q_1,q_2):=\sigma(t,x,\gamma,v)q_1+h(t,x,\gamma,v)lq_2-f(t,x,\gamma,v),
\end{equation}
for $(t,x,l,\gamma,v,q_1,q_2)\in[0,T]\times\mathbb{R}\times\mathbb{R}_+\times\mathcal{P}_2(\mathbb{R})\times U\times\mathbb{R}\times\mathbb{R}$, and the notations
\begin{equation*}
\begin{split}
\delta H(t)&=\delta\sigma(t)q_t^1+\delta h(t)L_t q_t^2-\delta f(t),\ \ H_{xx}(t)=\sigma_{xx}(t)q_t^1+h_{xx}(t)L_t q_t^2-f_{xx}(t),\\
H_{xl}(t)&=h_x(t)q_t^2,\ \qquad\qquad\qquad\qquad\quad \widetilde{H}_\mu^*(t)=\widetilde{\sigma}^*_\mu(t)\widetilde{q}_t^1+\widetilde{h}^*_\mu(t)\widetilde{L}_t \widetilde{q}_t^2
-\widetilde{f}_\mu^*(t),\\
\widetilde{H}_{z\mu}^*(t)&=\widetilde{\sigma}^*_{z\mu}(t)\widetilde{q}_t^1+\widetilde{h}^*_{z\mu}(t)\widetilde{L}_t \widetilde{q}_t^2-\widetilde{f}^*_{z\mu}(t),
\end{split}
\end{equation*}
where $\big((p^1,(q^1,\check{q}^1)),\, (p^2,(\check{q}^2,q^2))\big)$ is the solution of the first adjoint BSDE \eqref{Firadj}, we can rewrite \eqref{xingjing} as follows:
\begin{equation}\label{maincal}
\begin{split}
0\leq& J(u^\varepsilon)-J(u)=-E^Q\Big[\int_0^T\delta H(t)\mathbf{1}_{E_\varepsilon}(t)dt\Big]\\
&+\frac{1}{2}E^Q\Big[(Y_T^{1,\varepsilon})^2\Big(\Phi_{xx}(T)+L_T\varphi_x(X_T,Y_{\cdot \wedge T})^2\widetilde{E}^{Q}\big[\widetilde{\Phi}^*_{z\mu}(T)\big]\Big)\Big]\\
&+E^Q\Big[Y_T^{1,\varepsilon}K_T^{1,\varepsilon}\varphi_x(X_T,Y_{\cdot \wedge T})\widetilde{E}^{Q}\big[\widetilde{\Phi}^*_\mu(T)\big]\Big]\\
&-\frac{1}{2}E^Q\Big[\int_0^T(Y_t^{1,\varepsilon})^2\Big(H_{xx}(t)+L_t\varphi_x(X_t,Y_{\cdot \wedge t})^2\widetilde{E}^{Q}\big[\widetilde{H}^*_{z\mu}(t)\big]\Big)dt\Big]\\
&-E^Q\Big[\int_0^T Y_t^{1,\varepsilon}K_t^{1,\varepsilon}\Big(H_{xl}(t)+\varphi_x(X_t,Y_{\cdot \wedge t})\widetilde{E}^{Q}\big[\widetilde{H}^*_\mu(t)\big]\Big)dt\Big]+o(\varepsilon).
\end{split}
\end{equation}
This motivates to introduce the following second order adjoint equation which is quite different from (\ref{4.56starstarstar}):
{\small
\begin{equation}\label{Secadj}\left\{
\begin{split}
dP_t^1=&-\Big\{H_{xx}(t)+L_t\varphi_x(X_t,Y_{\cdot \wedge t})^2\widetilde{E}^{Q}\big[\widetilde{H}^*_{z\mu}(t)\big]+\sigma_x(t)^2 P_t^1+2\sigma_x(t)Q_t^{1,1}+2L_t h_x(t)Q_t^{2,2}\Big\}dt\\
&+Q_t^{1,1}dB_t^1+Q_t^{1,2}dY_t,\ t\in[0,T],\\
P_T^1=&-\Phi_{xx}(T)-L_T\varphi_x(X_T,Y_{\cdot \wedge T})^2\widetilde{E}^{Q}\big[\widetilde{\Phi}^*_{z\mu}(T)\big];\\
dP_t^2=&-\Big\{H_{xl}(t)+\varphi_x(X_t,Y_{\cdot \wedge t})\widetilde{E}^{Q}\big[\widetilde{H}^*_\mu(t)\big]+\sigma_x(t)Q_t^{2,1}+h(t)Q_t^{2,2}\Big\}dt\\
&+Q_t^{2,1}dB_t^1+Q_t^{2,2}dY_t,\ t\in[0,T],\\
P_T^2=&-\varphi_x(X_T,Y_{\cdot \wedge T})\widetilde{E}^{Q}\big[\widetilde{\Phi}^*_\mu(T)\big].
\end{split}\right.
\end{equation}}

The equation is a linear one and one proves with arguments similar to those used for the first adjoint equation that there exists a unique solution
 $\big((P^1,(Q^{1,1},Q^{1,2})),\,(P^2,(Q^{2,1},Q^{2,2}))\big)\in\big(S_{\mathbb{F}}^{\infty-}([0,T],Q)\times (L_{\mathbb{F}}^{\infty-}([0,T],Q))^2\times
  S_{\mathbb{F}}^{\infty-}([0,T],Q)\times (L_{\mathbb{F}}^{\infty-}([0,T],Q))^2\big)$. Let us now apply It\^{o}'s formula to $P_t^1(Y_t^{1,\varepsilon})^2$ and
  $P_t^2Y_t^{1,\varepsilon}K_t^{1,\varepsilon}$, and take then expectation. This yields
\begin{align}\label{Dua2}
&-\frac{1}{2}E^Q\!\Big[\!(Y_T^{1,\varepsilon})^2\Big(\Phi_{xx}(T)\!+\!L_T\varphi_x(X_T,Y_{\cdot \wedge T})^2\widetilde{E}^{Q}\big[\widetilde{\Phi}^*_{z\mu}(T)\big]\!\Big)\!\Big]\! -\!E^Q\Big[Y_T^{1,\varepsilon}K_T^{1,\varepsilon}\varphi_x(X_T,Y_{\cdot \wedge T})\widetilde{E}^{Q}\big[\widetilde{\Phi}^*_\mu(T)\big]\Big]\notag\\
=&E^Q\Big[\frac{1}{2}P_T^1(Y_T^{1,\varepsilon})^2+P_T^2Y_T^{1,\varepsilon}K_T^{1,\varepsilon}\Big]\notag\\
=&-E^Q\Big[\frac{1}{2}\int_0^T(Y_t^{1,\varepsilon})^2\Big(H_{xx}(t)+L_t\varphi_x(X_t,Y_{\cdot \wedge t})^2\widetilde{E}^{Q}\big[\widetilde{H}^*_{z\mu}(t)\big]\Big)dt\Big]\\
&-E^Q\Big[\int_0^TY_t^{1,\varepsilon}K_t^{1,\varepsilon}\Big(H_{xl}(t)+\varphi_x(X_t,Y_{\cdot \wedge t})\widetilde{E}^{Q}\big[\widetilde{H}^*_\mu(t)\big]\Big)dt\Big]\notag\\
&+E^Q\Big[\frac{1}{2}\int_0^T|\delta\sigma(t)|^2P_t^1\mathbf{1}_{E_\varepsilon}(t)dt\Big]+R'_\varepsilon, \notag
\end{align}
where
\begin{equation*}
\begin{split}
R'_\varepsilon=&E^Q\Big[\frac{1}{2}\int_0^T\Big\{P_t^1R_t^{1,\varepsilon}\big(R_t^{1,\varepsilon}+2\sigma_x(t)Y_t^{1,\varepsilon}+2\delta\sigma(t)\mathbf{1}_{E_\varepsilon}(t)\big)
+2P_t^1\sigma_x(t)Y_t^{1,\varepsilon}\delta\sigma(t)\mathbf{1}_{E_\varepsilon}(t)\\
&\qquad\qquad\qquad+2Q_t^{1,1}Y_t^{1,\varepsilon}\big(R_t^{1,\varepsilon}+\delta\sigma(t)\mathbf{1}_{E_\varepsilon}(t)\big)\Big\}dt\Big]\\
&+E^Q\Big[\int_0^T\Big\{Q_t^{2,1}K_t^{1,\varepsilon}\big(R_t^{1,\varepsilon}+\delta\sigma(t)\mathbf{1}_{E_\varepsilon}(t)\big)+Q_t^{2,2}Y_t^{1,\varepsilon}L_t\big(R_t^{2,\varepsilon}+\delta h(t)\mathbf{1}_{E_\varepsilon}(t)\big)\Big\}dt\Big],
\end{split}
\end{equation*}
with $\displaystyle R_t^{1,\varepsilon}=\widetilde{E}^Q\Big[\int_0^{\widetilde{U}_t}\sigma_\mu(t,y)dy\cdot\widetilde{K}_t^{1,\varepsilon}\Big]+\widetilde{E}^Q
\big[\widetilde{\sigma}_\mu(t)\widetilde{L}_t\widetilde{\overline{V}}_t^{1,\varepsilon}\big],$\ and $\displaystyle
R_t^{2,\varepsilon}=\widetilde{E}^Q\Big[\int_0^{\widetilde{U}_t}h_\mu(t,y)dy\cdot\widetilde{K}_t^{1,\varepsilon}\Big]+\widetilde{E}^Q\big[\widetilde{h}_\mu(t)\widetilde{L}_t
\widetilde{\overline{V}}_t^{1,\varepsilon}\big].
$
With the help of Proposition \ref{tech}, we obtain that $|R_t^{i,\varepsilon}|\leq \rho_t(\varepsilon)\sqrt{\varepsilon},\ \varepsilon>0,\
t\in[0,T]$, $i=1,2$, with $\rho_t(\varepsilon)\rightarrow0\ (\varepsilon\searrow0)$
 (Recall that $\overline{V}_t^{1,\varepsilon}=\varphi_x(X_t,Y_{\cdot\wedge t})Y_t^{1,\varepsilon}$). Moreover, it holds that
 $|\rho_t(\varepsilon)|\leq C,\ t\in[0,T]$. Then, we deduce from Proposition \ref{EstofXL} that $|R'_\varepsilon|\leq o(\varepsilon)$.

Now, substituting \eqref{Dua2} in \eqref{maincal} yields
$$ 0\leq J(u^\varepsilon)-J(u)=-E^Q\Big[\int_0^T\Big(\delta H(t)+\frac{1}{2}|\delta\sigma(t)|^2 P_t^1\Big)\mathbf{1}_{E_\varepsilon}(t)
dt\Big]+o(\varepsilon),\mbox{ as } \varepsilon\searrow0. $$
Then, as $v\in\mathcal{U}_{ad}$ has been fixed arbitrarily, we can use Lebesgue's differentiation theorem to deduce that for all $v\in
\mathcal{U}_{ad}$, $dtdQ$-a.e.,
\begin{equation*}
\begin{split}
E^Q\Big[H(t,X_t,L_t,\mu_t,v_t,q_t^1,q_t^2)-H(t,X_t,L_t,\mu_t,u_t,q_t^1,q_t^2)\\
+\frac{1}{2}P_t^1\big|\sigma(t,X_t,\mu_t,v_t)-\sigma(t,X_t,\mu_t,u_t)\big|^2\,\Big|\mathcal{F}_t^Y\,\Big]\leq 0.
\end{split}
\end{equation*}
This allows to formulate our stochastic maximum principle for the case $\mu_t = P_{U_t} = P_{\varphi(X_t,Y_{\cdot\wedge t})}.$
\begin{theorem}
Assume the Assumption (H2). Let $u\in\mathcal{U}_{ad}$ be optimal and $(X,L)$ be the associated solution of system \eqref{system}. Then,
for all $v\in U$, it holds that for $dtdQ$-a.e. $(t,\omega)\in[0,T]\times\Omega$,
\begin{equation}\label{Smpmain}
\begin{split}
E^Q\Big[H(t,X_t,L_t,\mu_t,v,q_t^1,q_t^2)-H(t,X_t,L_t,\mu_t,u_t,q_t^1,q_t^2)\\
+\frac{1}{2}P_t^1\big|\sigma(t,X_t,\mu_t,v)-\sigma(t,X_t,\mu_t,u_t)\big|^2\,\Big|\mathcal{F}_t^Y\,\Big]\leq 0,
\end{split}
\end{equation}
where $\big((p^1,(q^1,\check{q}^1)),\, (p^2,(\check{q}^2,q^2))\big)$ and $\big((P^1,(Q^{1,1},Q^{1,2}), (P^2,(Q^{2,1},Q^{2,2})\big)$ are the unique solutions
to \eqref{Firadj} and \eqref{Secadj}, respectively.
\end{theorem}
Comparing \eqref{Smpmain} with \eqref{finalSMP} in Theorem \ref{th4.2}, we see that for $U_t = \varphi(X_t,Y_{\cdot\wedge t})$ we rediscover,
in some sense, Peng's classical optimality condition (of course, with a more complicated Hamiltonian here involving mean-field terms),
while \eqref{finalSMP} is more general and contains new terms. Of course, also, the second order adjoint equations (\ref{4.56starstarstar}) and \eqref{Secadj} in
Theorem \ref{Wellof1stVar} and Theorem \ref{Smpmain}, respectively, are quite different.

\end{document}